\documentclass[a4paper,leqno,11pt]{article}
\usepackage{amssymb,amsfonts,amsmath,amsthm}
\usepackage{epsfig,epsf,curves}
\usepackage{bbm,epic,eepic}
\usepackage[mathscr]{eucal}
\theoremstyle{remark}

\begin{document}

\sloppy

%------------------------- Definitions --------------------------

\def\Mat#1#2#3#4{\left(\!\!\!\begin{array}{cc}{#1}&{#2}\\{#3}&{#4}\\ \end{array}\!\!\!\right)}

\def\tAA{{\Bbb A}}
\def\CC{{\Bbb C}}
\def\HH{{\Bbb H}}
\def\NN{{\Bbb N}}
\def\QQ{{\Bbb Q}}
\def\RR{{\Bbb R}}
\def\ZZ{{\Bbb Z}}
\def\FF{{\Bbb F}}
\def\SS{{\Bbb S}}
\def\GG{{\Bbb G}}
\def\PP{{\Bbb P}}
\def\LL{{\Bbb L}}

\title{Traces on Hecke algebras and $p$-adic families of modular forms}

\author{J. Mahnkopf}

\date{Universit{\"a}t Wien}

\maketitle

{\bf Abstract. } We prove that any modular eigenform $f$ of level $\Gamma_1(Np)$, finite slope $\alpha$ and weight $k_0$ can be placed into a $p$-adic family of modular eigenforms $f_k$ of the same level and slope and weight $k$ varying over all natural numbers which are sufficiently close to $k_0$ in the $p$-adic sense.  Here, the term $p$-adic family means that a $p$-adic congruence between two weights $k$ and $k'$ entails a certain $p$-adic congruence between the corresponding eigenforms $f_k$ and $f_{k'}$. We also prove that the dimension of the slope $\alpha$ subspace of the space of modular forms of weight $k$ does not depend on the weight as long as we consider weights $k$ which are sufficiently close to each other in the $p$-adic sense. Both these statements are predicted by the Mazur-Gouvea conjecture, which has been proven by Coleman (cf. [C]) using methods from rigid analytic geometry. In contrast, our proof of these statements is based on a comparison of trace formulas. %This is a technique which we adopted from the Langlands program.

\vspace{1.5cm}

\centerline{\Large\bf Introduction}

{\bf 0.1.}  We fix a prime $p\in {\Bbb N}$ and an integer $N\in{\Bbb N}$ such that $(N,p)=1$. We select a Dirichlet character $\chi:\,({\Bbb Z}/(Np))^*\rightarrow{\Bbb C}^*$ and we denote by ${\cal M}_k(\Gamma_1(Np),\chi)$ resp. ${\cal S}_k(\Gamma_1(Np),\chi)$ the space of holomorphic modular forms of level $\Gamma_1(Np)$ and nebentype $\chi$ resp. the subspace of ${\cal M}_k(\Gamma_1(Np),\chi)$ consisting of cusp forms.
For any prime $\ell\in{\Bbb N}$ we denote by $T_\ell$ the Hecke operator of level $\Gamma_1(Np)$.
%
%[[i.e. $T_\ell$ corresponds to the double coset $\Gamma_1(Np)\Mat{1}{}{}{\ell}\Gamma_1(Np)$]]
%
$T_\ell$ acts on ${\cal M}_k(\Gamma_1(Np),\chi)$ and we say that $f$ is a normalized eigenform if $f$ has leading Fourier coefficient equal to $1$ and is an eigenvector for all Hecke operators $T_\ell$, $\ell$ prime. We note that $T_p$ is the Atkin-Lehner operator sending $f=\sum_n a_nq^n$ to $\sum_na_{pn}q^{n}$. 

The Mazur-Gouvea conjecture is a statement about certain subspaces of the space of cusp forms and the eigenforms contained therein. To explain this, denote by ${\Bbb C}_p$ the completion of an algebraic closure of ${\Bbb Q}_p$ with valuation $w:\,{\Bbb C}_p\rightarrow{\Bbb Q}$ normalized by $w(p)=1$. We fix an isomorphism $i:\,{\Bbb C}\cong {\Bbb C}_p$; thus, for any $\lambda\in{\Bbb C}$ its "$p$-adic value" $v_p(\lambda)=w\circ i(\lambda)$ is defined. The slope $\alpha$ subspace ${\cal S}_k(\Gamma_1(Np),\chi)^\alpha$, $\alpha\in{\Bbb Q}$,of ${\cal S}_k(\Gamma_1(Np),\chi)$ is the image $p(T_p)\,{\cal S}_k(\Gamma_1(Np),\chi)$, where $p\in\bar{\Bbb Q}[X]$ is the factor of the characteristic polynomial $\tilde{p}$ of $T_p$ acting on ${\cal S}_k(\Gamma_1(Np),\chi)$, which contains all roots of $\tilde{p}$, whose $p$-adic value is different from $\alpha$ (cf. section 1.1).
%
%;[[ note that ${\Bbb C}$ is algebraically closed, hence, the roots exist]]
%
Let $\omega:\,({\Bbb Z}/(p))^*\rightarrow \mu_{p-1}$, $\mu_{p-1}$ the group of $(p-1)$-st roots of unity, denote the Teichmuller character. 

The {\bf Mazur-Gouvea conjecture} then reads

{\it
1.) The dimension of ${\cal S}_k(\Gamma_1(Np),\chi\omega^{-k})^\alpha$ is locally constant in the $p$-adic sense, i.e.
$$
{\rm dim}\,{\cal S}_k(\Gamma_1(Np),\chi\omega^{-k})^\alpha={\rm dim}\,{\cal S}_{k'}(\Gamma_1(Np),\chi\omega^{-k'})^\alpha
$$
if $k\equiv k'\pmod{p^{[\alpha]+1}}$ and $k,k'\ge 2\alpha+2$ ("$[\;]$" denotes the Gaussian bracket).

2.) Assume ${\rm dim}\,{\cal S}_{k_0}(\Gamma_1(Np),\chi\omega^{-k_0})^\alpha=1$. Then, any eigenform $f\in{\cal S}_{k_0}(\Gamma_1(Np),\chi\omega^{-k_0})$, can be placed into a $p$-adic family of eigenforms, i.e. for any weight $k$ contained in the arithmetic progression ${\cal K}=\{k_0+mp^{[\alpha]+1},\,m=0,1,2,3,\ldots\}$ there is an eigenform $f_k=\sum_n a_{k,n}q^n\in {\cal S}_k(\Gamma_1(Np),\chi\omega^{-k})^\alpha$ such that $f_{k_0}=f$ and $f_k\equiv f_{k'}\pmod{p^m}$ if $k\equiv k'\pmod{p^m}$. Here, the congruence $f_k\equiv f_{k'}\pmod{p^m}$ means that $a_{k,n}\equiv a_{k',n}\pmod{p^m}$ for all $n$.

3.) In general, i.e. if ${\rm dim}\,{\cal S}_k(\Gamma_1(Np),\chi\omega^{-k})^\alpha$ is arbitrary, $f$ can be placed into an R-family of modular forms $f_k\in {\cal S}_k(\Gamma_1(Np),\chi)$, i.e. there is a finite and free ${\Bbb Z}_p[[T]]$-algebra $R$, a power series $F=\sum_n r_n q^n\in R[[q]]$ with $r_1=1$ and a family of continuous morphisms $\eta_k:\,R\rightarrow \bar{{\Bbb Q}_p}$, $k$ running through the arithmetic progression ${\cal K}$, such that $f_k=\sum_n\eta_k(r_n)q^n$. Furthermore the ${\Bbb Z}_p[[T]]$-rank of $R$ is less than or equal to the dimension of the spaces ${\cal S}_k(\Gamma_1(Np))^\alpha$ for all $k\in{\cal K}$. (cf. [G-M], Conjectures 1,2,3.) 
}

%We note that in case ${\rm dim}\,{\cal S}_k(\Gamma_1(Np))^\alpha$ arbitrary the family $(f_k)_k$ depends on the initial weight $k_0$ because this is true of the ${\Bbb Z}_p[[T]]$-rank of $R$. 

In case $\alpha=0$ (i.e. the "ordinary" case) the conjecture has been settled by Hida using his theory of ordinary forms. 
%
%[[, which is an integral element projecting onto ${\cal S}_k(\Gamma_1(Np),\chi)^0$.]]
% 
In general, the conjecture has been proven by Coleman (cf. [C]) except for the determination of the range of weights over which the family exists. Wan has shown that the family exists over a domain whose size depends quadratically on the slope $\alpha$ (cf. [W]). Finally there is a counterexample to the conjecture in its strong form (cf. [B-C]). Coleman's proof uses methods from rigid analytic geometry but he also mentions work on $p$-adic properties of modular forms using the Selberg trace formula and says that this line of research seems to have stopped in the mid seventies.

{\bf 0.2.} In this article we describe an approach to the construction of $p$-adic families of modular forms using a comparison of trace formulas. %and {\it essentially} prove the Mazur-Gouvea conjecture using a {\it comparison} of Trace formulas (cf. section 0.3 for precise statements of our results; see section 0.4 for more details concerning the proofs). 
To this end we reinterpret the Mazur-Gouvea conjecture as follows. We denote by $\Lambda_{k,\chi}^\alpha$ the set of all sequences $(\lambda_\ell)_\ell$, where $\ell$ runs over all primes $\ell\in{\Bbb N}$, such that there is an eigenform form $f\in{\cal M}_k(\Gamma_1(Np),\chi\omega^{-k})^\alpha$ with $T_\ell$-eigenvalue equal to $\lambda_\ell$. Hence,  $\Lambda_{k,\chi}^\alpha$ is the set of systems of Hecke eigenvalues occuring in ${\cal M}_k(\Gamma_1(Np),\chi\omega^{-k})^\alpha$ and the elements $\lambda\in\Lambda_{k,\chi}^\alpha$ correspond to eigenforms $f\in{\cal M}_{k}(\Gamma_1(Np),\chi\omega^{-k})^\alpha$. The statement that any modular form $f\in{\cal M}_{k_0}(\Gamma_1(Np),\chi\omega^{-k_0})^\alpha$ can be placed into a family of modular forms $(f_k)_k$, $f_k\in{\cal M}_k(\Gamma_1(Np),\chi\omega^{-k})^\alpha$, then (nearly) is equivalent to the existence of a family of transfer maps
$$
\tilde{\Psi}_k:\,\Lambda_{k_0,\chi}^\alpha\rightarrow \Lambda_{k,\chi}^\alpha
$$
such that

$\bullet$ $\tilde{\Psi}_k$ exists for all $k$, which are sufficiently close to $k_0$ in the $p$-adic sense 

$\bullet$ $k\equiv k_0\pmod{p^m}$ implies that $\tilde{\Psi}_k(\lambda)_\ell$ is congruent to $\lambda_\ell$ for all primes $\ell$ modulo a certain power of $p$ (which depends on $m$). 

Thus, to prove the Mazur-Gouvea conjecture we have to establish the existence of the maps $\tilde{\Psi}_k$ and in doing so we shall use a comparison of Trace formulas. To this end we assume in addition that

$\bullet$ $\tilde{\Psi}_k$ is a bijection

$\bullet$ Denote by ${\cal M}_{k,\chi\omega^{-k}}^\alpha(\lambda)$ the generalized simultaneous eigenspace attached to $\lambda\in\Lambda_{k,\chi}^\alpha$. Then ${\rm dim}\,{\cal M}_{k_0,\chi\omega^{-k_0}}^\alpha(\lambda)={\rm dim}\,{\cal M}_{k,\chi\omega^{-k}}^\alpha(\tilde{\Psi}_k(\lambda))$.

We prove a Basic trace identity of the following type: let $T$ be an element of our Hecke algebra and denote by $\pi^{Np}_{Mp}$ the projection operator from level $\Gamma_1(Np)$ to level $\Gamma_1(Mp)$; then the trace of 
the operator  $\pi^{Np}_{Mp}\circ T$ on ${\cal M}_{k_0}(\Gamma_1(Np),\chi\omega^{-k_0})^\alpha$ is congruent to its trace on ${\cal M}_{k}(\Gamma_1(Np),\chi\omega^{-k})^\alpha$ modulo a certain power of $p$. The presence of the operator $\pi^{Np}_{Mp}$ enables us to detect in addition to the Hecke eigenvalues the level of a modular form. Clearly, this identity has been proven with the above two assumptions in mind. On the other hand, using the Basic trace identity we shall establish as an essential step the existence of the {\it unramified} transfer: in analogy with $\Lambda_{k,\chi}^\alpha$ we denote by $\Lambda_{k,\chi,Np}^\alpha$ the set of all sequences $(\lambda_\ell)_{\ell\not|Np}$ such that there is a modular form $f\in{\cal M}_k(\Gamma_1(Np),\chi\omega^{-k})^\alpha$ with $T_\ell$-eigenvalue equal to $\lambda_\ell$, $\ell\not|Np$; then, there is a family of transfer maps
$$
{\Psi}_k:\,\Lambda_{k_0,\chi,Np}^\alpha\rightarrow \Lambda_{k,\chi,Np}^\alpha
$$
such that

$\bullet$ ${\Psi}_k$ exists for all $k$, which are sufficiently close to $k_0$ in the $p$-adic sense 

$\bullet$ $k\equiv k_0\pmod{p^m}$ implies that for all $\ell\not|Np$, the element ${\Psi}_k(\lambda)_\ell$ is congruent to $\lambda_\ell$ modulo a certain power of $p$ (cf. the Transfer Theorem in 0.4 for a precise statement of this)

$\bullet$ the level of ${\Psi}_k(\lambda)$ divides the level of $\lambda$. 

The divisibility ${\rm Level}({\Psi}_k(\lambda))\,|\,{\rm Level}(\lambda)$ implies that the dimension of the generalized eigenspace ${\cal M}_{k,\chi\omega^{-k}}^\alpha(\Psi(\lambda))$ attached to $\Psi_k(\lambda)$ is greater than or equal to the dimension of the generalized eigenspace ${\cal M}_{k_0,\chi\omega^{-k_0}}^\alpha(\lambda)$ attached to $\lambda$ and together with the equality of dimensions ${\rm dim}\,{\cal M}_k(\Gamma_1(Np),\chi\omega^{-k})^\alpha={\rm dim}\,{\cal M}_{k_0}(\Gamma_1(Np),\chi\omega^{-k_0})^\alpha$ (note that this equality is the first part of the Mazur-Gouvea conjecture) we obtain that $\Psi_k$ is a bijection and even that ${\rm dim}\,{\cal M}_{k_0,\chi\omega^{-k_0}}^\alpha(\lambda)={\rm dim}\,{\cal M}_{k,\chi\omega^{-k}}^\alpha(\Psi_k(\lambda))$. Thus, in the course of establishing the existence of the map $\Psi_k$ we find that the above two assumptions hold, at first for the unramified transfer, but it is not difficult to derive from this the existence of the transfer maps $\tilde{\Psi}_k$. We note that we do not need the expansion of the geoemtric side of the (topological) Trace formula as a sum of orbital integrals.

{\bf 0.3. } We describe our results in more detail. We use the following convention: since we assume that $N$ and $p$ are fixed throughout, if we say that a constant $C$ depends (only) on the slope $\alpha$ this shall always mean that $C$ depends on $\alpha$ and on $N$ and $p$. Also, we set $\Gamma=\Gamma_1(Np)$. Our first main result then reads

{\bf Theorem B }(cf. section 4.1). {\it There are constants ${\sf B}_\alpha$ and ${\sf C}_\alpha$ only depending on $\alpha$ such that the following holds. Let $\alpha\in{\Bbb Q}_{\ge 0}$; for all pairs of integers $k,k'\in{\Bbb N}$ satisfying 

$\bullet$ $k,k'\ge({{\sf C}_\alpha}\alpha+1)^2+2$ and 

$\bullet$ $k\equiv k'\pmod{p^{{\sf B}_\alpha}}$

it holds that
$$
{\rm dim}\,{\cal M}_k(\Gamma_1(Np),\chi\omega^{-k})^\alpha={\rm dim}\,{\cal M}_{k'}(\Gamma_1(Np),\chi\omega^{-k'})^\alpha.
$$

}

Theorem B confirms the independence of weight of the dimension of the slope subspace, which is the first part of the Mazur-Gouvea conjecture. However, we have no control over the constants ${\sf C}_\alpha$ and ${\sf B}_\alpha$; they tend to be very large and in this respect we fall short of proving the conjecture. Of course, an analogous remark applies to the constant ${\sf L}(\alpha,k)$ appearing in Theorems D1, D2 below. 

To state our second main result we introduce some notation. For any sequence $\lambda=(\lambda_\ell)_\ell$ resp. $\lambda=(\lambda_\ell)_{\ell\not|Np}$ of complex numbers, where $\ell$ runs over all rational primes resp. all rational primes not dividing $Np$, we set
$$
{\cal M}_{k,\chi\omega^{-k}}^\alpha(\lambda)=\{f\in{\cal M}_k(\Gamma_1(Np),\chi\omega^{-k})^\alpha:\,\forall \ell\,\exists n\in{\Bbb N}: (T_\ell-\lambda_\ell)^nf=0\},
$$
resp.
$$
{\cal M}_{k,\chi\omega^{-k}}^\alpha(\lambda)=\{f\in{\cal M}_k(\Gamma_1(Np),\chi\omega^{-k})^\alpha:\,\forall \ell\not|Np\,\exists n\in{\Bbb N}: (T_\ell-\lambda_\ell)^nf=0\}.
$$
Thus, ${\cal M}_{k,\chi\omega^{-k}}^\alpha(\lambda)$ is the generalized simultaneous eigenspace for all Hecke operators $T_\ell$, resp. all Hecke operators $T_\ell$ with $\ell\not|Np$. Note that for all primes $\ell$, which do not divide $Np$, the Hecke operator $T_\ell$ is diagonalizable on ${\cal M}_k(\Gamma_1(Np),\chi\omega^{-k})$, whereas $T_\ell$, $\ell|Np$, need not be semisimple. 
%
%[[Note also that all $T_\ell$ commute]]
%
We also set 
$$
\Lambda_{k,\chi}^\alpha=\{\lambda=(\lambda_\ell)_\ell:\,{\cal M}_{k,\chi\omega^{-k}}^\alpha(\lambda)\not=0\}
$$
and
$$
\Lambda_{k,\chi,Np}^\alpha=\{\lambda=(\lambda_\ell)_{\ell\not|Np}:\,{\cal M}_{k,\chi\omega^{-k}}^\alpha(\lambda)\not=0\}.
$$
We note that $\lambda\in\Lambda_{k,\chi}^\alpha$ is equivalent to the existence of an eigenform $f$ in ${\cal M}_k(\Gamma_1(Np),\chi\omega^{-k})$
with $T_\ell$ eigenvalue equal to $\lambda_\ell$. 
%
%[[although $T_\ell$ need not be diagonalizable on ${\cal M}_k(\Gamma_1(Np),\chi\omega^{-k})]]
%
Since ${\rm dim}\,{\cal M}_k(\Gamma_1(Np),\chi\omega^{-k})$ is finite for any weight $k$, there is an integer $A_k$ such that elements $\lambda,\lambda'\in\Lambda_{k,\chi}^\alpha$ which are different, already are different modulo $p^{A_k}$, i.e. there is a prime $\ell=\ell(\lambda,\lambda')$ such that $\lambda_\ell\not\equiv \lambda_\ell'\pmod{p^{A_k}}$ (cf. equation (18) in section 4.2). Our second main result then reads

{\bf Theorem D1 }(cf. section 4.3). {\it For any $\alpha\in{\Bbb Q}_{\ge 0}$ and any $k>2$ there are constants ${\sf a}={\sf a}(\alpha)$ and ${\sf b}(\alpha)$, which only depend on $\alpha$, and ${\sf L}(\alpha,k)$ such that the following holds: Let $\alpha\in{\Bbb Q}_{\ge 0}$ and let $k_0\ge ({{\sf C}_\alpha}\alpha+1)^2+2$. Then, for any  $\lambda\in\Lambda_{k_0,\chi,Np}^\alpha$ there is a family $(\lambda(k))_k$, where $k$ runs over all integers satisfying $k\ge({{\sf C}_\alpha}\alpha+1)^2+2$ and $k\equiv k_0\pmod{p^{{\sf L}(\alpha,k_0)}}$, which satisfies the following properties

$\bullet$ $\lambda(k)\in\Lambda_{k,\chi,Np}^\alpha$

$\bullet$ $\lambda(k_0)=\lambda$

$\bullet$ if $k\equiv k'\pmod{p^m}$ then $\lambda(k)\equiv \lambda(k')\pmod{p^{{\sf a}m+{\sf b}-A_{k_0}}}$.

Moreover, ${\sf a}={\sf a}(\alpha)$ satisfies  
$$
{\sf a}\le \frac{1}{{\rm dim}\,{\cal M}_{k}(\Gamma,\chi\omega^{-k})^\alpha}
$$ 
for all $k\ge ({{\sf C}_\alpha}\alpha+1)^2+2$ and 
$$
{\sf a}={\rm min}\,\{\frac{1}{{\rm dim}\,{\cal M}_{k}(\Gamma,\chi\omega^{-k})^\alpha},\; k\ge ({{\sf C}_\alpha}\alpha+1)^2+2\}
$$ 
if $\alpha=0$. 

Finally, if $(\mu(k))_k$, $\mu(k)\in\Lambda_{k,\chi,Np}^\alpha$ is another family satisfying the above three conditions, then $\mu(k)=\lambda(k)$ for all $k$. 

}

%[[Does Corollary C3 imply that we can choose ${\sf L}(\alpha,k)$ such that it only depends on $k$ moduko $p^{{\sf L}(\alpha,k)}$ ???. That ${\sf a}$ only depends on $k$ modulo $p^{{\sf B}_\alpha}$ follows directly from the definition of ${\sf a}$ and Theorem B]]

We note that Theorem B implies that the set $\{\frac{1}{{\rm dim}\,{\cal M}_{k}(\Gamma,\chi\omega^{-k})^\alpha}, k\ge ({{\sf C}_\alpha}\alpha+1)^2+2\}$ is finite, in particular, ${\sf a}(\alpha)$ is strictly positive.

In the special case that ${\rm dim}\,{\cal M}_{k_0}(\Gamma,\chi\omega^{-k_0})^\alpha=1$ we obtain the following stronger result.

{\bf Theorem D2. } {\it Let $\alpha\in{\Bbb Q}_{\ge 0}$, let $k_0\ge ({{\sf C}_\alpha}\alpha+1)^2+2$ and assume in addition that ${\rm dim}\,{\cal M}_{k_0}(\Gamma,\chi\omega^{-k_0})^\alpha=1$. Then, for any $\lambda\in\Lambda_{k_0,\chi,Np}^\alpha$ there is a family $(\lambda(k))_k$, where $\lambda(k)\in\Lambda_{k,\chi,Np}^\alpha$ and $k$ runs over all integers satisfying $k\ge({{\sf C}_\alpha}\alpha+1)^2+2$ as well as $k\equiv k_0\pmod{p^{{\sf B}_\alpha}}$ which satisfies the properties

$\bullet$ $\lambda(k_0)=\lambda$

$\bullet$ if $k\equiv k'\pmod{p^m}$ then $\lambda(k)\equiv \lambda(k')\pmod{p^{{\sf a}m+{\sf b}}}$.

Here, ${\sf a}={\sf a}(\alpha)$ and ${\sf b}={\sf b}(\alpha)$ are the integers appearing in Theorem D1, hence, ${\sf a}\le 1$ and ${\sf a}= 1$ if $\alpha=0$. 

}

It remains to look at the primes $\ell$ dividing $Np$. To this end, let $\lambda\in\Lambda_{k_0,\chi}^\alpha$, where $k_0>({\sf C}_\alpha\alpha+1)^2+2$ and let $(\lambda_k)_k$ be the family of elements $\lambda_k\in\Lambda_{k,\chi,Np}^\alpha$ as in Theorem D1. We denote by $\Psi_{k,\chi,\ell}^\alpha$ the characteristic polynomial of $T_\ell$ acting on ${\cal M}_{k,\chi\omega^{-k}}^\alpha(\lambda(k))$. Our last result then reads

{\bf Theorem F }(cf. section 4.5). {\it Let $k_0>({\sf C}_\alpha\alpha+1)^2+2$ and let $k',k\ge ({\sf C}_\alpha\alpha+1)^2+2$ be weights which both are congruent to $k_0\pmod{p^{{\sf L}(\alpha,k)}}$. Then, with the above notations, for any $\ell|N$ we have
$$
v_p(\Psi_{k,\chi,\ell}^\alpha-\Psi_{k',\chi,\ell}^\alpha)
\ge 
{\sf a}m+{\sf b}-(d_{k_0,\chi}^\alpha+1)A_{k_0}-v_p(d_{k_0,\chi}^\alpha!)).
$$
where $d_{k_0,\chi}^\alpha={\rm dim}\,{\cal M}_{k_0}(\Gamma,\chi\omega^{-k_0})^\alpha$.

}

In conjunction with the proposition in section 4.5, Theorem F shows that for any prime $\ell$ dividing $N$ and any $k\equiv k_0\pmod{p^{{\sf L}(\alpha,k_0)}}$ we can select a root $\lambda(k)_\ell$ of $\Psi_{k,\chi,\ell}^\alpha$, i.e. an eigenvalue of $T_\ell$ acting on ${\cal M}_{k,\chi\omega^{-k}}^\alpha(\lambda(k))$, such that the resulting family $(\lambda(k)_\ell)_k$ satisfies the congruence
$$
v_p(\lambda(k)_\ell-\lambda(k')_\ell)\ge \frac{{\sf a}m+{\sf b}-(d_{k_0,\chi}^\alpha+1)A_{k_0}-v_p(d_{k_0,\chi}^\alpha!)}{d_{k_0,\chi}^\alpha}
$$ 
if $k\equiv k'\pmod{p^m}$. Thus, the congruences, which we obtain at the ramified places are somewhat weaker than those at the unramified places. Although we formulated Theorems D1, D2 only for systems of eigenvalues $(\lambda)_\ell$ with $\ell\not|Np$, an inspection of the proof shows that these Theorems also hold if we include the prime $p$, i.e. if we consider systems of eigenvalues $(\lambda_\ell)$ with $\ell\not|N$. Thus, altogether we obtain that any $\lambda\in\Lambda_{k_0,\chi}^\alpha$ can be placed into a $p$-adic family $(\lambda_k)_k$ of elements $\lambda_k\in\Lambda_{k,\chi}^\alpha$, where $k$ runs over all integers greater than $({\sf C}_\alpha\alpha+1)^2+2$ and congruent to $k_0$ $\pmod{p^{{\sf L}(\alpha,k_0)}}$; see our final Theorem G in section 4.5. Since any $\lambda\in \Lambda_{k,\chi}^\alpha$ corresponds to an eigenform $f\in {\cal M}_k(\Gamma_1(Np),\chi\omega^{-k})$ with $T_\ell$ eigenvalue equal to $\lambda_\ell$, Theorms D1 and D2 and F essentially confirm the existence of $p$-adic families of modular forms as predicted by the second part of the Mazur-Gouvea conjecture.

{\bf 0.4.} We explain the proof of the above results in some detail. We denote by $L_k$ the finite dimensional irreducible ${\rm GL}_2$ representation of highest weight $k$ and for any ${\Bbb Z}$-algebra $R$ we set $L_{k,R}=R\otimes_{\Bbb Z} L_k$. The space ${\cal M}_k(\Gamma_1(Np))$ of modular forms of weight $k$ and level $\Gamma_1(Np)$ is isomorphic as Hecke module to the group cohomology $H^1(\Gamma_1(Np),L_{k,{\Bbb C}})$ and we will replace the study of the Hecke module ${\cal M}_k(\Gamma_1(Np))$ by the study of the Hecke module $H^1(\Gamma_1(Np),L_{k,{\Bbb C}})$. In particular, we replace the Selberg trace formula by the Topological trace formula, which computes the Lefschetz number of a Hecke correspondence (cf. [B], [G-MacPh]). 

In  section 1.1 and 1.2, we discuss the slope decomposition of a finitely generated ${\cal O}$-module $V$, where ${\cal O}/{\Bbb Z}_p$ is a finite extension, with respect to an endomorphism $T$ of $V$. We look to the particular case of the space of modular forms with $T$ the Atkin-Lehner operator (cf. section 1.4) and we prove the following

{\bf Theorem A } (cf. section 1.5). {\it There is a constant $M(\alpha)$, which only depends on $\alpha$, such that
$$
\sum_{0\le \beta\le \alpha}{\rm dim}\,H^1(\Gamma_1(Np^{[\alpha]+1}),L_{k,{\Bbb C}})^\beta\le M(\alpha)
$$
for all weights $k\ge 2$.

}

To prove Theorem A we introduce a certain submodule of $L_{k,{\cal O}}$. We set $r=[\alpha]+1$, $\Gamma=\Gamma_1(Np^r)$ and we choose a finite, integral extension ${\cal O}/{\Bbb Z}_p$, which splits $T_p$. We then define the ${\cal O}$-submodule
$$
U_{k,{\cal O},r}=\bigoplus_{v\in L_{k,{\cal O}}\atop {\rm weight}(v)>r}{\cal O}v.
$$
Here, the weight of an element $v\in L_{k,{\cal O}}$ is defined with respect to the split torus of diagonal matrices in ${\rm GL}_2$. $U_{k,{\cal O},r}$ satisfies the following properties

- ${\cal O}/(p^r)\otimes U_{k,{\cal O},r}$ still is a $\Gamma_1(Np^r)$-module (cf. Lemma 1.) in section 1.4) and the long exact cohomology sequence will yield
$$
H^1(\Gamma,{\cal O}/(p^r)\otimes L_{k,{\cal O}})/H^1(\Gamma,{\cal O}/(p^r)\otimes U_{k,{\cal O},r})\cong H^1(\Gamma,{\cal O}/(p^r)\otimes (L_{k,{\cal O}}/U_{k,{\cal O},r})).
$$

- The isomorphism class of the quotient module ${\cal O}/(p^r)\otimes (L_{k,{\cal O}}/U_{k,{\cal O},r})$ only depends on $k\pmod{p^{2r}}$ (cf. the Proposition in section 1.4), hence, the ${\cal O}/(p^r)$-module $H^1(\Gamma,{\cal O}/(p^r)\otimes (L_{k,{\cal O}}/U_{k,{\cal O},r}))$ also depends only on $k\pmod{p^{2r}}$ and its ${\cal O}/(p^r)$-rank therefore can be bounded independently of $k$.

- The Atkin Lehner operator $T_p$ annihilates $H^1(\Gamma,{\cal O}/(p^r)\otimes U_{k,{\cal O},r})$ (cf. the Lemma in section 1.5).

The first two properties imply that the ${\cal O}/(p^r)$-rank of $H^1(\Gamma,{\cal O}/(p^r)\otimes L_{k,{\cal O}})/H^1(\Gamma,{\cal O}/(p^r)\otimes U_{k,{\cal O},r})$ is bounded independently of $k$. We denote by ${\rm ker}\,T_p$ the kernel of $T_p$ acting on $H^1(\Gamma,{\cal O}/(p^r)\otimes L_{k,{\cal O}})^\beta$. The third property implies that the quotient $H^1(\Gamma,{\cal O}/(p^r)\otimes L_{k,{\cal O}})^\beta/{\rm ker}\,T_p$ is a factor of $H^1(\Gamma,{\cal O}/(p^r)\otimes L_{k,{\cal O}})^\beta/H^1(\Gamma,{\cal O}/(p^r)\otimes U_{k,{\cal O},r})^\beta$. On the other hand, the key Lemma in section 1.2 allows us to bound the ${\cal O}$-rank of $H^1(\Gamma,L_{k,{\cal O}})^\beta$ for all $0\le \beta\le\alpha$ in terms of the cardinality of the quotient $H^1(\Gamma,{\cal O}/(p^r)\otimes L_{k,{\cal O}})^\beta/{\rm ker}\,T_p$ (note that $H^1(\Gamma,{\cal O}/(p^r)\otimes L_{k,{\cal O}})^\beta$ is a finite group). Thus, altogether
we arrive at the claim of Theorem A. Clearly, Theorem A implies the boundedness of $\sum_{0\le \beta\le \alpha}{\rm dim}\,H^1(\Gamma_1(Np),L_{k,{\Bbb C}})^\beta$ for all $k\ge 2$. We note that in the ordinary case, i.e. $\alpha=0$, the quotient $L_{k,{\cal O}}/U_{k,{\cal O},r}\cong {\cal O}$. Thus, we may replace $H^1(\Gamma,{\cal O}/(p^r)\otimes (L_{k,{\cal O}}/U_{k,{\cal O},r}))$ by $H^1(\Gamma,{\cal O}/(p^r))$ and if we do this, we essentially find Hida's proof of the boundedness of $H^1(\Gamma_1(Np),L_{k,{\Bbb C}})^0$. Since it is easily available using the methods introduced so far we also prove in section 1 that the Hecke operator $T_p$ acts nilpotently on the torsion part of the cohomology of $\Gamma_1(Np)$ with ${\mathfrak p}$-adic coefficients (cf. the Proposition in section 1.5).

{\bf Note.} After completion of the manuscript we have learned that a proof of the boundedness of the dimension of the slope spaces using a very similar idea has already appeared in [Bu] !

In section 2 we review the Topological trace formula. For the final formula see the Theorems in section 2.4 and 2.5. Since we only deal with the rank $1$ case our main reference will be [B]. A trace formula valid for groups of arbitrary rank can be found in [G-MacPh].

In section 3 we prove the basic trace identity from which we will deduce (in section 4) Theorems B, D1, D2 and F. We work with a certain subalgebra of the full Hecke algebra, which we define in sections 4.1 and 4.2 as follows. We denote by ${\cal H}_{\Bbb Z}=\langle T_\ell,\,\ell\;{\rm prime}\rangle_{\Bbb Z}$ the classical Hecke algebra generated by all Hecke operators $T_\ell$ of level $\Gamma_1(Np)$, $\ell$ prime, and all diamond operators $\langle\epsilon\rangle$, $(\epsilon,Np)=1$, of level $Np$; we then enlarge ${\cal H}$ by adjoining certain non-classical Hecke operators $S_{u,\delta}$, $u\in{\Bbb N}$, $(\delta,Np)=1$, which correspond to primes $\ell$ dividing $Np$ (cf. section 3.2). The idea is that the enlarged Hecke algebra contains the projection operator $\pi^{Np}_{Mp}:\,{\cal M}_k(\Gamma_1(Np),\chi\omega^{-k})\rightarrow {\cal M}_k(\Gamma_1(Mp),\chi\omega^{-k})$ from level $Np$ to level $Mp$ ($M|N$; cf. the Corollary in section 3.2). Thus, in addition to the eigenvalues of the Hecke operators $T_\ell$, $\ell$ prime, the enlarged Hecke algebra also is able to detect the level of a modular form $f\in {\cal M}_k(\Gamma_1(Np),\chi\omega^{-k})$. On the other hand, we do not know whether the Hecke operators $T_\ell$,$\ell|N$ commute with the operators $S_{u,\delta}$, hence, the enlarged Hecke algebra need not be commutative. This is the reason why at first we only consider Hecke eigenvalues for primes $\ell$ not dividing $Np$ (cf. Theorem D1, D2; $T_\ell$, $\ell\not|N$ and $S_{u,\delta}$ commute). In section 3.3 we review the slope decomposition of the space of modular forms and in section 3.4 we construct "approximative" idempotents attached to the slope decomposition of ${\cal M}_k(\Gamma_1(Np),\chi\omega^{-k})$: for any pair of weights $k,k'$ we inductively construct elements $e_{\alpha}=e_{k,k',\alpha}\in E[T_p]$, where $E$ is a splitting field for $T_p$ on ${\cal M}_k(\Gamma_1(Np),\chi\omega^{-k})$ as well as on ${\cal M}_{k'}(\Gamma_1(Np),\chi\omega^{-k'})$ and $\alpha$ runs over the slopes $\alpha_1<\alpha_2<\ldots<\alpha_r$ appearing in ${\cal M}_k(\Gamma_1(Np),\chi\omega^{-k})$ or in ${\cal M}_{k'}(\Gamma_1(Np),\chi\omega^{-k'})$ such that 
$$
{\rm tr}\,e_{k,k',\alpha}|_{{\cal M}_k(\Gamma_1(Np),\chi\omega^{-k})^\beta}\equiv\left\{
\begin{array}{ccc}
1\pmod{p}&{\rm if}&\beta=\alpha\\
0\pmod{p}&{\rm if}&\beta\not=\alpha\\
\end{array}
\right.
$$
and also
$$
{\rm tr}\,e_{k,k',\alpha}|_{{\cal M}_{k'}(\Gamma_1(Np),\chi\omega^{-k'})^\beta}\equiv\left\{
\begin{array}{ccc}
1\pmod{p}&{\rm if}&\beta=\alpha\\
0\pmod{p}&{\rm if}&\beta\not=\alpha\\
\end{array}
\right.
$$
(cf. the Proposition in section 3.4). 
%
%Moreover, $e_{\alpha,k,k'}$ does not depend "too much" on the weights $k$ and $k'$: it can be written $e_{\alpha,k,k'}=p_{\alpha,k,k'}(T_p)$, where the degree of the polynomial $p_{\alpha,k,k'}\in E[X]$ can be bounded independently of $k,k'$. The verification of this property needs the result of Theorem A. 
%
In section 3.5 we then prove the

{\bf Basic trace identity. } {\it There is a constant ${\sf C}_\alpha$, which only depends on $\alpha$, such that the following holds. Let $\alpha\in{\Bbb Q}_{\ge 0}$. Assume that $k,k'\in{\Bbb N}$ satisfy 

$\bullet$ $k,k'\ge ({{\sf C}_\alpha}\alpha+1)^2+2$  

$\bullet$ $k\equiv k'\pmod{p^m}$ with $m\ge {{\sf C}_\alpha}\alpha+1$. 

Let $\ell_1,\ldots,\ell_s\in{\Bbb N}$ be prime numbers and let $M\in{\Bbb N}$ be any divisor of $N$; then 
$$
{\rm tr}\,\pi^{Np}_{Mp} T_{\ell_1}^{r_1}\cdot\ldots\cdot T_{\ell_s}^{r_s}e_\alpha^{[\frac{m}{{\sf C}_\alpha\alpha}]}|_{{\cal M}_k(\Gamma,\chi\omega^{-k})^\alpha}
\equiv{\rm tr}\,\pi^{Np}_{Mp} T_{\ell_1}^{r_1}\cdot\ldots\cdot T_{\ell_s}^{r_s}e_\alpha^{[\frac{m}{{\sf C}_\alpha\alpha}]}|_{{\cal M}_{k'}(\Gamma,\chi\omega^{-k'})^\alpha}\pmod{p^{[\frac{m}{{{\sf C}_\alpha}\alpha}]-v_p(\varphi(N))}}.
$$
In case $\alpha=0$ we even obtain
$$
{\rm tr}\,\pi^{Np}_{Mp} T_{\ell_1}^{r_1}\cdot\ldots\cdot T_{\ell_s}^{r_s}|_{{\cal M}_k(\Gamma,\chi\omega^{-k})^\alpha}
\equiv{\rm tr}\,\pi^{Np}_{Mp} T_{\ell_1}^{r_1}\cdot\ldots\cdot T_{\ell_s}^{r_s}|_{{\cal M}_{k'}(\Gamma,\chi\omega^{-k'})^\alpha}\pmod{p^{m-v_p(\varphi(N))}}.
$$

%[[ die folgende Kongruenz brauchen wir nicht mehr

%2.) Assume that $k,k'\ge 2\alpha+2$; then then the following congruence, which does not involve the slope, holds true
%$$
%{\rm tr}\,\pi^{Np}_{Mp} T_{\ell_1}^{r_1}\cdot\ldots\cdot T_{\ell_s}^{r_s}T_p^{\left[\frac{m}{2\alpha}\right]+1}|_{{\cal M}_k(\Gamma,\chi\omega^{-k})^\alpha}
%\equiv{\rm tr}\,\pi^{Np}_{Mp} T_{\ell_1}^{r_1}\cdot\ldots\cdot T_{\ell_s}^{r_s}T_p^{\left[\frac{m}{2\alpha}\right]+1}|_{{\cal M}_{k'}(\Gamma,\chi\omega^{-k'})^\alpha}\pmod{p^{m-v_p(\varphi(N))}}.
%$$
%]]
}

Of course, the proof of this is an application of the trace formula. We note that the element $e_0$ has no denominators, which is the reason for the stronger statement in case $\alpha=0$. 

In section 4 we deduce Theorems B, D1, D2 and F from the basic trace identity. We note that in order to deduce these Theorems we need to know in advance that the dimension of ${\cal M}_k(\Gamma_1(Np),\chi\omega^{-k})^\alpha$ is bounded independently of $k$ (cf. Theorem A). 
%Without invoking Theorem A we are not able to prove anything. 
We first explain the proof of Theorem B. We denote by $\Psi=\sum_{j=0}^d (-1)^j a_j X^{d-j}$ resp. $\Psi'=\sum_{j=0}^{d'} (-1)^j a'_j X^{d'-j}$ the characteristic polynomial of $T_p$ acting on ${\cal M}_k(\Gamma_1(Np),\chi\omega^{-k})^\alpha$ resp. on ${\cal M}_{k'}(\Gamma_1(Np),\chi\omega^{-k'})^\alpha$. This defines the coefficients $a_j$ resp. $a'_j$ for $0\le j\le d$ resp. $0\le j\le d'$ and for $j>d$ resp. $j>d'$ we define $a_j$ resp. $a'_j$ to be equal to $0$. Clearly, $d={\rm dim}\,{\cal M}_k(\Gamma_1(Np),\chi\omega^{-k})^\alpha$ and $d'={\rm dim}\,{\cal M}_{k'}(\Gamma_1(Np),\chi\omega^{-k'})^\alpha$ and we have to show that $d=d'$ if $k$ and $k'$ are sufficiently close in the $p$-adic 
sense. Assume $d>d'$; we lead this assumption to a contradiction by computing the $p$-adic value of the constant coefficient $a_d$ of $\Psi$ in two different ways:

- since all eigenvalues of $T_p$ acting on ${\cal M}_k(\Gamma_1(Np),\chi\omega^{-k})^\alpha$ have $p$-adic value equal to $\alpha$ we find $v_p(a_d)=\alpha d$.

- By a classical formula (cf. equation (1) in section 4) $a_j$ is a linear expression in the terms ${\rm tr}\,T_p^h|_{{\cal M}_k(\Gamma_1(Np),\chi\omega^{-k})^\alpha}a_{j-h}$, where $h=1,\ldots,j$ and, quite analogous, $a'_j$ is a linear expression in the terms ${\rm tr}\,T_p^h|_{{\cal M}_{k'}(\Gamma_1(Np),\chi\omega^{-k'})^\alpha}a'_{j-h}$, where $h=1,\ldots,j$. This holds for all $j\in{\Bbb N}\cup\{0\}$, in particular it holds for $j=d$. The basic trace identity therefore will imply that $a_d\equiv a'_d\pmod{p^{v_p(\alpha d)+1}}$ if $k$ and $k'$ are sufficiently close. Since $d'<d$ we have $a'_d=0$, hence, $v_p(a_d)>\alpha d$. Contradiction !

%Thus, our assumption is wrong and we have proven that $d=d'$. 

In section 4.2 and 4.3 we show the existence of $p$-adic families of modular forms as in Theorems D1 and D2. To this end, in section 4.2 relying on the basic trace identity we prove the following

{\bf Theorem }({Transfer from weight $k$ to weight $k'$}).  {\it There are constants ${\sf B}_\alpha$, ${\sf a}={\sf a}(\alpha)$ and ${\sf b}={\sf b}(\alpha)$ such that the following holds. Let $\alpha\in{\Bbb Q}_{\ge 0}$ and assume that $k,k'\in{\Bbb N}$ satisfy

- $k,k'\ge ({{\sf C}_\alpha}\alpha+1)^2+2$ 

- $k\equiv k'\pmod{p^m}$ where $m\ge {{\sf B}_\alpha}$. 

Let $\lambda=(\lambda_\ell)_{\ell\not|Np}\in\Lambda_{k,\chi,Np}^\alpha$ and denote by $F$ the prime-to-$p$ part of the conductor of  $\lambda$. Then, there is a $\lambda'=(\lambda_\ell')_{\ell\not|Np}\in\Lambda_{k',\chi,Np}^\alpha$ such that

$\bullet$ ${\cal M}_{k',\chi\omega^{-k'}}^\alpha(\lambda')^{\Gamma_1(Fp)}\not=(0)$, i.e. the conductor of $\lambda'$ is a divisor of $Fp$

$\bullet$ $\lambda_\ell'\equiv\lambda_\ell\pmod{p^{{\sf a}m+{\sf b}-A_k}}$ for all $\ell\not|Np$. 

Moreover, ${\sf a}={\sf a}(\alpha)$ satisfies
$$
{\sf a}(\alpha,k)\le \frac{1}{{\rm dim}\,{\cal M}_k(\Gamma,\chi\omega^{-k})^\alpha}
$$ 
for all $k\ge ({\sf C}_\alpha\alpha+1)^2+1$ and even
$$
{\sf a}(\alpha,k)={\rm min}\,\{ \frac{1}{{\rm dim}\,{\cal M}_k(\Gamma,\chi\omega^{-k})^\alpha},\;k\ge ({\sf C}_\alpha\alpha+1)^2+1\}
$$
if $\alpha=0$. If in addition ${\rm dim}\,{\cal M}_k(\Gamma,\chi\omega^{-k})^\alpha=1$ holds we even obtain the congruence
$$
\lambda_\ell'\equiv\lambda_\ell\pmod{p^{{\sf a}m+{\sf b}}}.
$$
for all $\ell\not|Np$.
}

Here, the conductor of $\lambda$ is defined as the smallest (in the sense of divisibility) integer $M$ such that ${\cal M}_{k,\chi\omega^{-k}}^\alpha(\lambda)^{\Gamma_1(M)}\not=0$; clearly the conductor of $\lambda$ divides $Np$ and coincides with the conductor of the automorphic representation corresponding to $\lambda$. Thus, the transfer theorem not only yields a modular form satisfying the requested congruences, it shows in addition that the conductor of the transferred element $\lambda'$ divides the conductor of $\lambda$; this will be essential in the construction of $p$-adic families of modular forms (see below).  We note that ${\sf a}(\alpha)$ and ${\sf b}(\alpha)$ are the constants appearing in Theorem D1. 
%The idea of proof of the above theorem is somewhat similar to the proof of (cases of) the transfer of automorphic representations from a group $G$ to a group $G'$, which are predicted by the principle of functoriality
To prove the transfer theorem we assume there is no $\lambda'$ satisfying the two conditions of the theorem. The slope $\alpha$ subspace  decomposes ${\cal M}_k(\Gamma_1(Np),\chi\omega^{-k})^\alpha=\bigoplus_{\lambda\in\Lambda_{k,\chi,Np}^\alpha}{\cal M}_{k,\chi\omega^{-k}}^\alpha(\lambda)$ and there is a completely analogous decomposition for weight $k'$. We construct an element $e\in{\cal H}$ such that
$$
e\,{\cal M}_k(\Gamma_1(Np),\chi\omega^{-k})^\alpha={\cal M}_{k,\chi\omega^{-k}}^\alpha(\lambda)\quad\mbox{and}\quad
e\,{\cal M}_{k'}(\Gamma_1(Np),\chi\omega^{-k'})^\alpha=0.\leqno(*)
$$
To lead our assumption to a contradiction, again, we compare the trace of the operator $\pi^{Np}_{Fp}\circ e\circ e_\alpha^L$, where $L$ is a certain $p$-power (depending on $m$) on ${\cal M}_k(\Gamma_1(Np),\chi\omega^{-k})^\alpha$ and on ${\cal M}_{k'}(\Gamma_1(Np),\chi\omega^{-k'})^\alpha$ in two different ways:

- Using the basic trace identity we will find that the two traces are congruent to each other modulo a certain power of $p$

- Using equation (*) we will find that the two traces are not congruent to each other modulo the same power of $p$. Contradiction !

In section 4.3 we will derive from the transfer theorem the existence of $p$-adic families of modular forms as in Theorems D1, D2. We proceed as follows. Let $k_0,k>({\sf C}_\alpha\alpha+1)^2+2$ with $k\equiv k_0\pmod{p^m}$, $m>{\sf B}_\alpha$. For every $\lambda\in\Lambda_{k_0,\chi,Np}^\alpha$ there is a $\lambda'\in\Lambda_{k,\chi,Np}^\alpha$ satisfying the two conditions in the transfer theorem, hence, we obtain a map
$$
\begin{array}{cccc}
\Psi_k:&\Lambda_{k_0,\chi,Np}^\alpha&\rightarrow&\Lambda_{k,\chi,Np}^\alpha\\
&\lambda&\mapsto&\lambda'
\end{array}
$$
such that $\lambda'\equiv\lambda\pmod{p^{{\sf a}m+{\sf b}-A_{k_0}}}$.  Let $\lambda,\mu\in\Lambda_{k_0,\chi,Np}^\alpha$, $\lambda\not=\mu$. Then, $\lambda\not\equiv\mu\pmod{p^{A_{k_0}}}$ by the definition of $A_{k_0}$. Hence, if $k$ is close enough to $k_0$ such that ${\sf a}m+{\sf b}-A_{k_0}>A_{k_0}$, we see that $\lambda'\not\equiv\mu'\pmod{p^{A_{k_0}}}$. In particular, the $\lambda'$, $\lambda\in\Lambda_{k_0,\chi,Np}^\alpha$ are all distinct. On the other hand, since the conductor of $\lambda'$ divides the conductor of $\lambda$, we will find that ${\rm dim}\,{\cal M}_{k_0,\chi\omega^{-k_0}}^\alpha(\lambda)\le {\rm dim}\,{\cal M}_{k,\chi\omega^{-k}}^\alpha(\lambda')$ (cf. Corollary C1). Taking into account that ${\rm dim}\,{\cal M}_k(\Gamma_1(Np),\chi\omega^{-k})^\alpha={\rm dim}\,{\cal M}_{k_0}(\Gamma_1(Np),\chi\omega^{-k_0})^\alpha$ if $k$ is sufficiently close to $k_0$ by Theorem B we deduce that ${\cal M}_{k}(\Gamma_1(Np),\chi\omega^{-k})^\alpha=\bigoplus_{\lambda\in\Lambda_{k_0,\chi,Np}^\alpha} {\cal M}_{k,\chi\omega^{-k}}^\alpha(\lambda')$ and, hence, $\Lambda_{k,\chi,Np}^\alpha=\{\lambda',\,\lambda\in\Lambda_{k_0,\chi,Np}^\alpha\}$. Thus, altogether we obtain

{\bf Corollary C2. }{\it For any $\alpha\ge 0$ and $k>2$ there is a constant ${\sf L}(\alpha,k)$ such that the following holds. Assume that $k\equiv k_0\pmod{p^m}$ with $m\ge{\sf L}(\alpha,k_0)$; then for every $\lambda\in\Lambda_{k_0,\chi,Np}^\alpha$ there is precisely one $\lambda'\in\Lambda_{k,\chi,Np}^\alpha$ such that $\lambda'\equiv\lambda\pmod{p^{{\sf a}m+{\sf b}-A_{k_0}}}$ and the transfer map $\Psi_k:\Lambda_{k_0,\chi,Np}^\alpha\rightarrow \Lambda_{k,\chi,Np}^\alpha$ is a bijection }

Using the above Corollary we are now able to construct the $p$-adic family of modular forms as follows. Let $\lambda\in\Lambda_{k_0,\chi,Np}^\alpha$, $k_0>({\sf C}_\alpha\alpha+1)^2+1$. For any weight $k\equiv k_0\pmod{p^{{\sf L}(\alpha,k_0)}}$ there is a (unique) element  $\lambda_k=\Psi_k(\lambda)\in\Lambda_{k,\chi,Np}^\alpha$ such that $\lambda_k\equiv\lambda\pmod{p^{{\sf a}m+{\sf b}-A_{k_0}}}$, where $k\equiv k_0\pmod{p^m}$. To verify that this is a $p$-adic family, it remains to show that for {\it any} $k,k'\equiv k_0\pmod{p^{{\sf L}(\alpha,k_0}}$ the congruence $\lambda_k\equiv\lambda_{k'}\pmod{p^{{\sf a}m+{\sf b}-A_{k_0}}}$ holds, where now $m$ is given by $k\equiv k'\pmod{p^m}$. The transfer theorem yields an element $\lambda'\in\Lambda_{k',\chi,Np}^\alpha$ such that $\lambda'\equiv\lambda_k\pmod{p^{{\sf a}m+{\sf b}-A_{k}}}$. Since we may choose $A_k=A_{k_0}$ by Corollary C3 (note that $k\equiv k_0\pmod{p^{{\sf L}(\alpha,k_0)}}$) we obtain $\lambda'\equiv\lambda_k\pmod{p^{{\sf a}m+{\sf b}-A_{k_0}}}$. Since $m\ge {\sf L}(\alpha,k_0)$ this yields $\lambda'\equiv \lambda_k\pmod{p^{{\sf a}{\sf L}(\alpha,k_0)+{\sf b}-A_{k_0}}}$. Together with $\lambda_k\equiv \lambda\pmod{p^{{\sf a}{\sf L}(\alpha,k_0)+{\sf b}-A_{k_0}}}$ we obtain $\lambda'\equiv\lambda\pmod{p^{{\sf a}{\sf L}(\alpha,k_0)+{\sf b}-a_{k_0}}}$. Since also $\lambda_{k'}\equiv\lambda\pmod{p^{{\sf a}{\sf L}(\alpha,k_0)+{\sf b}-a_{k_0}}}$ we deduce from Corollary C2 applied to weights $k'$ and $k_0$ and $m={\sf L}(\alpha,k_0)$ that $\lambda'=\lambda_{k'}$. Clearly, this implies that $\lambda_k,\lambda_{k'}$ satisfy the requested congruence.

\newpage

\section {Group Cohomology and the slope decomposition}

{\bf 1.1. The slope decomposition. } 
%
%[[in diesem Abschnitt kommen keine Hecke Operatoren vor, d.h. $\Gamma_1(Np)$ spielt keine Rolle; $T$ ist irgendein Endomorphismus irgendeines ${\cal O}$-Moduls  ]]
%
We fix a prime $p\in{\Bbb N}$ and we denote by ${\Bbb C}_p$ the $p$-adic complex numbers with $p$-adic valuation $w$ normalized by $w(p)=1$.  
In addition we fix an isomorphism $i:\,{\Bbb C}\cong {\Bbb C}_p$; this immediately induces a "$p$-adic valuation" $v_p$ on ${\Bbb C}$ by defining $v_p(z)=w(i(z))$, $z\in {\Bbb C}$ and, hence, a compatible system of $p$-adic valuations $v_p=v_{p,E}$ on subfields $E$ of ${\Bbb C}$ via restriction. Here, compatible means that if $E/F$ is an extension of subfields of ${\Bbb C}_p$, the $p$-adic valuation $v_{p,E}$ obtained on $E$ restricts to the $p$-adic valuation $v_{p,F}$ obtained on $F$; in particular we may omit the index $E$ designating the field from $v_p$. 
%
%[[In fact, the choice of $i$ is equivalent to a choice of a (compatible) system of valuations on subfields $E\subset {\Bbb C}$.
%]]
% 
We shall use the valuation $v_p$ only as a convenient way to express congruences between elements which lie in arbitrary extension fields of ${\Bbb Q}$. To be more explicit, if $a,b$ are algebraic over ${\Bbb Q}$ and $c\in{\Bbb Q}$ we write 
$$
a\equiv b\pmod{p^c}
$$
to denote that $v_p(a-b)\ge c$. Let $E/{\Bbb Q}$ be a {\it finite} extension, which contains $a$ and $b$ and denote by ${\mathfrak p}={\mathfrak p}_E=\{x\in{\cal O}_E:\,v_p(x)>0\}$ the prime ideal in ${\cal O}_E$ corresponding to $v_p$ and by $e$ the ramification degree of ${\mathfrak p}|p$. The congruence $a\equiv b\pmod{p^c}$ then further is equivalent to 
$$
a\equiv b\pmod{{\mathfrak p}^{ec}}.
$$
%
%[[Auf $\bar{\Bbb Q}$ ist Vorgabe von $v_p$ aquivalent zur Vorgabe eines Systems von kompatiblen Bewertungen $v_{p,E}$ ist weiter aquivalent zur Vorgabe eines systems von kompatiblen Primidealen ${\mathfrak p}_E$ (d.h. $F\le E \Rightarrow {\mathfrak p}_E|{\mathfrak p}_F$): $E$ durchlauft alle endlöichen Erweiterungen von ${\Bbb Q}$.
%]]
%
The notation $a\equiv b\pmod{p^c}$ is weaker than the above congruence, because it leaves open in which field the congruence takes place. 

%[[Remark. 1.) If $E/{\Bbb Q}$ is algebraic, then $v_p$ is one extension of the valuation $v_p$ on ${\Bbb Q}$ to $E$, hence, to fix an isomorphism $i$ as above is equivalent to fixing valuations on all finite extensions $E/{\Bbb Q}$, which are compatible, i.e. if $F\subset E$, then the valuation on $E$ restricts to the one on $F$. If $E/{\Bbb Q}$ 
%is not algebraic then: let us assume $E={\Bbb Q}(T)$, $T$ an indeterminate; $v_p$ asigns a value $c$ to $T$ and the resulting valuation on ${\Bbb Q}(T)$ is given by
%$T^zP(T)\mapsto c^z$, $z\in{\Bbb Z}$, $(T,p(T))=1$.
%
%2.) All our finite extension $E/{\Bbb Q}$ that will show up come up as subfields of ${\Bbb C}$ (they are fields of definition of spaces of COMPLEX modular forms); thus we have to apply a map that sends subfields of ${\Bbb C}$ to subfields of ${\Bbb C}_p$ in order to apply $w$ $\leadsto$ $i:\,{\Bbb C}\cong {\Bbb C}_p$]]

Let $E/{\Bbb Q}$ be an arbitrary extension with ring of integers ${\cal O}_E$ 
and let $T$ be an operator on the finite dimensional $E$-vector space $V$. We define the slope $\alpha$ subspace $V^\alpha\le V$ as the image $p_\alpha(T)V$, where $p_\alpha$ is the factor of the characteristic polynomial $p$ of $T$, whose roots $\lambda$ (in a splitting field of $T$) have $p$-adic valuation $v_p(\lambda)$ different from $\alpha$. We note that $V^\alpha$ depends on the operator $T$ as well as on the choice of $v_p$. 
%
%[[Ist z.B. das Minimalpolynom von $T$ zerfallend in $E$, $E/{\Bbb Q}$ endlich, dann kann es sein, dass bezugl. einer Bewertung $v_1$ auf $E$ die Nullstellen von $\chi$ andere Werte haben als bezuglich einer anderen Bewertung $v_2$
%($v_1$ und $V_2$ setzen beide $V_p$ fort)]]
%
$V$ then decomposes as a direct sum of $E$-vector spaces
$$
V=\bigoplus_{\alpha\in{\Bbb Q}} V^\alpha\oplus V(0),\leqno(1)
$$
where $V(0)$ is the generalized eigenspace attached to $0$. 

%[[Note that $0\in{\Bbb Q}$, hence, the factor $T$ of the char. pol. of $T$ splits over ${\Bbb Q}$, i.e. the generalized eigenspace attached to $0$ is defined. This is not true if $\lambda$ is a root of char. Pol. of $T$ which is not contained in ${\Bbb Q}$, here the generalized eigenspace only appers over ${\Bbb Q}(\lambda)$]]

{\it Remark. } If we agree that $v_p(0)=\infty$ then $V(0)$ is the subspace of $V$ of slope $\infty$: $V(0)=V^\infty$. Hence, $V(0)$ is the subspace of $V$ of highest possible slope.

{\bf 1.2. Slope spaces over local fields. } We now assume that $E$ is a local field. More precisely, let $E/{\Bbb Q}_p$ be a finite extension, hence, the ring of integers ${\cal O}={\cal O}_E$ in $E$ is a principal ideal domain. We denote by ${\mathfrak p}\le{\cal O}$ the maximal ideal and we select a generator $\varpi$ of ${\mathfrak p}$, i.e. ${\mathfrak p}=(\varpi)$, hence, $v_p(\varpi)=1/e$, where $e$ is the ramification degree of ${\mathfrak p}|p$.  Furthermore, by $v_{\mathfrak p}$ we denote the valuation on $E$ normalized by $v_{\mathfrak p}(\varpi)=1$. We set $q=|{\cal O}/{\mathfrak p}|$.  Also, we assume that {\it $T$ is split over $E$}. This in particular implies that $V^\alpha=(0)$ unless $\alpha$ is contained in $\frac{1}{e}{\Bbb Z}$ 
%[[because the zeroes $\lambda$ of char pol. $p$ of $T$ are contained in $E$]]
and
$$
V^\alpha=\bigoplus_{\gamma\in E\atop v_p(\gamma)=\alpha}V_\gamma
$$
Here, $V_\gamma=\{v\in V:\,(T-\gamma)^kv=0\,\mbox{for some}\;k\in{\Bbb N}\}$ is the generalized eigenspace attached to $\gamma\in E$.

Working over a local field $E/{\Bbb Q}$ on the one hand has the advantage that ${\cal O}$ is a principal ideal domain and on the other hand enables us to construct topological idempotents attached to the slope decomposition (cf. section 3.4). 
As an example we want to look at the slope decomposition in case that $T$ stabilizes a lattice in $V$. To be more precise, let $V_{{\cal O}}\le V$ be a ${\cal O}$-lattice in $V$, i.e. $V_{\cal O}$ is a free 
${\cal O}$-submodule, which over $E$ generates $V$.
%[[note that ${\cal O}$ is a principal ideal domain, i.e. any torsion free module is free ($V_{\cal O}\subspace V\Rightarrow V_{\cal O}$ is torsion free)]]
For the remainder of this section we assume that $T$ leaves $V_{\cal O}$ stable, i.e. $T$ is an endomorphism of $V_{\cal O}$. This immediately implies that the eigenvalues of $T$ are integral over ${\cal O}$, hence, contained in ${\cal O}$
%[[Da die Nullstellen des Char polynom von $T$ ganz sind (fuhrender Koeffz $=1$, alle Koeffz in ${\cal O}$), ausserdem in $E$ liegen und ${\cal O}$ der ganzabschluss von ${\cal O}$ in $E$ ist]] 
%[[Ausfuhrlicher: Since $T$ stabilizes the lattice $V_{\cal O}^\alpha$ in $V^\alpha$ we deduce that the characteristic polynomial $\chi_T$ of $T$ has coefficients in ${\cal O}$
%(Beweis: Wahle Basis $(v_1,\ldots,v_d)$ des freien ${\cal O$-Moduls $V_{\cal O}^\alpha$; dann hat die Darstellungsmatrix von $T$ bezugl. $(v_1,\ldots,v_d)$ Koeffizienten in
%${\cal O$ (da $T(v_i)\in\langle v_1,\ldots,v_d\rangle=V_{\cal O}^\alpha$) $\Rightarrow$ alle Koeffizienten von $\chi_T$ liegen in ${\cal O}$.) with leading coefficient equal to $1$. Hence, the zero $\mu_i\in E$ of $\chi_T$ is integral over ${\cal O}$, which implies that $\mu_i\in {\cal O}$ because ${\cal O}$ is integrally closed]]. 
and we deduce that $V^\alpha\not=0$ only if $\alpha\in \frac{1}{e}{\Bbb N}\cup\{0\}$.
%
%[[Recall that we assume that $T$ is split over $E$ !]]

We put
$$
V_{\cal O}^\alpha=V_{\cal O}\cap V^\alpha.
$$ 
Obviously, $V^\alpha_{\cal O}$ is $T$-stable. Since $V_{\cal O}^\alpha\le V_{\cal O}$ and ${\cal O}$ is a principal ideal domain, $V_{\cal O}^\alpha$
is a finitely generated, free ${\cal O}$-module.

{\bf Remark. }{\it $V_{\cal O}^\alpha$ is a $T$-stable lattice in $V^\alpha$. In particular,
$$
{\rm rank}_{{\cal O}} V_{{\cal O}}^\alpha=\dim_{E} V^\alpha.
$$
}

{\it Proof. } Let $x\in V^\alpha$. Since $V_{{\cal O}}$ is a lattice in $V$ there is $n\in{\cal O}$ such that $nx\in V_{\cal O}$,
%[[Sei $\{m_i\}$ Erzeugendensystem fur $V_{{\cal O}}\Rightarrow $\{m_i\}$ ist auch EZS fur $V$ (d.h. uber $E$). Also existiert Darstellung $x=\sum_i a_i m_i, a_i\in E \Rightarrow nx=\sum_i (na_i) m_i\in V_{{\cal O}}$ wobei $n=$ produkt der nenner der $a_i$ (beachte: ${\cal O}$ ist faktoriell mit einem primelement da diskreter Bewertungsring.]] 
hence $nx\in V_{\cal O}\cap V^\alpha=V_{\cal O}^\alpha$. Thus, $V^\alpha/V_{\cal O}^\alpha$ is a torsion module,
which implies that $E\otimes_{\cal O}V_{\cal O}^\alpha=V^\alpha$. Since $V_{\cal O}^\alpha$ is a free ${\cal O}$-module the claim follows.

%[[Alternativ: Any basis for $V_{\cal O}$ remains linearly independent over $E$
%[[Allgemein: Ist $U\le V_E$ $\Rightarrow U\cap V_{\cal O}$ ist Gitter in $U$.]]

We call a submodule $U\le V_{\cal O}$ pure in $V_{\cal O}$ if $rv\in U$ for some $r\in {\cal O}$, $v\in V_{\cal O}$ implies that $v\in U$. This is equivalent to $U$ being complemented in $V_{\cal O}$ (note that $V_{\cal O}$ is free). If $U\le V_{\cal O}$ is any submodule we define 
the "purification of $U$" as
$$
U_{pure}=\{v\in V_{\cal O}:\,rv\in U\;\mbox{for some}\;r\in {\cal O}\}=(E\otimes_{\cal O}U)\cap V_{\cal O}.
$$
$U_{pure}$ is pure in $V_{\cal O}$ and has the same ${\cal O}$-rank as $U$.
%[[Da $U_{pure}/U$ ein Torsionsmodul ist $\Rightarrow$ $[U_{pure}:U]<\infty$ (beachte: $U,U_{pure}\le V_{\cal O}$ sind endlich erzeugt, da $V_{\cal O}$ endlich erzeugt)]] 
We call $v\in V_{\cal O}$ pure if $\langle v\rangle\le V_{\cal O}$ is pure. Obviously, $V^\alpha_{\cal O}\le V_{\cal O}$ is a pure submodule. 
%
%[[Because $V_{\cal O}^\alpha=V_{\cal O}\cap V^\alpha$ is the intersection of two pure modules.]]
%
In particular, we obtain
$$
V_{\cal O}^\alpha=(\bigoplus_{v_p(\gamma)=\alpha}V_{\cal O}(\gamma))_{pure},
$$ 
where $V_{\cal O}(\gamma)=\{v\in V_{\cal O}:\,(T-\gamma)^kv=0\,\mbox{for some}\;k\in{\Bbb N}\}=V_\gamma\cap V_{\cal O}$ ($\gamma\in {\cal O}$) 
%
%[[the sum is direct because the corresponding sum over $E$ is direct]]
%
(we note that $V_{\cal O}(\gamma)\le V_{\cal O}$ is a pure submodule). 
%[[Beweis: $V_{\cal O}^\alpha$ enthalt die rechte Seite und ist pure $\Rightarrow $ "$\supset$". Wegen ${\rm rk}_{\cal O} V_{\cal O}^\alpha={\rm rk}_{\cal O}$(rechte Seite folgt dann "$=$")]]
%[[Bemerkung: 1.) wir brauchen dieses Fakt nicht; $\sum_{v_{\mathfrak p}(\mu)=\alpha}V_{\cal O}(\mu)$ ist i. allg nicht pure !!
%2.) $V_{\cal O}^\alpha$ und $\sum_{v_{\mathfrak p}(\mu)=\alpha})V_{\cal O}(\mu)$ bestimmen sich gegenseitig: 
%$V_{\cal O}^\alpha=(\sum_{v_{\mathfrak p}(\mu)=\alpha})V_{\cal O}(\mu))_{pure}$ und $\sum_{v_{\mathfrak p}(\mu)=\alpha}V_{\cal O}(\mu)=\sum_{v_{\mathfrak p}(\mu)=\alpha} W(\mu)$, wobei $W(\mu)$ der $\mu$-Eigenraum von $T$ auf $V_{\cal O}^\alpha$ ist !
%]]
In particular, after reducing coefficients modulo ${\mathfrak p}^r$ we see that
$$
V^\alpha_{\cal O}\otimes {\cal O}/{\mathfrak p}^r\le {\cal O}/{\mathfrak p}^r\otimes V_{\cal O}
$$
is complemented. Moreover, $T$ induces an endomorphism of $V_{\cal O}^\alpha\otimes{\cal O}/{\mathfrak p}^r$ and we define the submodule
$$
V^\alpha_{\cal O}[0]=\{v\in V^\alpha_{\cal O}\otimes {\cal O}/{\mathfrak p}^r:\,T(v)=0\}\le {\cal O}/{\mathfrak p}^r\otimes V^\alpha_{\cal O}.
$$

\bigskip

We set $d={\rm rk}_{\cal O} V_{\cal O}^\alpha=\dim_E V^\alpha$. 
%[[note that $\dim_{E} V^\alpha={\rm rk}\, V_{\cal O}^\alpha=d$ by the above Remark]]
In view of the following Lemma we note that the cardinality of any ${\cal O}/{\mathfrak p}^r$-module is a power of $q$ and also that $e\alpha\in{\Bbb N}\cup\{0\}$.

\bigskip

{\bf Lemma. }{\it Let $r$ be any integer satisfying $r>e\alpha$. Then
$$
q^{d(r-e\alpha)}\big ||({\cal O}/{\mathfrak p}^r\otimes V^\alpha_{\cal O})/V^\alpha_{\cal O}[0]|.
$$
}

%[[Wir brauchen das Lemma nur fur den Fall $(p^r)$ anstelle von ${\mathfrak p}^r$, d.h. fur den Fall $e|r$.]]

%[[Das Lemma besagt, dass  $({\cal O}/{\mathfrak p}^r\otimes V^\alpha_{\cal O})/V^\alpha_{\cal O}[0]$ gross ist, d.h. dass $V^\alpha_{\cal O}[0]\le({\cal O}/{\mathfrak p}^r\otimes V^\alpha_{\cal O})$ ein kleiner Untermodul ist (idealerweise ist $V^\alpha_{\cal O}[0]=0$)]]

{\it Proof. } Since $T$ is split there is a basis $(v_1,\ldots,v_d)$ of $V^\alpha$  
such that the representing matrix of $T$ is upper triangular, i.e. $T(v_i)=\mu_iv_i+m$ for some $m\in \langle v_1,\ldots,v_{i-1}\rangle_{E}$ for all $1\le i\le d$. The characteristic polynomial of $T$ on $V^\alpha$ then reads $\chi_{T|_{V^\alpha}}=\prod_i (T-\mu_i)$; since $\chi_{T|_{V^\alpha}}=\prod_{v_p(\gamma)=\alpha\atop V_\gamma\not=0}(T-\gamma)^{n_\gamma}$ we deduce that $\mu_i=\gamma$ for some eigenvalue $\gamma$ of $T$, hence, 
$$
v_p(\mu_i)=\alpha\leqno(2)
$$ 
for all $i$. 
%[[anderes Argument: $\mu_i$ is EW von $T$ auf $V^\alpha$. Ware $\mu_i=\lambda$ mit einem EW $\lambda$, der $v_p(\lambda)\not=\alpha$ erfullt, dann folgte, dass die Summe
%$\sum_{v_p(\gamma)=\alpha}V_\gamma +V_\lambda$ nicht direkt ist Widerspruch !]]
%
Multiplying the basis vectors $v_i$ by a suitable scalar we may assume that $v_1,\ldots,v_d\in V_{\cal O}^\alpha$. 
%[[The relation $T(v_i)=\mu_i+m$, $m\in \langle v_1,\ldots,v_d\rangle_{E}$ then still holds]]
We set
$$
V_{\cal O}^i=\langle v_1,\ldots,v_i\rangle_{E}\cap V_{\cal O}^\alpha=\langle v_1,\ldots,v_i\rangle_{{\cal O},pure}.
$$
Obviously, $V_{\cal O}^i$ is $T$-stable. Taking into account that $\mu_i\in{\cal O}$ and that $T$ stabilizes $V_{\cal O}^\alpha$ we obtain 
$$
T(v_i)=\mu_i v_i+m, \leqno(3)
$$
where now $m\in V_{\cal O}^{i-1}$. 

%%%%%%%%%%%%%%%%%%%%%%%%%%%%%%%%%%%%%%%%%%%%%%%%%%%%%%%%%%%%%%%%%%%%%%%%%%%%%%%%%
We inductively construct a sequence of elements $v_1^0,\ldots,v_d^0\in V^\alpha_{\cal O}$ as follows. We define $v_1^0$ by demanding that
$\langle v_1^0\rangle=V_{\cal O}^1=\langle v_1\rangle_{pure}$. Assuming that $v_1^0,\ldots,v_i^0$ have been chosen we select any $v_{i+1}^0\in V^{i+1}_{\cal O}$ such that
$$
V^{i+1}_{\cal O}=V^i_{\cal O}\oplus {\cal O}v_{i+1}^0.\leqno(4)
$$
Notice here that $V_{\cal O}^i$ is pure in $V_{\cal O}^\alpha$ and thus in $V^{i+1}_{\cal O}$, hence, $V_{\cal O}^i$ is complemented in $V^{i+1}_{\cal O}$.
%[[das Komplement von $V^i_{\cal O}$ in $V_{\cal O}^{i+1}$ ist eindimnsional, denn 
%$\dim_{E} V^{i+1}_{\cal O}\otimes E=\dim_{E} V^{i+1}_{\cal O}\otimes E+1$]]
Obviously, $\{v_i^0\}$ is a basis of $V_{\cal O}^\alpha$.
%[[wegen (S)]]

%%%%%%%%%%%%%%%%%%%%%%%%%%%%%%%%%%%%%%%%%%%%%%%%%%%%%%%%%%%%%%%%%%%%%%%%%%%%%%%%%
By using the definition of $V_{\cal O}^{i+1}$ we may write
$$
v_{i+1}^0=\frac{1}{\varpi^n}(w+\epsilon v_{i+1})\leqno(5)
$$
for some $w\in V^i_{\cal O}$, $n\in{\Bbb N}\cup\{0\}$ and $\epsilon \in {\cal O}$.

Using (5) and (3) we find
\begin{eqnarray*}
T(v_{i+1}^0)&=&\frac{1}{\varpi^n}(\epsilon \mu_{i+1}v_{i+1}+\mu_{i+1}w-\mu_{i+1}w+T(w)+\epsilon m)\\
&=&\mu_{i+1} v_{i+1}^0+\frac{1}{\varpi^n}(T(w)-\mu_{i+1}w+\epsilon m)\\
\end{eqnarray*}

Since $v_{i+1}^0$ is contained in $V_{\cal O}^\alpha$ and $T$ leaves $V_{\cal O}^\alpha$ invariant
the last equation implies that $\frac{1}{\varpi^n}(T(w)-\mu_{i+1}w+\epsilon m)\in V_{\cal O}$, hence, $\frac{1}{\varpi^n}(T(w)-\mu_{i+1}w+\epsilon m)\in V_{\cal O}^i$. Thus, with respect to the basis $\{v_{i}^0\}$ $T$ has representaing matrix
$$
A=\left(\begin{array}{cccc}
\mu_1&*&\cdots&*\\
&\mu_2&\ddots&\vdots\\
&&\ddots&*\\
&&&\mu_d\\
\end{array}
\right)
$$ 
where all the $*$'s are contained in ${\cal O}$. The theorem on elementary divisors then shows that there are basis $\{c_i\}$ and $\{d_i\}$ of $V_{\cal O}^\alpha$ with respect to which $T$ has representing matrix in diagonal form
$$
B=\left(\begin{array}{cccc}
\lambda_1&&&\\
&\lambda_2&&\\
&&\ddots&\\
&&&\lambda_d\\
\end{array}
\right),\quad\lambda_i\in{\cal O}.\leqno(6)
$$ 
Moreover, $B$ is similar to $A$, hence, they share the same determinant
$$
\prod_{i=1}^d \lambda_i=\prod_{i=1}^d \mu_i.\leqno(7)
$$
Let $v=\sum_i \beta_i c_i\in V_{\cal O}^\alpha$ be arbitrary. (6) implies that 
$T(v)\in {\mathfrak p}^rV_{\cal O}^\alpha$ precisely if
$\beta_i\in {\mathfrak p}^{r-v_{\mathfrak p}(\lambda_i)}$ in case $v_{\mathfrak p}(\lambda_i)\le r$ and $\beta_i\in {\cal O}$ is arbitrary in case $v_{\mathfrak p}(\lambda_i)> r$. Since $T(v)=0$ in ${\cal O}/{\mathfrak p}^r\otimes V_{\cal O}^\alpha$ precisely if $T(v)\in{\mathfrak p}^r$ we obtain
$$
|V^\alpha_{\cal O}[0]|
=\prod_{i=1}^d q^{{\rm min}\{r,v_{\mathfrak p}(\lambda_i)\}}
=q^{\sum_{i=1}^d {\rm min}\{r,v_{\mathfrak p}(\lambda_i)\}}.
$$
Consequently, $|V^\alpha_{\cal O}[0]|$ divides $q^{\sum_i v_{\mathfrak p}(\lambda_i)}$. Using (7) and taking into account that $v_{\mathfrak p}(\mu_i)=e\alpha$ (cf. (2)) we thus obtain
$$
|V^\alpha_{\cal O}[0]|\, \big|\, q^{de\alpha}.
$$
Since $|{\cal O}/{\mathfrak p}^r\otimes V_{\cal O}^\alpha|=q^{dr}$ this proves the Lemma.

{\bf 1.3. Spaces of modular forms. } %[[Hier benutzen wir die ganze Hecke-Algebra $GL_2({\AA}_f)//K$, nicht die nur die von den Hecke Operatoren $T_p, T_{p,p}$ erzeugte Unteralgebra. Nur ganz zum Schluss taucht der Atkin Lehner Op auf]]
We denote by $L_{k,{\Bbb Z}}\le {\Bbb Z}[X,Y]$ the ${\Bbb Z}$-submodule consisting of all homogeneous polynomials in variables $X,Y$ of degree $k$. $L_{k,{\Bbb Z}}$ becomes a ${\rm GL}_2({\Bbb Z})$-module under the action 
$$
\gamma P(X,Y)=P({^t}(\gamma^\iota {X\choose Y})),\quad\gamma\in{\rm GL}_2({\Bbb Z}),
$$
where $\gamma^\iota=(\det \gamma) \gamma^{-1}$. 
%[[cf. [Hi], p. 165[]
For any ${\Bbb Z}$-algebra $R$ we put $L_{k,R}=L_{k,{\Bbb Z}}\otimes R$. $L_{k,{\Bbb Q}}$ is isomorphic to the unique irreducible finite dimensional representation of ${\rm GL}_2({\Bbb Q})$ of 
highest weight $(k,0)$ with respect to the torus $T_2$ of diagonal matrices in ${\rm GL}_2$.

We denote by $\Gamma_0(N)$ resp. $\Gamma_1(N)$, $N\in{\Bbb N}$, the set of all matrices  
$$
\gamma=\Mat{a}{b}{c}{d}\in{\rm SL}_2({\Bbb Z})
$$ 
satisfying $c\equiv 0\pmod{N}$ resp. $c\equiv 0, d\equiv 1\pmod{N}$. $\Gamma_0(N), \Gamma_1(N)$ are arithmetic subgroups of ${\rm SL}_2({\Bbb Z})$ and for any arithmetic subgroup $\Gamma\le{\rm SL}_2({\Bbb Z})$ we set
$$
S(\Gamma)=\Gamma\backslash {\rm SL}_2({\Bbb R})/{\rm SO}_2({\Bbb R}).
$$
We further denote by ${\cal M}_k(\Gamma)$ the space of holomorphic modular forms of level $\Gamma$ and weight $k$ and by ${\cal S}_k(\Gamma)$ the subspace of cuspidal forms. Moreover, by ${\cal M}_k(\Gamma_0(N),\chi)$ resp. ${\cal S}_k(\Gamma_0(N),\chi)$ we denote the subspaces of forms with nebentype $\chi$.

Similarly, we denote by $K_0(N)$ resp. $K_1(N)\le \prod_\ell {\rm GL}_2({\Bbb Z}_\ell)$ the subgroups consisting of all elements 
$$
k=\left(\Mat{a_\ell}{b_\ell}{c_\ell}{d_\ell}\right)_\ell
$$
satisfying $c_\ell\equiv 0\pmod{N}$ resp. $c_\ell\equiv 0, d\equiv 1\pmod{N}$. 
$K_0(N), K_1(N)\le \prod_\ell {\rm GL}_2({\Bbb Z}_\ell)$ are compact open subgroups satisfying $\det K_0(N)=\det K_1(N)=\prod_\ell {\Bbb Z}_\ell^*$. For any compact open subgroup $K\le \prod_\ell {\rm GL}_2({\Bbb Z}_\ell)$ we set
$$
S(K)={\rm GL}_2({\Bbb Q})\backslash {\rm GL}_2({\Bbb A}_f)/{\rm SO}_2({\Bbb R}){\Bbb R}^*_+.
$$
%[[note that ${\rm SO}_2({\Bbb R}){\Bbb R}^*_+={\rm SO}_2({\Bbb R}){\Bbb R}^*$; ${\Bbb R}^*$ ist das Zentrum in ${\rm SO}_2({\Bbb R})$]]
Assuming that $\det K=\hat{\Bbb Z}^*$ we obtain as analytic manifolds 
$$
S(K)\cong S(\Gamma),
$$
where $\Gamma=\Gamma_K={\rm GL}_2({\Bbb Q})\cap K\times {\rm GL}_2^+({\Bbb R})$ (in particular, $S(K)$ is connected). 
%
%[[man muss $g\in {\rm GL}_2({\Bbb A})$ mittels $\gamma\in{\rm GL}_2({\Bbb Q})$ nach $K\times{\rm GL}_2^+({\Bbb R})$ bringen; dann hat man noch die $\gamma\in {\rm GL}_2^+({\Bbb R})\cal K$ frei (solche Elemente andern $\det g_\infty$ nicht, d.h. sie lassen ${\rm GL}_2^+({\Bbb R})/{\rm SO}_2({\Bbb R}){\Bbb R}_+^*$ invarinat) mit denen man auf $g_\infty$ wiken kann.]]

For the remainder of section 1.3 we set $\Gamma=\Gamma_1(N)$ and $K=K_1(N)$, where $N\in{\Bbb N}$ is arbitrary. Following [H], Theorem 2, p. 77, the space of modular forms of weight $k$ with respect to this congruence subgroup decomposes as module under the Hecke algebra ${\cal H}={\rm GL}_2({\Bbb A}_f)//K$ as follows. We define the unitarily induced representation
$$
I_k={\rm Ind}_{B_2({\Bbb R})}^{{\rm GL}_2({\Bbb R})} ({\rm sgn}^k\otimes|\cdot|_\infty^{1/2},|\cdot|_\infty^{3/2-k}),\quad k\ge 2.
$$
$I_k$ is reducible and contains a unique proper subrepresentation $D_k$, which has lowest ${\rm SO}_2({\Bbb R})$-type $k$ and central character $\Mat{a}{}{}{a}\mapsto {\rm sgn}^k(a)a^{2-k}$, $a\in{\Bbb R}^*$.
%[[d.h. genauer $I_k$ hat lowest weight=lowest ${\rm SO}_2({\Bbb R})$-type $\chi_k:\,e^{i\Theta}\mapsto e^{ik\Theta}$]] 
%[[denn: $\chi_1\chi_2^{-1}={\rm sgn}^k |t|_\infty^{1-k}={\rm sgn}(t)t^{k-1}$ $\Rightarrow$ reducible and lowest ${\rm SO}_2({\Bbb R})$-Type $k$ nach [Gelbart, Automorphic forms on adele groups, Theorem 4.4]]]

The space of holomorphic modular forms with respect to $\Gamma=\Gamma_1(N)$ then decomposes
$$
{\cal M}_{k+2}(\Gamma)\cong H^1(S(K),L_{k,{\Bbb C}})=H^1_{\rm Eis}(S(K),L_{k,{\Bbb C}})\oplus H^1_{\rm cusp}(S(K),L_{k,{\Bbb C}}).
$$
Here,
$$
H^1_{\rm cusp}(S(K),L_{k-2,{\Bbb C}})=\bigoplus_{\pi\in {\cal A}_0(K,k)} \pi_f^K\leqno(8)
$$
with ${\cal A}_0(K,k)$ denoting the space of all cuspidal automorphic representation with infinity component isomorphic to $D_k$ and $\pi_f^K\not=0$,
$$
H^1_{\rm Eis}(S(K),L_{k-2,{\Bbb C}})=\bigoplus_{\chi\in {\cal T}(k)} ({\rm Ind}\chi_f)^K\leqno(9)
$$
if $k>2$ and 
$$
H^1_{\rm Eis}(S(K),L_{k-2,{\Bbb C}})=\bigoplus_{\chi\in {\cal T}(k)\atop \chi_1\chi_2^{-1}\not=|\cdot|_{\Bbb A}^2} ({\rm Ind}\chi_f)^K\oplus \bigoplus_{\chi\in {\cal T}(k)\atop\chi_1\chi_2^{-1}=|\cdot|_{\Bbb A}^2}(\tilde{\rm Ind}\chi_f)^K\bigoplus\chi_1|\cdot|_{\Bbb A}^{-1/2}\circ\det\leqno(10)
$$
if $k=2$, where ${\cal T}(k)$ is the set of all characters $\chi=(\chi_1,\chi_2):\,T_2({\Bbb A})\rightarrow {\Bbb C}^*$ satisfying $\chi_{1,\infty}|_{{\Bbb R}^*_+}=x^{1/2}$, $\chi_{2,\infty}|_{{\Bbb R}^*_+}=x^{3/2-k}$, $\chi_{1,\infty}\chi_{2,\infty}(-1)=(-1)^k$ and $({\rm Ind}\chi_f)^K\not=0$. 
%
%[[beachte: $\chi_\infty|_{{\Bbb R}^*_+}=x^a \Leftrightarrow \chi_\infty|_{{\Bbb R}^*_+}=|\cdot|_\infty^a$ und $\cih_1\chi_2(-1)=(-1)^k\Leftrightarrow \chi_{1,\infty}\chi_{2,\infty}(-1)=(-1)^k \Leftrightarrow \chi_{1,\infty}=\chi_{2,\infty}{\rm sgn}^k$ (hier bezeichnen $\chi_1,\chi_2$ einen Dirichletcharakter und gleichzeitig auch die zu $\chi_1,\chi_2$ gehorigen Ideleklassencharaktere)]]
%
Moreover, $\tilde{\rm Ind}\chi_f$ is the subspace of ${\rm Ind}\chi_f=\otimes_{v\not=\infty}{\rm Ind}\chi_v$ generated by those tensors $\otimes_v \varphi_v$, for which at least one finite component is contained in the unique proper submodule of the reducible representation $\tilde{\rm Ind}\chi_v$.

We deduce that any irreducible subrepresentation $\rho$ appearing in $H^1(S(K),L_{k,{\Bbb C}})$ either coincides with the isotypical
component $H^1_{\rm cusp}(S(K),L_{k,{\Bbb C}})(\pi_f)$ for some cuspidal $\pi$ (notice multiplicity-1) or is contained in the isotypical component $H^1_{\rm Eis}(S(K),L_{k,{\Bbb C}})({\rm Ind}\chi_f)$ for some $\chi\in{\cal T}(k)$. Hence, $\rho$ is isomorphic to $\pi_f^K$ or to $({\rm Ind}\chi_f)^K$. We denote by ${\Bbb Q}(\pi)$ the field of definition of the representation $\pi_f$, $\pi\in {\cal A}_0(K,k)$. The corresponding representation $\pi_f^K$ of the Hecke algebra ${\cal H}$ also is defined over ${\Bbb Q}(\pi)$. On the other hand, ${\rm Gal}(\bar{\Bbb Q}/{\Bbb Q})$ acts on the set
${\cal X}(k)$ of all characters $\chi:\,{\Bbb Q}\backslash {\Bbb A}\rightarrow {\Bbb C}^*$ satisfying $\chi|_{{\Bbb R}^*_+}=|\cdot|_\infty^k$ by $({^\sigma\chi})(x)=\sigma(\chi(x))$. To see that this action is well defined note that
${^\sigma\chi}={^\sigma}(\chi|\cdot|^{-k})|\cdot|^k$. As a consequence, ${\rm GL}_2({\Bbb Q})$ acts on ${\cal T}(k)$ by
${^\sigma}(\chi_1,\chi_2)=({^\sigma}\chi_1,{^\sigma\chi_2})$ and we deduce that ${\rm Ind}\chi_f$ as well as the representation $({\rm Ind}\chi_f)^K$ of the Hecke algebra ${\cal H}$ are defined over ${\Bbb Q}(\chi)={\Bbb Q}(\chi_1,\chi_2)$. 
We denote by $E=E_{\Gamma,k}$ the composite of the fields ${\Bbb Q}(\pi)$, $\pi\in{\cal A}_0(K,k)$ and the fields ${\Bbb Q}(\chi_1,\chi_2)$, $\chi\in{\cal T}(K,k)$. In particular, ${\cal M}_k(\Gamma)$ and the decomposition (8) and (8) are defined over $E$. Note that the decomposition (8) is finite, hence, $E$ is a finite extension of ${\Bbb Q}$.

Using group cohomology we define an integrality structure on $H^1(S(\Gamma),L_{k,E})$. To this end we denote by ${\cal O}$ the ring of integers of $E$. We fix a prime $p\in{\Bbb N}$ and we select a prime ideal ${\mathfrak p}\le {\cal O}$ lying above $p$. $E_{\mathfrak p}$ resp. ${\cal O}_{\mathfrak p}$ is the completion of $E$ resp. ${\cal O}$ with respect to ${\frak p}$. The inclusion
$$
i:\,L_{k,{\cal O}_{\frak p}}\hookrightarrow L_{k,{E}_{\mathfrak p}}
$$
induces a map
$$
i^*:\,H^1(\Gamma,L_{k,{\cal O}_{\frak p}})\rightarrow H^1(\Gamma,L_{k,E_{\mathfrak p}}).
$$
We define
$$
H^1(\Gamma,L_{k,{\cal O}_{\frak p}})_{\rm int}=i^*(H^1(\Gamma,L_{k,{\cal O}_{\frak p}}).
$$

{\bf Lemma. }{\it $H^1(\Gamma,L_{k,{\cal O}_{\frak p}})_{\rm int}$ is an ${\cal O}_{\mathfrak p}$-lattice in $H^1(\Gamma,L_{k,E_{\mathfrak p}})$. 
%
%[[i.e. $H^1(S(K),L_{k,{\cal O}_{\frak p}})_{\rm int}$ is a free ${\cal O}_{\mathfrak p}$-{\bf sub}module that generates $H^1(S(K),L_{k,E_{\frak p}})$ over $E$.]]
% 
Moreover, $i:\,L_{k,{\cal O}_{\mathfrak p}}\rightarrow L_{k,E_{\mathfrak p}}$ induces an isomorphism
$$
H^1(\Gamma,L_{k,{\cal O}_{\frak p}})
/H^1(\Gamma,L_{k,{\cal O}_{\frak p}})_{\rm tor}
\stackrel{i^*}{\cong}
H^1(\Gamma,L_{k,{\cal O}_{\frak p}})_{\rm int}, 
$$
where $V_{\rm tor}$ denotes the torsion submodule of the module $V$. In particular,
$$
\dim_{E_{\mathfrak p}} H^1(\Gamma,L_{k,E_{\mathfrak p}})={\rm rank}_{{\cal O}_{\mathfrak p}} H^1(\Gamma,L_{k,{\cal O}_{\mathfrak p}})_{\rm int}.
$$
}

{\it Proof.} The short exact sequence 
$$
0\rightarrow L_{k,{\cal O}_{\mathfrak p}}\stackrel{i}{\rightarrow} L_{k,E_{\mathfrak p}}\stackrel{\pi}{\rightarrow} L_{k,E_{\mathfrak p}}/L_{k,{\cal O}_{\mathfrak p}}\rightarrow 0
$$
yields an exact sequence
$$
H^1(\Gamma,L_{k,{\cal O}_{\mathfrak p}})\stackrel{i^*}{\rightarrow} H^1(\Gamma,L_{k,E_{\mathfrak p}})\stackrel{\pi^*}{\rightarrow}
H^1(\Gamma,L_{k,E_{\mathfrak p}}/L_{k,{\cal O}_{\mathfrak p}})\stackrel{\delta}{\rightarrow} H^2(\Gamma,L_{k,{\cal O}_{\mathfrak p}}).
$$
Taking into account that $H^2(\Gamma,L_{k,{\cal O}_{\mathfrak p}})=0$ 
%[[da $H^2(\Gamma,M)=0$ fur alle $\Gamma$-Moduln $M$; s. [Hida], Proposition 6.1, p. 162]] 
we deduce that
$$
H^1(\Gamma,L_{k,E_{\mathfrak p}})/H^1(\Gamma,L_{k,{\cal O}_{\mathfrak p}})_{\rm int}\stackrel{\pi^*}{\rightarrow}  H^1(\Gamma,L_{k,E_{\mathfrak p}}/L_{k,{\cal O}_{\mathfrak p}})
$$
is an isomorphism. Since $L_{k,E_{\mathfrak p}}/L_{k,{\cal O}_{\mathfrak p}}\cong (E_{\mathfrak p}/{\cal O}_{\mathfrak p})^{k+1}$ is a torsion ${\cal O}_{\mathfrak p}$-module 
%[[Multipliziere $(x_i)\in E_{\mathfrak p}^{k+1}$ mit dem Hauptnenner aller $x_i$ durch, dann kommmt es in ${\cal O}_{\mathfrak p}^{k+1}$ zu liegen]]
we deduce that $H^1(\Gamma,L_{k,E_{\mathfrak p}}/L_{k,{\cal O}_{\mathfrak p}})$ and, hence, 
$$
H^1(\Gamma,L_{k,E_{\mathfrak p}})/H^1(\Gamma,L_{k,{\cal O}_{\mathfrak p}})_{\rm int}\leqno(11)
$$ 
is a torsion ${\cal O}_{\mathfrak p}$-module. On the other hand,
$H^1(\Gamma,L_{k,{\cal O}_{\mathfrak p}})_{\rm int}$ is a torsion free, finitely generated ${\cal O}_{\mathfrak p}$-module and since ${\cal O}_{\mathfrak p}$ is a principal ideal domain we deduce that
$H^1(\Gamma,L_{k,{\cal O}_{\mathfrak p}})_{\rm int}$ is a free module. Since (11) is torsion we deduce that $H^1(\Gamma,L_{k,{\cal O}_{\frak p}})_{\rm int}$ is a lattice in $H^1(S(K),L_{k,E_{\frak p}})$. 

%[[Details: Schreibe 
%$$
%H^1(\Gamma,L_{k,{\cal O}_{\mathfrak p}})_{\rm int}=\bigoplus {\cal O}_{\mathfrak p}m_i 
%$$
%Since (4) is torsion we deduce that
%$H^1(\Gamma,L_{k,E})=\sum_i Em_i$
%( Proof: $x\in H^1(E)\Rightarrow nx\in H^1({\cal O}_{\mathfrak p})_{\rm int} \Rightarrow nx=\sum \alpha_i m_i\Rightarrow x=\sum_i (\alpha_i/n) m_i $), hence, $H^1(\Gamma,L_{k,{\cal O}_{\mathfrak p}})_{\rm int}$ is an 
%${\cal O}_{\mathfrak p}$-lattice in $H^1(\Gamma,L_{k,E})$. Note that since $\{m_i\}$ is ${\cal O}_{\mathfrak p}$-free, it is also $E$-free. Thus, $H^1(\Gamma,L_{k,E})=\bigoplus_i Em_i$. In particular, we see that
%$$
%\dim_{E} H^1(\Gamma,L_{k,E})={\rm rank}_{{\cal O}_{\mathfrak p}} H^1(\Gamma,L_{k,{\cal O}_{\mathfrak p}})_{\rm int}.
%$$
%]] 

To prove the second statement, we note that the short exact sequence above yields another piece of the long exact cohomology sequence
$$
H^0(\Gamma,L_{k,E_{\mathfrak p}}/L_{k,{\cal O}_{\mathfrak p}})\stackrel{\delta}{\rightarrow} H^1(\Gamma,L_{k,{\cal O}_{\mathfrak p}})\stackrel{i^*}{\rightarrow} H^1(\Gamma,L_{k,E_{\mathfrak p}}).
$$
$H^0(\Gamma,L_{k,E_{\mathfrak p}}/L_{k,{\cal O}_{\mathfrak p}})$ is a torsion module, which implies that ${\rm ker}\,i^*={\rm im}\,\delta$ is torsion, hence, ${\rm ker}\,i^*=H^1(\Gamma,L_{k,{\cal O}_{\mathfrak p}})_{\rm tor}$, because $H^1(\Gamma,L_{k,E_{\mathfrak p}})$ is torsion free.
%[[Detail: wir haben ${\rm ker}\,i^*\subset H^1(\Gamma,L_{k,{\cal O}_{\mathfrak p}})_{\rm tor}$; andererseits ist $=H^1(\Gamma,L_{k,{\cal O}_{\mathfrak p}})_{\rm tor}$ naturlich in ${\rm ker}\,i^*$ enthalten, da $H^1(\Gamma,L_{k,E})$ torsionsfrei ist]]
Thus, we obtain an isomorphism
$$
H^1(S(K),L_{k,{\cal O}_{\frak p}})/
H^1(S(K),L_{k,{\cal O}_{\frak p}})_{\rm tor}
\stackrel{i^*}{\rightarrow} 
H^1(S(K),L_{k,{\cal O}_{\frak p}})_{\rm int},\leqno(12)
$$ 
which proves the Lemma.

%[[note that $H^1(S(K),L_{k,{\cal O}_{\frak p}})_{\rm int}={\rm im}\,i^*$]]

%[[
%{\it Remark. } The above proof works for any pair ${\cal O}\subset E$ instead of ${\cal O}_{\mathfrak p}\subset E_{\mathfrak p}$ as long as ${\cal O}$ is a principal ideal domain. For example we may choose $E={\Bbb Q}$ and ${\cal O}={\Bbb Z}$. For an arbitrary field $E$ the belonging ring of integrers ${\cal O}_E$ is not a principal ideal domain, i.e. a torsion free ${\cal O}_E$-module is not necessarily free. Note that we cannot work over ${\Bbb Q}$ as soon as we consider Hecke operators, because as Hecke module the cohomology is not defined over ${\Bbb Q}$. (This probably can be saved be generalizing the definition of lattice over non PID (projective, finitely generated and ${\rm Quot}{\cal O}=E$)
%
%For later purpose we note that ${\rm dim}_{\Bbb C} {\cal M}_k(K)={\rm dim}_{E_{\mathfrak p}} H^1(S(K),L_{k,E_{\mathfrak p}})$
%
%{\it Proof. } We calculate
%
%\begin{eqnarray*}
%{\rm dim}_{\Bbb C} {\cal M}_k(\Gamma)&=&\dim_{\Bbb C} H^1(\Gamma,L_{k,{\Bbb C}})\\
%&\stackrel{(2)}{=}&\dim_{\Bbb C} H^1(\Gamma,L_{k,E})\otimes{\Bbb C}\\
%&=&\dim_E H^1(\Gamma,L_{k,E})\\
%&=&\dim_{E_{\mathfrak p}} H^1(\Gamma,L_{k,E})\otimes E_{\mathfrak p}\\
%&\stackrel{(3)}{=}&\dim_{E_{\mathfrak p}} H^1(\Gamma,L_{k,E_{\mathfrak p}})\\
%\end{eqnarray*}
%Equations (2) resp. (3) hold because ${\Bbb C}$ resp. $E_{\mathfrak p}$ are flat $E$-modules.

%]]

{\bf 1.4. Mod $p^r$ reduction of irreducible ${\rm GL}_2$-modules. } 
%
%[[Auch in diesem Abschnitt kommen keine Hecke Operatoen vor, d.h. $\Gamma_1(Np, Np^r, p^r)$ spielt keine Rolle. Es kommt auch der Slope $\alpha$ nicht vor und der globale Korper $E=$ Definitionskorper aller automorphen darstellungen zu $\Gamma$ kommt auch nicht vor]]
%
In this section $E/{\Bbb Q}$ is a finite extension with ring of integers ${\cal O}$ and ${\mathfrak p}\le {\cal O}$ is a prime ideal lying above $p$. We denote by $E_{\mathfrak p}$ and ${\cal O}_{\mathfrak p}$ the completions with respect to ${\mathfrak p}$. 

Recall that $L_{k,{\cal O}_{\mathfrak p}}$ is ${\cal O}_{\mathfrak p}$-free with basis $X^{k-i}Y^i$, $i=0,\ldots,k$. For any $r\in{\Bbb N}$
%[[i.e. $r\ge 1$]] 
we define a submodule
$$
U_{k,{\cal O}_{\mathfrak p},r}=\langle  X^{k-r-1}Y^{r+1},\ldots,Y^k\rangle\le L_{k,{\cal O}_{\mathfrak p}}.
$$

%[[Spater wird $r=[\alpha]+1$ sein]]

We note that ${\cal O}_{\mathfrak p}/(p^r)\otimes L_{k,{\cal O}_{\mathfrak p}}\cong L_{k,{\cal O}_{\mathfrak p}}/p^rL_{k,{\cal O}_{\mathfrak p}}$ and ${\rm SL}_2({\Bbb Z})$ acts on ${\cal O}_{\mathfrak p}/(p^r)\otimes L_{k,{\cal O}_{\mathfrak p}}$ via $\gamma (P(X,Y)+p^rL_{k,{\cal O}_{\mathfrak p}})=(\gamma P(X,Y))+p^rL_{k,{\cal O}_{\mathfrak p}}$.

\bigskip

{\bf Lemma. } {\it 1.) The submodule ${\cal O}_{\mathfrak p}/(p^r)\otimes U_{k,{\cal O}_{\mathfrak p},r}\le {\cal O}_{\mathfrak p}/(p^r)\otimes L_{k,{\cal O}_{\mathfrak p}}$ is $\Gamma_1(p^r)$-stable.

2.) Any $\gamma^\iota$, $\gamma\in \Gamma_1(p^r)\Mat{1}{}{}{p}\Gamma_1(p^r)$ annihilates ${\cal O}_{\mathfrak p}/(p^r)\otimes U_{k,{\cal O}_{\mathfrak p},r}$.
}
%[[Alter Ansatz: betrachte ${\Bbb Z}/p^r{\Bbb Z}\otimes L_{k,{\cal O}_{\mathfrak p}}$. Dabei betrachten wir den
%${\cal O}_{\mathfrak p}$-modul $L_{k,{\cal O}_{\mathfrak p}}$ als ${\Bbb Z}$_modul, d.h. als abelsche Gruppe

{\it Proof. } 1.) Let $\gamma=\Mat{a}{b}{c}{d}\in\Gamma_1(p^r)$ and let $X^{k-i}Y^i\in U_{k,{\cal O}_{\mathfrak p},r}$, i.e. $i>r$. Then, $\gamma^\iota=\Mat{d}{-b}{-c}{a}\in \Gamma_1(p^r)$ and we find
\begin{eqnarray*}
\gamma X^{k-i}Y^i&=&(dX-bY)^{k-i}(-cX+aY)^i\\
&\equiv&\underbrace{(dX-bY)^{k-i}(Y)^i}_{\in U_{k,{\cal O}_{\mathfrak p},r}}\pmod{p^rL_{k,{\cal O}_{\mathfrak p}}}.\\
\end{eqnarray*}

Since $i>r$, the term $(dX-bY)^{k-i}(Y)^i$ only contains powers of $Y$ with exponent strictly larger than $r$, hence, it is contained in $U_{k,{\cal O}_{\mathfrak p},r}$. Thus, the coset $\gamma X^{k-i} Y^i+p^r L_{k,{\cal O}_{\mathfrak p}}$ again is contained in ${\cal O}_{\mathfrak p}/(p^r)\otimes U_{k,{\cal O}_{\mathfrak p},r}(\le {\cal O}_{\mathfrak p}/(p^r)\otimes L_{k,{\cal O}_{\mathfrak p}})$.

2.) Since $\Gamma_1(p^r)^\iota=\Gamma_1(p^r)$ and ${\cal O}_{\mathfrak p}/(p^r)\otimes U_{k,{\cal O}_{\mathfrak p},r}$ is $\Gamma_1(p^r)$-stable, it suffices to prove that $\Mat{1}{}{}{p}^\iota$ annihilates ${\cal O}_{\mathfrak p}/(p^r)\otimes U_{k,{\cal O}_{\mathfrak p},r}$. But this is obvious, since for any $i>r$ we find
$$
\Mat{1}{}{}{p}^\iota X^{k-i}Y^i=p^iX^{k-i}Y^i\in p^r U_{k,{\cal O}_{\mathfrak p},r}.
$$
This completes the proof of the Lemma.

We note that ${\cal O}_{\mathfrak p}/(p^r)\otimes L_{k,{\cal O}_{\mathfrak p}}/U_{k,{\cal O}_{\mathfrak p},r}
\cong 
{\cal O}_{\mathfrak p}/(p^r)\otimes L_{k,{\cal O}_{\mathfrak p}}/{\cal O}_{\mathfrak p}/(p^r)\otimes U_{k,{\cal O}_{\mathfrak p},r}$; in particular, the above Lemma shows that this quotient is a $\Gamma_1(p^r)$-module.

{\bf Proposition. }{\it Fix $r\in{\Bbb N}$. If $k,k'\in{\Bbb N}$ satisfy $k,k'\ge r$ and $k\equiv k'\pmod{p^{2r}}$, then 
$$
{\cal O}_{\mathfrak p}/(p^r)\otimes L_{k,{\cal O}_{\mathfrak p}}/U_{k,{\cal O}_{\mathfrak p},r}
\cong
{\cal O}_{\mathfrak p}/(p^r)\otimes L_{k',{\cal O}_{\mathfrak p}}/U_{k',{\cal O}_{\mathfrak p},r}
$$
as $\Gamma_1(p^r)$-modules 

%[[and hence, their belonging cohomology groups are isomorphic! this is what we will need below]]
}

{\it Proof. } We assume that $k'\ge k$. We then may define a map
$$
\begin{array}{cccc}
i:&{\cal O}_{\mathfrak p}/(p^r)\otimes L_{k,{\cal O}_{\mathfrak p}}/U_{k,{\cal O}_{\mathfrak p},r}&\rightarrow&{\cal O}_{\mathfrak p}/(p^r)\otimes L_{k',{\cal O}_{\mathfrak p}}/U_{k',{\cal O}_{\mathfrak p},r}\\
&X^{k-i}Y^i+p^rL_{k,{\cal O}_{\mathfrak p}}/U_{k,{\cal O}_{\mathfrak p},r}&\mapsto&X^{k'-i}Y^i+p^rL_{k',{\cal O}_{\mathfrak p}}/U_{k',{\cal O}_{\mathfrak p},r},
\end{array}\leqno(13)
$$
where $i=0,\ldots,r$. Since $L_{k,{\cal O}_{\mathfrak p}}=\bigoplus_{i=0}^k {\cal O}_{\mathfrak p}X^{k-i}Y^i$ and recalling the definition of $U_{k,{\cal O}_{\mathfrak p},r}$ we find that
$$
{\cal O}_{\mathfrak p}/(p^r)\otimes L_{k,{\cal O}_{\mathfrak p}}/U_{k,{\cal O}_{\mathfrak p},r}=\bigoplus_{i=0}^{r} {\cal O}_{\mathfrak p}/(p^r)\, X^{k-i}Y^i.
$$
(Here we need the asumption $k\ge r$.) The same holds true with $k$ replaced by $k'$ and we immediately deduce that (13) is an isomorphism of ${\cal O}_{\mathfrak p}/(p^r)$-modules.

We want to show that (13) is $\Gamma_1(p^r)$-equivariant. To this end let $\gamma=\Mat{a}{b}{c}{d}\in\Gamma_1(p^r)$ be
arbitrary. Since $\gamma^\iota=\Mat{d}{-b}{-c}{a}$ and $a\equiv d\equiv 1\pmod{p^r}$, $c\equiv 0\pmod{p^r}$ we obtain for any $0\le i\le r$
\begin{eqnarray*}
&&\gamma X^{k'-i}Y^i\\
&&=\,(dX-bY)^{k'-i}(-cX+aY)^i\\
&&\equiv\,(X-bY)^{k'-i}(Y)^i\pmod{p^r L_{k',{\cal O}_{\mathfrak p}}/U_{k',{\cal O}_{\mathfrak p},r}}\\
&&=\,\sum_{j=0}^{k'-i} {k'-i\choose j} \,(-b)^{k'-i-j}\, X^j Y^{k'-j}\\
&&\equiv\,\sum_{j=k'-r}^{k'-i} {k'-i\choose j} \, (-b)^{k'-i-j}\, X^j Y^{k'-j}\pmod{p^r L_{k',{\cal O}_{\mathfrak p}}/U_{k',{\cal O}_{\mathfrak p},r}}\\
&&=\,(-b)^0{k'-i\choose k'-i}X^{k'-i}Y^{i}+(-b)^1{k'-i\choose k'-i-1}X^{k'-i-1}Y^{i+1}+\cdots+(-b)^{r-i}{k'-i\choose k'-r}X^{k'-r}Y^{r}\\
&&=\,(-b)^0{k'-i\choose 0}X^{k'-i}Y^{i}+(-b)^1{k'-i\choose 1}X^{k'-i-1}Y^{i+1}+\cdots+(-b)^{r-i}{k'-i\choose r-i}X^{k'-r}Y^{r}.
%&\equiv&\sum_{j=k'-[\alpha]}^{k-i} {k'-i\choose j} \,d^j (-b)^{k'-i-j}a^i\, X^j Y^{k'-j}\\
\end{eqnarray*}
%[[beachte: $r\ge i \Rightarrow {k'-i\choose r-i}$ wohldefiniert; das ist die Ungleichung fur $i$, die also benutzt wird]]

In the same way we find for any $0\le i\le r$
\begin{eqnarray*}
\gamma X^{k-i}Y^i=(-b)^0{k-i\choose 0}X^{k-i}Y^{i}&+&(-b)^1{k-i\choose 1}X^{k-i-1}Y^{i+1}+\cdots\\
&&\cdots+(-b)^{r-i}{k-i\choose r-i}X^{k-r}Y^{r}\pmod{p^rL_{k,{\cal O}_{\mathfrak p}}/U_{k,{\cal O}_{\mathfrak p},r}}.\\
\end{eqnarray*}
%$$
%[[=\sum_{j=k-[\alpha]}^{k-i} {k-i\choose j} \,d^j (-b)^{k-i-j}a^i\, X^j Y^{k-j}\pmod{p^r L_{k,{\cal O}_{\mathfrak %p}}/U_{k,{\cal O}_{\mathfrak p},r}}]]
%$$

We claim that for all $0\le j\le r-i$      
%
%[[das ist die Ungleichung fur $j$, die also auch benutzt wird]]        
%
the following congruence holds
$$
{k'-i\choose j}\equiv{k\choose j} \pmod{p^r}.\leqno(14)
$$
The above equations together with (14) immediately imply that for all $0\le i\le r$, 
\begin{eqnarray*}
i\gamma (X^{k-i}Y^i+p^rL_{k,{\cal O}_{\mathfrak p}}/U_{k,{\cal O}_{\mathfrak p},r})&=&\gamma X^{k'-i}Y^i+p^rL_{k',{\cal O}_{\mathfrak p}}/U_{k',{\cal O}_{\mathfrak p},r}\\
&=&\gamma i(X^{k-i}Y^i+p^rL_{k,{\cal O}_{\mathfrak p}}/U_{k,{\cal O}_{\mathfrak p},r}),
\end{eqnarray*}
hence, $i$ is $\Gamma_1(p^r)$-equivariant.

It remains to prove (14). Since $j\le r-i\le r$ we obtain 
$$
v_p(j!)\le v_p(r!).
$$ 
%[[beachte: $j'\le r \Rightarrow j'!|r! \Rightarrow v_p(j'!)\le v_p(r!)$]]
Let $\sum _h a_hp^h$ be the $p$-adic expansion of $r$; the $p$-adic valuation of $r!$ is then given as
$$
v_p(r!)=\frac{1}{p-1}\sum_h a_h(p^h-1)
$$
(cf. [N] (5.6) Lemma, p. 143). Hence,
$$
v_p(j!)\le v_p(r!)\le \sum_h a_hp^h=r.
$$
Taking into account that
$$
{k'-i\choose j}=\frac{(k'-i)\cdots(k'-i-j'+1)}{j'!}
$$
we see that $k\equiv k'\pmod{p^{2r}}$ and $v_p(j!)\le r$ imply 
$$
{k'-i\choose j}\equiv {k-i\choose j}\pmod{p^r}.
$$
%[[beachte: $j\le r-i\le k'-i,k-i$, da $k,k'\ge r$; dies ist die Ungleichung fur $k,k'$, die also auch benutzt wird]].

This finishes the proof of equation (14) and, hence, the proposition is proven.

{\it Remark. } The statements of the preceeding Lemma and Proposition in particular hold if we replace $\Gamma_1(p^r)$ by the (smaller) group $\Gamma_1(Np^r)$.

{\bf 1.5. The slope $\alpha$ subspace of ${\cal M}_k(\Gamma_1(Np))$}.
%
%[[Im Folgenden taucht zum ersten Mal die Gruppe $\Gamma_1(Np)$ und der Atkin-Lehner Hecke Operator auf (sowie der zugehorige slope $\alpha$-Unterraum) und das Theorem uber die Beschranktheit der Dimenion von $H^1(\Gamma_1(Np^r),L_{k,{\Bbb C}})^\alpha$ bewiesen.]]
%]]
%
%
We fix a prime $p\in{\Bbb N}$, $p>2$ and a "tame level" $N\in{\Bbb N}$, i.e. $(N,p)=1$. We set $\Gamma=\Gamma_1(Np)$ and we denote by $T_p$ the Hecke Operator
$$
T_p=\Gamma_1(Np)\Mat{1}{0}{0}{p}\Gamma_1(Np).
$$
We choose a finite extension $E/{\Bbb Q}$ such that $T_p$ is split on $H^1(\Gamma,L_{k,E})$ and we denote by ${\cal O}$ the integers of $E$. Furthermore, we select a prime ideal ${\mathfrak p}\le{\cal O}$ dividing $p$ and we denote by $E_{\mathfrak p}$ and ${\cal O}_{\mathfrak p}$ the completions. 

We recall the definition of the action of the Hecke operators on group cohomology. Let
$$
T=\Gamma\alpha\Gamma=\dot{\bigcup}_{i=1,\ldots,m} \Gamma\alpha_i
$$
be the decomposition into $\Gamma$ left cosets. For any $\gamma\in\Gamma$ the element $\alpha_i\gamma$ is 
contained in $\Gamma\alpha\Gamma$ and we therefore can write  $\alpha_i\gamma=\rho_i(\gamma)\alpha_{\pi(i)}$ for some $1\le\pi(i)\le m$ and some $\rho_i(\gamma)\in \Gamma$. Let $\omega:\,\Gamma^{d+1}\rightarrow M$ be a homogeneous cocycle representing a class in $H^d(\Gamma,M)$; then
$$
(T\omega)(\gamma_0,\ldots,\gamma_d)=\sum_{i=1}^m \alpha_i^\iota \omega(\rho_i(\gamma_0),\ldots,\rho_i(\gamma_d)).\leqno(15)
$$
Note, that with this definition of the action of Hecke operators on group cohomology the Eichler-Shimura isomorphism becomes Hecke equivariant (cf. [Hi], p. 176, 177)

%[[cf. Hida], p. 175,176, [Kuga-Parry-Sah, p. 227]; beachte dabei folgendes: Der Hecke operator $T$ operiert von RECHTS (s. [Hida], p. 184, (6b)), [Kuga-Parry-Sah], p. 226, Abschnitt 1.2), deshalb mussen wir setzen
%$$
%(T\omega)(\gamma_0,\ldots,\gamma_d)=\sum_i \alpha_i^{-1} \omega(\rho_i(\gamma_0),\ldots,\rho_i(\gamma_d))
%$$
%oder
%$$
%(T\omega)(\gamma_0,\ldots,\gamma_d)=\sum_i \alpha_i^\iota \omega(\rho_i(\gamma_0),\ldots,\rho_i(\gamma_d)).
%$$
%der Unterschied ist wieder, dass wir im ersten Fall den Ismomorphietyp useres Koeffizientensystems $L_{k,{\Bbb Q}}$ abandern: der zentrale Charakter geht nun in sein INverses uber !. Wir wollen aber eine rechts Aktion des Hecke Operators auf einem Koeffizientensystem das isomorph ist zu $L_{k,{\Bbb Q}}$ (und nicht isomorph zum dualen modul $L_{k,{\Bbb Q}}^\wedge$ !). Deshalb mussen wir die zweite Variante wahlen. Blodsinn: der Grund ist naturlich, dass nur mit der Definition der Hecke Aktion mit $\alpha^\iota$ der Eichler-Shimura isomorphismus ein Isomorphismus von Hecke Moduln wird ! (Bei verwendung der defintion mit $\alpha^{-1}$ muss man die Hecke Aktion noch twisten damit Eichler-Shimura Hecke equivariant wird !)
%]]

Since $T_p$ acts on cohomology we obtain a corresponding slope decomposition
$$
H^1(\Gamma,L_{k,E})=\bigoplus_{\alpha}  H^1(\Gamma,L_{k,E})^\alpha.
$$
Moreover, $T_p$ acts on cohomology with integer coefficients, hence, $T_p$ leaves the lattice $H^1(\Gamma,L_{k,{\cal O}_{\mathfrak p}})_{\rm int}$ stable. We deduce that 

$\bullet$ the slope $\alpha$ subspace is non-trivial only if $\alpha\in\frac{1}{e}{\Bbb N}\cup\{0\}$, 
where $e$ is the ramification index of $p|{\mathfrak p}$ and

$\bullet$ $$
H^1(\Gamma,L_{k,{\cal O}_{\mathfrak p}})_{\rm int}^\alpha=H^1(\Gamma,L_{k,{\cal O}_{\mathfrak p}})_{\rm int}
\cap
H^1(\Gamma,L_{k,E})^\alpha
$$
is an ${\cal O}_{\mathfrak p}$-lattice in $H^1(\Gamma,L_{k,E})^\alpha$ (cf. Remark in section 1.2). 

In particular
$$
{\rm rank}_{{\cal O}_{\mathfrak p}} H^1(\Gamma,L_{k,{\cal O}_{\mathfrak p}})_{\rm int}^\alpha
=
{\rm dim}_{E} H^1(\Gamma,L_{k,E})^\alpha.\leqno(16)
$$

We denote by $T_{p,r}=\Gamma_1(Np^r)\Mat{1}{}{}{p}\Gamma_1(Np^r)$ the Hecke operator of level $Np^r$, i.e. $T_p=T_{p,1}$.

{\bf Lemma. }{\it  $T_{p,r}$ annihilates $H^\bullet(\Gamma_1(Np^r),U_{k,{\cal O}_{\mathfrak p},r}\otimes {\cal O}_{\mathfrak p}/(p^r))$.
}

{\it Proof. } We write
$$
T_{p,r}=\dot{\bigcup}_{i} \Gamma_1(Np^r)\beta_i
$$
The Lemma then is an immediate consequence of section 1.4. Lemma 2.) and (15) taking into account that $\beta_i\in \Gamma_1(p^r)\Mat{1}{}{}{p}\Gamma_1(p^r)$,  $i=1,\ldots,m$.

Although we will not need the following result we record it because it is of indepenedent interest.

{\bf Proposition. }{\it $T_p$ acts nilpotently on the Torsion submodule $H^1(\Gamma_1(Np),L_{k,{\cal O}_{\mathfrak p}})_{\rm tor}$ of $H^1(\Gamma_1(Np),L_{k,{\cal O}_{\mathfrak p}})$.
}

{\it Proof.} We continue to set $\Gamma=\Gamma_1(Np)$ and we denote by
$$
H^i(\Gamma_1(Np),L_{k,{\cal O}_{\mathfrak p}})[{\mathfrak p}^a]
$$ 
the submodule of $H^i(\Gamma_1(Np),L_{k,{\cal O}_{\mathfrak p}})$ consisting of all elements which are annihilated by ${\mathfrak p}^a$. The short exact sequence
$$
0\rightarrow L_{k,{\cal O}_{\mathfrak p}}\stackrel{\cdot\varpi^a}{\rightarrow} L_{k,{\cal O}_{\mathfrak p}}\rightarrow L_{k,{\cal O}_{\mathfrak p}}/{\mathfrak p}^a L_{k,{\cal O}_{\mathfrak p}}\rightarrow 0
$$
yields the exact sequnce
$$
(L_{k,{\cal O}_{\mathfrak p}}/{\mathfrak p}^aL_{k,{\cal O}_{\mathfrak p}})^\Gamma=H^0(\Gamma,L_{k,{\cal O}_{\mathfrak p}}/{\mathfrak p}^aL_{k,{\cal O}_{\mathfrak p}})\stackrel{\delta}{\rightarrow} H^1(\Gamma,L_{k,{\cal O}_{\mathfrak p}})\stackrel{\cdot\varpi^a}{\rightarrow}
H^1(\Gamma,L_{k,{\cal O}_{\mathfrak p}}).
$$
Thus,
$$
H^1(\Gamma,L_{k,{\cal O}_{\mathfrak p}})[{\mathfrak p}^a]=\delta((L_{k,{\cal O}_{\mathfrak p}}/{\mathfrak p}^aL_{k,{\cal O}_{\mathfrak p}})^\Gamma)
$$
and since $\delta$ commutes with the action of Hecke Operators 
%[[Referenz: [Ash-Stevens] "Cohomology of arithmetic groups and congruences between systems of Hecke eigenvalues" in Crelle Band 365, 1985, Lemma 1.1.1, p. 194]]]
it is sufficient to prove that $T_p$ acts nilpotently on $(L_{k,{\cal O}_{\mathfrak p}}/{\mathfrak p}^aL_{k,{\cal O}_{\mathfrak p}})^\Gamma$. An easy computation yields
$$
\Gamma_1(Np)\Mat{1}{}{}{p}\Gamma_1(Np)=\bigcup_{u=0}^{p-1} \Gamma_1(Np)\Mat{1}{u}{}{p},
$$
%[[(cf. [M] Lemma 4.5.6)]] 
hence, for any $v\in(L_{k,{\cal O}_{\mathfrak p}}/{\mathfrak p}^aL_{k,{\cal O}_{\mathfrak p}})^\Gamma=H^0(\Gamma,L_{k,{\cal O}_{\mathfrak p}}/{\mathfrak p}^aL_{k,{\cal O}_{\mathfrak p}})$ we obtain 
$$
Tv=\sum_{u=0}^{p-1}\Mat{1}{u}{}{p}^\iota v.\leqno(17)
$$ 
%[[$\bullet$ Sei $M$ ein {\cal O}-Modul auf dem $\Gamma$ operiert.
%
%$\bullet$ Der Komplex: $C^i(\Gamma,M)=\{c:\,\Gamma^{i+1}\rightarrow M|\,c(\gamma\gamma_0,\ldots,\gamma\gamma_i)=\gammac(\gamma_0,\ldots,\gamma_i)}$, $i=0,1,2,\ldots$. $d_0: C^0\rightarrow C^1 $ ist gegeben durch $dc(\gamma_0,\gamma_1)=c(\gamma_1)-c(\gamma0)$
%(vgl. [Kuga-Parry-Sah]. p. 227 oder [Cassels-Frohlich], Aufsatz von Atiyah-Wall, p. 96,97)
%
%$\bullet$ $H^0(\Gamma,M)$: $H^0(\Gamma,M)={\rm ker}\,d_0$. Sei $c\in C^0(\Gamma,M)$; dann ist $dc=0$ $\Leftrightarrow$ $c(\gamma_1)=c(\gamma_0)$ fur alle $\gamma_0,\gamma_1\in\Gamma$, es folgt also $c(1)=c(\gamma)=\gamma c(1)$ $\Rightarrow
%c(1)\in M^\Gamma$. Wir erhalten Isomorphismus von ${\cal O}$-Moduln
%$$
%H^0(\Gamma,M)\cong M^\Gamma,\quaf c\mapsto c(1).
%$$
%
%$\bullet$ Hecke Operator $T_p$: Sei $c\in H^0(\Gamma,M)$; nach Definition ist 
%$T_pc(\gamma_0)=\sum_i \alpha_i^{\iota} c(\rho_i(\gamma_0)$ $\Rightarrow$
%$T_pc(1)=\sum_i \alpha_i^{\iota} c(\underbrace{\rho_i(1)}_{=1})$ $\Rightarrow$
%$T_pc(\gamma_0)=\sum_i \alpha_i^\iota c(1)$ $\Rightarrow$ $T_p$ operiert auf $M^\Gamma$ durch $T_pv=\sum_i \alpha_i^\iota v$. (Beachte, dass $\alpha_i 1=1 \alpha_i$, d.h. $\rho_i(1)=1$)
%]]
Equation (17) defines an Operator on all of $L_{k,{\cal O}_{\mathfrak p}}/{\mathfrak p}^aL_{k,{\cal O}_{\mathfrak p}}$, which we also denote by $T_p$.
We then will show that $T_p$ acts nilpotently on $L_{k,{\cal O}_{\mathfrak p}}/{\mathfrak p}^aL_{k,{\cal O}_{\mathfrak p}}$, which implies that it does so on $\Gamma$-invariants.
%[[this implies that $T_p$ acts nilpotently on the $\Gamma$-invariants $(L_{k,{\cal O}_{\mathfrak p}}/{\mathfrak p}^aL_{k,{\cal O}_{\mathfrak p}})^\Gamma$, because $\Gamma$-invariants is a ${\cal O}$-submodule]]
Recall that $L_{k,{\cal O}_{\mathfrak p}}=\bigoplus_{i=0}^k {\cal O}_{\mathfrak p} X^{k-i}Y^i$. For any $j=k,\ldots,0$ we denote by $L_j$ the image of
$$
\bigoplus_{i=0}^j {\cal O}_{\mathfrak p} X^iY^{k-i}\le L_{k,{\cal O}_{\mathfrak p}}
$$
in $L_{k,{\cal O}_{\mathfrak p}}/{\mathfrak p}^aL_{k,{\cal O}_{\mathfrak p}}$. (17) implies that $L_j$ is invariant under $T_p$ and we show that for all $j=k,\ldots,1$ the operator $T_p^a$ annihilates $L_j/L_{j-1}$; this in particular implies that $T_p$ acts nilpotently on $L_{k,{\cal O}_{\mathfrak p}}/{\mathfrak p}^aL_{k,{\cal O}_{\mathfrak p}}$. We first assume $j=k$; then $L_k/L_{k-1}=\langle X^k\rangle$ and $T_p X^k=\sum_{u=0}^{p-1} (X+uY)^k=\sum_u X^k=pX^k\pmod{L_{k-1}}$,
%[[note that $\alpha^\iota P(X,Y)=P((\alpha^\iota)^\iota (X,Y)^t)=P(\alpha(X,Y))$ $\Rightarrow \Mat{1}{u}{}{p}^\iota X^jY^{k-j}=(X+uY)^j(pY)^{k-j}$.]]
hence, $T_p^a X^k=p^a X^k=0\pmod{L_{k-1}}$.
%[[beachte, dass $(p^{aj})=(\varpi^{eaj})={\mathfrak p}^{eaj}$, wobei $e=$ Verzweigungindex von ${\mathfrak p}$ uber $p$.]] 
If $j<k$ then
$L_j/L_{j-1}=\langle X^jY^{k-j}\rangle$ and $T_pX^jY^{k-j}=\sum_{u=0}^{p-1} (X+uY)^j(pY)^{k-j}=\sum_u p^{k-j} X^jY^{k-j}\pmod{L_{j-1}}$, hence, $T_p^a X^jY^{k-j}=p^{a(k-j)}X^jY^{k-j}=0\pmod{L_{j-1}}$. Thus, $T_p^a$ acts nilpotently on $L_{k,{\cal O}_{\mathfrak p}}/{\mathfrak p}^aL_{k,{\cal O}_{\mathfrak p}}$ and therefore also on $H^1(\Gamma,L_{k,{\cal O}_{\mathfrak p}})[{\mathfrak p}^a]$, which proves the Proposition.

\bigskip

{\bf Theorem A. }{\it Fix any $\alpha\in{\Bbb Q}_{\ge 0}$. Then, there is $M(\alpha)\in{\Bbb N}$ such that 
$$
\dim \bigoplus_{0\le\beta\le\alpha} H^1(\Gamma_1(Np^{[\alpha]+1}),L_{k,{\Bbb C}})^\beta\le M(\alpha)
$$
for all $k\ge 2$. In different words, 
$\sum_{0\le\beta\le\alpha} \dim H^1(\Gamma_1(Np^{[\alpha]+1}),L_{k,{\Bbb C}})^\beta$ is bounded independently of $k$.

}

In the above theorem, the slope $\alpha$ subspace $H^1(\Gamma_1(Np^{[\alpha]+1}),L_{k,{\Bbb C}})^\alpha$ is
defined with respect to the Hecke operator $T_{p,[\alpha]+1}=\Gamma_1(Np^{[\alpha]+1})\Mat{1}{}{}{p}\Gamma_1(Np^{[\alpha]+1})$ of level $Np^{[\alpha]+1}$.

{\it Proof. } We set $r=[\alpha]+1$, hence, $r>\alpha$, and {\it only in this proof} we put $\Gamma=\Gamma_1(Np^r)$. We choose a finite extension $E/{\Bbb Q}$ such that $T_{p,r}$ splits on $H^1(\Gamma,L_{k,E})$ (notice that $E$ depends on $k$).  
We denote by ${\cal O}$ the ring of integers of $E$ and we select a prime ideal ${\mathfrak p}\le {\cal O}$ above $p$. As before, ${\cal O}_{\mathfrak p}$ resp. $E_{\mathfrak p}$ is the completion of ${\cal O}$ resp. $E$ with respect to ${\mathfrak p}$, $e$ is the degree of ramification of ${\mathfrak p}|p$, hence, $\alpha\in\frac{1}{e}{\Bbb N}\cup\{0\}$ and $q=p^f=|{\cal O}_{\mathfrak p}/{\mathfrak p}|$ is the cardinality of the residue field. Finally, if $L$ is any coefficient system, we shall sometimes use $H^1(L)$ as a short hand notation for $H^1(\Gamma,L)$.

Since the slope decomposition of $H^1(\Gamma,L_{k,{\Bbb C}})$ is defined over ${\Bbb Q}$ we have $H^1(\Gamma,L_{k,{\Bbb C}})^\beta=H^1(\Gamma,L_{k,{E}})^\beta\otimes{\Bbb C}$,
%[[Insbesondere ist slope decomposition also uber $E$ definiert, woraus die obige Gleichung unmittelbar folgt]],
thus, after tensoring with $E_{\mathfrak p}$ we are reduced to showing that
$$
\sum_{0\le\beta\le\alpha} \dim_{E_{\mathfrak p}} H^1(\Gamma,L_{k,E_{\mathfrak p}})^\beta
$$
is bounded independently of $k$. By using (15) this in turn is equivalent to proving the boundedness of $\sum_{0\le\beta\le\alpha}d_k^\beta$ for varying $k$, where
$$
d_k^\beta:={\rm rank}_{{\cal O}_{\mathfrak p}} H^1(\Gamma,L_{k,{\cal O}_{\mathfrak p}})_{\rm int}^\beta.
$$
%[[$H^1(L_{k,{\cal O}_{\mathfrak p}})_{\rm int}^\alpha\subset H^1(L_{k,{\cal O}_{\mathfrak p}})_{\rm int}$ is a pure submodule, hence, complemented and we deduce that
%$$
%{\rm rk}_{{\cal O}_{\mathfrak p}/(p)^r} H^1(L_{k,{\cal O}_{\mathfrak p}})_{\rm int}^\alpha\otimes{\cal O}_{\mathfrak p}/(p^r)={\rm rk}_{{\cal O}_{\mathfrak p}} H^1(L_{k,{\cal O}_{\mathfrak p}})_{\rm int}^\alpha.
%$$
%($\rightarrow$ Beweis fur diese Gleichheit !)
%]]
We set ${\cal TOR}=H^1(\Gamma,L_{k,{\cal O}_{\mathfrak p}})_{\rm tor}$. The isomorphism
$$
i^*:\,H^1(\Gamma,L_{k,{\cal O}_{\mathfrak p}})/{\cal TOR}\stackrel{\sim}{\rightarrow}H^1(\Gamma,L_{k,{\cal O}_{\mathfrak p}})_{\rm int} ,
$$
(cf. (12))
%[[beachte: 1.) Es gilt $H^1(\Gamma_1(Np^{[\alpha]+1}),L_{k,E})=H^1(\Gamma_1(Np^{[\alpha]+1}),L_{k,E})\otimes_E E$ und  $T_{p,r}=T_{p,r}'\otimes E$, wobei $T_{p,r}'$ der Hecke Operator auf $H^1(\Gamma_1(Np^{[\alpha]+1}),L_{k,E})$ ist; deshalb gilt $H^1(\Gamma_1(Np^{[\alpha]+1}),L_{k,E})^\alpha=H^1(\Gamma_1(Np^{[\alpha]+1}),L_{k,E})^\alpha\otimes_E E$
%woraus folgt dass $\dim_E H^1(\Gamma_1(Np^{[\alpha]+1}),L_{k,E})\alpha=\dim_{E }H^1(\Gamma_1(Np^{[\alpha]+1}),L_{k,E})$ ist. 2.) $H^1(...,L_{\cal O})_{\rm int}^\alpha$ ist Gitter in $H^1(...,L_E)^\alpha$, vgl. (lattice) 3.) der Isomorphismus $i^*:\,H^1(\Gamma_1(Np^{[\alpha]+1}),L_{k,{\cal O}_{\mathfrak p}})_{\rm int}^\alpha
%\cong (H^1(\Gamma_1(Np^{[\alpha]+1}),L_{k,{\cal O}_{\mathfrak p}})/{\rm Tor})$ ist induziert durch die Inklusion $i:\,L_{k,{\cal O}_{\mathfrak p}}\hookrightarrow L_{k,E}$, d.h. er kommutiert mit Hecke Operatoren (Kohomologie ist ein Funktor), insbesondere also mit $T_{p,r}$, also bildet $i^*$ den slope $\alpha$ auf den slope $\alpha$ Unterraum ab.]]
induces an isomorphism
$$
i^*:\,(H^1(\Gamma,L_{k,{\cal O}_{\mathfrak p}})/{\cal TOR})\otimes{\cal O}_{\mathfrak p}/(p^r)\stackrel{\sim}{\rightarrow} H^1(\Gamma,L_{k,{\cal O}_{\mathfrak p}})_{\rm int}\otimes{\cal O}_{\mathfrak p}/(p^r).
$$
%[[(explicitly, $i^*:\,v\otimes 1\mapsto i^*(v)\otimes 1$).]]
%[[After reducing modulo $(p^r)$, $i$ induces isomorphisms
%$$
%\begin{array}{ccc}
%H^1(L_{k,{\cal O}_{\mathfrak p}})_{\rm int}^\alpha&\stackrel{i}{\cong}&i(H^1(L_{k,{\cal O}_{\mathfrak p}})_{\rm int}^\alpha)\\
%\cap&&\cap\\
%H^1(L_{k,{\cal O}_{\mathfrak p}})_{\rm int}&\stackrel{i}{\cong}&H^1(L_{k,{\cal O}_{\mathfrak p}})/{\cal TOR}\\
%\downarrow&&\downarrow\\
%H^1(L_{k,{\cal O}_{\mathfrak p}})_{\rm int}\otimes{\cal O}_{\mathfrak p}/(p^r)&\stackrel{i}{\cong}&
%H^1(L_{k,{\cal O}_{\mathfrak p}})/{\cal TOR}\otimes {\cal O}_{\mathfrak p}/(p^r)\\
%\cup&&\cup\\
%H^1(L_{k,{\cal O}_{\mathfrak p}})_{\rm int}^\alpha\otimes{\cal O}_{\mathfrak p}/(p^r)&\stackrel{i}{\cong}&i(H^1(L_{k,{\cal O}_{\mathfrak p}})_{\rm int}^\alpha)\otimes{\cal O}_{\mathfrak p}/(p^r)\\
%\end{array}
%$$
%
%Note that $H^1(L_{k,{\cal O}_{\mathfrak p}})_{\rm int}^\alpha$ and $H^1(L_{k,{\cal O}_{\mathfrak p}})_{\rm int}^\alpha\otimes{\cal O}_{\mathfrak p}/(p^r)$ are complemented submodules, hence, the same is true for their images under $i$ (because $i$ is an isomorphism). in particular, we have 
%${\rm rk}_{{\cal O}_{\mathfrak p}}H^1(L_{k,{\cal O}_{\mathfrak p}})_{\rm int}^\alpha={\rm rk}_{{\cal O}_{\mathfrak p}/(p^r)H^1(L_{k,{\cal O}_{\mathfrak p}})_{\rm int}^\alpha\otimes{\cal O}_{\mathfrak p}/(p^r)$ (BEWEIS !!!; brauchen wir das uberhaupt ????) SIEHE DAZU DIE PRINZIPIELLE ANMERKUNG GANZ ZUM SCHLUSS
%]]
Moreover, the short exact sequence
$$
0\rightarrow L_{k,{\cal O}_{\mathfrak p}}\stackrel{\cdot p^r}{\rightarrow} L_{k,{\cal O}_{\mathfrak p}}
\stackrel{\pi}{\rightarrow} {\cal O}_{\mathfrak p}/(p^r)\otimes L_{k,{\cal O}_{\mathfrak p}}\rightarrow 0
$$
induces an exact sequence (via the long exact cohomology sequence)
$$
H^1(\Gamma,L_{k,{\cal O}_{\mathfrak p}})
\stackrel{\cdot p^r}{\rightarrow} 
H^1(\Gamma,L_{k,{\cal O}_{\mathfrak p}})
\stackrel{\pi^*}{\rightarrow} H^1(\Gamma,{\cal O}_{\mathfrak p}/(p^r)\otimes L_{k,{\cal O}_{\mathfrak p}})
\rightarrow 
0.
$$
Note that $H^2(\Gamma,M)=0$ for all $\Gamma$-modules $M$ and $(\cdot p^r)^*$ again just is multiplication by $p^r$. 
%[[which need not be injective]] 
We obtain a $T_{p,r}$-equivariant isomorphism
$$
H^1(\Gamma,L_{k,{\cal O}_{\mathfrak p}})/p^rH^1(\Gamma,L_{k,{\cal O}_{\mathfrak p}})
\stackrel{\pi^*}{\cong} 
H^1(\Gamma,{\cal O}_{\mathfrak p}/(p^r)\otimes L_{k,{\cal O}_{\mathfrak p}})\leqno(18)
$$
induced by reduction mod $p^r$ of the values of cohomology classes.

The short exact sequence
$$
0\rightarrow U_{k,{\cal O}_{\mathfrak p},r}\otimes{\cal O}_{\mathfrak p}/(p^r)\stackrel{j}{\rightarrow} L_{k,{\cal O}_{\mathfrak p}}\otimes{\cal O}_{\mathfrak p}/(p^r)\stackrel{p}{\rightarrow} L_{k,{\cal O}_{\mathfrak p}}/U_{k,{\cal O}_{\mathfrak p},r}\otimes{\cal O}_{\mathfrak p}/(p^r)\rightarrow 0,
$$
where $j$ is the inclusion map and $p:\,v\otimes\alpha\mapsto (v+U_{k,{\cal O}_{\mathfrak p},r})\otimes\alpha$ the canonical projection (note that $U_{k,{\cal O}_{\mathfrak p},r}$ is complemented in $L_{k,{\cal O}_{\mathfrak p}}$, which yields the exactness of the above sequence), 
yields via the long exact cohomology sequence an exact sequence
$$
\begin{array}{ccccc}
H^1(U_{k,{\cal O}_{\mathfrak p},r}\otimes{\cal O}_{\mathfrak p}/(p^r))
&\stackrel{j^*}{\rightarrow}&
H^1(L_{k,{\cal O}_{\mathfrak p}}\otimes{\cal O}_{\mathfrak p}/(p^r))
&\stackrel{p^*}{\rightarrow}&
H^1(L_{k,{\cal O}_{\mathfrak p}}/U_{k,{\cal O}_{\mathfrak p},r}\otimes{\cal O}_{\mathfrak p}/(p^r))\\
&&\parallel\;\pi^*&&\\
&&H^1(L_{k,{\cal O}_{\mathfrak p}})\otimes{\cal O}_{\mathfrak p}/(p^r)).&&\\
\end{array}\leqno(19)
$$
We denote by ${\cal TOR}\otimes {\cal O}_{\mathfrak p}/(p^r)=\{t\otimes 1,\,t\in{\cal TOR}\}$ the image of ${\cal TOR}$ in $H^1(L_{k,{\cal O}_{\mathfrak p}})\otimes{\cal O}_{\mathfrak p}/(p^r)$; the assignment 
$$
\begin{array}{cccc}
\varphi:&H^1(L_{k,{\cal O}_{\mathfrak p}})\otimes{\cal O}_{\mathfrak p}/(p^r)\big/{\cal TOR}\otimes {\cal O}_{\mathfrak p}/(p^r)&\cong&(H^1(L_{k,{\cal O}_{\mathfrak p}})/{\cal TOR})\otimes{\cal O}_{\mathfrak p}/(p^r)\\
&x\otimes 1+{\cal TOR}\otimes{\cal O}_{\mathfrak p}/(p^r)&\mapsto&(x+{\cal TOR})\otimes 1\\
\end{array}
$$
then defines an isomorphism
and after dividing by ${\cal T}={\cal TOR}\otimes{\cal O}_{\mathfrak p}/(p^r)$ we obtain from (19) the following diagramm with exact rows
$$
\begin{array}{ccccc}
&&(H^1(L_{k,{\cal O}_{\mathfrak p}})\otimes{\cal O}_{\mathfrak p}/p^r)/{\cal T}&&\\

&\stackrel{j^*}{\nearrow}&\cong\;\downarrow\; i^*\circ\varphi&\stackrel{p^*}{\searrow}&\\

H^1(U_{k,{\cal O}_{\mathfrak p},r}\otimes{\cal O}_{\mathfrak p}/p^r)/j^{-1}({\cal T})&\stackrel{(j^*)'}{\rightarrow}&H^1(L_{k,{\cal O}_{\mathfrak p}})_{\rm int}\otimes{\cal O}_{\mathfrak p}/(p^r)&\stackrel{(p^*)'}{\rightarrow}&H^1\big(\frac{L_{k,{\cal O}_{\mathfrak p}}}{U_{k,{\cal O}_{\mathfrak p},r}}\otimes{\cal O}_{\mathfrak p}/p^r\big)\big/p^*({\cal T})\\
&&\cup&&\\
&&\bigoplus_{0\le\beta\le\alpha} H^1(L_{k,{\cal O}_{\mathfrak p}})_{\rm int}^{\beta}\otimes{\cal O}_{\mathfrak p}/(p^r).&&\\
\end{array}
$$
Here, $(j^*)'=i^*\circ\varphi\circ j^*$ and $(p^*)'=p^*\circ (i^*\circ\varphi)^{-1}$.
%
%[[(I.) ${\cal TOR}$ ist $T_{p,r}$-invariant (da $T_{p,r}$ ${\cal O}_{\mathfrak p}$-linear ist) und da $j$ und $p$ $T_{p,r}$-equivariant sind folgt: auch  $j^{-1}({\cal TOR}), p({\cal TOR})$ sind $T_{p,r}$-invariant.
%
%
%[[(II.) Sei $A\stackrel{f}{\rightarrow} B\stackrel{g}{\rightarrow} C$ exakt und sei $U\le B$ Untermodul; dann ist auch 
%$A/f^{-1}(U)\stackrel{\bar{f}}{\rightarrow} B/U\stackrel{\bar{g}}{\rightarrow} C/g(U)$ exakt, wobei $\bar{f}(a+f^{-1}(U))=f(a)+U$ und $\bar{g}(b+U)=g(b)+g(U)$.
%
%Beweis. 1.) $\bar{f},\bar{g}$ sind wohldefiniert: $f(f^{-1}(U)\subset U \Rightarrow \bar{f}$ wohldefiniert. $g(U)\subset g(U) \Rightarrow \bar{g}$ wohldefiniert.
%
%2.) $\bar{g}\circ\bar{f}=0$: $\bar{g}\circ\bar{f}(a+f^{-1}(U))=\bar{g}(f(a)+U)=g(f)(a)+g(U)=0+U=U=0$
%
%3.) Exaktheit: Sei $\bar{g}(b+U)=0 \Rightarrow g(b)+g(U)=0 \Rightarrow g(b)\in g(U) \Rightarrow g(b)=g(u), u\in U \Rightarrow g(b-u)=0 \Rightarrow b-u=f(a), a\in A \Rightarrow b+U=f(a)+U \Rightarrow b+U=\bar{f}(a+f^{-1}(U))$.
%]]
%
%[[Die Exaktheit der horizontalen Reihe in obigem Diagramm ist der kern des Beweises des Theorems: aus ihr folgt die Beschranktheit des ranges von $H^1(...)_{\rm int}^\alpha$]]
%
By the above Lemma, $T_{p,r}=\Gamma\Mat{1}{}{}{p}\Gamma$ annihilates $H^1(U_{k,{\cal O}_{\mathfrak p},r}\otimes{\cal O}_{\mathfrak p}/p^r)$, hence, $T_{p,r}$ also annihilates
$$
X=H^1(U_{k,{\cal O}_{\mathfrak p},r}\otimes{\cal O}_{\mathfrak p}/(p^r))/j^{-1}({\cal TOR}\otimes{\cal O}_{\mathfrak p}/p^r)
$$ 
as well as the image of $X$ under $(j^*)'$. 
%
%[[$(j^*)'$ is $T_{p,r}$-equivariant, da $j^*$, $i^*$ via Funktorialitat von Morphismen auf den Koeffzsyst definiert sind, also $T_{p,r}$-equivariant sind (s. [Ash-Stevens] in Crelle J.) und $\varphi$ $T_{p,r}$-equivariant ist, wie wir oben gesehen haben]]
%
For any $\beta\in{\Bbb Q}$ satisfying $0\le \beta\le \alpha$ we set 
$$
X^\beta:=(j^*)'(X)\cap H^1(L_{k,{\cal O}_{\mathfrak p}})_{\rm int}^\beta\otimes{\cal O}_{\mathfrak p}/(p^r)
$$
and we deduce that $X^\beta\subset H^1(L_{k,{\cal O}_{\mathfrak p}})_{\rm int}^\beta[0]$ for all $0\le \beta\le \alpha$ (for the definition of $H^1(L_{k,{\cal O}_{\mathfrak p}})_{\rm int}^\beta[0]$ see section 1.2; recall that  $T_{p,r}$ annihilates $X^\beta$). Taking into account that $(p^r)={\mathfrak p}^{er}$ and $er>e\alpha\ge e\beta$, the Lemma in section 1.2 yields
$$
q^{e d_k^\beta(r-\beta)}  \,\big|\, |H^1(L_{k,{\cal O}_{\mathfrak p}})_{\rm int}^\beta\otimes {\cal O}_{\mathfrak p}/(p^r)/X^\beta|.
$$

%[[Bei Anwenden des Lemmas aus section 1.2: $H^1(L_{k,{\cal O}_{\mathfrak p}})_{\rm int}$ entspricht $V_{\cal O}$]]

%[[denn es ist $\alpha\in 1/e{\Bbb N}$, also schreibe schreibe $[\alpha]=\alpha-f/e$ mit $f<e$ $\Rightarrow  e([\alpha]+1)-e\alpha=-f+e\ge 1$; beachte auch dass $[\alpha]+1-\alpha\in{\Bbb N}$ eine ganze Zahl ist]]
%
%[[Genauer: 
%\begin{eqnarray*}
%q^{{\rm rk}_{{\cal O}_{\mathfrak p}} H^1(L_{k,{\cal O}_{\mathfrak p}})_{\rm int}^\alpha(e-r\alpha)} &|& |H^1(L_{k,{\cal O}_{\mathfrak p}})_{\rm int}^\alpha\otimes {\cal O}_{\mathfrak p}/(p^r)/ H^1(L_{k,{\cal O}_{\mathfrak p}})_{\rm int}^\alpha[0]|\\ 
%&|&|(H^1(L_{k,{\cal O}_{\mathfrak p}})_{\rm int}^\alpha\otimes {\cal O}_{\mathfrak p}/p^r)/X^\alpha|\\
%\end{eqnarray*}
%]]
%

Since $r-\beta\ge r-\alpha(\ge 0)$ this immediately yields
$$
q^{e d_k^\beta(r-\alpha)}  \,\big|\, |H^1(L_{k,{\cal O}_{\mathfrak p}})_{\rm int}^\beta\otimes {\cal O}_{\mathfrak p}/(p^r)/X^\beta|.\leqno(20)
$$ 

On the other hand, $(p^*)'$ induces an injection
$$
\bigoplus_{0\le\beta\le\alpha}(H^1(L_{k,{\cal O}_{\mathfrak p}})_{\rm int}^\beta\otimes{\cal O}_{\mathfrak p}/p^r)/X^\beta\stackrel{(p^*)'}{\hookrightarrow} H^1(L_{k,{\cal O}_{\mathfrak p}}/U_{k,{\cal O}_{\mathfrak p},r}\otimes{\cal O}_{\mathfrak p}/(p^r))\big/p^*({\cal T}).\leqno(21)
$$
To see that this holds let $\sum_{0\le \beta\le \alpha} h^\beta\in\bigoplus_{0\le \beta\le\alpha} H^1(L_{k,{\cal O}_{\mathfrak p}})_{\rm int}^\beta\otimes{\cal O}_{\mathfrak p}/p^r$ be contained in the kernel of $(p^*)'$. Then $(p^*)'(h^\beta)$ is contained in the slope-$\beta$ subspace of the target space of $(p^*)'$ and since the slope decomposition is a direct sum decomposition, $(p^*)'(\sum_\beta h^\beta)=0$ implies that $(p^*)'(h^\beta)=0$ for all $\beta$. Thus, $h^\beta\in {\rm ker}\,(p^*)'\cap H^1(L_{k,{\cal O}_{\mathfrak p}})_{\rm int}^\beta\otimes{\cal O}_{\mathfrak p}/p^r$, which by the exactness of the rows of the above diagram equals $X^\beta$. Therefore,  ${\rm ker}\,(p^*)'=\bigoplus_{0\le\beta\le\alpha} X^\beta$.

As an immediate consequence of (21) we find
$$
\prod_{0\le\beta\le\alpha}  |(H^1(L_{k,{\cal O}_{\mathfrak p}})_{\rm int}^\beta\otimes{\cal O}_{\mathfrak p}/(p^r))/X^\beta|\,\big|\, |H^1(L_{k,{\cal O}_{\mathfrak p}}/U_{k,{\cal O}_{\mathfrak p},r}\otimes{\cal O}_{\mathfrak p}/(p^r))|.\leqno(22)
$$
Since $r=[\alpha]+1$ equations (20) and (22) together yield
$$
q^{e \sum_{0\le\beta\le\alpha} d_k^\beta ([\alpha]+1-\alpha)}\,\big|\,|H^1(L_{k,{\cal O}_{\mathfrak p}}/U_{k,{\cal O}_{\mathfrak p},r}\otimes{\cal O}_{\mathfrak p}/(p^r))|.\leqno(23)
$$

We want to replace ${\cal O}_{\mathfrak p}/(p^r)$-modules appearing in (23) by ${\Bbb Z}/(p^r)$-modules. Taking into account that ${\cal O}_{\mathfrak p}$ is a flat ${\Bbb Z}_p$-module
%[[da ${\Bbb Z}_p$ Hauptidealring und ${\cal O}_{\mathfrak p}$ torsionsfrei ist, ist ${\cal O}$ eine flache ${\Bbb Z}_p$-Algebra; vgl. [Curtis-Reiner, p. 35, Corollary (2.34)]; Insbesondere gilt also $H^1(\Gamma,M\otimes_{{\Bbb Z}_p}{\cal O})=H^1(\Gamma,M)\otimes_{{\Bbb Z}_p}{\cal O}$; s. [Hida], p. 168, Formel (1a).]]
we find
\begin{eqnarray*}
(flat:24)\quad H^1(L_{k,{\Bbb Z}_{p}}/U_{k,{\Bbb Z}_{p},r}\otimes_{{\Bbb Z}_p}{\Bbb Z}_{p}/(p^r))\otimes_{{\Bbb Z}_p} {\cal O}_{\mathfrak p}
&=&H^1(L_{k,{\Bbb Z}_{p}}/U_{k,{\Bbb Z}_{p},r}\otimes_{{\Bbb Z}_p}{\Bbb Z}_{p}/(p^r) \otimes_{{\Bbb Z}_p} {\cal O}_{\mathfrak p})\\
%[[this holds because ${\cal O}_{\mathfrak p}$ is flat]]
&=&H^1(L_{k,{\Bbb Z}_{p}}/U_{k,{\Bbb Z}_{p},r}\otimes_{{\Bbb Z}_p} {\cal O}_{\mathfrak p}/(p^r))\\
&=&H^1(L_{k,{\Bbb Z}_{p}}/U_{k,{\Bbb Z}_{p},r}\otimes_{{\Bbb Z}_p}{\cal O}\otimes_{{\cal O}_{\mathfrak p}} {\cal O}_{\mathfrak p}/(p^r))\\
&=&H^1(L_{k,{\cal O}_{\mathfrak p}}/U_{k,{\cal O}_{\mathfrak p},r}\otimes_{{\cal O}_{\mathfrak p}}{\cal O}_{\mathfrak p}/(p^r))\\ 
%[[hier haben wir die Assoziativitat des Tensorproduktes benutzt wobei jeweils uber verschiedenen Ringen tensoriert wird; Referenz dafur: [Curtis-Reiner], p. 25, Proposition (2.18); beachte: unsere Ringe sind kommutativ, also Rechts=Links Moduln]]
\end{eqnarray*}

Since $H^1(L_{k,{\Bbb Z}_p}/U_{k,{\Bbb Z}_{p},r}\otimes_{{\Bbb Z}_p}{\Bbb Z}_{p}/(p^r))$ is a torsion ${\Bbb Z}_p$-module we have
$$
H^1(L_{k,{\Bbb Z}_p}/U_{k,{\Bbb Z}_{p},r}\otimes_{{\Bbb Z}_p}{\Bbb Z}_{p}/(p^r))\cong\bigoplus_{i=1}^n {\Bbb Z}_p/(p^{n_i}).
$$
Using (24) we deduce that
$$
H^1(L_{k,{\Bbb Z}_p}/L_{k,{\cal O}_{\mathfrak p},r}\otimes_{{\Bbb Z}_p}{\cal O}_{\mathfrak p}/(p^r))\cong \bigoplus_{i=1}^n {\cal O}_{\mathfrak p}/(p^{n_i}).
$$
Furthermore, using that $|{\cal O}_{\mathfrak p}/(p^n)|=p^{fen}=|{\Bbb Z}_p/(p^n)|^{[E_{\mathfrak p}:{\Bbb Q}_p]}$ we obtain
$$
|H^1(L_{k,{\cal O}_{\mathfrak p}}/L_{k,{\cal O}_{\mathfrak p},r}\otimes_{{\Bbb Z}_p}{\cal O}_{\mathfrak p}/(p^r))|=|H^1(L_{k,{\Bbb Z}_p}/U_{k,{\Bbb Z}_{p},r}\otimes_{{\Bbb Z}_p}{\Bbb Z}_{p}/(p^r))|^{[E_{\mathfrak p}:{\Bbb Q}_p]}
$$
Together with (Div) and taking into acount that $q^e=p^{ef}=p^{[E_{\mathfrak p}:{\Bbb Q}_p]}$ we thus obtain
$$
p^{\sum_{0\le\beta\le\alpha}d_k^\beta[E_{\mathfrak p}:{\Bbb Q}_p]([\alpha]+1-\alpha)}\,|\, |H^1(L_{k,{\Bbb Z}_p}/U_{k,{\Bbb Z}_{p},r}\otimes_{{\Bbb Z}_p}{\Bbb Z}_{p}/(p^r))|^{[E_{\mathfrak p}:{\Bbb Q}_p]},\leqno(25)
$$
%[[note that $ef=[E_{\mathfrak p}:{\Bbb Q}_p]$. Referenz: [N], p. 157, (6.8) Satz; beachte $E_{\mathfrak p}/{\Bbb Q}_p$ muss nur separable (nicht galoissch) sein, was automatisch erfullt ist, da $char\,{\Bbb Q}_p=0$]]
which immediately yields that
$$
\sum_{0\le\beta\le\alpha} d_k^\beta\le s(r,k):=\frac{1}{[\alpha]+1-\alpha}\,\log_p |H^1(L_{k,{\Bbb Z}_p}/L_{k,{\Bbb Z}_{p},r}\otimes_{{\Bbb Z}_p}{\Bbb Z}_{p}/(p^r))|.
$$
%[[beachte, dass $[\alpha]+1-\alpha>0$, d.h. beim Dividieren durch $[\alpha]+1-\alpha$ bleibt "$\le$" erhalten]]

The Proposition in section 1.4 implies that $s(r,k)$ for $k>r=[\alpha]+1$ only depends on the residue class of $k$ modulo $p^{2r}{\Bbb Z}$, hence, we obtain
$$
\sum_{0\le\beta\le\alpha} d_k^\beta \le M(\alpha):={\rm max}\{s(r,2),\ldots,s(r,r),s(r,k),\,k=r+1,\ldots,p^{2r}\}.
$$

Since $r=[\alpha]+1$ we see that $M(\alpha)$ only depends on $\alpha$, i.e. it does not depend on $k$. 
Thus, $\sum_{0\le\beta\le\alpha} d_k^\beta$ is bounded by a constant, which only depends on $\alpha$ and the proof of the theorem therefore is complete.

We fix an arbitrary $\alpha\in{\Bbb Q}$. The embedding
$$
{\cal M}_k(\Gamma_1(Np))\hookrightarrow {\cal M}_k(\Gamma_1(Np^{[\alpha]+1}))
$$
commutes with the action of the Hecke Operator $T_p$ on ${\cal M}_k(\Gamma_1(Np))$ and $T_{p,r}$ on  ${\cal M}_k(\Gamma_1(Np^{[\alpha]+1}))$, because the level of both spaces is divisible by $p$ (cf. [M], Theorem 4.5.10, p. 143). 
Hence, we deduce an embedding
$$
{\cal M}_k(\Gamma_1(Np))^\alpha\hookrightarrow {\cal M}_k(\Gamma_1(Np^{[\alpha]+1}))^\alpha
$$
which immediately yields that
\begin{eqnarray*}
\sum_{0\le \beta\le \alpha}\dim_{\Bbb C} {\cal M}_k(\Gamma_1(Np))^\beta&\le &\sum_{0\le \beta\le \alpha} \dim_{\Bbb C} {\cal M}_k(\Gamma_1(Np^{[\alpha]+1}))^\beta\\
&\le&\sum_{0\le\beta\le\alpha} {\rm dim}_{\Bbb C}\,H^1(\Gamma_1(Np^{[\alpha]+1}),L_{k,{\Bbb C}})^\beta\\
&\le& M(\alpha).
\end{eqnarray*}
for all $k\ge 2$. Thus, we have proven

\bigskip

{\bf Corollary. }{\it Let $\alpha\in{\Bbb Q}$ be arbitrary; for all $k\ge 2$ the following inequality is true:
$$
\sum_{0\le \beta\le\alpha} \dim_{\Bbb C} {\cal M}_k(\Gamma_1(Np))^\beta\le M(\alpha).
$$

}

\bigskip

In particular for any Dirichlet character $\chi:\,({\Bbb Z}/Np{\Bbb Z})^*\rightarrow {\Bbb C}^*$ we obtain 
$$
{\rm dim}\, {\cal M}_k(\Gamma_1(Np),\chi)^\alpha\le M(\alpha),
$$
where ${\cal M}_k(\Gamma_1(Np),\chi)$ is the space of modular forms of weight $k$, level $\Gamma_1(Np)$ and nebentype $\chi$. We set
%$$
%d(\chi,\alpha):=\max_{k\ge 2}\; \dim_{\Bbb C} {\cal M}_k(\Gamma_1(Np),\chi)^\alpha.
%$$
%and
$$
d(\alpha):=\max_{k\ge 2,\chi}\; \dim_{\Bbb C} {\cal M}_k(\Gamma_1(Np),\chi)^\alpha,
$$
where $\chi$ runs over all Dirichlet characters of level $Np$. The Corollary may then be reformulated by saying that $d(\alpha)<\infty$.

\newpage

\section{The topological trace formula}

{\bf 2.1. The Borel-Serre compactification. }  We set $G={\rm GL}_2$ and we denote by $Z_\infty$ resp. $Z_\infty^0$ the center of $G({\Bbb R})$ resp. the connected component of the center of $G({\Bbb R})$, which contains the identity. We set $K_\infty={\rm SO}_2({\Bbb R})$ and
denote by $\Theta=\Theta_{K_\infty}$ the Cartan involution attached to $K_\infty$.
%
%[[i.e. $\Theta(X)=(X^{-1})^t$]]
% 
We set $P_0=\{\Mat{*}{*}{}{*}\}\le G({\Bbb R})$ and $U_0=\{\Mat{1}{*}{}{1}\}$. We let $P\le G({\Bbb R})$ be any minimal parabolic subgroup, which is defined over ${\Bbb Q}$. $P({\Bbb R})$ contains a unique real torus $S_P$, 
%
%[[d.h. definiert uber ${\Bbb R}$]]
%
which is stable under $\Theta$ and satisfies $P({\Bbb R})={^0P}({\Bbb R})\times A_P$, where $A_P=S_P({\Bbb R})^0$ is the connected component of $1\in S_P({\Bbb R})$. 
In case $P=P_0$ we find $A_P=\{\Mat{a}{}{}{d}, a,d\in{\Bbb R}^*_+\}$ and ${^0P}({\Bbb R})=\{\Mat{\pm 1}{u}{}{\pm 1},\,u\in{\Bbb R}\}$. We denote by $\rho_P\in{\rm Hom}(S_P({\Bbb R}),{{\Bbb R}^*})$ the weight of the action of $S_P({\Bbb R})$ on $U_P({\Bbb R})$. Clearly, 
$$
\rho_{P_0}\left(\Mat{a}{}{}{d}\right)=a/d.
$$
Moreover, if $Q={^k}P$, $k\in K_\infty$, then $S_Q={^k}S_P$ 
%
%[[da $\Theta_K(k S_P k^{-1})=\theta_K(k)\Theta_K(S_P)\Theta_K(k^{-1})=\Theta_K(S_P)=S_P$ (da $K_\infty$ und $S_P$ $\Theta$-invariant sind), also ist ${^k}S_P=S_Q$]]
%
and, hence, 
$$
\rho_Q(t)=\rho_P(k^{-1}tk)\leqno(1)
$$ 
for all $t\in S_Q({\Bbb R})$. $\rho_P$ induces an isomorphism of groups
$$
\rho_P:\,A_P/Z_\infty^0\rightarrow (0,\infty).\leqno(2)
$$ 
%
%[[We extend $\rho_P$ to $P({\Bbb R})$ by setting 
%$$
%\rho_P(p)\rho_P(a)
%$$
%where $p=ua\in U_P({\Bbb R})\times S_P$. Beachte:
%$$
%\rho_P(p)=\rho_P(a)
%$$
%wobei $p={^0}pa\in{^0}P({\Bbb R})A_P$, ist nicht wohldefiniert. Aber %$|\rho_P(p)|:=|\rho_P(a)|$ ist so noch wohldefiniert]]
%
We define the distance function attached to $P$ as
$$
\begin{array}{cccc}
\ell_P:&G({\Bbb R})/Z_\infty^0&\rightarrow&(0,\infty)\\
&g&\mapsto&|\rho_P(a)|_\infty,\\
\end{array}
$$
where $g=uakZ_\infty^0\in {^0P}({\Bbb R}) A_PK_\infty/Z_\infty^0$. 
%
%[[beachte: der faktor "$a$" ist eindeutig bestimmt ducrh $g$ !]]
%
Clearly, $\ell_P$ factorizes over $K_\infty$, i.e. $\ell_P$ 
descends to a function on $G({\Bbb R})/K_\infty Z_\infty^0$.
%
%[[$\ell$ does not depend on the choice of $K_\infty$, da $\ell$ auf kompakten Untergruppen verschwindet.]]
%
$(1)$ immediately implies that for $k\in K_\infty$ and $p\in P({\Bbb R})$
$$
\ell_{{^k}P}(kpk^{-1})=\ell_{P}(p).\leqno(3)
$$
If $g\in G({\Bbb Q})$ is arbitrary with Iwasawa decomposition $g=kb\in K_\infty P({\Bbb R})$
%
%[[Sei $g^{-1}=pk$, $p\in P({\Bbb R})$, $k\in K_\infty$ die Iwasawazerlegung von $g^{-1}$, dann folgt $g=k^{-1}p^{-1}$, $k^{-1}\in K_\infty$, $p^{-1}\in P({\Bbb R})$, d.h. eine Iwasawazerlegung mit vertauschter Reihenfolge]]
%
we obtain 
$$
\ell_{{^g}P}({^g}p)=\ell_{{^k}P}(^{kb}p)=\ell_P({^b}p)=\ell_P(p),\leqno(3')
$$ 
i.e. (3) even holds for arbitrary $g\in G({\Bbb R})$ instead of $k\in K_\infty$.

%[[beachte: $b,p\in P({\Bbb R})$ $\Rightarrow \ell_P({^b}p)=\ell_P(p)$]]

The face associated to $P$ is defined as
$$
e(P)={^0P}({\Bbb R})/K_P=P({\Bbb R})/A_PK_P,
$$
where $K_P=K_\infty\cap P({\Bbb R})$. 
%
%[[beachte: $P({\Bbb R})={^0}P({\Bbb R})\times A_P$.]]
%
%[[Grundsatzliche Bemerkung: $g\mapsto ga$, $a\in S_P$ ist Aktion durch Linkstranslation von $S_P$ auf $G$. 
%$gK_\infty Z_\infty^0\mapsto gaK_\infty Z_\infty^0$ ist die {\it geodatische Aktion} von $S_P$ auf $X=G({\Bbb R})/K_\infty Z_\infty^0$, d.h. fur geodatische Aktion "ist das $K_\infty$ dabei", so dass man im symmetrischen Raum landet]]
%
For any $t\in (0,\infty)$ we set $U_{P,t}=e(P)\times(t,\infty)$. The map
$$
\begin{array}{cccc}
\varphi:&e(P)\times A_P/Z_\infty^0=P({\Bbb R})/K_PZ_\infty^0&\stackrel{\sim}{\rightarrow}&G({\Bbb R})/K_\infty Z_\infty^0\\
&(pK_PA_P, aZ_\infty^0)=paK_PZ_\infty^0&\mapsto&paK_\infty Z_\infty^0,\\
\end{array}\leqno(4)
$$
where $p\in {^0}P({\Bbb R})=P({\Bbb R})/A_P$, 
%
%In the opposite direction the map is given by  $g\mapsto(uK_PA_P,aZ_infty^0)$, where  $g=uaK_\infty Z_\infty^0\in U({\Bbb R})A_PK_\infty$]]
%
%
%[[Beachte (WICHTIG): in der Iwasawazerlegung von $g=pk$ von $g\in G({\Bbb R})$ ist $p$ modulo $K_P=\{\pm 1\}$ EINDEUTIG bestimmt ! d.h. $p$ ist fast EINDEUTIG BESTIMMT DURCH $g$ !! D.h. wenn man $g$ durch $p$ ersetzt beim Rechnen in $X=G({\Bbb R})/K_\infty Z_\infty^0$ tut man nichts willkurliches, denn $p$ ist durch $g4 schon festgelegt !]]
%
defines a homeomorphism 
%
%[[Note: 1.) in $e(P)={^0}P({\Bbb R})/K_P$ wird rechts nach $K_P$ faktorisiert; trotzdem 
%ist die obige Abbildung wohldefiniert, denn: $a(K_P)a^{-1}\subset K_P$, da $K_P=\{\pm \Mat{1}{}{}{1}\}$ im Falle $P=P_0$. Fur beliebiges %$P={^\gamma}P_0$ besteht $K_\infty\cap P$ aus diagonalisierbaren Matrizen $g$ mit Eigenwerten in $\{\pm 1\}$ and ${\rm det}\,g=1$ %$\Rightarrow$ Eigenwerte von $g$ sind entweder $(1,1)$ oder $(-1,-1)$ $\Rightarrow g=\pm 1$ ($\in Z_\infty$)
%
%2.) $\varphi$ macht aus $K_\infty$- rechts Nebenklassen $K_P$-rechts Nebenklassen; hier passiert das zum ersten Mal !!]]
%
and we obtain the following diagram
$$
\begin{array}{ccc}
&U_{P,t}=e(P)\times (t,\infty)&\\
\qquad\quad"\subset"\swarrow&&\searrow\varphi\quad \qquad\\
e(P)\times(0,\infty]&&G({\Bbb R})/K_\infty Z_\infty^0.\\
\end{array}
$$
%
%[[Beachte: 1.) Via der Abbildung "$\swarrow$" wird die Anheftung von $e(P)$ definiert, via der Abbildung "$\searrow$" wird die $G({\Bbb Q})$-Aktion auf dem Rand $\bigcup e(P)$ definiert
%2.) (WICHTIG !!!) $varphi$ bildet rechts $K_P$-Nebenklassen auf rechts %$K_\infty$-Nebenklassen ab !!!!!, d.h. aus $K_P$ wird $K_\infty$ !!!]]
%

We identify the face $e(P)$ with $e(P)\times\{\infty\}$ and we define a topological space by gluing $e(P)\times (0,\infty]$ to $G({\Bbb R})/K_\infty Z_\infty^0$ along $U_{P,t}$, where $t\in{\Bbb R}_+$ is chosen sufficiently large. The resulting space is the disjoint union $G({\Bbb R})/K_\infty Z_\infty^0\cup e(P)$. 
Repeating this gluing process for all $P\in{\cal P}$, where ${\cal P}$ is the set of all proper, parabolic subgroups $P\le G$, which are defined over ${\Bbb Q}$, we finally obtain the Borel-Serre compactification
$$
\overline{G({\Bbb R})/K_\infty Z_\infty^0}=G({\Bbb R})/K_\infty Z_\infty^0\dot{\cup}\dot{\bigcup}_{P\in{\cal P}} e(P).
$$

%[[ Sei $x=pK_PA_P\in e(P)$. Was sind die typischen konvergenten Folgen $x_n\in X$, die gegen ein $x\in e(P)\subset \partial \overline{X}$ konvergieren ? Antwort: die typische Folge, die gegen $x$ konvergiert ist 
%$x_n=pa_nK_\infty Z_\infty^0$ mit $a_n\in A_P/Z_\infty^0$, $a_n\rightarrow \infty$, denn: $pa_nK_\infty Z_\infty^0=\varphi(pK_PA_P,a_nZ_\infty^0)$ und $(pK_PA_P,a_nZ_\infty^0)\rightarrow pK_PA_P\in e(P)\times\{\infty\}=e(P)$. Da wir Elemente in $X$ via $varphi$ identifizieren (d.h. $x$ und $\varphi(x)$ sind in $\overline{X}$ dasselbe Element !), folgt $pa_nK_\infty Z_\infty^0\rightarrow pK_PA_P=x$.
%
%REGEL: wollen wir die Seite $e(P)$ untersuchen via ihrer Anheftung an $X$ untersuchen (d.h. als Limes von Elementen in $X$), dann mussen wir die Elemente $x\in X$ in der Form $paK_\infty Z_\infty=(pK_PA_P,aZ_\infty^0)$ darstellen (insbesondere benutzen wir dann $K_P$-rechts nebenklassen !!). Dadurch kommen wir in die Umgebung $U_{P,t}$ von $e(P)$ mit ihrer Topologie: Formal druckt sich dies im Verkleben uber $U_{P,t}$ aus.]]

We set 
$$
X=G({\Bbb R})/K_\infty Z_\infty^0.
$$ 
The action of $G({\Bbb Q})$ on $X$ extends continuously to the boundary
$\partial \overline{X}=\overline{X}-X$, i.e. for an arbitrary element $x\in\partial\bar{X}$ we set $\gamma x=\lim_{n\rightarrow \infty}\gamma x_n$, where $x_n\in X$ is any sequence such that $\lim_{n\rightarrow \infty}=x$. To make this more explicit, let 
$$
x=pK_PA_P\in e(P)\subset \partial\overline{X},
$$
i.e. $p\in{^0P}({\Bbb R})$, and let $\gamma\in G({\Bbb Q})$. For simplicity, in the remainder of this section we assume $P=P_0$. Then,  
$x=\lim_n x_n$, where $x_n=(pK_PA_P,a_nZ_\infty^0)=pa_nK_\infty Z_\infty^0\in X$ with $a_n\in A_P$ (cf. (4)) such that $a_nZ_\infty^0\rightarrow \infty$ in $A_P/Z_\infty^0$ for $n\rightarrow \infty$ (cf. (2)). Also, write $\gamma =uak$ with $u\in {^\gamma}{^0}P({\Bbb R})$, $a\in A_{{^\gamma}P}$ and $k\in K_\infty$; we obtain
\begin{eqnarray*}
(5)\qquad\gamma x_n =\gamma pa_n K_\infty Z_\infty^0&=&(u akp k^{-1}a^{-1}) \, (a ka_nk^{-1})\, K_\infty Z_\infty^0\\
&\stackrel{\varphi}{\leftrightarrow}&(u akp k^{-1}a^{-1} K_{{^\gamma}P}A_{{^\gamma}P}, a ka_nk^{-1}Z_\infty^0)\\
&\in&{^\gamma {^0P}({\Bbb R})}/K_{^\gamma P}A_{^\gamma P}\times A_{{^\gamma P}}/Z_\infty^0.
\end{eqnarray*}
%
%[[1.) note that we identify points via $\varphi$, i.e. $x$ and $\varphi(x)$ are the same point in $\overline{X}$ !
%2.) note that ${^\gamma}P}={^(ua)^{-1}}({^\gamma}P)$ (da $ua\in {^\gamma}P$)
%$={^{(ua)^{-1}}\gamma}P={^k}P$ i.e. ${^k}P={^\gamma}P$. Thus, the coset $u akp k^{-1}a^{-1} K_{{^\gamma}P}A_{{^\gamma}P}$ is well defined]]
%
Since ${^\gamma}P={^k}P$ 
%
%[[denn: wegen $ua\in{^\gamma P}$ ist auch $(ua)^{-1}\in{^\gamma}P$ woraus folgt: ${^{(ua)^{-1}\gamma}P}={^\gamma P}$ $\Leftrightarrow$ ${^kP}={^\gamma P}$]]
%
we find $uak p(ak)^{-1}\in {^\gamma} {^0}P({\Bbb R})$ and since $\Theta(kA_Pk^{-1})=kA_Pk^{-1}$ and $kA_Pk^{-1}\subset {^\gamma}{P}$ we deduce that $kA_Pk^{-1}=A_{^\gamma P}$, i.e. $a ka_nk^{-1}\in A_{{^\gamma}P}$. 
%
%[[d.h. (T) ist die Iwasawa Zerlegung von $\gamma x_n$ bezugl. ${^\gamma}P$. Beachte: ${^\gamma P}={^k P}$ sowie $\Theta(X)=(X^{-1})^t$ und $\Theta(k)=k$ fur $k\in K_\infty$ sowie $\Theta(A_P)=A_P$. Beachte auch: $a{^0P}({\Bbb R})a^{-1}\subset {^0P}({\Bbb R})$ (einfach direkt nachrechnen mit ${^0P}({\Bbb R})=\Mat{\pm 1}{*}{}{\pm 1}$).]]
%
In different words, the diagram
$$
\begin{array}{ccc}
G({\Bbb R})/K_\infty Z_\infty^0&\stackrel{x\mapsto\gamma x}{\rightarrow}&G({\Bbb R})/K_\infty Z_\infty^0\\
\varphi\uparrow&&\uparrow\varphi\\
{^0P}({\Bbb R})\times A_P&\stackrel{(p,a')\mapsto (uakpk^{-1}a^{-1}, aka'k^{-1})}{\rightarrow}&{^0{^\gamma P}}({\Bbb R})\times A_{^\gamma P}\\
\end{array}\leqno(5')
$$
commutes. In particular, using (3) and taking into account that ${^\gamma P}={^kP}$, we see that 
\begin{eqnarray*}
\ell_{^\gamma P}(\gamma x_n)
&\stackrel{(5)}{=}&\ell_{^\gamma P}(aka_nk^{-1}K_\infty Z_\infty^0)\\
&=&\ell_{^\gamma P}(a)\ell_{^\gamma P}(ka_nk^{-1})\\
&=&\ell_{{^\gamma}P}(a)\ell_{^kP}(ka_nk^{-1})\\
&=&\ell_{{^\gamma}P}(a)\ell_P(a_n)\\
&=&\ell_{{^\gamma}P}(\gamma)\ell_P(x_n).\\
\end{eqnarray*}
Hence, we obtain 
$$
\ell_{{^\gamma}P}(\gamma x_n)=\ell_{{^\gamma}P}(\gamma)\ell_P(x_n).\leqno(6)
$$ 
(6) immediately shows that $\ell_{^\gamma P}(\gamma x_n)\rightarrow \infty$ for $n\rightarrow \infty$, because $\ell_P(x_n)\rightarrow \infty$. Hence, $\gamma x=\lim_n \gamma x_n\in e({^\gamma}P)$, i.e. the action of $\gamma\in G({\Bbb Q})$ leaves the boundary invariant: it maps the face $e(P)$ to the face $e({^\gamma P})$. More precisely, passing to the limit for $n\rightarrow \infty$ in (5) we obtain 
\begin{eqnarray*}
(7)\qquad\gamma x=\lim_n \gamma x_n&=&\lim_n (u akp k^{-1}a^{-1} K_{^\gamma P}A_{^\gamma P}, a ka_nk^{-1}Z_\infty^0)\\
&=&(uakp k^{-1}a^{-1}K_{{^\gamma}P}A_{{^\gamma}P},\infty)\\
&=&uak p k^{-1}a^{-1} K_{{^\gamma}P}A_{{^\gamma}P}\in e({^\gamma P}).\\
\end{eqnarray*}
Thus, we obtain the following: let $\gamma=uak\in{^\gamma{^0 P}({\Bbb R})}A_{^\gamma P}K_\infty$; the map 
$$
\begin{array}{cccc}
f_\gamma:&X&\rightarrow&X\\
&x&\mapsto&\gamma x\\
\end{array}
$$
extends continuously to a map 
$$
f_\gamma:\,\bar{X}\rightarrow\bar{X},
$$
whose restriction to a single face is given as
$$
\begin{array}{ccccc}
\gamma:&e(P)&\rightarrow&e({^\gamma}P)&\\
&x=pK_PA_P&\mapsto&uak pk^{-1}a^{-1}K_{{^\gamma}P}A_{{^\gamma}P}&(p\in{^0P}({\Bbb R})).\\
\end{array}\leqno(8)
$$
%[[0.) Insbesondere ist $\gamma(e(P))=e({^\gamma}P)$ und nicht nur $\gamma(e(P))\subset e({^\gamma}P)$ !

%1.) Beachte: $k K_P k^{-1}=K_{{^k}P}=K_{{^\gamma}P}$ und $kA_Pk^{-1}=A_{{^k}P}=A_{{^\gamma}P}$, d.h. die Abbildung ist wohldefiniert.

%2.) Beachte: die folgenden Elemente entsprechen sich unter den beiden Darstellungen von $e(P)$ (beachte: definiert war $e(P)$ als $P({\Bbb R})/A_PK_P$)
%$$
%\begin{array}{ccc}
%{^0}P({\Bbb R})&\leftrightarrow&P({\Bbb R})/A_PK_P\\
%ua va^{-1}&\mapsto&uav \\
%\end{array}
%$$

%3.) Note: formula (Z) explicitely gives the action of $\gamma\in G({\Bbb Q})$ on elements in the boundary. This formula shows: if $\gamma$ is contained in the unipotent radical of the parabolic $P$ to which the face is associated, then the action of $\gamma$ is simple: it is just multiplication by $\gamma$. In general this is not true !! (If $\gamma=ua$ is contained in $P({\Bbb Q})=U({\Bbb Q})A_P$ the action still is simple: multiply by $u$ and conjugate by $a$. But if there is a $k$-component for $\gamma$, i.e. $\gamma\not\in P({\Bbb Q})$ then the action of $\gamma$ is not easy to describe]]

We note that in order to extend $f_\gamma$ to the face $e(P)$ it would have been sufficient to consider the map $f_\gamma$ on a neighbourhood
$[t,\infty]\times e(P)$ of $e(P)$, i.e. to consider the map $f\gamma:\,[t,\infty]\times e(P)\rightarrow [s,\infty]\times e(P)$ where $t$ is sufficiently large as compared to $s$.
%
%[[so dass $f_\gamma(U_t)\subset U_s$]]

We look closer at the case $\gamma \in P({\Bbb Q})$. Since $P$ equals its own normalizer we see that $\gamma e(P)\subset e(P)$ if and only if $\gamma\in P({\Bbb Q})$. There is a homeomorphism
$$
\begin{array}{ccc}
\{\pm 1\}\times {\Bbb R}&\cong& e(P)\\
(\epsilon,u)&\mapsto&\Mat{\epsilon}{}{}{1} \Mat{1}{u}{}{1} A_PK_P.\\
\end{array}
$$
Notice that $e(P)_\epsilon=\{\Mat{\epsilon}{}{}{1}\Mat{1}{u}{}{1}K_PA_P,\;u\in {\Bbb R}\}$, $\epsilon=\pm 1$, are the connected components of $e(P)$. Let 
$$
x=(\epsilon,u,a)=\Mat{\epsilon}{}{}{1}\Mat{1}{u}{}{1}\in e(P)
$$ 
and let 
$$
\gamma=\Mat{\sigma_1}{}{}{\sigma_2}\Mat{1}{\beta}{}{1}\Mat{\alpha}{}{}{\delta}=\Mat{1}{\frac{\sigma_1}{\sigma_2}}{}{1}\Mat{\sigma_1\alpha}{}{}{\sigma_2\delta}\in P({\Bbb Q})
$$
with $\sigma_1,\sigma_2\in\{\pm 1\}$, $\beta\in{\Bbb Q}$ and $\alpha,\delta\in{\Bbb Q}^*_+$. We then obtain from (8)
$$
\gamma x=\Mat{\sigma_1\epsilon}{}{}{\sigma_2}\Mat{1}{\frac{\alpha u}{\delta}+\frac{\beta}{\epsilon}}{}{1}\quad(\in{^0P}({\Bbb R})),
$$ 
i.e.
$$
\begin{array}{cccc}
\gamma:&e(P)&\rightarrow&e(P)\\
&\Mat{\epsilon}{}{}{1}\Mat{1}{u}{}{1}K_PA_P&\mapsto&\Mat{\sigma_1\epsilon}{}{}{\sigma_2}\Mat{1}{\frac{\alpha u}{\delta}+\frac{\beta}{\epsilon}}{}{1}K_PA_P.\\
\end{array}
$$
%
%[[beachte: die matrizen $\Mat{\epsilon}{}{}{1},\Mat{\epsilon\sigma_1}{}{}{\sigma_2}$ gehoren zu ${^0P}({\Bbb R})$, nicht Zu $A_P$ !!]]
%
We deduce the following:

(9) Let $x\in e(P)$; then, $\gamma x=x$ holds precisely if either 

$\bullet$ $\alpha=\delta$ in which case $\beta=0$, i.e. $\gamma\in Z({\Bbb Q})$ 

or 

$\bullet$ $\alpha\not=\delta$ in which case $\sigma_1=\sigma_2$ (to ensure that $\gamma$ respects the connected components of $e(P)$) and $u=\frac{\alpha u}{\delta}+\frac{\beta}{\epsilon}$. 

In the first case the whole face $e(P)$ is fixed by $\gamma$, whereas in the second case there is only one fixed point $p=(\epsilon,\frac{\beta}{\epsilon(1-\alpha \delta^{-1})})$ in each of the two connected components $e(P)$ corresponding to $\epsilon=\pm 1$. Moreover, in the first case $\gamma$ is central and in the second case $\gamma$ is hyperbolic, thus, $\gamma$ never is unipotent. Note that we call an element $\gamma$ hyperbolic if its characteristic polynomial decomposes $\chi_\gamma=(T-\alpha_1)(T-\alpha_2)$ with $\alpha_1\not=\alpha_2$.
%
%[[i.e. hyperbolic excludes central]]

We recall that equation (6) holds for elements $x_n=(pK_PA_p,a_nZ_\infty^0)\leftrightarrow pa_nK_\infty Z_\infty^0\in X$ with arbitrary $a_n\in A_{P}/Z_\infty^0$ and $p\in{^0P}({\Bbb R})$, i.e. equation (6) holds for arbitrary $x=paK_\infty Z_\infty^0\in X$. Since we assume that  ${^\gamma P}=P$ this yields
$$
\ell_P(\gamma x)=\ell_P(\gamma)\ell_P(x).\leqno(10)
$$

We look at the quotient space. For the remainder of section 2 $\Gamma\le G({\Bbb Q})$ denotes an arbitrary arithmetic subgroup, although later on we shall only need the case $\Gamma=\Gamma_1(Np)$, $(N,p)=1$. Since $\Gamma$ acts on $X$ and on $\bar{X}$ we can form the quotient spaces $\Gamma\setminus X$ and $\Gamma\setminus\bar{X}$. We set $\Gamma_P=\Gamma\cap P({\Bbb Q})$; the map $\varphi$ then induces a map
$$
\varphi:\,\Gamma_P\setminus \big(e(P)\times[t,\infty)\big)\rightarrow\Gamma\setminus X.
$$
Here, $\Gamma_P$ acts on $e(P)={^0P}({\Bbb R})/K_P$ by left multiplication. It is a consequence of reduction theory that for $t$ sufficiently large this map is an injection, i.e. $\varphi$ extends to an injective map
$$
\varphi:\,\Gamma_P\setminus e(P)\times[t,\infty]\hookrightarrow\Gamma\setminus \bar{X}.
$$ 
In particular, we obtain an injection
$$
\varphi:\,\Gamma_P\setminus e(P)\rightarrow\Gamma\setminus \partial\bar{X}
$$
We shall identify the space $\Gamma_P\setminus e(P)$ with its image under $\varphi$. 
%
%[[d.h. ein Ausdruck auf der linken Seite steht in Wahrheit fur den (unter $\psi$) entsprechenden Ausdruck auf der rechten Seite]]
%
Clearly, the boundary is invariant under the action of $\Gamma$
%
%[[$\gamma\in\Gamma$ maps $e(P)$ to $e({^\gamma P})$]]
%
and we obtain
$$
\partial(\Gamma\backslash \overline{X})=\Gamma\backslash \partial\overline{X}=\bigcup_{P\in{\cal P}/\sim_\Gamma} \Gamma_P\backslash e(P).
$$
Here and in the following, "$a\sim_G b$" denotes that $a$ and $b$ are conjugate by an element $\gamma\in G$, i.e. $P\sim_\Gamma Q$ means that $Q={^\gamma P}=\gamma P \gamma^{-1}$ for some $\gamma\in \Gamma$.

The map $f_\alpha$, $\alpha\in G({\Bbb Q})$, also extends continuously to the boundary: we set $\Gamma'=\Gamma\cap \alpha^{-1}\Gamma\alpha$. $f_\alpha$ induces a map
$$
\begin{array}{cccc}
f_\alpha:&\Gamma'\setminus\bar{X}&\rightarrow&\Gamma\setminus\bar{X}\\
&x&\mapsto&\alpha x.\\
\end{array}
$$
For any parabolic subgroup $P\le G({\Bbb Q})$ we set $\Gamma_P=\Gamma\cap P({\Bbb Q})$ and $\Gamma'_P=\Gamma'\cap P({\Bbb Q})(=\Gamma_P\cap \alpha^{-1}\Gamma\alpha)$. The restriction of $f_\alpha$ to a face then reads
$$
f_\alpha:\,\Gamma_P'\setminus e(P)\rightarrow \Gamma_{^\alpha P}\setminus e({^\alpha P})\leqno(11)
$$ 
and it is the unique continuous extension of the map
$$
f_\alpha:\,\Gamma_P'\setminus e(P)\times[t,\infty)\rightarrow \Gamma_{^\alpha P}\setminus e({^\alpha P})\leqno(12)
$$
(note that $\alpha\Gamma_{P}'\alpha^{-1}\subset \Gamma_P$). As a consequence, the map (11) is uniquley determined by the map (12), which has domain of definition resp. target inside the interior of $\Gamma'\setminus\bar{X}$ resp. $\Gamma\setminus\bar{X}$.

{\bf 2.2. The Hecke Correspondence. } We fix an arithmetic subgroup $\Gamma\le G({\Bbb Q})$. Let $\alpha\in G({\Bbb Q})$ and set $\Gamma'=\Gamma\cap \alpha^{-1}\Gamma \alpha$. The pair of maps
$$
\begin{array}{ccccc}
&&\Gamma\backslash X\\
&\stackrel{f_\alpha}{\nearrow}&\\
\Gamma'\backslash X&\\
&\stackrel{f_1}{\searrow}&\\
&&\Gamma\backslash X,\\
\end{array}\leqno(13)
$$
where $f_\alpha:\,\Gamma' x\mapsto \Gamma \alpha x$ and $f_1:\,\Gamma' x\mapsto \Gamma x$, induces a correspondence on $\Gamma\backslash X$ via
$$
\Gamma x\mapsto \sum_{y\in f_\alpha f_1^{-1}(\Gamma x)} y\in C^0(\Gamma\backslash X).
$$

%[[beachte: $C^0(\Gamma\backslash X)\cong F_{\Gamma\backslash X}=$ die freie abelsche Gruppe uber $\Gamma\backslash X$, d.h. jedes $x\in \Gamma\backslash X$ ist ein Basisvektor.  $(C^0(\Gamma\backslash X))=$ Gruppe der 0-kozykel $=$ Linearkombinationen von Punkten aus $\Gamma\backslash X$]]

Since the maps $f_\alpha$ and $f_1$ extend to the boundary, the correspondence (13) too extends continuously to a correspondence on the boundary of $\Gamma\backslash X$:
$$
\begin{array}{ccccc}
&&\Gamma\backslash \partial\overline{X}\\
&\stackrel{f_\alpha}{\nearrow}&\\
\Gamma'\backslash \partial\overline{X}&\\
&\stackrel{f_1}{\searrow}&\\
&&\Gamma\backslash \partial \overline{X}.\\
\end{array}\leqno(14)
$$
Clearly, this correspondence is uniquely determined by the correspondence (13) and an explicit formula for the restriction of the correspondence (14) to the boundary follows from equation (8). 
%[[Again, this yields a correspondence on $\partial\Gamma\backslash \overline{X}$
%$$
%\Gamma x\mapsto f_\alpha f_1^{-1}(\Gamma x).
%$$
%]]

We denote by ${\cal L}_k$ the locally constant sheaf on $\Gamma\backslash X$, which is attached to the representation $L_k$ of $G({\Bbb Q})$. ${\cal L}_k$ extends to a sheaf on the Borel-Serre compactification $\Gamma\backslash \overline{X}$, which we again denote by ${\cal L}_k$. The Hecke correspondence 
$(f_1,f_\alpha):\,\Gamma'\backslash\overline{X} \stackrel{\rightarrow}{\rightarrow} \Gamma\backslash \overline{X}$ induces an endomorphism of the cohomology $H^\bullet(\Gamma\backslash \overline{X},{\cal L}_k)$ as follows. Denote by ${\cal L}'_k$ the locally constant sheaf on $\Gamma'\backslash \overline{X}$ attached to $L_k$. The topological morphism $f_\alpha:\,\Gamma'\backslash X\rightarrow \Gamma\backslash X$ induces a morphism of sheafs
$$
f_\alpha^*:\,{\cal L}_k\rightarrow{\cal L}_k'.
$$ 
Explicitely, if $U\subset \Gamma\backslash \overline{X}$ is an open subset, then
$$
\begin{array}{cccc}
f_\alpha^*:&{\cal L}_k(U)&\rightarrow&{\cal L}_k'(f_\alpha^{-1}(U))\\
&s&\mapsto&s\circ f_\alpha.\\
\end{array}
$$

%[[Beachte: Es ist (Definitionsbereich von $s\circ f_\alpha)=f_\alpha^{-1}(U)$.]]

Fix a decomposition $\Gamma\alpha\Gamma=\dot{\bigcup}_{i=1}^s\Gamma\alpha\gamma_i$, $\gamma_i\in \Gamma$. We define a morphism of sheafs
$$
f_{1*}:\,{\cal L}_k'\rightarrow {\cal L}_k
$$
as follows: let $U\subset \Gamma\backslash \overline{X}$ be arbitrary, then
$$
\begin{array}{cccc}
f_{1*}:&{\cal L}_k'(f_1^{-1}(U))&\rightarrow&{\cal L}_k(U)\\
&s&\mapsto&\{\Gamma x\mapsto\sum_{i=1}^s \gamma_i^{-1}s(\gamma_i x)\}\quad (x\in X).\\
\end{array}
$$

%[[Beachte: 1.) $f_{1*}$ ist {\it nicht} Mitteln uber die Faser $f_1^{-1}(\Gamma x)$.
%2.) $f_\alpha^*$ und $f_{1*}$ hangen beide von $\alpha$ ab !!!]]

$f_{1*}$ and $f_\alpha^*$ induce morphisms in cohomology and we define the Hecke operator $T(\alpha)$ on cohomology as the composition
$$
T(\alpha):\,H^\bullet(\Gamma\backslash \overline{X},{\cal L}_k)\stackrel{f_\alpha^*}{\rightarrow} H^\bullet(\Gamma'\backslash \overline{X},{\cal L}_k')\stackrel{f_{1*}}{\rightarrow} H^\bullet(\Gamma\backslash \overline{X},{\cal L}_k).
$$

{\bf 2.3. Restriction of the Hecke correspondence to a single face. } We say that a correspondence $(c_1,c_2):\,\partial \Gamma'\backslash \overline{X}\stackrel{\rightarrow}{\rightarrow} \partial \Gamma\backslash \overline{X}$ restricts to the face $\Gamma_P\backslash e(P)\subset \partial \Gamma\backslash \overline{X}$ if there is a face $\Gamma'_Q\backslash e(Q)\subset \Gamma'\backslash \partial\overline{X}$ such that 

$$
\Gamma'_Q\backslash e(Q)\subset c_1^{-1}(\Gamma_P\backslash e(P))\leqno(15)
$$
and
$$
c_2(\Gamma'_Q\backslash e(Q))\subset\Gamma_P\backslash e(P).\leqno(16)
$$
We recall that $\Gamma_Q'=\Gamma'\cap Q({\Bbb Q})$. 
%
%[[Im Falle einer Hecke Korrespondenz gilt in (G2) stets $c_2(\Gamma'_Q\backslash e(Q))=\Gamma_P\backslash e(P)$, denn eine Seite wird unter $f_\alpha$ stets surjektiv auf eine andere Seite abgebildet]]
%
We define the restriction of $(c_1,c_2)$ to $\Gamma_P\backslash e(P)$ as the union of the correspondences
$$
\begin{array}{ccc}
&&\Gamma_P\backslash e(P)\\
&\stackrel{c_2}{\nearrow}&\\
\Gamma'_Q\backslash e(Q)&&\\
&\stackrel{c_1}{\searrow}&\\
&&\Gamma_P\backslash e(P),\\
\end{array}
$$
where $\Gamma'_Q\backslash e(Q)$ runs over all faces in $\Gamma'\backslash \partial\overline{X}$ satisfying (15,16). 

%[[Note that the above union of correspondences is disjoint, which means that the faces $\Gamma'_Q\backslash e(Q)$ are disjoint. This is a consequence of the trivial fact that different faces have empty intersection !
%Moreover, note that $\Gamma'\backslash X$ in general has more (inequivalent, i.e. different) boundary components than $\Gamma\backslash X$ !!]]

We specialize to the case of a Hecke correspondence 
$$
(f_1,f_\alpha):\,\Gamma'\backslash X\stackrel{\rightarrow}{\rightarrow}\Gamma\backslash X,\leqno(17)
$$ 
for a fixed $\alpha\in G({\Bbb Q})$. Of course, in this case $\Gamma'=\Gamma\cap\alpha^{-1}\Gamma\alpha$. We have seen that (AE) extends continuously to a correspondence on the boundary
$$
(f_1,f_\alpha):\,\partial\Gamma'\backslash \overline{X}\stackrel{\rightarrow}{\rightarrow}\partial\Gamma\backslash \overline{X}.\leqno(18)
$$ 
We want to determine the restriction of (18) to the face $\Gamma_P\backslash e(P)$ and proceed as follows. First, we note that there are bijections between the following sets of objects:

\begin{description}
\item[(i)] faces $\Gamma'_Q\backslash e(Q)$ in $\Gamma'\backslash \partial\overline{X}$ satisfying (15,16),

\item[(ii)] $\Gamma'$-conjugacy classes of parabolic subgroups $Q\le G$ satisfying
$$
Q\sim_\Gamma P\leqno(15')
$$
$$
\alpha Q\alpha^{-1}\sim_\Gamma P,\leqno(16')
$$

\item [(iii)] double cosets $\Gamma_P\beta\Gamma_P\subset \Gamma\alpha\Gamma\cap P({\Bbb Q})$ (note that necessarily $\beta\in P({\Bbb Q})$).
\end{description}

The bijection $(ii)\rightarrow(i)$ is given by
$$
Q\mapsto \Gamma'_Q\backslash e(Q).
$$
Furthermore, let $\Gamma_P\beta\Gamma_P\subset \Gamma\alpha\Gamma\cap P({\Bbb Q})$ (cf. (iii)). $\beta$ then may be written as  $\beta=\gamma_2\alpha\gamma_1$ and the bijection $(iii)\rightarrow(ii)$ is given by  
$$
\Gamma_P\beta\Gamma_P\mapsto \gamma_1 P\gamma_1^{-1}
$$
(cf. [G-MacPh], 7.4 Lemma). We write
$$
\Gamma\alpha\Gamma\cap P({\Bbb Q})=\bigcup_{i=1}^m \Gamma_P\alpha_i\Gamma_P\leqno(19)
$$
with $\alpha_i\in P({\Bbb Q})$ (note that $\Gamma\alpha\Gamma\cap P({\Bbb Q})$ is $\Gamma_P$-invarinat from the left and from the right, hence, (19) is the $\Gamma_P\times\Gamma_P$-orbit decomposition of $\Gamma\alpha\Gamma\cap P({\Bbb Q})$ and that $\Gamma\alpha\Gamma\cap P({\Bbb Q})$ is contained in $P({\Bbb Q})$, hence, $\alpha_i\in P({\Bbb Q})$). Moreover, $\alpha_i$ may be written as $\alpha_i=\gamma_{i,2}\alpha\gamma_{i,1}$, where $\gamma_{i,1},\gamma_{i,2}\in \Gamma$, and we set 
$$
Q_i=\gamma_{i,1}P\gamma_{i,1}^{-1}.
$$
%[[Dass $\alpha_i\in P({\Bbb Q})$ steht in [G-M], 7.3 Proposition; ausserdem: wenn $\gamma$ $\Gamma\alpha\Gamma\cap P({\Bbb Q})$ normalisiert, dann normalisiert $\gamma$ auch $P({\Bbb Q})$, also ist $\gamma\in P({\Bbb Q})$ $\Rightarrow \gamma\in \Gamma\cap P({\Bbb Q})=\Gamma_P$.]]
The above bijections then show that $Q_1,\ldots,Q_m$ are precisely a set of representatives of the $\Gamma'$-conjugacy classes of parabolic subgroups satisfying (15') and (16'). The restriction of $(UE)$ to the face $\Gamma_P\backslash e(P)$ therefore is the (disjoint) union of the correspondences
$$
\begin{array}{ccc}
&&\Gamma_P\backslash e(P)\\
&\stackrel{f_\alpha}{\nearrow}&\\
\Gamma'_{Q_i}\backslash e(Q_i)&&\\
&\stackrel{f_1}{\searrow}&\\
&&\Gamma_P\backslash e(P),\\
\end{array}\leqno(20)
$$
where $i=1,\ldots,m$. 

%[[Disjoint union because: $\bullet$ $e(P)\cap e(Q)=e(P)=e(Q)$ if $P\sim_\Gamma Q$, $e(P)\cap e(Q)=\emptyset$ if $P\not\sim_\Gamma Q$ $\bullet$ $f_\alpha:\,e(Q)\rightarrow e(P)$ ist surjektiv
%]]

%[[Note that if $(f_1,f_\alpha)$ does not restrict to $\Gamma_P\backslash e(P)$, then $\Gamma\alpha\Gamma\cap P({\Bbb Q})=\emptyset$ (i.e. the restriction of the correspondence to the boundary is the empty correspondence, sending any $x$ to $\emptyset$) ]]

We want to analyze more closely the correspondence (20). Quite generally we say that two correspondences $(f_{\alpha_1},f_{\alpha_2}):\,\Gamma'\setminus X\stackrel{\rightarrow}{\rightarrow}\Gamma\setminus X$ and $(f_{\beta_1},f_{\beta_2}):\,\Gamma''\setminus X\stackrel{\rightarrow}{\rightarrow}\Gamma\setminus X$ are isomorphic if there is $\gamma\in G({\Bbb Q})$ such that the diagram
$$
\begin{array}{cccc}
(f_{\alpha_1},f_{\alpha_2}):&\Gamma'\backslash X&\stackrel{\rightarrow}{\rightarrow}&\Gamma\backslash X\\
&&&\\
&f\,\uparrow&&\parallel\\
&&&\\
(f_{\beta_1},f_{\beta_2}):&\Gamma''\backslash X&\stackrel{\rightarrow}{\rightarrow}&\Gamma\backslash X\\
\end{array}\leqno(21)
$$
commutes, i.e. $f_{\beta_i}\circ f=f_{\alpha_1}$, $i=1,2$. Here, $\Gamma'=\Gamma\cap \alpha^{-1}\Gamma\alpha$ and $\Gamma''=\Gamma\cap \beta^{-1}\Gamma\beta$ and $f$ denotes the map which sends $x\mapsto\gamma x$.  We obtained the correspondence (20) on the boundary $\partial\Gamma\setminus\bar{X}$ as the continuous extension of a Hecke correspondence on $X$. There is a second type of correspondences 
$$
\begin{array}{ccccc}
&&\Gamma_P\backslash  e({P})\\
&\stackrel{f_\alpha}{\nearrow}&\\
\Gamma_P'\backslash e(P)&\\
&\stackrel{f_1}{\searrow}&\\
&&\Gamma_P\backslash e(P),\\
\end{array}\leqno(22)
$$
defined on the boundary: for any $\alpha\in P({\Bbb Q})$ these are obtained as the continuous extension of the correspondence defined by the pair of maps
$$
\begin{array}{ccccc}
&&\Gamma_P\backslash ([s,\infty)\times e({P}))\\
&\stackrel{f_\alpha}{\nearrow}&\\
\Gamma_P'\backslash ([t,\infty)\times e(P))&\\
&\stackrel{f_1}{\searrow}&\\
&&\Gamma_P\backslash ([s,\infty)\times e(P)),\\
\end{array}\leqno(23)
$$
where $f_\alpha:\,\Gamma_P' x\mapsto \Gamma_P \alpha x$ and $f_1:\,\Gamma_P' x\mapsto \Gamma_P x$ (note that ${^\alpha P}=P$). 
The correspondence (20) then is isomorphic to the correspondence
$$
\begin{array}{ccccc}
&&\Gamma_P\backslash e({P})\\
&\stackrel{f_{\alpha_i}}{\nearrow}&\\
\Gamma_{P,i}'\backslash  e(P)&\\
&\stackrel{f_1}{\searrow}&\\
&&\Gamma_P\backslash e(P),\\
\end{array}\leqno(24)
$$
where $\Gamma_{P,i}'=\Gamma_P\cap \alpha_i^{-1}P\alpha_i$ (cf. [G-MacPh], section 7.6).

{\bf 2.4. The trace formula. } We will apply the result from the previous section to determine the contribution of the boundary to the Lefschetz number of the Hecke correspondence. We fix an arithmetic subgroup $\Gamma\le G({\Bbb Q})$ and an element $\alpha\in G({\Bbb Q})$. As in section 2.3 we write $\Gamma\alpha\Gamma\cap P({\Bbb Q})=\dot{\bigcup}_{i=1,\ldots,m}\Gamma_P\alpha_i\Gamma_P$ with $\alpha_i\in P({\Bbb Q})$ and we also write $\alpha_i=\gamma_{i,2}\alpha\gamma_{i,1}$, $\gamma_{i,1},\gamma_{i,2}\in\Gamma$ and set $Q_i=\gamma_{i,1}P\gamma_{i,1}^{-1}$.

We denote by 
$$
F_\partial=\{x\in\partial \Gamma'\backslash \overline{X}: f_1(x)=f_\alpha(x)\}
$$  
the set of fixed points of $(f_1,f_\alpha)$ on the boundary. Furthermore, for any $P\in{\cal P}/\sim_\Gamma$ we set 
$$
F_{\partial,P}=\{x\in F_\partial:\,f_1(x)=f_\alpha(x)\in \Gamma_P\backslash e(P)\},
$$ 
i.e. $F_\partial$ is the disjoint union
$$
F_\partial=\bigcup_{P\in{\cal P}/\sim_\Gamma} F_{\partial,P}.
$$
The above description of the restriction of the Hecke correspondence to the face $e(P)$ implies that
$$
F_{\partial,P}=\dot{\bigcup}_{i=1,\ldots,m} F_{Q_i},
$$
where 
$$
F_{Q_i}=\{x\in\Gamma'_{Q_i}\backslash e(Q_i):\,f_\alpha(x)=f_1(x)\}
$$ 
is the set of fixed points of the correspondence $(20)$. 

We denote by 
$$
F_{P,i}=\{x\in \Gamma'_{P,i}\backslash e(P):\,f_1(x)=f_{\alpha_i}(x)\},
$$ 
$i=1,\ldots,m$, the set of fixed points of the correspondence $(24)$. Since the correspondences (24) and (20) are isomorphic the definition of isomorphism of correspondences shows that the fixed point sets of $(20)$ and $(24)$ are isomorphic (under the map $f:\,x\mapsto\gamma x$ defining the isomorphism of (20) and (24); cf. equation (21)). Hence, for our purposes we may replace $F_{Q_i}$ by $F_{P,i}$. 
We look closer at $F_{P,i}$ and set 
$$
X_{P,i}=\{x\in e(P): \eta x=x \;\mbox{for some}\;\eta\in\Gamma_P\alpha_i\}.
$$ 
Moreover, for any $\eta\in\Gamma_P\alpha_i$ we put 
$$
X_{P,\eta}=\{x\in e(P): \eta x=x\}\leqno(25)
$$
and we define $F_{P,\eta}=\pi'(X_{P,\eta})$. Here, $\pi':\,e(P)\rightarrow \Gamma'_{P,i}\backslash e(P)$ is the canonical projection. $\Gamma'_{P,i}x\in F_{P,i}$ is equivalent to the existence of a $\gamma\in \Gamma_P$ such that $\gamma\alpha_ix=x$, hence $x\in X_{P,\eta}$ with $\eta=\gamma\alpha_i\in\Gamma_P\alpha_i$. Thus, $\Gamma'_{P,i}x\in F_{P,\eta}$ and we deduce that $F_{P,i}$ is the union of the $F_{P,\eta}$, $\eta\in\Gamma_P\alpha_i$. Moreover, if $F_{P,\eta}\cap F_{P,\mu}\not=\emptyset$, $\eta,\mu\in \Gamma_P\alpha_i$, then there are $x\in X_{P,\eta}$ and $y\in X_{P,\mu}$ such that $\gamma x=y$ for some $\gamma\in \Gamma'_{P,i}$. 
%
%[[$\Gamma'_{P,i}x=\Gamma'_{P,i}y\Leftrightarrow \gamma x=y$ for some $\gamma\in\Gamma'_{P,i}$]] 
%
By using the definition of $X_{P,\mu}$ we obtain $\mu\gamma x=\gamma x$, i.e. $\gamma^{-1}\mu\gamma x=x$. Since $\eta x=x$ and $\Gamma$ is assumed to be torsion free we deduce that $\eta=\gamma^{-1}\mu\gamma$. If, on the other hand, $\eta=\gamma^{-1}\mu\gamma$, $\gamma\in \Gamma'_{P,i}$ then, for any $x\in X_{P,\eta}$ we find $\mu\gamma x=\gamma\gamma^{-1}\mu\gamma x=\gamma \eta x=\gamma x$, hence, $\gamma x\in X_{P,\mu}$. This implies that $\gamma X_{P,\eta}\subset X_{P,\mu}$ and, by symmetry, $\gamma X_{P,\eta}= X_{P,\mu}$. Thus $F_{P,\eta}\cap F_{P,\mu}\not=\emptyset\Leftrightarrow $ $\eta$ and $\mu$ are $\Gamma'_{P,i}$-conjugate and in this case we even have $F_{P,\eta}=F_{P,\mu}$.
%
%[[beachte: wegen $\gamma\in\Gamma'_{P,i}$ ist $\gamma \Gamma'_{P,i}\gamma\subset \Gamma'_{P,i}$, d.h. $\Gamma'_{P,i}\backslash e(P)\rightarrow \Gamma'_{P,i}\backslash e(P)$, $\Gamma'_{P,i}x\mapsto \Gamma'_{P,i}\gamma x$ ist wohldefiniert]]
%
Hence, altogether we have seen that
$$
F_{P,i}=\dot{\bigcup}_{\eta\in\Gamma_{P}\alpha_i/\sim_{\Gamma_{P,i}'}} F_{P,\eta}.
$$

Taking into account that $\gamma\alpha_i{\Gamma_{P,i}'}\mapsto \gamma\alpha_i{\Gamma_P}$ defines a bijection $\Gamma_P\alpha_i/\sim_{\Gamma_{P,i}'}\stackrel{1-1}{\leftrightarrow} \Gamma_P\alpha_i\Gamma_P/\sim_{\Gamma_P}$ (cf. [B], Lemma 3.3 (iii), p. 46), we finally obtain
\begin{eqnarray*}
(26)\qquad F_\partial&=&\bigcup_{P\in{\cal P}/\sim_\Gamma} \bigcup_{i=1}^m \bigcup_{\eta\in \Gamma_P\alpha_i\Gamma_P/\sim_{\Gamma_P}} F_{P,\eta}\\
&=&\bigcup_{P\in{\cal P}/\sim_\Gamma} \bigcup_{\eta\in (\Gamma\alpha\Gamma\cap P({\Bbb Q}))/\sim_{\Gamma_P}} F_{P,\eta}.\\
\end{eqnarray*}
Here, the second equality is a consequence of (19). We note that (19) implies that $\eta\in \Gamma\alpha\Gamma\cap P({\Bbb Q})$ is contained in a uniquely determine double coset 
$\eta\in \Gamma_P\alpha_i\Gamma_P/\sim_{\Gamma_P}=\Gamma_P\alpha_i/\sim_{\Gamma_{P,i}'}$  
for a some $i\in\{1,\ldots,m\}$. Thus, $\eta$ uniquely determines the index $i\in\{1,\ldots,m\}$ and, hence, the set  $F_{P,\eta}$, which is the set of fixed points of the correspondence (24). 

We say that the face $\Gamma_P\backslash e(P)$ is attracting resp. equidistant resp. repelling if the correspondence (23) is attracting resp. equidistant resp. repelling, i.e. if $\ell_P(f_{\alpha_i}(x))>\ell_P(f_1(x))$ resp. $\ell_P(f_{\alpha_i}(x))=\ell_P(f_1(x))$ resp. $\ell_P(f_{\alpha_i}(x))<\ell_P(f_1(x))$ for some $x\in \Gamma'_{P,i}\backslash U_{P,t}$. Since $\alpha_i\in P({\Bbb Q})$, equation (F1) shows that $\Gamma_P\backslash e(P)$ is attracting resp. equidistant resp. repelling precisely if $\ell_P(\alpha_i)>1$ resp. $\ell_P(\alpha_i)=1$ resp. $\ell_P(\alpha_i)<1$. 

{\bf Remark.} 1.) In particular, we see that if $\ell_P(f_{\alpha_i}(x))?\ell_P(f_1(x))$, $?=">,=,<"$, holds for one $x\in \Gamma'_{P,i}\backslash U_{P,t}$ then it holds for all $x\in \Gamma'_{P,i}\backslash U_{P,t}$ (because then $\ell_P(\alpha_i)?1)$. 

2.) The attracting property of the face $\Gamma'_{P,i}\backslash e(P)$ is well defined, i.e. it only depends on the $\Gamma'$-conjugacy class of $P$. To see this, let $\gamma\in \Gamma'$; then, conjugating with $\gamma$ and taking into account that $\Gamma'\subset \Gamma$ we find that $\Gamma_p\alpha_i\Gamma_p\subset \Gamma\alpha\Gamma\cap P({\Bbb Q})$ is equivalent to $\Gamma_{{^\gamma}P}{^\gamma}\alpha_i\Gamma_{{^\gamma}P}\subset\Gamma\alpha\Gamma\cap {^\gamma}P({\Bbb Q})$.
%
%[[Beachte: ${^\gamma}\Gamma_P=\Gamma_{{^\gamma}P}$ und $\Gamma\alpha\Gamma=\Gamma\gamma\alpha\gamma^{-1}\Gamma$, da $\gamma\in \Gamma$]]
%
Thus, the assignment $\Gamma_P\alpha_i\Gamma_P\mapsto\Gamma_{^\gamma P}{^\gamma}\alpha_i\Gamma_{^\gamma P}$ defines a bijection between $\Gamma_P$-double cosets appearing in $\Gamma\alpha\Gamma\cap P({\Bbb Q})$ and $\Gamma_{^\gamma P}$-double cosets appearing in $\Gamma\alpha\Gamma\cap {^\gamma}P({\Bbb Q})$. Let $x=pK_PZ_\infty^0\in \Gamma'_{P,i}\backslash U_{P,t}$ be arbitrary; using the commutativity of the diagram
$$
\begin{array}{ccc}
\Gamma'_{^{\gamma}P,i}\backslash U_{{^\gamma}P,t}&\stackrel{f_{\gamma\alpha_i\gamma^{-1}}}{\rightarrow}&\Gamma_{^{\gamma\alpha_i}P}\backslash U_{^{\gamma\alpha_i}P,s}\\
&&\\
\uparrow\;\gamma&&\uparrow\;\gamma\\
&&\\
\Gamma'_{P,i}\backslash U_{P,t}&\stackrel{f_{\alpha_i}}{\rightarrow}&\Gamma_P\backslash U_{P,s}\\
\end{array}\leqno(27)
$$
%
%[[Es gilt ${^\gamma}\Gamma'_{P,i}={^\gamma}(\Gamma\cap \alpha_i^{-1}\Gamma\alpha_i\cap P({\Bbb Q})=\Gamma\cap \gamma\alpha_i^{-1}\Gamma\alpha_i\gamma^{-1}\cap {^\gamma}P({\Bbb Q})$ und $\Gamma'_{{^\gamma}P,i}=\Gamma\cap 
%\gamma\alpha_i\gamma^{-1}\Gamma \gamma^{-1}\alpha_i^{-1}\gamma \cap {^\gamma}P({\Bbb Q})$ $\Rightarrow {^\gamma}\Gamma'_{P,i}=\Gamma'_{^{\gamma}P,i}$. D.h. die linke vertikale Abbildung ist wohldefiniert
%]]
%
we may then compute (note that ${^{\gamma\alpha_i}}P={^\gamma}P$, because $\alpha_i\in P({\Bbb Q})$)
$$
\ell_{{^\gamma}P}(f_{\gamma\alpha_i\gamma^{-1}}(\gamma x))=\ell_{^{\gamma} P}(\gamma \alpha_i p)
\stackrel{(6)}{=}\ell_{^{\gamma} P}(\gamma)\ell_{P}(\alpha_i p).
$$
Similarly, upon setting $\alpha_i=1$ in (27), we find $\ell_{^\gamma P}(f_1(x))=\ell_{{^\gamma}P}(\gamma x)\stackrel{(6)}{=}\ell_{{^\gamma}P}(\gamma)\ell_P(x)$. Thus,
$$
\frac{\ell_{^\gamma P}(f_{\gamma\alpha_i\gamma^{-1}}(\gamma x))} {\ell_{{^\gamma}P}(f_1(\gamma x))}=\frac{\ell_P(f_{\alpha_i} x)} {\ell_P(f_1(x))}.\leqno(28)
$$
Since $x\in \Gamma'_{P,i}\backslash U_{P,t}$ and $\gamma x\in \Gamma'_{{^\gamma} P,i}\backslash U_{{^\gamma}P,t}$, equation (28) implies that $\Gamma'_{P,i}\backslash U_{P,t}$
is attracting (repelling, equidistant) precisely if $\Gamma'_{{^\gamma P},i}\backslash U_{{^\gamma}P,t}$ is attracting (repelling, equidistant).

We note that only if $\Gamma'_{P,i}\backslash e(P)$ is attracting or equidistant, the correspondence (23) (or, more precisely, the correspondence (24)) will contribute to the trace formula. In particular, a fixed point component $F_{P,\eta}\subset \Gamma'_{P,i}\backslash e(P)$, $\eta\in \Gamma_P\alpha_i$, only contributes if $e(P)$ is attracting or equidistant. Since $f_{\alpha_i}(\Gamma'_{P,i}x)=\Gamma_P\eta x$ for $x\in F_{P,\eta}$ and $\ell_P(\eta x)=\ell_P(\eta)\ell_P(x)$ because $\eta\in P({\Bbb Q})$ (cf. equation (10)), this is equivalent to $\ell_P(\eta)> 1$ resp. $\ell_P(\eta)= 1$. We say that the pair $(P,\eta)$ is attracting, resp. repelling resp. equidistant if $\alpha_P(\eta)>1$ resp. $\alpha_P(\eta)<1$ resp. $\alpha_P(\eta)=1$, where $\alpha_P$ is the weight of the action of $S_P$ on $U_P({\Bbb R})$. We define $(\Gamma\alpha\Gamma\cap P({\Bbb Q}))^+$ as the set of those $\eta\in \Gamma\alpha\Gamma\cap P({\Bbb Q})$ such that $(P,\eta)$ is attracting or equidistant. Using (26) and following the arguments in the proof of [B], Lemma 2.5, p. 39 ff, we find that the contribution of the boundary $\partial\Gamma\backslash \bar{X}$ to the Lefschetz number of the Hecke corspondence on $H^\bullet({\Gamma\backslash X,\cal L}_{k,{\Bbb C}})$ is given as
$$
L((f_1,f_\alpha),H^\bullet(\partial \Gamma\backslash\bar{X},{\cal L}_{k,{\Bbb C}}))=\sum_{P\in{\cal P}/\sim_\Gamma}\sum_{\eta\in(\Gamma\alpha\Gamma\cap P({\Bbb Q}))^+/\sim_{\Gamma_P}}{\rm tr}(\eta^\iota|L_k)\chi(F_{P,\eta}).
$$
Together with the contribuition from the inner terms, which is given in [B], Satz 5.3, p. 65 (the index "ell" is missing in [B], 5.3 Satz) we therefore obtain the

{\bf Trace formula. }

\begin{eqnarray*}
(29)&&L(T(\alpha)|H^\bullet(\Gamma\backslash \overline{X},{\cal L}_{k,{\Bbb C}}))\\
&&\\
&&=\sum_{\xi\in(\Gamma\alpha\Gamma)_{\rm ell.}/\sim_{\Gamma}}{\rm tr}(\xi^\iota|L_k)\chi(F_\xi)+\sum_{P\in{\cal P}/\sim_\Gamma}\sum_{\eta\in(\Gamma\alpha\Gamma\cap P({\Bbb Q}))^+/\sim_{\Gamma_P}}{\rm tr}(\eta^\iota|L_k)\chi(F_{P,\eta}).\\
\end{eqnarray*}

Here, $F_\xi$ is the fixed point component corresponding to $\xi\in (\Gamma\alpha\Gamma)_{\rm ell.}$, i.e. $F_\xi=\pi(X_\xi)$, where $X_\xi=\{x\in X:\,\xi x=x\}$ and $\pi:\,X\rightarrow \Gamma\backslash X$ is the canonical projection.  

%[[Beachte: Da wir die Aktion des Heckeoperators $\Gamma\alpha\Gamma$ auf der GruppenKohomologie durch $T\omega=\alpha_i^\iota\omega(\gamma_i)$ definiert war (cf. [Hida, El. Theory of Eis.series..], p. 176) und nicht durch $T\omega=\alpha_i^{-1}\omega(\gamma_i)$ (cf. [B], (1.2.2), (1.2.3), p. 12 wo $\alpha^{-1}$ und $\gamma_i^{-1}$ in auftauchen; s. auch [Kuga-Parray-Sah]) erscheint in der Spurformel ${\rm tr}\,(\xi^\iota|L_k)$ und nicht ${\rm tr}\,(\xi^{-1}|L_k)$ wie in [B], (2.6) Satz, (5.3) Satz
%
%2.) Wenn wir die Aktion der Hecke algebra durch $T\omega=\alpha_i^\iota\omega(\gamma_i)$ definieren dann wird der Eichler-Shimura Isomorphismus Hecke equivariant (cf. [Hida, El. theory of Eis.series...], p. 177), d.h. ${\cal M}_k(\Gamma)$ und $H^1(\Gamma,L_k)$ sind isomorphe Heckemoduln. (Mit der Definition $T\omega=\alpha_i^{-1}\omega(\gamma_i)$ muss die Aktion des heckeoperators noch twisten, damit Eichler-Shimura Isomorphismus Hecke equivariant wird !)]]
%
%3.) Die Defintion der Hecke Aktion auf Gruppenkohomologie via $T\omega=\alpha_i^\iota\omega(\gamma_i)$ entspricht der Hecke Aktion auf Garbenkohomologie wie wir sie oben definiert haben (vgl. [Hida, El. theory of Eis.series], p. 181- 183)
%]]

{\bf 2.5. The boundary term. } We want to rewrite the boundary term of the trace formula. We proceed in two steps. 

1.) First we note that for $\eta\in\Gamma_P\alpha_i$
$$
X_{P,\eta}=
\left\{
\begin{array}{cl}
\mbox{a point}&\mbox{if $\eta$ is hyperbolic}\\
\emptyset&\mbox{if $\eta$ is unipotent}\\
e(P)&\mbox{if $\eta$ is central.}\\
\end{array}
\right.
$$
(cf. (9) and the definition of $X_{P,\eta}$ in (25)). In particular $\chi(F_{P,\eta})=1$ if $\eta$ is hyperbolic and $\chi(F_{P,\eta})=0$ if $\eta$ is unipotent or central.
%
%[[beachte: ist $\eta$ zentral, dann ist $F_{P,\eta}=\Gamma'_{P,i}\backslash e(P)\cong m{\Bbb Z}\backslash {\Bbb R}\cong S^1$ und $\chi(S^1)=0$.]]
%
We say that a fixed point component $F_{P,\eta}$ is hyperbolic if $\eta$ is hyperbolic. Since $F_{P,\eta}= F_{P,\eta'}$ 
%[[$\Leftrightarrow F_{P,\eta}\cap F_{P,\eta'}\not=\emptyset \Leftrightarrow$ 
is equivalent to $\eta$ and $\eta'$ being $\Gamma'_{P,i}$-conjugate, this notion is well defined. The above together with (26) then shows that only 
$$
F_{\partial,{\rm hyp.}}=\bigcup_{P\in{\cal P}/\sim_\Gamma}\bigcup_{\eta\in(\Gamma\alpha\Gamma\cap P({\Bbb Q}))_{\rm hyp.}/\sim_{\Gamma_P}} F_{P,\eta}
$$
contributes to the boundary term of the trace formula.

%[[$M_{\rm hyp}:=$ Menge der hyperbolischen Elemente in $M$]]

2.) To proceed further, let $\eta\in G({\Bbb Q})$ be a ${\Bbb Q}$-hyperbolic element. We denote by ${\Bbb Q}[\eta]\subset M_2({\Bbb Q})$ the ${\Bbb Q}$-subalgebra of the algebra of $2\times 2$ matrices generated by $\eta$ and we set $C_\eta={\Bbb Q}[\eta]^*\subset G({\Bbb Q})$. Since $\eta$ is hyperbolic, which in particular implies that $\eta$ is regular (i.e. non-central), $C_\eta$ is a maximal ${\Bbb Q}$-split torus in $G({\Bbb Q})$. We denote by $W_\eta={\cal N}_{G({\Bbb Q})}(C_\eta)/C_\eta$ the Weyl group of $C_\eta$ and we select any $w \in {\cal N}_{G({\Bbb Q})}(C_\eta)$ such that $w\not=1$ in $\in W_\eta$. There are precisely two different parabolic subgroups $P,P'\le G$ containing $\eta$ and they are conjugate to each other by $w_\eta$:  $P'={^{w_\eta}}P$. 
%
%[[Beweis fur diese Behauptung: Es ex. $\tau$ mit ${^\tau}C_\eta=C_0$, wobei $C_0$ der Torus der Diagonalmatrizen ist. Seien $P_0,P_0^o$ die beiden parabolischen untergruppen, die $C_0$ enthalten, d.h. $P_0, P_0^0$ sind die Gruppe der oberen bzw. unteren Dreiecksmatrizen und es gilt $P_0^0={^{w_0}}P_0$. Insbesondere gilt dann ${^\tau }P$ enthalt $C_0$, also ist  ${^\tau}P=P_0$ oder $=P_O^0$, d.h. es gilt: Es gibt nur zwei parabolische Untergruppen, die $C_\eta$ enthalten, namlich $P={^{\tau^{-1}}}P_0$ und $P'={^{\tau^{-1}}}P_O^0={^{\tau^{-1}w_0}}P_0$. Weiter: ist $\sigma\in{\rm Aut}(C_\eta)$ $\Rightarrow \tau\sigma\tau^{-1}\in {\rm Aut}(C_0)$, d.h. $W_eta={^{\tau^{-1}}}W_0$ und $w_\eta=\tau^{-1} w_0\tau$. Damit folgt $P'={^{\tau^{-1}w_0}}P_0={^{\tau^{-1}w_0}\tau} {^{\tau^{-1}}}P_0={^{\tau^{-1}w_0}\tau} {^{\tau^{-1}}}P_0={^{w_\eta}}P$. Fazit: Es gibt nur zwei parabolische Untergruppen, die $\eta$ enthalten und diese sind unter $w_\eta$ konjugiert. 
%
%Erganzung: Beweis dafur dass $C_\eta$ oder $C_0$ nur in zwei parabolischen untergruppen enthalten ist (Beweis auf Liealgebrenniveau): die minimalen parabolischen Unteralgebren, die die Cartanunteralgebra $Lie(C_\mu)$ enthalten sind von der Form $Lie(C_\mu)\oplus\bigoplus_{\alpha>0}g_\alpha$, d.h. sie entsprechen den Auswahlen 
%einer Basis des Wurzelsystems, die aber sind alle unter der Weylgruppe konjugiert. 
%]]
%
Moreover, for any $x\in X$ we have
$$
\ell_{P'}(x)=\ell_P(x)^{-1}.\leqno(30)
$$
This can easily be seen by writing ${^\tau}C_\eta=C_0$, where $\tau\in G({\Bbb Q})$ and $C_0$ is the torus of diagonal matrices in $G({\Bbb Q})$. Hence ${^\tau}P=P_0$, ${^\tau}P'=P_0^{\rm opp}$ and the claim follows taking into account that $\ell_{{^\tau}P}({^\tau}x)=\ell_P(x)$ (cf. equation (3')) and $\ell_{P^{\rm opp}}(x)=\ell_P(x)^{-1}$. Moreover we fix any system of representatives $P_1,\ldots,P_k$ for ${\cal P}/\sim_\Gamma$ and $P\in{\cal P}/\sim_\Gamma$ will mean that $P\in\{P_1,\ldots,P_k\}$.
%
%[[note that ${\rm SL}_2({\Bbb Z})$ acts transitively on ${\cal P}$, da alle Spitzen ${\Bbb Q}\cup\{\infty\}$ der oberen Halbebene unter ${\rm SL}_2({\Bbb Z})$ aquivalent sind]]
%
We set 
$$
{\cal M}=\{(P,[\eta]_{\Gamma_P}):\,P\in{\cal P}/\sim_\Gamma,\,[\eta]_{\Gamma_P}\in (\Gamma\alpha\Gamma\cap P({\Bbb Q}))_{\rm hyp.}/\sim_{\Gamma_P}\}.
$$

{\bf Lemma. }{\it The map
$$
\begin{array}{cccc}
\psi:&{\cal M}&\rightarrow&(\Gamma\alpha\Gamma)_{\rm hyp.}/\Gamma\\ 
&(P,[\eta]_{\Gamma_P})&\mapsto&[\eta]_\Gamma\\
\end{array}
$$
is surjective. Moreover, the following holds: if $P\not\sim_\Gamma P'$ then the fibre $\psi^{-1}([\eta]_\Gamma)=\{(P,[\eta]_{\Gamma_P}),(Q,[\mu]_{\Gamma_Q})\}$ consists of two elements such that $(P,\eta)$ is repelling if $(Q,\mu)$ is attracting and $(P,\eta)$ is attracting if $(Q,\mu)$ is repelling. If $P\sim_\Gamma P'$ then $\psi^{-1}([\eta]_\Gamma)=\{(P,[\eta]_{\Gamma_P})\}$ consists of a single element. 
}

{\it Proof. } Obviously, $\psi$ is well defined. We prove surjectivity. If $[\eta]_\Gamma\in(\Gamma\alpha\Gamma)_{\rm hyp.}/\Gamma$ then $\eta$ is contained in a parabolic subgroup $P({\Bbb Q})\le G({\Bbb Q})$, hence, after conjugating $\eta$ by some $\gamma\in \Gamma$ we may assume that $\eta\in P({\Bbb Q})$ for some $P\in{\cal P}/\sim_\Gamma$. 
In particular, $\eta\in (\Gamma\alpha\Gamma\cap P({\Bbb Q}))_{\rm hyp.}$ and $(P,[\eta]_{\Gamma_P})$ obviously maps to $[\eta]_{\Gamma}$, i.e. $\psi$ is surjective.

It remains to examine the fibres of $\psi$. To this end we fix $(P,[\eta]_{\Gamma_P})\in\psi^{-1}([\eta]_\Gamma)$ and we let $(Q,[\mu]_{\Gamma_Q})\in{\cal M}$ ($Q\in{\cal P}/\sim_\Gamma$, $\mu\in (\Gamma\alpha\Gamma\cap Q({\Bbb Q}))_{\rm hyp.}$) be any element, which is contained in the same fibre $\psi^{-1}([\eta]_\Gamma)$ as $(P,[\eta]_{\Gamma_P})$. Hence, $[\mu]_\Gamma=[\eta]_\Gamma$ and there is $\gamma\in \Gamma$ such that ${^\gamma}\eta=\mu$. 

We set $C_\eta={\Bbb Q}[\eta]^*$ and $C_\mu={\Bbb Q}[\mu]^*$. $C_\eta\le P({\Bbb Q})$ and $C_\mu\le Q({\Bbb Q})$ are ${\Bbb Q}$-split tori, which satisfy $\gamma C_\eta\gamma^{-1}=C_\mu$. We denote by $P$ and $P'$ resp. $Q$ and $Q'$ the two parabolic subgroups containing $C_\eta$ resp. $C_\mu$. Since ${^\gamma}P$ contains $C_\mu$ we deduce that ${^\gamma}P=Q$ or ${^\gamma}P=Q'$. We now distinguish cases.

Case A: $P'\not\sim_\Gamma P$. We distinguish further:

If ${^\gamma}P=Q$ we obtain $P=Q$, because $P,Q$ are contained in ${\cal P}/\sim_\Gamma$. Hence, $P={^\gamma P}$ is invariant under conjugation by $\gamma$. Since $P({\Bbb Q})$ equals its own normalizer we obtain $\gamma\in \Gamma\cap P({\Bbb Q})=\Gamma_P$. Thus, $[\mu]_{\Gamma_P}=[\eta]_{\Gamma_P}$, i.e. $(P,[\eta]_{\Gamma_P})=(Q,[\mu]_{\Gamma_Q})$. We conclude that $\psi^{-1}([\eta]_\Gamma)=\{(P,[\eta]_{\Gamma_P})\}$.

If ${^\gamma}P=Q'$ we immediately deduce $Q={^\gamma}P'$. Moreover, we have ${^\gamma \eta}=\mu\in Q({\Bbb Q})={^\gamma P'}({\Bbb Q})$. Thus, 
${^\gamma \eta}\in {^\gamma P'}({\Bbb Q})\cap\Gamma\alpha\Gamma$ and $({^\gamma}P',[{^\gamma}\eta]_{{^\gamma}P'})$ therefore is a well defined element of ${\cal M}$ (note that ${^\gamma}P'=Q\in{\cal P}/\sim_\Gamma$) and it clearly holds that $(Q,\mu)=({^\gamma}P',{^\gamma}\eta)$. (Note that the condition ${^\gamma}P'=Q$ uniquely determines the coset $\gamma\Gamma_{P'}$,
% 
%[[Sei auch ${^\tau}P'=Q$ $\Rightarrow {^{\tau^{-1}\gamma}}P'=P' \Rightarrow \tau^{-1}\gamma\in \Gamma_{P'} \Rightarrow \gamma\Gamma_{P'}=\tau\Gamma_{P'}$]]
%
hence, the element $({^\gamma}P',[{^\gamma}\eta]_{{^\gamma}P'})$ is uniquely determined, i.e. it does not depend on the choice of a representative in $\gamma\Gamma_{P'}$.) Thus, we have seen that necessarily $(Q,\mu)=({^\gamma}P',{^\gamma}\eta)$ and $\psi^{-1}([\eta]_\Gamma)$ therefore consists of at most two elements. We will show that the fibre consists of $2$ elements. To this end we select $\tau\in\Gamma$ such that ${^\tau}P'=P_i$ for some $P_i\in{\cal P}/\sim_\Gamma$. Since $\eta\in P'({\Bbb Q})$ we obtain ${^\tau}\eta\in {^\tau}P'({\Bbb Q})$, hence, ${^\tau}\eta\in {^\tau}P'({\Bbb Q})\cap \Gamma\alpha\Gamma$. Thus, $(P,[\eta]_{\Gamma_P})$ and 
$({^\tau}P',[{^\tau}\eta]_{{^\tau}\Gamma_{P'}})$ are elements of ${\cal M}$, which are both contained in the fibre of $[\eta]_{\Gamma}$ and which are different, because, by assumption, $P'\not\sim_\Gamma P$. Hence, 
$$
\psi^{-1}([\eta]_\Gamma)=\{(P,[\eta]_{\Gamma_P}),({^\tau}P',[{^\tau}\eta]_{{^\tau}\Gamma_{P'}})\}.
$$
%
%[[Does it hold that $[\mu]_{\Gamma_P}$ and
%$[{^{\tau^{-1}}}\mu]_{\Gamma_P}$ aredifferent, i.e. $\mu\not\sim_{\Gamma_P}{^{\tau^{-1}}}\mu$. To this end we assiume that 
%${^{\tau^{-1}}}\mu={^\sigma\mu}$ for some $\sigma\in \Gamma_P$. We obtain ${^{\tau\sigma}}\mu=\mu$, hence, $\tau\sigma\in
%{\cal C}_{G({\Bbb Q})}(C_\mu)=C_\mu$ (note that $\mu$ is hyperbolic, i.e. non-central and therefore regular). Thus, 
%$\tau\sigma\in C_\mu\cap \Gamma=\{\pm 1\}$.
%(Fur $Q=\{\Mat{*}{*}{}{*}\}$ ist das klar, der allgemeine Fall folgt wieder nach Konjugation; beachte, dass $\pm 1$ zentral ist, also mit allen Matrizen kommutiert. Benutze eventuell, dass $P^{\rm opp}={^w}P$. Here, we denote by $W_\eta={\cal N}_{G({\Bbb Q})}(C_\eta)/C_\eta$ the Weyl group of $C_\eta$. We let $w$ be any element in $N_{G({\Bbb Q})}(C_\eta)-C_\eta$, i.e. $w$ is a representative of the non-trivial element of $W_\eta$.)
%) and we find $\tau=\pm \sigma\in \Gamma_P$. Mit anderen Worten: $\Gamma$ operiert via Konjugation {\it treu} auf $\mu$ !]]
%
Finally, we calculate
$$
\ell_{{^\tau}P'}({^\tau}\eta)\stackrel{(3')}{=}\ell_{P'}(\eta)\stackrel{(30)}{=}\ell_{P}(\eta)^{-1}.
$$
%
%[[Alternativer Beweis fur $(Kon)$: Ist $\eta=\Mat{a}{b}{}{d}k$, $k\in K_\infty$, dann folgt ${^\gamma}\eta=\gamma \Mat{a}{b}{}{d}\gamma^{-1} \gamma k\gamma^{-1}$. Weiter: schreibe ${^\gamma}k=\Mat{a'}{b'}{}{d'}k'$, $k'\in K_\infty$. Dann gilt $\ell_P(\Mat{a'}{b'}{}{d'})=1$, denn ${^\gamma}k\in {^\gamma}K_\infty$ und ${^\gamma}K_\infty$ ist kompakt $\Rightarrow \ell_P({^\gamma}k)=1$. Insgesamt folgt $\ell_{{^\gamma}P}({^\gamma}x)  
%=\ell_{{^\gamma}P}(\gamma\Mat{a}{b}{}{d}\gamma^{-1}\Mat{a'}{b'}{}{d'}k')\stackrel{???}{=}a/d=\ell_P(x)^{-1}$.
%
%2.) Ist $x=\Mat{a}{b}{}{d}k$, $k\in K_\infty$, dann folgt $x=w\Mat{a}{b}{}{d}w^{-1} wkw^{-1}$ $\Rightarrow \ell_P(x)=\ell_{P^{\rm opp}}(x)$ (beachte $P^{\rm opp}={^w}P$ und $wkw^{-1}\in K_\infty)$. Beachte auch: damit 
%$\ell_P(x)$ und $\ell_{P^{\rm opp}}(x)$ definiert sind muss nicht $x\in P\cap P^{\rm opp}$ gelten !]]
%
In particular, if $(P,\eta)$ is attracting then $({^\tau}P',[{^\tau}\eta]_{{^\tau}\Gamma_{P'}})$ is repelling and if $(P,\eta)$ is repelling then $({^\tau}P',[{^\tau}\eta]_{{^\tau}\Gamma_{P'}})$ is attracting. This completes the proof of the Lemma in case A.

Case B: $P\sim_\Gamma P'$. If ${^\gamma}P=Q$ we deduce as above that $P=Q$ and hence $(P,[\eta]_{\Gamma_P})=(Q,[\mu]_{\Gamma_P}])$. If ${^\gamma}P=Q'$ then ${^\gamma}P'=Q$ 
%
%[[klar fur $\Mat{*}{*}{}{*}$, allgemeiner Fall wieder durch Konjugieren. ]]
%
and we find that $P\sim_\Gamma P'\sim_\Gamma Q$, i.e. again we have $P=Q$ and therefore  $(P,[\eta]_{\Gamma_P})=(Q,[\mu]_{\Gamma_P}])$. Thus, $\psi^{-1}([\eta]_\Gamma)=(P,[\eta]_{\Gamma_P})$. This completes the proof of the lemma in case B and, hence, the proof of the lemma is complete.

%[[Beachte: $\Gamma$ oder $\Gamma_P$-Konjugation erhalt Hyperbolizitat]]

%[[Anmerkung: $[\eta]_\Gamma$ zerfalllt bei Einschrankung auf $\Gamma_P$ in sehr viele verschiedene $\Gamma_P$-Konjugationsklassen. Warum liegt dann in $\psi^{-1}([\eta]_\Gamma)$ nur eine bzw. zwei $\Gamma_P$ Konjugationsklassen $[\mu]_{\Gamma_P}$? der Grund ist, dass $\mu\in P$ fur ein $P\in{\cal P}/\sim_\Gamma=\{P_1,\ldots,P_k\}$b sein muss ! fur beliebiges $\gamma\in \Gamma$ ist ${^\gamma}\eta\in{^\gamma}P\not\in{\cal P}/\sim_\Gamma$. Nur fur ein $\gamma\in \Gamma$ kommt fur ein beliebiges $P$ ${^\gamma}P$ nach ${\cal P}/\sim_\Gamma$ hinein. Allerdings erhalten wir auch noch etwas wenn ${^\gamma}P^{\rm opp}$
%nach ${\cal P}/\sim_\Gamma$ hinein kommt, das gibt dann eventuell ein zweites Element.]] 

The Lemma immediately implies that
\begin{eqnarray*}
&&\bigcup_{P\in{\cal P}/\sim_\Gamma} \bigcup_{\eta\in(\Gamma\alpha\Gamma\cap P({\Bbb Q}))_{\rm hyp}/\sim_{\Gamma_P}} F_{P,\eta}\\
&&=\bigcup_{[\eta]_\Gamma\in(\Gamma\alpha\Gamma_{\rm hyp})/\sim_{\Gamma}} 
\bigcup_{(P,[\mu]_{\Gamma_P})\in\psi^{-1}([\eta]_\Gamma)} F_{P,\eta}\\
&&=\bigcup_{[\eta]_\Gamma\in(\Gamma\alpha\Gamma_{\rm hyp})/\sim_{\Gamma}\atop \psi^{-1}([\eta]_\Gamma)=\{(P,[\eta]_{\Gamma_P})\}} F_{P,\eta}
\cup \bigcup_{[\eta]_\Gamma\in(\Gamma\alpha\Gamma_{\rm hyp})/\sim_{\Gamma}\atop \psi^{-1}([\eta]_\Gamma)=\{(P,[\eta]_{\Gamma_P}),(P',[\eta']_{\Gamma_{P'}})\}} F_{P,\eta}\cup F_{P',\eta'}.\\
\end{eqnarray*}
Here, according to the above Lemma, $(P,\eta)$ is attracting if $(P',\eta')$ is repelling and vice versa. 

%[[In particular, there is no $[\eta]_\Gamma\in(\Gamma\alpha\Gamma_{\rm hyp})\Gamma$, which is attracting for more than one parabolic $P\in{\cal P}/\sim_\Gamma$.]]

We define $(\Gamma\alpha\Gamma_{\rm hyp})^+$ as the set of those $\eta\in\Gamma\alpha\Gamma_{\rm hyp}$ for which there is a parabolic subgroup
$P\in{\cal P}$ such that $\eta\in P({\Bbb Q})$ and $(P,\eta)$ is attracting. Notice that, if $(P,\eta)$ is equidistant, $\eta$ is either central or unipotent, hence the fixed point component $F_{P,\eta}$ does not contribute to the Lefschetz number of the Hecke correspondence. Hence, only attracting fixed point sets contribute to the Lefschetz number and we finally see that the contribution of the boundary to the trace formula is given by
$$
\sum_{[\eta]_\Gamma\in(\Gamma\alpha\Gamma)_{\rm hyp}^+/\sim_{\Gamma}} {\rm tr}(\eta^\iota|L_k).
$$

Since the inner fixed point components $F_\xi$ are connected and the action of $G({\Bbb Q})$ on $X$ only has isolated fixed points we obtain  $\chi(F_\xi)=1$ for all $\xi\in (\Gamma\alpha\Gamma)_{\rm ell}$. Thus, the trace formula finally becomes

{\bf Theorem. }{\it
$$
L(T(\alpha)|H^\bullet(\Gamma\backslash X,{\cal L}_{k,{\Bbb C}}))
=\sum_{\xi\in(\Gamma\alpha\Gamma)_{\rm ell.}/\sim_\Gamma} {\rm tr}(\xi^\iota|L_k)
+\sum_{\eta\in(\Gamma\alpha\Gamma)_{\rm hyp.}^+/\sim_\Gamma}{\rm tr}(\eta^\iota|L_k).
$$

}

\newpage

\section{Hecke algebra and the basic trace identity.}

{\bf 3.1. The Hecke algebra. } Throughout section 3, we fix a tame level $N$ and a prime $p\in{\Bbb N}$ satisfying $p>2$ and $(p,N)=1$ and we write $\Gamma=\Gamma_1(Np)$ and $\Gamma_0=\Gamma_0(Np)$. We set 
$$
\Delta_0=\{\alpha=\Mat{a}{b}{c}{d}\in M_2({\Bbb Z }):\,c\equiv 0\pmod{Np},\,(a,Np)=1,\,\det \alpha>0\}
$$ 
and 
$$
\Delta_1=\{\alpha=\Mat{a}{b}{c}{d}\in M_2({\Bbb Z }):\,c\equiv 0\pmod{Np},\,a\equiv 1\pmod{Np},\,\det \alpha>0\}.
$$ 
We define the Hecke algebras $\tilde{\cal H}=\Gamma\backslash \Delta_0/\Gamma$, ${\cal H}_1=\Gamma\backslash \Delta_1/\Gamma$ and ${\cal H}_0=\Gamma_0\backslash \Delta_0/\Gamma_0$.
We look closer at these Hecke algebras and their interrelation. We define elements in $\tilde{\cal H}$ as follows

$\bullet$ $T(l,m)=\Gamma\Mat{l}{}{}{m}\Gamma$, where $l,m\in{\Bbb Z}$ such that $l|m,\, (l,Np)=1$

$\bullet$ $\langle\epsilon\rangle=\Gamma\sigma_\epsilon\Gamma=\Gamma\sigma_{\epsilon}=\sigma_\epsilon\Gamma$, 
where $\epsilon\in{\Bbb Z}$ is such that $(\epsilon,Np)=1$ and 
$$
\sigma_\epsilon=\Mat{a}{b}{c}{d}\in\Gamma_0(Np)
$$ 
is any matrix satisfying $a\equiv \epsilon, c\equiv 0, d\equiv \epsilon^{-1}\pmod{Np}$ (note that $\sigma_\epsilon$ normalizes $\Gamma_0$ and $\Gamma$). 

%[[Da $(\epsilon,Np)=1$ existieren $a,b\in{\Bbb N}$ mit $a\epsilon+bNp=1$; dann ist $\Mat{\epsilon}{-b}{Np}{a}\in\Gamma_0(Np)$ eine Matrix wie fur $\sigma_\epsilon$ gefordert.]]

$\bullet$ For any Dirichlet character $\chi:\,({\Bbb Z}/(Np))^*\rightarrow{\Bbb C}^*$ we set
$$
e_\chi=\frac{1}{\varphi(Np)}\,\sum_{\epsilon\in({\Bbb Z}/Np{\Bbb Z})^*} \chi(\epsilon)\langle\epsilon\rangle \in\tilde{\cal H}_{{\Bbb Z}[Np]}.
$$

%[[Beachte: In der Definition von $e_\chi$ kommt $\chi$, nicht $\chi^{-1}$ vor. Grund: die Definition von $\sigma_\epsilon$: der "a"-Eintrag von $\sigma_\epsilon$ ist kongruent zu $\epsilon$, nicht der "d" Eintrag.]]

Here, for an arbitrary ${\Bbb Z}$-algebra $E$ we have set $\tilde{\cal H}_E=\tilde{\cal H}\otimes E$. Note that $T(l,m)$ is not contained in ${\cal H}_1$. Moreover, we define the Hecke operators
$$
T(l,m)_0=\Gamma_0\Mat{l}{}{}{m}\Gamma_0, \quad l,m\in{\Bbb Z},\; l|m,\,(l,Np)=1,
$$
which are contained in ${\cal H}_0$.

We note the following facts. 

$\bullet$ ${\cal H}_0=\{T(l,m):\;l|m,\,(l,Np)=1\}$ (cf. [M], Lemma 4.5.2, p. 132) and ${\cal H}_0$ is generated by the operators $T(1,\ell)_0$, $\ell$ any prime, and $T(\ell,\ell)_0$ with $(\ell,Np)=1$.

$\bullet$ There is a canonical isomorphism of ${\Bbb Z}$-algebras
$$
\begin{array}{ccc}
{\cal H}_1&\cong&{\cal H}_0\\
\Gamma\alpha\Gamma&\mapsto&\Gamma_0\alpha\Gamma_0\\
\Gamma\sigma_a^{-1}\alpha\Gamma&\leftarrow&\Gamma_0\alpha\Gamma_0.
\end{array}
$$
Here, $a$ denotes the upper left entry of $\alpha\in\Delta_0$ (cf. [M], Theorem 4.5.18, p. 151).

$\bullet$ We deduce that $\tilde{\cal H}=\{\sigma_aT(l,m):\;(a,Np)=1,\,l|m,\,(l,Np)=1\}$ and $\tilde{\cal H}$ is generated by the Hecke operators $\langle\epsilon\rangle$ with $(\epsilon,Np)=1$, $T(1,\ell)$ with $\ell$ any prime and $T(\ell,\ell)$ with $(\ell,Np)=1$.

$\bullet$ Since the operators $T(\ell,\ell)$ commute with any Hecke operator and since $\langle\epsilon\rangle$ commutes with $T(1,\ell)$ (cf. [D-S], p. 169), we see that $\tilde{\cal H}$ is commutative. In particular,  $\{e_\chi\}_\chi$ is a complete set of pairwise orthogonal (central) idempotents in $\tilde{\cal H}$ and we obtain a decomposition of the Hecke algebra
$$
\tilde{\cal H}_{{\Bbb Z}[Np]}=\bigoplus_\chi e_\chi\tilde{\cal H}_{{\Bbb Z}[Np]}.
$$
We note that $e_{1}\tilde{\cal H}_{{\Bbb Z}[Np]}={{\cal H}_0}_{{\Bbb Z}[Np]}$.

From now on we set ${\cal H}={\cal H}_1$ and we shall identify ${\cal H}_1$ with ${\cal H}_0$, i.e. a Hecke operator $\Gamma_0\alpha\Gamma_0\in{\cal H}_0$ with $\alpha\in\Delta_0$ also stands for the operator $\Gamma\sigma_a^{-1}\alpha\Gamma$ in ${\cal H}$ ($a$ denotes the upper left entry of $\alpha$). In particular, $\Gamma_0 \Mat{1}{}{}{\ell} \Gamma_0$ resp. $\Gamma_0 \Mat{\ell}{}{}{\ell} \Gamma_0$ also denotes the Hecke operator $T(1,\ell)$ resp. $\sigma_\ell^{-1}T(\ell,\ell)$. We set 
$$
T_n=T(1,n),
$$ 
$n\in{\Bbb N}$. Note that $T_\ell=T(1,\ell)=T(\ell)$, while in general $T_n$ differs from the Hecke operator $T(n)=\sum_{ab=n\atop a|b}T(a,b)$.

%[[Wir werden nur ${\cal H}_1$ verwenden, nicht ${\cal H}_0$, d.h. wir sind bestrebt Elemente in ${\cal H}$ als Elemente in ${\cal H}_1$ darzustellen, d.h. also in der Form $\Gamma\alpha\Gamma$, wobei $a\equiv 1\pmod{Np}$. Grund: ${\cal H}_1$ operiert naturlich auf ${\cal M}_k(\Gamma,\chi)$ (Elemente in ${\cal M}_k(\Gamma,\chi)$ sind nur $\Gamma$-invariant, nicht $\Gamma_0$ invariant).]]

%[[We note that ${\cal H}$ is commutative.]] 

%[[Wir haben also $T_n$, $T(1,n)$, $T(n)$; die ersten beiden sind identisch]]

{\bf Remark. } The following relation is an immediate consequence of [M], Lemma 4.5.7 (1):
$$
T_\ell^m=T(1,\ell^m)+f_{m-2}(T(\ell,\ell))T(1,\ell^{m-2})+f_{m-4}(T(\ell,\ell))T(1,\ell^{m-4})+\cdots,\leqno(1)
$$
where $f_i\in{\Bbb Z}[X]$ (the last summand appearing on the right hand side either is $f_1(T(\ell,\ell))T(1,\ell)$ or $f_0(T(\ell,\ell))$) and all $f_i=0$ in case $\ell|Np$. 
%
%[[benutze vollstandige Induktion nach $m$ ]]
%
Since $T(1,m)T(1,n)=T(1,mn)$ if $(m,n)=1$, 
%
%[[cf. [M], Lemma 4.5.8, p. 141]]
%
using (x) we find that
$$
T_{\ell_1}^{r_1}\cdot\ldots\cdot T_{\ell_s}^{r_s}=\sum_{h=(h_1,\ldots,h_s)} g_h T_{\ell_1^{h_1}\cdot\ldots\cdot\ell_{s}^{h_s}},\leqno(2)
$$
where the sum runs over all $h$ such that $h_i\equiv r_i\pmod{2}$ and $0\le h_i\le r_i$ and $g_h\in{\Bbb Z}[T(\ell_1,\ell_1),\ldots,T(\ell_s,\ell_s)]$. Since $T(a,a)T(b,b)=T(ab,ab)$ we deduce that $g_h$ is a ${\Bbb Z}$-linear combination of operators $T(n',n')$, where $n'|\ell_1^\infty\cdot\ldots\cdot \ell_s^\infty$ and $(n',Np)=1$. In different words, $T_{\ell_1}^{r_1}\cdot\ldots\cdot T_{\ell_s}^{r_s}$ is a ${\Bbb Z}$-linear combination, which obviously does not depend on $k$, of the operators $T(n',n')T_n$, where $n$ runs over (certain) divisors of $\ell_1^{r_1}\cdot\ldots\cdot\ell_s^{r_s}$ and $n'$ runs over (certain) divisors of $\ell_1^\infty\cdot\ldots\cdot\ell_s^\infty$ satisfying $(n',Np)=1$. 

This reduces the computation of ${\rm tr}\, T_{\ell_1}^{r_1}\cdot\ldots\cdot T_{\ell_s}^{r_s}|_{{\cal M}_k(\Gamma,\chi)^\alpha}$ to the computation of ${\rm tr}\, T(n',n')\, T_n|_{{\cal M}_k(\Gamma,\chi)^\alpha}$, where $n$ runs over divisors of $\ell_1^{r_1}\cdot\ldots\cdot\ell_s^{r_s}$ and $n'$ runs over divisors of $\ell_1^\infty\cdot\ldots\cdot\ell_s^\infty$ satisfying $(n',Np)=1$(note that $T_1={\rm id}$, i.e. the Hecke Operator $T(n',n')$ is included).

%[[Beachte: Alle bei uns vorkommenden Hecke operatoren kommutieren, d.h. es gibt keine Probleme mit Polynomen (Ist $F(X,Y)$ ein Polynom, dann ist $F(T_1,T_2)$ nur wohldefiniert wenn $T_1,T_2$ kommutieren).]]

{\bf 3.2. An extension of the Hecke algebra. }

{\bf Lemma. }{\it Let $N\in{\Bbb N}$ be any integer with prime decomposition $N=p_1^{n_1}\cdot\ldots\cdot n_s^{p_s}$ and let $M$ be any divisor of $N$ with prime decomposition $m=p_1^{m_1}\cdot\ldots\cdot p_s^{m_s}$. Then $\Gamma_1(M)$ is the (not necessarily disjoint) union of double cosets
$$
\Gamma_1(M)=\bigcup_{0\le m\le N \atop M|m}\bigcup_{\bar{a}\in({\Bbb Z}/(N))^*\atop a\equiv 1\pmod{M}} \Gamma_1(N)\Mat{1}{}{m}{1}\sigma_a\Gamma_1(N).\leqno(3)
$$
Here, $\sigma_a$ is defined in 3.1, i.e. $\sigma_a\in\Gamma_0(N)$ is any matrix whose upper left entry is congruent to $a$ modulo $N$. 

%[[i.e. $m$ runs over a system of representatives of ${\Bbb Z}/(N)$]]
}

{\it Proof. } We denote by $\Gamma(N)=\{\Mat{a}{b}{c}{d}: \,a\equiv b\equiv c\equiv d\equiv1\pmod{N} \}$ the principal congruence subgroup of level $N$. Since $\Gamma(N)$ is contained in both sides of (3) equation (3) is equivalent to
$$
\Gamma_1(M)/\Gamma(N)=\bigcup_{0\le m\le N\atop\,M|m}\bigcup_{\bar{a}\in({\Bbb Z}/(N))^*\atop a\equiv 1\pmod{M}} \Gamma_1(N)\Mat{1}{}{m}{1}\sigma_a\Gamma_1(N)/\Gamma(N).\leqno(4)
$$
We apply the isomorphism 
$$
\begin{array}{ccccc}
i:&{\rm SL}_2({\Bbb Z})/\Gamma(N)&\cong&{\rm SL}_2({\Bbb Z}/(N))&\\
&\Mat{a}{b}{c}{d}&\mapsto&\Mat{\bar{a}}{\bar{b}}{\bar{c}}{\bar{d}}&(\bar{x}=x+(N))\\
\end{array}
$$
to both sides of (4), i.e. we reduce the entries of all matrices appearing in (4) modulo $(N)$: since $\Gamma_1(N)/\Gamma(N)\stackrel{i}{\cong} N_2({\Bbb Z}/(N))$, where $N_2\le {\rm GL}_2$ is the subgroup of unipotent upper triangular matrices and $\Gamma_1(M)/\Gamma(N)\stackrel{i}{\cong}\Gamma_1(M;{\Bbb Z}/(N))$, where
$$
\Gamma_1(M;{\Bbb Z}/(N))=\{\Mat{\bar{a}}{\bar{b}}{\bar{c}}{\bar{d}}\in{\rm SL}_2({\Bbb Z}/(N)),\,a\equiv d\equiv 1, c\equiv 0\pmod{M}\},
$$
equation (4) becomes equivalent to 
$$
\Gamma_1(M;{\Bbb Z}/(N))=\bigcup_{\bar{m}\in{\Bbb Z}/(N)\atop M|m}\bigcup_{\bar{a}\in({\Bbb Z}/(N))^*\atop a\equiv 1\pmod{M}}N_2({\Bbb Z}/(N))\Mat{1}{}{\bar{m}}{1}\sigma_{\bar {a}} N_2({\Bbb Z}/(N)).\leqno(5)
$$
where $\sigma_{\bar{a}}=i(\sigma_a)$.

The Chinese remainder theorem yields an isomorphism $j:\,{\Bbb Z}/(N)\cong \prod_i{\Bbb Z}/(p_i^{n_i}),\,x\mapsto(x\pmod{p_i^{n_i}})_i$ under which 

$\bullet$ the ideal $(M)=\{\bar{m}\in{\Bbb Z}/(N):\,M|m\}$ maps to $\prod_i (p_i^{m_i})$; here $(p_i^{m_i})=\{\bar{x}\in{\Bbb Z}/(p_i^{n_i}):\,p_i^{m_i}|x\}$

%[[klar: $j((M))\subset\prod_i (p_i^{m_i})$; ausserdem: ist $j(x)\in \prod_i (p_i^{m_i})$ dann folg $p_i^{m_i}|x\Rightarrow M|x$ ]]

$\bullet$ $\Gamma_1(M;{\Bbb Z}/(N))$ maps to $\prod_i \Gamma_1(p_i^{m_i};{\Bbb Z}/(p_i^{n_i}))$.

Hence, (5) is equivalent to showing that
$$
\Gamma_1(p_i^{m_i};{\Bbb Z}/(p_i^{n_i}))=\bigcup_{\bar{x}\in{\Bbb Z}/(p_i^{n_i})\atop p_i^{m_i}|x} \bigcup_{\bar{a}\in({\Bbb Z}/(p_i^{n_i}))^*\atop a\equiv 1\pmod{p_i^{m_i}}}N_2({\Bbb Z}/(p_i^{n_i}))\Mat{1}{}{\bar{x}}{1}\sigma_{\bar{a}}N_2({\Bbb Z}/(p_i^{n_i}))\leqno(6)
$$
for all all primes $p_1,\ldots,p_s$. Straightforward computation shows that the right hand side of (6) is contained in the left hand side. To prove the reverse inclusion let $\Mat{a}{b}{c}{d}\in{\rm SL}_2({\Bbb Z})$ represent an arbitrary element of $\Gamma_1(p_i^{m_i};{\Bbb Z}/(p_i^{n_i}))$. Hence, $v_{p_i}(c)\ge m_i$ and $a\equiv 1\pmod{p_i^{m_i}}$. We distinguish cases

Case 1. $p\not|a$, i.e. $\bar{a}\in{\Bbb Z}/(p_i^{n_i})^*$. The decomposition
$$
\Mat{a}{b}{c}{d}=\Mat{1}{}{c/a}{1}\sigma_a  \left(\sigma_{a^{-1}}\Mat{a}{b}{}{-\frac{bc}{a}+d}\right)\leqno(7)
$$
then shows that 
$$
\Mat{\bar{a}}{\bar{b}}{\bar{c}}{\bar{d}}\in\Mat{1}{}{\overline{c/a}}{1}\sigma_{\bar {a}} N_2({\Bbb Z}/(p_i^{n_i})).
$$
Here, the bar denotes reduction modulo $p_i^{n_i}$, i.e. $\bar{x}=x+(p_i^{n_i})$
%
%[[Beachte: die letzte Matrix auf der rechten Seite hat Determinante $=1$]]
%
Since $v_p(c)\ge m_i$, i.e. $c/a\in(p_i^{m_i})$, and $a\equiv 1\pmod{p_i^{m_i}}$ this is contained in the right hand side of (6). 

Case 2: $p|a$, i.e. $\bar{a}\not\in{\Bbb Z}/(p_i^{n_i})^*$. Since $a\equiv 1\pmod{p_i^{m_i}}$ this case can only occur if $m_i=0$. We obtain $p\not|c$, i.e. $\bar{c}\in({\Bbb Z}/(p_i^{n_i}))^*$ and $a+c\in ({\Bbb Z}/(p_i^{n_i}))^*$. We look at the decomposition
$$
\Mat{a}{b}{c}{d}=\Mat{1}{-1}{}{1}\Mat{1}{}{\frac{c}{a+c}}{1}\sigma_{a+c}\left(\sigma_{(a+c)^{-1}}\Mat{a+c}{b+d}{}{-\frac{(b+d)c}{a+c}+d}\right).\leqno(8)
$$
Since $m_i=0$ and $v_p(c/(a+c))=0$, the above decomposition again shows that $\Mat{a}{b}{c}{d}\pmod{p_i^{m_i}}$ is contained in the right hand side of (6) (note that the condition "$a\equiv 1\pmod{p_i^{m_i}}$" is empty because $m_i=0$). 
%Finally, since the right hand side of (3) contains the principal congruence subgroup $\Gamma(N)$ it is sufficient if the indices $m$ and $a$ only run over cosets modulo $N$. 
Thus, the proof of the Lemma is complete. 

%[[Dividiert man die rechte seite von (3) durch $\Gamma(N)$, dann laufen alle Matrixeintrage modulo $N$]]

We return to the case of relevance to us. We let $N\in{\Bbb N}$ and $p\in{\Bbb N}$ be any prime such that $(N,p)=1$ and we let $M\in{\Bbb N}$ be any divisor of $Np$. We continue to set $\Gamma=\Gamma_1(Np)$ and we let the group $\Gamma\times \Gamma$ act on $\Gamma_1(M)$ by
$$
(\alpha,\beta) x=\alpha x\beta^{-1}.
$$
The sets $\Gamma \Mat{1}{}{m}{1}\sigma_a \Gamma$ are orbits under this action, hence, any two of these sets either are equal or disjoint. In particular, there are subsets 
$$
{\bf M}\subset \{m\in{\Bbb N}\cup\{0\}:\,0\le m\le Np,\,M|m\}
$$ 
and 
$$
{\bf A}\subset\{ a\in({\Bbb Z}/(Np))^*,\,a\equiv 1\pmod{M}\}
$$ 
such that
$$
\Gamma_1(M)=\dot{\bigcup}_{m\in{\bf M}}\dot{\bigcup}_{a\in{\bf A}} \Gamma\Mat{1}{}{m}{1}\sigma_a\Gamma\quad(\mbox{disjoint union}).\leqno(9)
$$
We look closer at the special case $p|M$. Let $\gamma=\Mat{a}{b}{c}{d}\in\Gamma_1(M)$; since $p|M$ we deduce that $p|c$ and $a\equiv 1\pmod{p}$. In particular, $p\not|a$ and equation (7) yields
$$
\gamma\in\Gamma \Mat{1}{}{u}{1}\sigma_x \Gamma
$$
with $x\equiv 1\pmod{p}$ and $p|u$. Thus, if $p|M$ we obtain
$$
\Gamma_1(M)=\dot{\bigcup}_{m\in{\bf M}}\dot{\bigcup}_{a\in{\bf A}} \Gamma\Mat{1}{}{m}{1}\sigma_a\Gamma\qquad(\mbox{disjoint union})\leqno(10)
$$
with any $u\in{\bf M}$ divisible by $p$ and ${\bf A}\subset 1+p{\Bbb Z}$.

For any integer $u\in{\Bbb N}$ and any $\delta\in({\Bbb Z}/(Np))^*$ we define Hecke operators 
$$
S_{u,\delta}=\Gamma\Mat{1}{}{u}{1}\sigma_\delta\Gamma.
$$

%[[Beachte $\sigma_\delta$ normalisiert $\Gamma(=\Gamma_1(N))$, aber $\sigma_\delta$ kommutiert nicht mit $\Mat{1}{}{m}{1}$]]

The operators $S_{u,\delta}$ are contained in the larger Hecke algebra $\Gamma\backslash M_2({\Bbb Z})^+/\Gamma$, where $M_2^+({\Bbb Z})$ is the (multiplicative) semigroup of integral matrices with positive determinant. 

%[[1.) Positive Determinante heisst insbesondere die Determinante ist ungleich $0$. 2.) Halbgruppe versteht sich bezuglich der Multiplikation von Matrizen]]

We define a projection operator to level $M$ as follows. Let $V$ be any right $M_2({\Bbb Z})^+$-module and let $M,N\in{\Bbb Z}$ such that $M|Np$; we then define
$$
\begin{array}{cccc}
\pi^{Np}_M:&V^{\Gamma_1(Np)}&\rightarrow&V^{\Gamma_1(M)}\\
&v&\mapsto&\sum_{\gamma\in \Gamma_1(Np)\backslash\Gamma_1(M)}v\gamma.\\
\end{array}
$$

{\bf Corollary. }{\it As operators on ${\cal M}_k(\Gamma)$ the following identity holds
$$
\pi^{Np}_M=\sum_{u\in{\bf M}, \,\delta\in{\bf A}} S_{u,\delta}.
$$
Furthermore, if $p|M$ then any $u\in{\bf M}$ is divisible by $p$ and any $a\in{\bf A}$ is congruent to $1$ modulo $p$.
}

{\it Proof. } We write
$$
S_{u,\delta}=\bigcup_{j} \Gamma x_{u,\delta,j}
$$
with $x_{u,\delta,j}\in M_2({\Bbb Z})^+$. $S_{u,\delta}$ then acts on elements $f\in{\cal M}_k(\Gamma)$ as
$$
S_{u,\delta}(f)=\sum_j f|_{x_{u,\delta,j}}.
$$
In particular, we obtain
$$
\sum_{u\in{\bf M},\delta\in{\bf A}} S_{u,\delta}(f)=\sum_{u\in{\bf M},\,\delta\in{\bf A}} \sum_j f|_{x_{u,\delta,j}}.\leqno(11)
$$
On the other hand, the previous Lemma shows that
\begin{eqnarray*}
\Gamma\backslash\Gamma_1(M)&=&\bigcup_{u\in{\bf M},\,\delta\in{\bf A}} \Gamma\backslash \Gamma\Mat{1}{}{u}{1}\sigma_\delta\Gamma\\
&=&\bigcup_{u\in{\bf M},\,\delta\in{\bf A}} \Gamma\backslash S_{u,\delta}\\
&=&\bigcup_{u\in{\bf M},\,\delta\in{\bf A}}\bigcup_j \Gamma x_{u,\delta,j}.\\
\end{eqnarray*}
and we deduce that
$$
\pi^{Np}_M(f)=\sum_{u\in{\bf M},\,\delta\in{\bf A}}\sum_j f|_{x_{u,\delta,j}}.\leqno(12)
$$
Comparing (11) and (12) yields the first claim of the Corollary. The second claim is immediate by equation (10). This completes the proof of the Corollary.

%[[Moral/Fazit: The level lowering operator $\pi^N_M$ is a Hecke operator and we may compute its trace via the trace formula. 
%]]

{\bf Convention.}  Since the full Hecke algebra $\Gamma\backslash M_2({\Bbb Z})^+/\Gamma$ is not commutative $\Gamma\backslash M_2({\Bbb Z})^+/\Gamma$ does no longer act from the left, it only acts fom the right on spaces of modular forms. Nevertheless we will write the action of $\Gamma\backslash M_2({\Bbb Z})^+/\Gamma$ from the left, noting that this now is an anti action. In doing so we use the following convention: if $X_1,\ldots,X_n\in \Gamma\backslash M_2({\Bbb Z})^+/\Gamma$ then 
$X_n\cdots X_1 f$ means $((((f X_1)X_2)\cdots )X_{n-1})X_n$ i.e. we first apply $X_1$ then $X_2$ and so on. It does not mean to apply the whole product form the right to $f$, i.e. it does not mean $f(X_n\cdots X_1)$. The term  $f(X_n\cdots X_1)$ would be denoted by $(X_n\cdots X_1)f$, but this will not appear.

{\bf Remark. } Using the above Corollary and Remark in section 3.1 we see that the product of operators $\pi^N_M T_{\ell_1}^{r_1}\cdot\ldots\cdot T_{\ell_s}^{r_s}$ is a ${\Bbb Z}$-linear combination of the operators
$S_{u,\delta} T(n',n') T_n$, where $n$ runs over certain divisors of $\ell_1^{r_1}\cdot\ldots\cdot \ell_s^{r_s}$, $n'$ runs over certain divisors of $\ell_1^\infty\cdot\ldots\cdot \ell_s^\infty$, which are prime to $Np$ and $u\in{\bf M}$, $\delta\in{\bf A}$. Since this is an identity in the abstract Hecke algebra and not in ${\rm End}({\cal M}_k(\Gamma))$, this linear combination clearly does not depend on $k$.

{\bf Example.} We compute the product of Hecke operators $S_{u,\delta}T_nT_p^h\langle\epsilon\rangle$. We first note that $T_nT_p^h=T_nT_{p^h}=T_{p^hn}$ (notice that $p$ divides the level $Np$ and write $n=n_1n_2$ with $(n_1,p)=1$ and $n_2$ a power of $p$). 
%
%[[Zunachst gilt $T_p^e=T_{p^e}$, denn: $T_p^e=T(p)^e=T(p^e)$ (s. [M], Lemma 4.5.7 (2), p. 140]; beachte $p$ teilt den Level) $=\sum_{a+b=e\atop (p^a,Np)=1} T(p^a,p^b)=T(1,p^e)=T_{p^e}$ (Erinnere: wir haben $T_n:=T(1,n)$ gesetzt). Schreibe nun $n=n'p^a$ mit $(n',p)=1$. Dann folgt
%$T_nT_p^h=T_{n'}T_{p}^aT_p^h$ (wegen [M], Lemma 4.5.8 (1), p. 141 und obigem) $=T_{n'}T_{p}^{a+h}=T_{n'}T_{p^{a+h}}$ (wegen obigem) $=T_{n'p^{a+h}}$ (wieder wegen [M], Lemma 4.5.8 (1), p. 141) $=T_{np^h}=\Gamma\Mat{1}{}{}{np^h}\Gamma$. Insgesamt also haben wir
%$$
%T_n T_p^h=\Gamma \Mat{1}{}{}{np^h}\Gamma
%$$
%]]
%
We choose a left coset decomposition
$$
\Gamma \Mat{1}{}{}{p^hn}\Gamma=\bigcup_{\alpha\in{\cal V}} \Gamma\Mat{1}{}{}{p^hn}\alpha\quad(\mbox{disjoint union}),
$$
i.e. $\alpha$ runs over a system of representatives ${\cal V}$ of $\Mat{1}{}{}{p^hn}^{-1}\Gamma\Mat{1}{}{}{p^hn}\cap \Gamma\backslash \Gamma$ 
%
%[[cf. [Bewersdorff], p. 11]]
%
We then find
\begin{eqnarray*}
\Gamma\sigma_\epsilon\Gamma \Mat{1}{}{}{p^hn}\Gamma \Mat{1}{}{u}{1}\sigma_\delta\Gamma
&=&\sigma_\epsilon \Gamma\Mat{1}{}{}{p^hn}\Gamma \Mat{1}{}{u}{1}\Gamma\sigma_\delta\\
&=&\bigcup_{\alpha\in{\cal V}}\sigma_\epsilon \Gamma \Mat{1}{}{}{p^hn}\alpha\Mat{1}{}{u}{1}\Gamma\sigma_\delta.\\
\end{eqnarray*}
Replacing ${\cal V}$ by a certain subset ${\cal V}'\subset {\cal V}$ we obtain a disjoint union
$$
\Gamma\sigma_\epsilon\Gamma \Mat{1}{}{}{p^hn}\Gamma \Mat{1}{}{u}{1}\sigma_\delta\Gamma
=\dot{\bigcup}_{\alpha\in{\cal V}'}\Gamma \sigma_\epsilon \Mat{1}{}{}{p^hn}\alpha\Mat{1}{}{u}{1}\sigma_\delta\Gamma.
$$
We then caculate using the definition of the product of Hecke operators
\begin{eqnarray*}
S_{u,\delta} T_{p^hn}\langle\epsilon\rangle
&\stackrel{def}{=}&\sum_{x\in\Gamma\backslash \Gamma \sigma_\epsilon\Gamma \Mat{1}{}{}{p^hn}\Gamma\Mat{1}{}{u}{1}\sigma_\delta\Gamma/\Gamma} m_x \,x\\
&=&\sum_{x\in \Gamma\backslash\dot{\bigcup}_{\alpha\in{\cal V}'}\Gamma \sigma_\epsilon \Mat{1}{}{}{p^hn}\alpha\Mat{1}{}{u}{1}\sigma_\delta\Gamma/\Gamma} m_x\, x,\\ 
\end{eqnarray*}
where the multiplicity $m_x\in{\Bbb N}\cup\{0\}$ is a certain integer. 
%
%[[cf. [Hida, modular forms and galois cohomology, p. 79, die erste Gleichung auf dieser seite]]]
%
Obviously, the last equation may be rewritten as follows
$$
S_{u,\delta} T_{p^hn}\langle\epsilon\rangle
=\sum_{\alpha\in{\cal V}'} m_\alpha \,\Gamma  \sigma_\epsilon \Mat{1}{}{}{p^hn}\alpha\Mat{1}{}{u}{1}\sigma_\delta\Gamma,\leqno(13)
$$
where we have set $m_\alpha=m_x$ with $x=\Gamma  \sigma_\epsilon \Mat{1}{}{}{p^hn}\alpha\Mat{1}{}{u}{1}\sigma_\delta\Gamma$. Thus, altogether we have found:

$S_{u,\delta}T_nT_p^h\langle\epsilon\rangle$ is a ${\Bbb Z}$-linear combination of the Hecke operators 
$$
\Gamma  \sigma_\epsilon \Mat{1}{}{}{p^hn}\alpha\Mat{1}{}{u}{1}\sigma_\delta\Gamma,
$$ 
where $\epsilon,\delta\in({\Bbb Z}/(Np))^*$, $\alpha\in{\cal V}'$.

%[[Zur definition von $\Gamma x\Gamma \cdot \Gamma y\Gamma$ s. [Hida, Modular forms and Galois cohomology], chapter 3.1, die Formel auf S. 79 oben]]

{\bf 3.3. The slope decomposition of the space of modular forms. } We denote by ${\cal M}_k(\Gamma,\chi)$ resp. ${\cal S}_k(\Gamma,\chi)$ the space of complex modular forms resp. of complex cusp forms of level $\Gamma$ and nebentype $\chi$. For any subring $A\subset {\Bbb C}$ we define ${\cal M}_k(\Gamma;A)$ resp. ${\cal M}_k(\Gamma,\chi;A)$
as the subspace of ${\cal M}_k(\Gamma)$ resp. ${\cal M}_k(\Gamma,\chi)$ consisting of forms, whose Fourier coefficients are contained in $A$; in a completely analogous manner we define ${\cal S}_k(\Gamma;A)$ and ${\cal S}_k(\Gamma,\chi;A)$. $\tilde{\cal H}$ and, hence, its subalgebra ${\cal H}$ act on ${\cal M}_k(\Gamma,\chi)$  
%
%[[${\cal H}$ acts on ${\cal M}_k(\Gamma,\chi)$, da $\langle\epsilon\rangle$ und $T_{\ell}$ kommutieren.]]
%
%[[Alternatives Argument: Dass $T_{\ell,\epsilon}$ auf ${\cal M}_k(\Gamma)$ operiert ist klar, da $T_{\ell,\epsilon}$ ein $\Gamma$-double coset ist. $T_{\ell,\epsilon}$ acts on ${\cal M}_k(\Gamma,\chi)$, da 
%$$
%\Gamma\langle\epsilon\rangle \Mat{1}{}{}{\ell}=\langle\epsilon\rangle \Gamma\Mat{1}{}{}{\ell}=\langle\epsilon\rangle\bigcup_i\Gamma \alpha_i,
%$$
%woraus folgt, dass $f|T_{\ell,\epsilon}=f|\langle\epsilon\rangle T_{\ell}=\chi(d)f|T_{\ell}\in {\cal M}_k(\Gamma,\chi)$, da $T_{\ell}$ ${\cal M}_k(\Gamma,\chi)$ invariant lasst.
%]]
%
and the submodules ${\cal M}_k(\Gamma;A)$ and ${\cal M}_k(\Gamma,\chi;A)$ then are invariant under the action of $\tilde{\cal H}$. 
% 
%[[kann man an der Formel fur die Aktion der Hecke Operatoren auf der $q$-Entwicklung ablesen; BEACHTE: ${\cal M}_k(\Gamma)$ ist als Hecke modul uber $A$ definiert (d.h. ${\cal M}_k(\Gamma)={\cal M}_k(\Gamma;A)\otimes {\Bbb C}$ und ${\cal M}_k(\Gamma;A)$ ist Hecke invariant, ABER: ${\cal M}_k(\Gamma)$ zerfallt nicht uber $A$ als Hecke Modul !!! auch nicht wenn $A=E$ ein Korper ist.)]] 
%
Moreover, we denote ${\Bbb Z}[Np]$ resp. ${\Bbb Q}(Np)$ the ring resp. the field obtained from ${\Bbb Z}$ resp. ${\Bbb Q}$ by adjoining all $\varphi(Np)$-th roots of unity. We note that ${\cal M}_k(\Gamma)={\cal M}_k(\Gamma;{\Bbb Z})\otimes_{\Bbb Z}{\Bbb C}$ as well as ${\cal M}_k(\Gamma,\chi)={\cal M}_k(\Gamma,\chi;{\Bbb Z}[Np])\otimes_{{\Bbb Z}[Np]}{\Bbb C}$ are defined over ${\Bbb Z}[Np]$ as Hecke modules.
%
%[[s. [Hida], p. 132 Theorem 1 und p. 141, Zeile -15 sowie Theorem 6.3.2 p. 178]]
%
%
%[[Beachte (WICHTIG !!): es gilt ${\cal M}_k(\Gamma,\chi)={\cal M}_k(\Gamma,\chi;{\Bbb Z}[Np])\otimes_{{\Bbb Z}[Np]}{\Bbb C}$, aber die {\it Zerlegung}
%$${\cal M}_k(\Gamma,\chi;{\Bbb Q}(Np))=\bigoplus_\chi  {\cal M}_k(\Gamma,\chi;{\Bbb Q}[Np])\otimes_{{\Bbb Z}[Np]}{\Bbb C}$$ ist nicht uber %${\Bbb Z}[Np]$ definiert, da die Idempotenten $e_\chi$ den nenner $\varphi(Np)$ haben !!
%]]
%
For later purpose we note that 
$$
{\cal M}_k(\Gamma,{\Bbb Q}(Np))=\bigoplus_\chi e_\chi{\cal M}_k(\Gamma,{\Bbb Q}(Np))\leqno(14)
$$ 
and $e_\chi{\cal M}_k(\Gamma;{\Bbb Q}(Np))={\cal M}_k(\Gamma,\chi;{\Bbb Q}(Np))$.

Since ${\cal H}$ acts on ${\cal M}_k(\Gamma,\chi)$ and using the canonical isomorphism ${\cal H}\cong {\cal H}_0$ the Hecke algebra 
${\cal H}_0$ also acts on ${\cal M}_k(\Gamma,\chi)$. In particular, the Hecke operator $T(n',n')\in{\cal H}_0$, $(n',Np)=1$ acts via $\Gamma\sigma_{n'}^{-1}\Mat{n'}{}{}{n'}\Gamma$ on ${\cal M}_k(\Gamma,\chi)$ and for $f$ in ${\cal M}_k(\Gamma,\chi)$ we obtain
$$
T(n',n')f=\chi(n') {n'}^{k-2}\,f.\leqno(15) 
$$  
%
%[[d.h. $T(n',n')$ operiert auf ${\cal M}_k(\Gamma,\chi)$ obwohl $T(n',n')$ nicht in der hecke algebra zum Level $\Gamma$ liegt ($T(n',n')$ liegt in hecke algebra zum level $\Gamma_0(Np)$)
%]]

%[[siehe [M], (4.5.27), p. 143. Wir rechnen das aber nocheinmal nach: 

%$T(n',n')_0f=(\det \Mat{n'}{}{}{n'})^{k/2-1}\,f|_{\sigma_{n'}^{-1}\Mat{n'}{}{}{n'}}=(\det \Mat{n'}{}{}{n'})^{k/2-1}\,\chi(n') %f|_{\Mat{n'}{}{}{n'}}=
%(\det \Mat{n'}{}{}{n'})^{k/2-1}\,\chi(n') (\det \Mat{n'}{}{}{n'})^{k/2} j(\Mat{n'}{}{}{n'},z)^{-k} f(\Mat{n'}{}{}{n'}z)=\chi(n') {n'}^{(k-2)} f(z)$ 

%(fur Defintion von $f|_{\alpha}$ siehe [M], (2.1.5), p. 37, zur Definition von $f|_{\Gamma\alpha\Gamma}$ siehe [M], (2.8.2), p. 74]]

Since the operators $\langle \epsilon\rangle$, $(\epsilon,Np)=1$ commute with the Hecke operators $T_\ell$, the space of modular forms ${\cal M}_k(\Gamma,\chi)$ is a ${\cal H}$-module 
%
%[[${\cal H}$ is generated by $T_\ell$ and $T(\ell,\ell)$ $\Rightarrow$ any $T\in{\cal H}$ commutes with $\langle\epsilon\rangle$]]
%
and we obtain a slope decomposition with respect to $T_p$. More precisely, as in section 1, we write ${\cal M}_{k}(\Gamma,\chi;{\Bbb Q}(Np))^\alpha$ resp. 
${\cal S}_k(\Gamma,\chi;{\Bbb Q}(Np))^\alpha$ to denote the slope $\alpha$ subspace, i.e. ${\cal M}_{k}(\Gamma,\chi;{\Bbb Q}(Np))^\alpha$ is the image of ${\cal M}_{k}(\Gamma,\chi;{\Bbb Q}(Np))$ under $p_\alpha(T_p)$, where $p_\alpha$ is the factor of the characteristic polynomial of $T_p$ on ${\cal M}_k(\Gamma,\chi;{\Bbb Q}(Np))$, whose roots $\lambda$ (in a splittig field of $T_p$) have $p$-adic valuation $v_p(\lambda)$ different from $\alpha$. We note that ${\cal M}_k(\Gamma,\chi)^\alpha$ is an ${\cal H}$-module. 

%[[Grund: ${\cal M}_k(\Gamma,\chi)^\alpha=e_\alpha e_\chi{\cal M}_k(\Gamma)$ und da $\tilde{\cal H}$ kommutativ ist, vertauschen $e_\chi, e_\alpha$ mit allen $T(1,\ell), T(\ell,\ell)\in{\cal H}$, aber $T_\ell, T(\ell,\ell)$ erzeugen
%${\cal H}$ $\Rightarrow$ $e_\chi, e_\alpha$ vertauschen mit allen $T\in{\cal H}$ $\Rightarrow$ ${\cal H}$ lasst $e_\alpha e_\chi{\cal %M}_k(\Gamma)={\cal M}_k(\Gamma,\chi)^\alpha$ invariant.
%]]

Since the slope decomposition is defined over ${\Bbb Q}$, we obtain 
$$
{\cal M}_k(\Gamma,\chi;{\Bbb Q}(Np))=\bigoplus_{\alpha\in{\Bbb Q}_{\ge 0}} {\cal M}_k(\Gamma,\chi;{\Bbb Q}(Np))^\alpha.\leqno(16)
$$

{\bf Lemma. }{\it The characteristic polynomial $\Psi_{k,\chi,\ell}^\alpha$ of $T_\ell\in{\cal H}$ acting on ${\cal M}_k(\Gamma,\chi)^\alpha$ is contained in ${\Bbb Z}[Np][X]$, i.e. its coefficients are integers in ${\Bbb Q}(Np)$.
}

{\it Proof. } Since ${\cal M}_k(\Gamma,\chi;{\Bbb Z}[Np])$ is a lattice in ${\cal M}_k(\Gamma,\chi;{\Bbb Q}(Np))$, the ${\Bbb Z}[Np]$-submodule ${\cal M}_k(\Gamma,\chi;{\Bbb Z}[Np])$ contains a ${\Bbb Q}(Np)$-basis $\{b_i\}$ of ${\cal M}_k(\Gamma,\chi;{\Bbb Q}(Np))$. Thus, $\{p_\alpha(T_p) b_i\}\subset p_\alpha(T_p) {\cal M}_k(\Gamma,\chi;{\Bbb Z}[Np])$ is a system of generators for the ${\Bbb Q}(Np)$-vector space
${\cal M}_k(\Gamma,\chi;{\Bbb Q}(Np))^\alpha$ and we see that the ${\Bbb Z}[Np]$-submodule $p_\alpha(T_p){\cal M}_k(\Gamma,\chi;{\Bbb Z}[Np])$ contains a ${\Bbb Q}(Np)$-basis of ${\cal M}_k(\Gamma,\chi;{\Bbb Q}(Np))^\alpha$. Since furthermore $p_\alpha(T_p) {\cal M}_k(\Gamma,\chi;{\Bbb Z}[Np])$ as submodule of a vector space necessarily is a free ${\Bbb Z}[Np]$-module we deduce that $p_\alpha(T_p) {\cal M}_k(\Gamma,\chi;{\Bbb Z}[Np])$ is a lattice in ${\cal M}_k(\Gamma,\chi;{\Bbb Q}(Np))$. Moreover, $p_\alpha(T_p)\in{\cal H}$,
%
%[[da $T_p\in{\cal H}$]]
%
hence, $p_\alpha(T_p)$ commutes with $T_\ell$ for all primes $\ell$ 
%
%[[da ${\cal H}$ kommutativ ist]]
%
and $p_\alpha(T_p) {\cal M}_k(\Gamma,\chi;{\Bbb Z}[Np])$ therefore is ${\cal H}$-stable, hence, $T_\ell$-stable for all $\ell$. This immediately implies that $\Psi_{k,\chi,\ell}^\alpha$ is contained in ${\Bbb Z}[Np][X]$.

Let $E/{\Bbb Q}(Np)$ be a finite extension, which splits  $T_p$ on ${\cal M}_k(\Gamma,\chi)$. We denote by $e$ the ramification index of the prime ideal ${\mathfrak p}\le {\cal O}_E$ corresponding to $v_p=w\circ i$ on $E$ (cf. section 1.1). As a consequence of the above Lemma, all eigenvalues of $T_p$ are integers, hence their $p$-adic value is contained in $\frac{1}{e}{\Bbb N}\cup\{0\}$. In different words, ${\cal M}_k(\Gamma,\chi,{\Bbb Q}(Np))^\alpha$ vanishes unless
$$
\alpha\in \frac{1}{e}{\Bbb N}\cup\{0\}.
$$

{\bf 3.4. Quasi idempotents attached to the slope decomposition. } We introduce one more piece of notation: if $\alpha\in{\Bbb Q}$ we set
$$
{\cal M}_k(\Gamma,\chi\omega^{-k})^{\le\alpha}=\bigoplus_{\beta:\,0\le\beta\le\alpha} {\cal M}_k(\Gamma,\chi\omega^{-k})^{\beta}.
$$

We fix arbitrary weights $k,k'>2$ and we fix an arbitrary slope $\alpha\ge 0$. We denote by 
$$
0\le\beta_1<\beta_2<\ldots<\beta_r\le\alpha+1\leqno(17)
$$ 
the slopes, which appear in ${\cal M}_k(\Gamma,\chi\omega^{-k})^{\le\alpha+1}$ or in ${\cal M}_{k'}(\Gamma,\chi\omega^{-k'})^{\le\alpha+1}$, i.e. for any $i$ we have ${\cal M}_k(\Gamma,\chi\omega^{-k})^{\beta_i}\not=0$ or ${\cal M}_{k'}(\Gamma,\chi\omega^{-k'})^{\beta_i}\not=0$. In particular, we may write the slope decomposition as
$$
{\cal M}_k(\Gamma,\chi\omega^{-k})^{\le\alpha+1}=\bigoplus_{i=1}^r {\cal M}_k(\Gamma,\chi\omega^{-k})^{\beta_i},\leqno(18)
$$
where some of the slope spaces ${\cal M}_k(\Gamma,\chi\omega^{-k})^{\beta_i}$ now may be trivial. Of course, the same decomposition with $k$ replaced by $k'$ holds. Since ${\rm dim}\,{\cal M}_k(\Gamma,\chi\omega^{-k})^{\le\alpha+1}\le M(\alpha+1)$ and ${\rm dim}\,{\cal M}_{k'}(\Gamma,\chi\omega^{-k'})^{\le\alpha+1}\le M(\alpha+1)$ we know that $r\le 2M(\alpha+1)$. 

We denote by $\Xi_k^{\alpha+1}$ resp. $\Xi_{k'}^{\alpha+1}$ the characteristic polynomial of $T_p$ acting on ${\cal M}_k(\Gamma,\chi\omega^{-k})^{\le\alpha+1}$ resp. on ${\cal M}_{k'}(\Gamma,\chi\omega^{-k'})^{\le\alpha+1}$ and we denote by
$E_k^{\alpha+1}$ resp. $E_{k'}^{\alpha+1}$ the field obtained from ${\Bbb Q}$ by adjoining the roots of $\Xi_k^{\alpha+1}$ resp. of 
$\Xi_{k'}^{\alpha+1}$; we set $E_{k,k'}^{\alpha+1}=E_k^{\alpha+1}E_{k'}^{\alpha+1}$. Thus, $E=E_{k,k'}^{\alpha+1}$ is a splitting field for $T_p$ on ${\cal M}_k(\Gamma,\chi\omega^{-k})^{\le\alpha+1}$ and for $T_p$ on ${\cal M}_{k'}(\Gamma,\chi\omega^{-k'})^{\le\alpha+1}$ and since ${\rm deg}\,\Xi^{\alpha+1}_k={\rm dim}\,{\cal M}_k(\Gamma,\chi\omega^{-k})^{\le\alpha+1}\le M(\alpha+1)$ we find
$$
[E:{\Bbb Q}]\le [E_k^{\alpha+1}:{\Bbb Q}][E_{k'}^{\alpha+1}:{\Bbb Q}]\le M(\alpha+1)^2.
$$
We denote by ${\mathfrak p}={\mathfrak p}_{k,k'}$ the prime ideal in the ring of integers ${\cal O}_E$ of $E$, which corresponds to the valuation $v_p=i\circ w$ on $E$ (cf. section 1.2), $e=e_{k,k'}$ is the ramification index of ${\mathfrak p}|p$ and $\varpi=\varpi_{k,k'}\in E$ any element satisfying $v_p(\varpi)=1/e$. In particular,
$$
e=e_{k,k'}\le [E:{\Bbb Q}]\le M(\alpha+1)^2\leqno(19)
$$
is bounded independently of $k,k'$. Moreover, we may write
$$
\beta_i=b_i/e
$$
with $b_i\in{\Bbb N}\cup\{0\}$ for all $i=1,\ldots,r$; in particular,
$$
\beta_j-\beta_i> 1/e\leqno(20)
$$
for all $1\le i,j\le r$, $i\not=j$.

The space ${\cal M}_k(\Gamma,\chi\omega^{-k})$ decomposes
$$
{\cal M}_k(\Gamma,\chi\omega^{-k})=\bigoplus_\mu V(\mu),
$$
where $\mu\in E$ runs over the eigenvalues of $T_p$ and $V(\mu)$ is the generalized eigenspace of $T_p$ with respect to $\mu$. Hence, there is a basis ${\cal B}_\mu$ of $V(\mu)$ such that $T_p$ is represented by the upper triangular matrix
$$
{\cal D}_{{\cal B}_\mu}(T_p|_{V(\mu)})=\left(\begin{array}{ccc} \mu&&*\\ &\ddots&\\ &&\mu\\ \end{array}\right).
$$
Of course an analogous decomposition holds:
$$
{\cal M}_{k'}(\Gamma,\chi\omega^{-k'})=\bigoplus_{\mu'} V(\mu'),
$$
where $\mu'\in E$ runs over the eigenvalues of $T_p$ and
$$
{\cal D}_{{\cal B}_{\mu'}}(T_p|_{V(\mu')})=\left(\begin{array}{ccc} \mu'&&*\\ &\ddots&\\ &&\mu'\\ \end{array}\right).
$$
with respect to some basis ${\cal B}_{\mu'}$ of $V(\mu')$.

We attach to $T_p$ and $\varpi$ a sequence of operators $\tilde{e}_i=\tilde{e}_{i,k,k'}$, $i=1,2,3,\ldots$, on ${\cal M}_k(\Gamma,\chi\omega^{-k})$ in the following way. We set
$$
\tilde{e}_1=\tilde{e}_{1,k,k'}=\left(\frac{T_p}{\varpi^{b_1}}\right)^{q-1}.\leqno(21)
$$
Assuming that $\tilde{e}_1,\ldots,\tilde{e}_{i-1}$ have been defined we set
$$
\tilde{e}_{i}=\tilde{e}_{i,k,k'}=\left(\frac{T_p}{\varpi^{b_i}}\right)^{q-1}\,(1-\tilde{e}_1^{q^{b_i(q-1)+1}}-\tilde{e}_2^{q^{b_i(q-1)+1}}-\cdots-\tilde{e}_{i-1}^{q^{b_i(q-1)+1}}).\leqno(22)
$$
If we want to emphasize that $\tilde{e}_{i}$ also depends on $\varpi$ we write $\tilde{e}_{i,k,k',\varpi}$. Clearly, $\tilde{e}_i\in E[T_p]$, 
because $\varpi\in E$. We determine the effect of the elements $\tilde{e}_i$ for those $i$, which satisfy $\beta_i\le\alpha$, on ${\cal M}_k(\Gamma,\chi\omega^{-k})$:

%Alternativ: Definiere 
%$$
%\tilde{e}_{i}=\left(\frac{T}{\varpi^i}\right)^{q-1}\,(1-\tilde{e}_0^{q^{b_i(q-1)+1}}-\tilde{e}_1^{q^{b_{i-1}(q-1)+1}}-\cdots-\tilde{e}_i^{q^{b_1(q-1)+1}}).
%$$
%das funktioniert auch und benotigt weniger Potenzen von $T$, was aber spatestens beim Grenzubergang $N\rightarrow\infty$ keien Unterschied %mehr macht.
%]]

{\bf Lemma. }{\it 1.) Let $\mu$ be any eigenvalue of $T_p$ on ${\cal M}_k(\Gamma,\chi\omega^{-k})$. We set $F=E_{k,k'}^{\alpha+1}(\mu)$ 
%
%[[beachte: es ist nicht notwendig $v_p(\mu)\le \alpha+1$, d.h. $\mu$ muss nicht in $E$ liegen]]
%
and we denote by $(p)$ the prime ideal in ${\cal O}_F$ generated by $p$. For all $1\le i\le r$ satisfying $\beta_i\le\alpha$ the element $\tilde{e}_i$ restricted to $V(\mu)$ looks as follows:
$$
{\cal D}_{{\cal B}_\mu}(\tilde{e}_i|_{V(\mu)})
=\left(\begin{array}{ccc} \gamma&&*\\ &\ddots&\\ &&\gamma\\ \end{array}\right),\leqno(23)
$$
where
$$
\gamma\in\left\{\begin{array}{ccl}
{\mathfrak p}&\mbox{if}&v_p(\mu)\le\alpha+1\quad\mbox{and}\quad v_p(\mu)\not=\beta_i\\
1+{\mathfrak p}&\mbox{if}&v_p(\mu)=\beta_i\\
(p)&\mbox{if}&v_p(\mu)>\alpha+1
\end{array}\right.
$$  
(recall that ${\mathfrak p}={\mathfrak p}_{k,k'}$). Of course, the same statement holds true if we replace ${\cal M}_k(\Gamma,\chi\omega^{-k})$ by ${\cal M}_{k'}(\Gamma,\chi\omega^{-k'})$, $V(\mu)$ by $V(\mu')$ and ${\cal B}_\mu$ by ${\cal B}_{\mu'}$.

2.) For any $i=1,2,3,\ldots$ there is a polynomial $p_i=p_{i,\varpi}\in{\cal O}_E[X]$ ($E=E_{k,k'}^{\alpha+1}$) such that $X|p_i$ and
$$
\tilde{e}_i=p_i\left(\frac{T}{\varpi^{b_i}}\right).
$$
Furthermore, ${\rm deg}\,p_i$ only depends on $i$.
}

%[[Uberblick: Wirkung von $\tilde{e}_i$ auf den einzelnen SlopeRaumen (Hier: Annahme, dass $T$ auf $V$ diagonalisierbar)
%
%$$
%\begin{array}{ccccc ccccc}
%\qquad {\rm Slope Raum}&V^{\beta_1}&V^{\beta_2}&\cdots&V^{\beta_{i-1}}&V^{\beta_i}&V^{\beta_{i+1}}&\cdots\\
%&&&&&&&&\\
%{\rm Operator}\backslash {\rm EW}&\varpi^{b_1}&\varpi^{b_2}&\cdots&\varpi^{b_{i-1}}&\varpi^{b_i}&\varpi^{b_{i+1}}&\cdots&\\
%&&&&&&&&\\
%\tilde{e}_0^{q^{c_i}}+\cdots+\tilde{e}_{i-1}^{q^{c_i}}&1+{\mathfrak p}^{c_i}&1+{\mathfrak p}^{c_i}&\cdots&1+{\mathfrak p}^{c_i}&{\mathfrak %p}^{q^{c_i}}&{\mathfrak p}^{q^{c_i}}&\cdots\\
%&&&&&&&&\\
%1-\tilde{e}_0^{q^{c_i}}-\cdots-\tilde{e}_{i-1}^{q^{c_i}}&{\mathfrak p}^{c_i}&{\mathfrak p}^{c_i}&\cdots&{\mathfrak p}^{c_i}&1+{\mathfrak p}^{q^{c_i}}&1+{\mathfrak p}^{q^{c_i}}&\cdots\\
%&&&&&&&&\\
%(\frac{T}{\varpi^i})^{q-1}\circ(1-\cdots-\tilde{e}_{i-1}^{q^{c_i}})&{\mathfrak p}&{\mathfrak p}&\cdots&{\mathfrak p}&(1+{\mathfrak p}^{q^{c_i}})\cdot\underbrace{(\frac{\lambda}{\varpi^i})^{q-1}}_{\equiv 1 \,(\varpi)}&{\varpi}(1+{\mathfrak p}^{q^{c_i}})&\cdots\\
%\end{array}
%$$
%($\lambda=$ EW von $T$ auf $V^{i/e}$, $\varpi=$ lokales Primelement in $E$, i.e. $v_{\mathfrak p}(\varpi)=1/e$, $c_i=i(q-1)+1$, $(\varpi)$ bedeutet "modulo $\varpi$".
%Beachte: Im letzen Schritt ($=$Anwenden von $(\frac{T}{\varpi^i})^{q-1}$) benotigen wir den Exponenten $c_i=i(q-1)+1$ im Faktor vor dem slope $0$-Vektor; weniger Exponent wurde nach Anwenden von $(\frac{T}{\varpi^i})^{q-1}$ keine $\varpi$-Potenz ubrig lassen!)
%]]

{\it Proof. } 1.) We abbreviate $c_i=b_i(q-1)+1$ and ${\cal D}={\cal D}_{{\cal B}_\mu}$. We use induction on $i=1,2,3,\ldots$. We note in advance that if $\mu$ is an eigenvalue of $T_p$ on ${\cal M}_k(\Gamma,\chi\omega^{-k})$ then precisely one of the following two holds (cf. (18)): 

- $v_p(\mu)=\beta_j$ for some $1\le j\le r$ in which case $\mu\in E=E_{k,k'}^{\alpha+1}$, i.e. $F=E$ 

or 

- $v_p(\mu)>\alpha+1$. 

We first look at the case $i=1$. We assume $\beta_1\le \alpha$ since otherwise there is nothing to prove. The element $\tilde{e}_1$ equals $(\frac{T_p}{\varpi^{b_1}})^{q-1}$, thus we find
$$
{\cal D}(\tilde{e}_1|_{V(\mu)})=\left(
\begin{array}{ccc}
\left(\frac{\mu}{\varpi^{b_1}}\right)^{q-1}&&*\\
&\ddots&\\
&&\left(\frac{\mu}{\varpi^{b_1}}\right)^{q-1}\\
\end{array}
\right).
$$
If $v_p(\mu)=\beta_1$ then $(\frac{\mu}{\varpi^{b_1}})^{q-1}\in 1+{\mathfrak p}\subset E$ (note that $\mu\in E=E_{k,k'}^{\alpha+1}$); if $v_p(\mu)=\beta_j$ with $j>1$ then $v_p(\mu)-\beta_1=\beta_j-\beta_1>1/e$ ($e=e_{k,k'}$; cf. (20)), hence, $(\frac{\mu}{\varpi^{b_1}})^{q-1}\in{\mathfrak p}$; if $v_p(\mu)> \alpha+1$ then $v_p(\mu/\varpi^{b_1})=v_p(\mu)-\beta_1>\alpha+1-\beta_1$ 
and since we assume that $\beta_1\le \alpha$ we find $v_p(\mu/\varpi^{b_1})>1$, hence, $(\frac{\mu}{\varpi^{b_1}})^{q-1}\in (p)\subset F$. Altogether, the claim about $\tilde{e}_1$ follows.

%[[$\left(\frac{\mu}{\varpi^{b_1}}\right)^{q-1}$ even is contained in ${\mathfrak p}^{q-1}$ if $v_p(\mu)>\beta_1$]]

%[[sei $x\in F$ mit $v_p(x)\ge 1$. Schreibe $x\tilde {\mathfrak P}^a\prod_i{\mathfrak q}_i^{a_i}$, wobeo ${\mathfrak P}$ das zu $v_p$ gehorige Primideal in ${\cal O}_F$ ist und ${\mathfrak q}_i$ primideale in ${\cal O}_F$ sind. Dann folgt $v_p(x)=a/e\ge 1\Rightarrow a\ge e$; dabei ist $e$ der Verzweigungsindex von ${\mathfrak P}|p$. Andererseits gilt $(x)={\mathfrak P}^a\prod_i{\mathfrak q}_i^{a_1}\Rightarrow x\in(x)\subset {\mathgfrak P}^a$, d.h. $x\in{\mathfrak P}^a$ mit $a\ge e$. $a\ge e$ heisst aber dass ${\mathfrak P}^a\subset (p)$, also insgesamt
%$x\in (p)$ qed.]]

We assume that the Lemma holds for $i-1$. We assume $\beta_i\le\alpha$ since otherwise there is nothing to prove. In particular, we obtain $\beta_1,\ldots,\beta_i\le\alpha$. Equation (23), which by our induction assumption holds for all $\tilde{e}_j$, $j\le i-1$, therefore yields
$$
{\cal D}((1-\tilde{e}_1^{q^{c_i}}-\ldots-\tilde{e}_{i-1}^{q^{c_i}})|_{V(\mu)})
=\left(\begin{array}{ccc} \gamma&&*\\ &\ddots&\\ &&\gamma\\ \end{array}\right),\leqno(24)
$$
where 
$$
\gamma\in \left\{\begin{array}{ccc}
{\mathfrak p}^{c_i}& {\rm if} &v_p(\mu)=\beta_j\;\mbox{for some $j\le i-1$}\\
1+{\mathfrak p}^{c_i}&{\rm if} &v_p(\mu)=\beta_j\;\mbox{for some $j> i-1$}\\
1+(p)^{c_i}&{\rm if}&v_p(\mu)>\alpha+1.\\
\end{array}\right.\leqno(25)
$$

%[[Note that $*\not\in{\cal O}_E$ !!!!!!!!!!!!]]

We mention that even $\gamma\in 1+{\mathfrak p}^{q^{c_i}}$ in case that $v_p(\mu)=\beta_j$ for some
$j> i-1$ and $\gamma\in 1+(p)^{q^{c_i}}$ in case that $v_p(\mu)> \alpha+1$. 
%
%[[Note that $q^{c_i}\ge c_i$, $q=p^f\ge 3$, da $p>2$ ]]
%
The operator $\left(\frac{T_p}{\varpi^{b_i}}\right)^{q-1}|_{V(\mu)}$ is represented by the matrix
$$
{\cal D}\left(\left(\frac{T_p}{\varpi^{b_i}}\right)^{q-1}|_{V(\mu)}\right)
=\left(\begin{array}{ccc} (\frac{\mu}{\varpi^{b_i}})^{q-1}&&*\\ &\ddots&\\ &&\left(\frac{\mu}{\varpi^{b_i}}\right)^{q-1}\\ \end{array}\right),\leqno(26)
$$
where $\mu/\varpi^{b_i}$ is as follows:
$$
\left(\frac{\mu}{\varpi^{b_i}}\right)^{q-1}\in\left\{\begin{array}{ccc}
{\mathfrak p}^{-b_i(q-1)}&{\rm if}&v_p(\mu)=\beta_j,\;j\le i-1\\
1+{\mathfrak p}&{\rm if}&v_p(\mu)=\beta_i\\
{\mathfrak p}&{\rm if}&v_p(\mu)=\beta_j,\,j>i\\
(p)&{\rm if}&v_p(\mu)>\alpha+1.\\
\end{array}
\right.\leqno(27)
$$
In case $v_p(\mu)=\beta_j$, $j\le i-1$, this is obvious, because $\mu\in E$ is integral; in case $v_p(\mu)=\beta_i$ the ratio $\mu/\varpi^{b_i}$ is a unit, whence the claim; in case $v_p(\mu)=\beta_j$, with $j>i$ equation (20) implies that $v_p(\mu/\varpi^{b_i})=v_p(\mu)-\beta_i=\beta_j-\beta_i>1/e$, hence, $\mu/\varpi^{b_i}\in{\mathfrak p}$; finally, if $v_p(\mu)\ge\alpha+1$ 
then $v_p(\mu/\varpi^{b_i})\ge \alpha+1-\beta_i$ and since we assume that $i$ is such that $\beta_i\le \alpha$, we furthermore obtain $v_p(\mu/\varpi^{b_i})\ge \alpha+1-\alpha=1$, i.e. $\mu/\varpi^{b_i}\in(p)\subset F$. Thus, equation (27) is proven.

Multiplying (24) by the matrix (26) we obtain
$$
{\cal D}(\tilde{e}_i|_{V(\mu)})
={\cal D}\left(\left(\frac{T}{\varpi^i}\right)^{q-1}\circ (1-\tilde{e}_1^{q^{c_i}}-\ldots-\tilde{e}_{i-1}^{q^{c_i}})|_{V(\mu)}\right)
=\left(\begin{array}{ccc} \gamma&&*\\ &\ddots&\\ &&\gamma\\ \end{array}\right),
$$
where equations (25) and (27) imply that $\gamma$ is as follows: $\gamma\in {\mathfrak p}$ if $v_p(\mu)=\beta_j$ with $j\le i-1$, $\gamma\in 1+{\mathfrak p}$
if $v_p(\mu)=\beta_i$, $\gamma\in {\mathfrak p}$ if $v_p(\mu)=\beta_j$ with $j\ge i+1$ and $\gamma\in(p)$ if $v_p(\mu)\ge \alpha+1$. But this is precisely the statment about $\tilde{e}_i$, which therefore is proven. Clearly, the above proof still holds if we replace ${\cal M}_k(\Gamma,\chi\omega^{-k})$ by ${\cal M}_{k'}(\Gamma,\chi\omega^{-k'})$, $V(\mu)$ by $V(\mu')$ and ${\cal B}_\mu$ by ${\cal B}_{\mu'}$. Thus, the first part of the Lemma is proven.

2.) Again, we use induction on $i$. If $i=1$ we know that $\tilde{e}_1=(T_p/\varpi^{b_1})^{q-1}$, hence, we choose 
$$
p_1=X^{q-1}.\leqno(28)
$$ 
Let $i$ be arbitrary and assume that the claim holds for $j=1,\ldots,i-1$. Then,
$\tilde{e}_j=p_j(T/\varpi^{b_j})=p_j^*(T/\varpi^{b_i})$, where $p_j^*\in{\cal O}_E[X]$ is
obtained from $p_j=\sum_h a_hX^h\in{\cal O}_E[X]$ by multiplying the coefficient $a_h$ with $\varpi^{h(b_i-b_j)}$. Recalling that $c_i=b_i(q-1)+1$ we obtain 
$$
\tilde{e}_i=(T_p/\varpi^{b_i})^{q-1}\circ(1-\sum_{j=1}^{i-1}p_j^*(T_p/\varpi^{b_i})^{q^{c_i}}.
$$ 
Thus, the polynomial
$$
p_i=X^{q-1}\left(1-\sum_{j=1}^{i-1} (p_j^*)^{q^{c_i}}\right)\in{\cal O}_E[X]\leqno(29)
$$ 
satisfies the claim. Equations (28) and (29) immediately yield that $X|p_i$ and since ${\rm deg}\,p_j={\rm deg}\,p_j^*$ we inductively see from equations (28) and (29) that ${\rm deg}\,p_i$ does not depend on $\varpi$, i.e. it only depends on $i$.
%
%2.) Fur unsere Zwecke reicht es sogar wenn wenn ${\rm deg}\, p_i={\rm deg}\,p_{i,\varpi}$ nur durch eine von $\varpi$ unabhangige Zahl $t_i$ beschrankt ist, d.h. ${\rm deg}\,p_{i,\varpi}\le t_i$ fur alle $\varpi$.
%
%]]
This completes the proof of the Lemma.

Recall that we fixed arbitrary weights $k,k'$ and a slope $\alpha\in{\Bbb Q}_{\ge 0}$ and $\beta_i$ are the non-trivial slopes appearing in 
${\cal M}_k(\Gamma,\chi\omega^{-k})^{\le\alpha+1}$ or in ${\cal M}_{k'}(\Gamma,\chi\omega^{-k'})^{\le\alpha+1}$. We now assume that $\alpha$ is a non-trivial slope, i.e.
${\cal M}_k(\Gamma,\chi\omega^{-k})^\alpha\not=0$ or ${\cal M}_{k'}(\Gamma,\chi\omega^{-k'})^\alpha\not=0$. Without loss of generality we may assume that ${\cal M}_k(\Gamma,\chi\omega^{-k})^\alpha\not=0$. In particular, there is  $s=s_{\alpha,k,k'}$ such that $1\le s\le r$ and  $\alpha=\beta_s$. Moreover, part 2.) of the above Lemma tells us that $\tilde{e}_{i,k,k'}=p_i\left(\frac{T_p}{\varpi_{k,k'}^{b_i}}\right)$. We define
$$
e_{k,k',\alpha}:=\tilde{e}_{s,k,k'}^{p^{M(\alpha+1)^2}}
$$
and 
$$
p_{\alpha,k,k'}:=p_{s}^{p^{M(\alpha+1)^2}}\in{\cal O}_E[X],
$$
where $s=s_{\alpha,k,k'}$ and $E=E_{k,k'}^{\alpha+1}$. Moreover, we define $i=i_{\alpha,k,k'}$ by $\alpha=i/e_{k,k'}$.

{\bf Proposition. } {\it 1.) Let $V(\mu)$ denote the generalized eigenspace attached to the eigenvalue $\mu$ of $T_p$ acting on ${\cal M}_k(\Gamma,\chi\omega^{-k})$. Then,
$$
{\cal D}_{{\cal B}_\mu}(e_{k,k',\alpha}|_{V(\mu)})
=\left(\begin{array}{ccc} \gamma&&*\\ &\ddots&\\ &&\gamma\\ \end{array}\right),\leqno(30)
$$
where $v_p(\gamma)\ge 1$ if $v_p(\mu)\not=\alpha$ and $v_p(\gamma-1)\ge 1$ if $v_p(\mu)=\alpha$. In different notation,
$$
\gamma\equiv\left\{\begin{array}{ccc}
0\pmod{p}&\mbox{if}&v_p(\mu)\not=\alpha\\
1\pmod{p}&\mbox{if}&v_p(\mu)=\alpha.\\
\end{array}
\right.
$$
The same statement holds if we replace $V(\mu)$ by the generalized eigenspaces $V(\mu')$ attached to an eigenvalue $\mu'$ of $T_p$ acting on ${\cal M}_{k'}(\Gamma,\chi\omega^{-k'})$.

2.)
$$
e_{k,k',\alpha}=p_{\alpha,k,k'}\left(\frac{T_p}{\varpi_{k,k'}^{i}}\right).
$$ 
Moreover, the polynomial $p_{\alpha,k,k'}\in {\cal O}_E[X]$ satisfies the following properties: 

$\bullet$ $X|p_{\alpha,k,k'}$ 

$\bullet$ There is a constant $N(\alpha)\in{\Bbb N}$ such that ${\rm deg}\,p_{\alpha,k,k'}\le N(\alpha)$ for all $k$, $k'$. (In different words, although $p_{\alpha,k,k'}$ depends on $k,k'$ its degree is bounded independently of $k,k'$.)
}

{\it Proof. }  1.) This is an immediate consequence of part 1.) of the above Lemma, taking into account that $e_{k,k'}\le M(\alpha+1)^2$ (cf. (19)), which implies that 

- ${\mathfrak p}^{p^{M(\alpha+1)^2}}\subset {\mathfrak p}^{{M(\alpha+1)^2}}\subset (p)$

and

- $(1+{\mathfrak p})^{p^{M(\alpha+1)^2}}\le 1+(p)$ (note that $e_{k,k'}$ is the ramification index of ${\mathfrak p}={\mathfrak p}_{k,k'}$and that $(1+{\mathfrak p}^a)/(1+{\mathfrak p}^{a+1})\cong({\cal O}/{\mathfrak p},+)$, which is annihilated by multiplication with $p$).

2.) All statements are obvious from part 2.) of the above Lemma except for the claim about the degree of $p_{\alpha,k,k'}$. But 
equation (29) shows that ${\rm deg}\,p_a<{\rm deg}\,p_b$ if $a<b$. Furthermore, equation (20) implies that $s=s_{\alpha,k,k'}\le i=i_{\alpha,k,k'}$. Since $i=\alpha e_{k,k'}\le \alpha M(\alpha+1)^2$ we therefore obtain
$$
{\rm deg}\,p_s\le {\rm deg}\,p_i\le {\rm deg}\,p_{\alpha M(\alpha+1)^2},
$$
hence, 
$$
{\rm deg}\,p_{\alpha,k,k'}={\rm deg}\,p_s^{p^{M(\alpha+1)^2}}\le {\rm deg}\,p_{\alpha M(\alpha+1)^2}^{p^{M(\alpha+1)^2}}=:N(\alpha),\leqno(31)
$$
which is the claim.

We note that equations (28) and (29) imply that ${\rm deg}\, p_j\ge q-1$ for all $j\in{\Bbb N}$. Hence, the definition of $N(\alpha)$ (cf. (31)) implies
$$
N(\alpha)\ge 1.\leqno(32)
$$

{\bf Remark.} 1.) In view of equation (30) we call $e_{k,k',\alpha}$ a quasi idempotent (or approximative idempotent) attached to the slope $\alpha$ subspace of ${\cal M}_k(\Gamma,\chi\omega^{-k})$ and ${\cal M}_{k'}(\Gamma,\chi\omega^{-k'})$.

{\bf 3.5. The basic trace identity. } We denote by $\omega:\,{\Bbb Z}_p^*\rightarrow\mu_{p-1}$ the Teichmuller character, i.e. $\omega(x)\equiv x\pmod{p}$ for all $x\in{\Bbb Z}_p^*$. In particular, $\omega$ factorizes over $1+p{\Bbb Z}_p$ and therefore induces a character $({\Bbb Z}/p{\Bbb Z})^*\rightarrow\mu_{p-1}$ by sending $\epsilon+p{\Bbb Z}\mapsto\omega(\epsilon)$, which, again, we denote by $\omega$.
%
%[[Detail des Wohldefiniertheitsbeweises fur $\omega:\,({\Bbb Z}/p{\Bbb Z})^*\rightarrow\mu_{p-1}$: Z.z. ist $\omega(\epsilon)=\omega(\epsilon+ph)$ fur alle $\epsilon\in{\Bbb Z}$, $(\epsilon,p)=1$ $h\in{\Bbb Z}$. Dazu: es gilt $\epsilon+ph=\epsilon\underbrace{(1+\epsilon^{-1}ph)}_{\in 1+p{\Bbb Z}_p}$
%$\Rightarrow \omega(\epsilon)=\omega(\epsilon+ph)$.
%]]
%
%[[Beachte: Wir konnen $\mu_{p-1}$ auffassen als Teimenge von $\bar{\Bbb Q}$, von ${\Bbb Z}_p$, ${\Bbb C}$, ${\Bbb C}_p$]]

{\bf Lemma 1. }{\it Let $k,k'\in{\Bbb N}$ be integers, which satisfy $k\equiv k'\pmod{p^m}$. Let furthermore
$$
[\xi]_\Gamma\in\Gamma\sigma_{\epsilon}\Mat{1}{}{}{p^hn}\alpha\Mat{1}{}{u}{1}\sigma_\delta\Gamma/\sim_\Gamma,
$$ 
be a $\Gamma$-conjugacy class, where $h,n,u\in{\Bbb N}\cup\{0\}$, $\epsilon,\delta\in({\Bbb Z}/(Np))^*$ with $h\ge 1$, $n\ge 1$, $u$ divisible by $p$ and $\alpha\in \Gamma$. Assume that $\xi$ is elliptic or hyperbolic; then the following congruence holds 
$$
\omega^{-k}(\epsilon\delta)\,{\rm tr}\,(\xi^\iota|L_k)\equiv \omega^{-k'}(\epsilon\delta)\,{\rm tr}\,(\xi^\iota|L_{k'})\pmod{p^{{\rm min}(m,hk,hk')}}.
$$
}

We note that the expressions $\omega^{-k}(\epsilon\delta)\,{\rm tr}\,(\xi^\iota|L_k)$ and 
$\omega^{-k'}(\epsilon\delta)\,{\rm tr}\,(\xi^\iota|L_{k'})$ take values in the field ${\Bbb Q}(\mu_{p-1})$, because the representation $L_k$ is defined over ${\Bbb Q}$ .

%[[Beachte: $p\in {\Bbb Q}$ liegt naturlich auch in $E={\Bbb Q}(\mu_{p-1})$, deshalb macht die Kongruenz "$\pmod{p}$" auch in $E$ Sinn ("Differenz ist durch $p$ teilbar"). Sie ist aquivalent zu "$\pmod{(p)}$", wobei $(p)$ das von $p$ in ${\cal O}_E$ erzeugte Hauptideal ist. Es gilt $(p)={\mathfrak p}_0$, in $E$; s. den folgenden Beweis. ]]

{\it Proof. } Conjugating $\xi$ by a suitable element in $\Gamma$ we may assume that $\xi\in\Gamma\sigma_{\epsilon}\Mat{1}{}{}{p^hn}\alpha\Mat{1}{}{u}{1}\sigma_\delta$. We write 
$$
\xi=\Mat{a}{b}{c}{d};
$$ 
straightforward computation shows that
$$
\xi\equiv\Mat{\epsilon\delta}{*}{0}{0}\pmod{p},
$$
where $*\in{\Bbb Z}$.
%
%[[Wir prufen das nach:
%$$
%\zeta=\Mat{a}{b}{c}{d}\Mat{1}{}{}{p^hn}\Mat{\alpha}{\beta}{\gamma}{\delta}\Mat{1}{}{u}{1}
%=\Mat{\alpha a+a\beta u+p^hnb(\gamma+\delta u)}{a\beta+\delta p^hnb}{c\alpha+c\beta u+p^hnd(\gamma+\delta u)}{c\beta+\delta p^hnd},
%$$
%wobei
%$$
%\alpha=\Mat{\alpha}{\beta}{\gamma}{\delta}\in \Gamma=\Gamma_1(Np)\Rightarrow \alpha,\delta\equiv 1, \gamma\equiv 0\pmod{Np}
%$$
%und
%$$
%\Gamma_0(Np)\ni\Mat{a}{b}{c}{d}=\sigma_{\epsilon}\Rightarrow a\equiv \epsilon,d\equiv \epsilon, c\equiv 0\pmod{Np}
%$$
%(klar: $\alpha,\beta,\gamma,\delta,a,b,c,d\in{\Bbb Z}$). Daraus ergibt sich wegen $p|u$ und $h\ge 1$
%$$
%\zeta\equiv \Mat{\epsilon}{*}{0}{0}\pmod{p}.
%$$ 
%Multiplizieren wir $\zeta$ noch von rechts mit $\sigma_\delta\in\Gamma_0(Np)$ und von links mit einem beliebigen $\gamma\in\Gamma$, dann ergibt sich sofort 
%$$
%\xi\equiv\Mat{1}{*}{0}{1}\Mat{\epsilon}{*}{0}{0}\Mat{\delta}{*}{0}{\delta^{-1}}=\Mat{\epsilon\delta}{*}{0}{0}\pmod{p}
%$$
%
%Beachte: wir haben Kongruenzen modulo $Np$, die von der Defintiion von $\Gamma=\Gamma_1(Np)$ herkommen und Kongruenzen modulo $p^hn$, die daher kommen, dass $\xi$ den faktor $\Mat{1}{}{}{p^hn}$ enthalt. Im Schnitt bekommen wir nur Kongruenzen modulo $p$ !
%]]
%
Thus, the characteristic polynomial of $\xi$ reads 
$$
\chi_\xi=T^2-{c_1}T+c_2\in{\Bbb Z}[T]\leqno(33)
$$ 
with $c_1=a+d\equiv \epsilon\delta\pmod{p}$ and $c_2=\det \xi=p^hn$. We distinguish:

%[[Beachte: 1.) $\Mat{\alpha}{\beta}{\gamma}{\delta},\Mat{a}{b}{c}{d},\sigma_{[\epsilon,\delta]},\alpha$ haben alle Determinante $=1$, da in ${\rm, SL}_2({\Bbb Z})$ enthalten $\Rightarrow {\rm det}\,\xi={\rm det}\,\Mat{1}{}{}{p^hn}=p^hn$.

%2.) (N) gilt sowohl im elliptischen wie im hyperbolischen Fall; wir unterscheiden erst spater diese beiden Falle.]]

Case A: $\xi$ is elliptic. We denote by $F/{\Bbb Q}$ a minimal splittig field for $\chi_\xi$ and $\tau:\,\alpha\mapsto{^\tau}\alpha$ is the non-trivial automorphism of $F/{\Bbb Q}$. We further denote by $\lambda,{^\tau}\lambda\in F$ the zeros of $\chi_\xi$. Since $\chi_\xi\equiv (T-\epsilon\delta)T\pmod{p}$ we deduce that $p$ splits in $F$, i.e. $(p)={\mathfrak p}^{\tau}{\mathfrak p}$. On the other hand, $(p)_{{\cal O}_{{\Bbb Q}_{\zeta_{p-1}}}}={\mathfrak p}_0$ is totally inert in ${\Bbb Q}(\zeta_{p-1})/{\Bbb Q}$, hence, $F$ is not contained in ${\Bbb Q}(\zeta_{p-1})$. 
%
%[[da die eine Erweiterung uber $p$ rein zerlegt ist und die andere vollig unzerlegt ist]]
%
In particular, $\lambda, {^\tau}\lambda$
are not contained in ${\Bbb Q}(\zeta_{p-1})$. Since the minimal polynomial $\chi_\xi^*$ of $\lambda\in F(\zeta_{p-1})$ over ${\Bbb Q}(\zeta_{p-1})$ is a divisor of $\chi_\xi$ 
%
%[[da $\chi_\xi\in{\Bbb Q}(\zeta_{p-1})[T])$ und $\chi_\xi$ $\lambda$ anulliert $\Rightarrow \chi_\xi$ teilt das minmalpolynom von $\lambda$ uber ${\Bbb Q}(\zeta_{p-1})$]]
%
we thus obtain that $[F(\zeta_{p-1}):{\Bbb Q}(\zeta_{p-1})]= 2$ and $\chi_\xi$ equals $\chi_\xi^*$. We also denote by $\tau$ the non-trivial
automorphism of $F(\zeta_{p-1})/{\Bbb Q}(\zeta_{p-1})$. Since $\chi_\xi^*=\chi_\xi\equiv (T-\epsilon\delta)T\pmod{{\mathfrak p}_0}$ the prime ${\mathfrak p}_0={\mathfrak P}{^\tau}{\mathfrak P}$ is split in $F(\zeta_{p-1})$:
$$
\begin{array}{ccrc cclcc}
{\mathfrak P}\bar{\mathfrak P}&&&&F(\zeta_{p-1})&&&\qquad&{\mathfrak P}\bar{\mathfrak P}\\
\parallel&&&/&&\backslash&&\qquad&\parallel\\
{\mathfrak p}_0&\qquad&{\Bbb Q}(\zeta_{p-1})&&&&F\ni\lambda,{^\tau\lambda}&\qquad&{\mathfrak p}\bar{\mathfrak p}\\
\parallel&&&\backslash&&/&&\qquad&\parallel\\
(p)&&&&{\Bbb Q}&&&\qquad&(p).\\
\end{array}
$$  

Since ${\rm tr}\,(\xi^\iota|L_k)={\rm tr}\,(\xi^\iota|L_k\otimes F(\zeta_{p-1}))$ we may as well compute the trace of $\xi^\iota$ on $L_k\otimes F(\zeta_{p-1})$. But $\chi_\xi$ splits over $F$, hence, $\xi$ is diagonalizable over $F$, i.e.
$$
\xi=g\Mat{\lambda}{}{}{^\tau\lambda}g^{-1}
$$ 
with $g\in {\rm GL}_2(F)$. Thus, applying the involution $\iota$, we obtain
$$
\xi^\iota=(g^\iota)^{-1}\Mat{^\tau\lambda}{}{}{\lambda}g^\iota.\leqno(34)
$$
Using $(34)$ we immediately derive (e.g. use Weyl's character formula)
\begin{eqnarray*}
(35)\quad\omega^{-k}(\epsilon\delta)\,{\rm tr}\,(\xi^\iota|L_k)
%&=&\omega^{-k}(\epsilon\delta)\,\frac{\lambda^{k+1}-\bar{\lambda}^{k+1}}{\lambda-\bar{\lambda}}\\
&=&\frac{(\omega^{-1}(\epsilon\delta)\lambda)^{k+1}-(\omega^{-1}(\epsilon\delta)^\tau{\lambda})^{k+1}}{\omega(\epsilon\delta)\lambda-\omega(\epsilon\delta)^\tau{\lambda}}.\\
%&=&(\omega^{-1}(\epsilon\delta\lambda)^k+(\omega^{-1}(\epsilon\delta)\bar{\lambda})^k+\sum_{i=1}^{k-1}
%(\omega^{-1}(\epsilon\delta)\lambda)^i\,(\omega^{-1}(\epsilon\delta\bar{\lambda})^{k-i}.\\
\end{eqnarray*}
%
%[[Details dazu: Es gilt $(\Mat{a}{}{}{b}^\iota)^\iota=\Mat{b}{}{}{a}$
%$\Mat{b}{}{}{a}X^{k-i}Y^i=X^{k-i}Y^i(\Mat{b}{}{}{a}(X,Y)})=X^{k-i}Y^i(b X,a Y)=b^{k-i}a^i X^{k-i}Y^i$; damit folgt ${\rm tr}\,(\Mat{a}{}{}{b}|L_k)=\sum_{i=0}^k
%a^ib^{k-i}$
%]]

To evaluate $(35)$ we look closer at $\lambda,{^\tau}\lambda$. Since $\lambda, {^\tau}\lambda$ in particular satisfy $\chi_\xi$ modulo ${\mathfrak p}$ and $\chi_\xi\equiv (T-\epsilon\delta)T\pmod{\mathfrak p}$ we find (after eventually permuting $\lambda$ and ${^\tau}\lambda$)
$$
(\lambda-\epsilon\delta){^\tau}\lambda\equiv 0\pmod{\mathfrak p}.
$$
Since ${\cal O}_F/{\mathfrak p}\cong {\Bbb Z}/(p)$ is a field 
%
%[[(i.e. there are no divisors of $0$)]] 
%
we deduce $\lambda\equiv \epsilon\delta$, ${^\tau}\lambda\equiv 0\pmod{\mathfrak p}$ and since $\lambda{^\tau}\lambda=p^hn\sim {\mathfrak p}^h{^\tau}{\mathfrak p}^h(n)$ we finally obtain
$$
\lambda\equiv \epsilon\delta\pmod{{\mathfrak p}},\;\;{^\tau}\lambda\equiv 0\pmod{{\mathfrak p}^h}.\leqno(36)
$$
We now evaluate $(35)$ as follows.

$\bullet$ Since $\omega(\epsilon\delta)\equiv \epsilon\delta\pmod{p}$ and ${\mathfrak P}$ divides $p$ as well as ${\mathfrak p}$, equation (36) implies $\omega^{-1}(\epsilon\delta)\lambda\equiv 1\pmod{\mathfrak P}$ hence, taking into account that $k\equiv k'\pmod{p^m}$ we see that $(\omega^{-1}(\epsilon\delta)\lambda)^{k+1}\equiv (\omega^{-1}(\epsilon\delta)\lambda)^{k'+1}\pmod{{\mathfrak P}^m}$. (Note that because ${\mathfrak P}|(p)$ multiplication by $p$ annihilates $1+{\mathfrak P}^i/1+{\mathfrak P}^{i+1}\cong{\mathfrak P}^i/{\mathfrak P}^{i+1}$)

$\bullet$ Moreover, $\omega^{-1}(\epsilon\delta){^\tau}\lambda\equiv 0\pmod{{\mathfrak p}^h}$, hence, $(\omega^{-1}(\epsilon\delta){^\tau}\lambda)^{k+1}\equiv 0\pmod{{\mathfrak p}^{hk}}$. Analogously, we have $(\omega^{-1}(\epsilon\delta){^\tau}\lambda)^{k'+1}\equiv 0\pmod{{\mathfrak p}^{hk'}}$. 

%[[Hier verschenken wir ein Potenzen von $p$ d,h, es gilt sogar $(\omega^{-1}(\epsilon\delta){^\tau}\lambda)^{k+1}\equiv 0\pmod{{\mathfrak p}^{h(k+1)}}$ und $(\omega^{-1}(\epsilon\delta){^\tau}\lambda^{k'+1})^{k'}\equiv 0\pmod{{\mathfrak p}^{h(k'+1)}}$
%]]

$\bullet$ $\omega(\epsilon\delta)^{-1}\lambda\equiv 1\pmod{\mathfrak p}$ and $\omega^{-1}(\epsilon\delta){^\tau\lambda}\equiv 0\pmod{\mathfrak p}$ imply that $\omega^{-1}(\epsilon\delta)\lambda-\omega^{-1}(\epsilon\delta){^\tau}\lambda$ is a ${\mathfrak p}$-adic unit.

Using these facts and taking into account that the ramification index of ${\mathfrak P}|{\mathfrak p}$ equals $1$ it is immediate to see that
$$
\frac{(\omega^{-1}(\epsilon\delta)\lambda)^{k+1}-(\omega^{-1}(\epsilon\delta){^\tau\lambda})^{k+1}}{\omega^{-1}(\epsilon\delta)\lambda-\omega^{-1}(\epsilon\delta){^\tau\lambda}}
\equiv
\frac{(\omega^{-1}(\epsilon\delta)\lambda)^{k'+1}-(\omega^{-1}(\epsilon\delta){^\tau\lambda})^{k'+1}}{\omega^{-1}(\epsilon\delta)\lambda-\omega^{-1}(\epsilon\delta){^\tau\lambda}}
\pmod{{\mathfrak P}^{{\rm min}\,(m,hk,hk')}}.
$$
Hence, by using (35), which also holds with $k$ replaced by $k'$ we find
$$
\omega^{-k}(\epsilon\delta)\,{\rm tr}\,(\xi^\iota|L_k)\equiv \omega^{-k'}(\epsilon\delta)\,{\rm tr}\,(\xi^\iota|L_{k'})\pmod{{\mathfrak P}^{{\rm min}(m,hk,hk')}}.\leqno(37)
$$ 
Applying $\tau\in{\rm Gal}(F(\zeta_{p-1})/{\Bbb Q}(\zeta_{p-1}))$ to (37) and noting that ${\rm tr}\,(\xi^\iota|L_k),{\rm tr}\,(\xi^\iota|L_{k'})\in{\Bbb Q}$ because the representations $L_k,L_{k'}$ are defined over ${\Bbb Q}$,
%
%[[and $\xi^\iota\in G({\Bbb Q})$ ($L_k,L_{k'}$ are even defined over ${\Bbb Z}$)]]
%
we obtain
$$
\omega^{-k}(\epsilon\delta)\,{\rm tr}\,(\xi^\iota|L_k)\equiv \omega^{-k'}(\epsilon\delta)\,{\rm tr}\,(\xi^\iota|L_{k'})\pmod{^\tau{\mathfrak P}^{{\rm min}(m,hk,hk')}}.\leqno({^\tau}37)
$$ 
Since ${\mathfrak P}{^\tau}{\mathfrak P}=(p)$, $(37)$ and $({^\tau 37})$ finally yield the claim in the elliptic case.

Case B: $\xi$ hyperbolic. Again, we denote by $\lambda,\lambda'\in{\Bbb Q}$ the zeroes of $\chi_\xi$. Equation (33) shows that over ${\Bbb Z}/p{\Bbb Z}$ the characteristic polynomial of $\xi$ decomposes $\chi_\xi=T(T-\epsilon\delta)$,
hence, $\lambda\equiv\epsilon\delta\pmod{p}$ and $\lambda'\equiv 0\pmod{p}$. 
%
%[[beachte: ${\Bbb Z}/p{\Bbb Z}$ ist KORPER, d.h. ist Nullteilerfrei bzw. es gibt hochstens $2$ Nullstellen modulo $p$]]
%
In particular, $\lambda$ is a $p$-adic unit and $\lambda\lambda'=p^hn$ even implies $\lambda'\equiv 0\pmod{p^h}$. As above, we deduce that
\begin{itemize}

\item $\omega^{-1}(\epsilon\delta)\lambda\equiv 1\pmod{p}$, hence, $(\omega^{-1}(\epsilon\delta)\lambda)^{k+1}\equiv (\omega^{-1}(\epsilon\delta)\lambda)^{k'+1}\pmod{p^m}$

\item $(\omega^{-1}(\epsilon\delta)\lambda')^{k+1}\equiv 0\pmod{p^{hk}}$ and $(\omega^{-1}(\epsilon\delta)\lambda')^{k'+1}\equiv 0\pmod{p^{hk'}}$

\item $\omega^{-1}(\epsilon\delta)\lambda-\omega^{-1}(\epsilon\delta)\lambda'\equiv 1\pmod{p}$ is a $p$-adic unit.

\end{itemize}

$\xi$ is $G({\Bbb Q})$-conjugate to $\Mat{\lambda}{}{}{\lambda'}$, hence, $\xi^\iota$ is $G({\Bbb Q})$-conjugate to $\Mat{\lambda'}{}{}{\lambda}$ and 
$\omega^{-1}(\epsilon\delta)^k\,{\rm tr}\,(\xi|L_k)$ therefore is given by the expression in (35) with
$^\tau{\lambda}$ replaced by $\lambda'$. As in the case $\xi$ elliptic we thus obtain 
$$
\omega^{-k}(\epsilon\delta)\,{\rm tr}\,(\xi^\iota|L_k)\equiv \omega^{-k'}(\epsilon\delta)\,{\rm tr}\,(\xi^\iota|L_{k'})\pmod{p^{{\rm min}(m,hk,hk')}},
$$ 
which is the claim in the hyperbolic case. Thus, the proof of the Lemma is complete.

We recall that in section 3.4 we introduced the constant $N(\alpha)$, which only depends on $\alpha$ and which bounds the degree of the polynomials $p_{\alpha,k,k'}$.

{\bf Proposition 1. }{\it Fix $\alpha\in{\Bbb Q}_{\ge 0}$ and let $C\in{\Bbb N}$ be any integer. Assume that $k,k'\in{\Bbb N}$ satisfy 

- $k,k'\ge (C\alpha+1)^2 +2$

- $k\equiv k'\pmod{p^m}$ with $m\ge C\alpha+1$. 

Set $L=[\frac{m}{C\alpha}]$ if $\alpha>0$ and $L=p^m$ if $\alpha=0$ and let $n,u\in{\Bbb N}$ with $u|Np$; then the following congruence holds true: 

$$
\omega^{-k}(\delta)\,{\rm tr}\,S_{u,\delta} T_n e_\alpha^{L}|_{{\cal M}_k(\Gamma,\chi\omega^{-k})}\equiv \omega^{-k'}(\delta)\,{\rm tr}\,S_{u,\delta} T_n e_\alpha^{L}|_{{\cal M}_{k'}(\Gamma,\chi\omega^{-k'})}\pmod{p^{(1-\frac{N(\alpha)}{C})m-v_p(\varphi(N))}}
$$
if $\alpha>0$ and

$$
\omega^{-k}(\delta)\,{\rm tr}\,S_{u,\delta} T_n e_\alpha^{L}|_{{\cal M}_k(\Gamma,\chi\omega^{-k})}\equiv \omega^{-k'}(\delta)\,{\rm tr}\,S_{u,\delta} T_n e_\alpha^{L}|_{{\cal M}_{k'}(\Gamma,\chi\omega^{-k'})}\pmod{p^{m-v_p(\varphi(N))}}
$$
if $\alpha=0$. Here, $e_\alpha=e_{k,k',\alpha}$ is the quasi idempotent defined in section 3.4.

}

%[[Beachte: $n$ ist nicht notwendig prim zu $Np$; diese Bedingung wird erst bei der konstruktion der $p$-adischen Familien von Modulformen gebraucht damit $\pi^N_M$ und $T_\ell$ FUR $\ell\not|Np$ kommutieren, also SIMULTAN diagobar sind $\Rightarrow$ ${\rm tr}\,\pi^N_MT_\ell$
%ist bestimmbar]]

{\it Proof. } Using the definition of $e_\chi$ (cf. section 3.1) and the isomorphism ${\cal M}_k(\Gamma)\cong H^1(S(\Gamma),L_{k-2})$ it is immediate that for $k\ge 2$ 
$$
{\rm tr}\,S_{u,\delta} T_n \,e_\alpha^{L}|_{{\cal M}_k(\Gamma,\chi\omega^{-k})}
={\rm tr}\,S_{u,\delta} T_n\,e_\alpha^{L} e_{\chi\omega^{-k}}|_{{\cal M}_k(\Gamma)}
={\rm tr}\, S_{u,\delta} T_n\,e_\alpha^{L} e_{\chi\omega^{-k}}|_{H^1(S(\Gamma),L_{k-2,{\Bbb C}})}.
$$ 
(recall that we have written the action of ${\cal H}$ on ${\cal M}_k(\Gamma)$ from the left, i.e. we first apply $e_{\chi\omega^{-k}}$).
%
%[[beachte: ${\cal H}$ operiert auf ${\cal M}_k(\Gamma,\chi;)$, deshalb ist $e_i\circ e_\chi$ wohldefiniert (d.h. $e_i$ lasst das bild von $e_\chi$ invariant)
%2.) ${\cal H}$ operiert von RECHTS auf ${\cal M}_k(\Gamma)$, d.h. $T_1T_2$ bedeutet zuerst $T_1$, dann $T_2$ anwenden
%
%
Since $H^i(S(\Gamma),L_{k-2,{\Bbb C}})=0$ in degrees $i=0,2$ if $k> 2$, 
%
%[[Im Grad $i=2$ verschwindet die Kohomologie immer; im Grad $i=0$ ist $H^0(...)=0$ falls $k>2$ und $H^0(...)={\Bbb C}$ falls $k=2$, d.h. wenn residuelles Spektrum existiert.]]
%
the above equation further implies that
$$
{\rm tr}\,S_{u,\delta} T_n\,e_\alpha^L |_{{\cal M}_k(\Gamma,\chi\omega^{-k})} 
=L(S_{u,\delta} T_n\,e_\alpha^{L}\,e_{\chi\omega^{-k}}|_{H^\bullet(S(\Gamma),L_{k-2,{\Bbb C}})}).\leqno(38)
$$

We recall some of the notations that we introduced in connection with the quasi idempotent $e_{k,k',\alpha}$: 

- $E=E_{k,k'}^{\alpha+1}$ is obtained from ${\Bbb Q}$ by adjoining the roots of the characteristic polynomial of $T_p$ acting on ${\cal M}_k(\Gamma,\chi\omega^{-k})^{\le\alpha+1}$ and of the characteristic polynomial of $T_p$ acting on ${\cal M}_{k'}(\Gamma,\chi\omega^{-k'})^{\le\alpha+1}$

- ${\mathfrak p}={\mathfrak p}_{k,k'}$ is the prime ideal in ${\cal O}_{E}$ corresponding to $v_p$ 

- $e=e_{k,k'}$ is the ramification degree of ${\mathfrak p}|p$ 

- $i=i_{\alpha,k,k'}$ is defined by $\alpha=i/e_{k,k'}$ 

- $\varpi=\varpi_{k,k'}\in {\cal O}_E$ is a local prime, i.e. $v_p(\varpi)=1/e$. 

- $e_\alpha=e_{k,k',\alpha}=p_\alpha(T_p/\varpi^i)$ ($i=i_{\alpha,k,k'}$), where  $p_\alpha=p_{\alpha,k,k'}=\sum_{h=1}^{t_\alpha} b_h X^h$ is a polynomial with integer coefficients $b_h\in{\cal O}_E$ (cf. Proposition 2.) in section 3.4) and degree $t_\alpha=t_{\alpha,k,k'}$. Hence, we obtain
$$
e_\alpha^{L}=\sum_{h={L}}^{{L}t_\alpha} \frac{a_h}{\varpi^{ih}} T_p^h\leqno(39)
$$ 
with $a_h\in{\cal O}_E$. 
%
%[[$a_h$ are given by $p_\alpha^L=\sum_{h=L}^{Lt_\alpha} a_hX^h$]]
%
%[[Es gilt sogar $X^{s(q-1)}|p_\alpha$ nicht nur $X|p_\alpha$; bringt das eine Verbesserung ??
%]]
%
%
Using equation (13) and taking into account that $T_nT_p^h=T_{p^hn}$ we then find
\begin{eqnarray*}
S_{u,\delta} T_n\,e_\alpha^L\, e_{\chi\omega^{-k}}
&=&\frac{1}{\varphi(Np)}\,\sum_\epsilon \sum_{h={L}}^{Lt_\alpha} \chi\omega^{-k}(\epsilon) \frac{a_h}{\varpi^{ih}}\,S_{u,\delta} T_n T_p^h \langle \epsilon\rangle\\
&=&\frac{1}{\varphi(Np)}\,\sum_\epsilon \sum_{h=L}^{Lt_\alpha} \sum_j \chi\omega^{-k}(\epsilon) m_j \frac{a_h}{\varpi^{ih}}\, \Gamma\sigma_{\epsilon}\Mat{1}{}{}{p^{h}n}\alpha_j\Mat{1}{}{u}{1}\sigma_\delta\Gamma,\\
\end{eqnarray*}
where $m_j\in{\Bbb N}\cup\{0\}$ and $\alpha_j\in \Gamma$. The trace formula (cf. the Theorem in section 2.5) together with equation (38) therefore yields

\begin{eqnarray*}
(40)&&\omega^{-k}(\delta)\,{\rm tr}\, S_{u,\delta} T_n \,e_\alpha^{L}|_{{\cal M}_k(\Gamma,\chi\omega^{-k})}\\
&&=\frac{1}{\varphi(Np)}\,\sum_\epsilon\sum_{h={L}}^{Lt_\alpha}\sum_j\chi(\epsilon) m_j \frac{a_h}{\varpi^{ih}}\,\left(\sum_{[\xi]\in(\Gamma\sigma_\epsilon\Mat{1}{}{}{p^hn}\alpha_j\Mat{1}{}{u}{1}\sigma_\delta\Gamma)_{\rm ell}/\sim_\Gamma}\omega^{-k}(\epsilon\delta)\,{\rm tr}\,(\xi^\iota|L_{k-2})\right.\\
&&\qquad\left.+\sum_{[\eta]\in(\Gamma\sigma_\epsilon\Mat{1}{}{}{p^hn}\alpha_j\Mat{1}{}{u}{1}\sigma_\delta\Gamma)_{\rm hyp}^+/\sim_\Gamma}\omega^{-k}(\epsilon\delta)\,{\rm tr}\,(\eta^\iota|L_{k-2})\right).
\end{eqnarray*}

Of course, equation (40) also holds true if we replace $k$ with $k'$. If 
$$
[\zeta]\in(\Gamma\sigma_\epsilon\Mat{1}{}{}{p^hn}\alpha_j\Mat{1}{}{u}{1}\sigma_\delta\Gamma)_{\rm ell}/\sim_\Gamma\quad\mbox{or}\quad [\zeta]\in(\Gamma\sigma_\epsilon\Mat{1}{}{}{p^hn}\alpha_j\Mat{1}{}{u}{1}\sigma_\delta\Gamma)_{\rm hyp}^+/\sim_\Gamma
$$
then Lemma 1 implies that
$$
m_j\,\frac{a_h}{\varpi^{ih}}\,\omega^{-k}(\epsilon\delta)\,{\rm tr}\,(\zeta^\iota|L_{k-2})\equiv m_j\,\frac{a_h}{\varpi^{ih}}\,\omega^{-k'}(\epsilon\delta)\,{\rm tr}\,(\zeta^\iota|L_{k'-2})\pmod{{\mathfrak p}^{e\clubsuit}}.\leqno(41)
$$
Here 
$$
\clubsuit\ge{\rm min}(m,h(k-2),h(k'-2))-\frac{i}{e}h\leqno(42)
$$ 
(note that $m_j\in{\Bbb N}$ and $a_h\in{\cal O}_{E_{k,k'}^\alpha}$ are integers). 

We now distinguish cases and first assume that $\alpha>0$, hence, $L=\left[\frac{m}{C\alpha}\right]$. An easy calculation shows that $m\ge C\alpha+1$ implies  $\frac{m}{m-C\alpha}\le C\alpha+1$. Hence, we obtain
$$
k-2\ge (C\alpha+1)^2\ge(C\alpha+1)\frac{m}{m-C\alpha}>C\alpha\frac{m}{m-C\alpha}>\frac{m}{\frac{m}{C\alpha}-1}\ge \frac{m}{L};
$$
here, the last inequality holds because $L\ge m/(C\alpha)-1$. 
%
%[[note that $m-(C\alpha+1)\ge 0$ and use that $C\alpha+1>C\alpha$]]
%
Thus, we find
$$
(k-2)h\ge (k-2)L\ge \frac{m}{L} L=m.
$$
Of course, since also $k'-2\ge (C\alpha+1)^2$, the same equation with $k$ replaced by $k'$ holds. 
In particular, ${\rm min}(m,h(k-2),h(k'-2))=m$ and 
$$
\clubsuit\ge m-\frac{i}{e} h=m-\alpha h.
$$ 
Since further $h\le t_\alpha L$ and $t_\alpha=t_{\alpha,k,k'}\le N(\alpha)$ ( cf. Proposition 2.) in section 3.4) we obtain
$$
\clubsuit\ge m-\alpha t_\alpha L\ge m-\alpha N(\alpha)L\ge m-\alpha N(\alpha)(\frac{m}{C\alpha})=(1-\frac{N(\alpha)}{C})m\leqno(43)
$$
(note that $L\le m/(C\alpha)$). Equations (41) and (43) imply that
$$
m_j\frac{a_h}{\varpi^{ih}}\,\omega^{-k}(\epsilon\delta)\,{\rm tr}\,(\zeta^\iota|L_{k-2})\equiv m_j\frac{a_h}{\varpi^{ih}}\,\omega^{-k'}(\epsilon\delta)\,{\rm tr}\,(\zeta^\iota|L_{k'-2})\pmod{{\mathfrak p}^{e(1-\frac{N(\alpha)}{C})m}}.\leqno(44)
$$
Equations (40) and (44) taken together imply the claim in case $\alpha>0$. 

Finally, we look at the case $\alpha=0$, i.e. $L=p^m$. In this case $\beta_1=0=\alpha$ (cf. equation (17)); moreover, $i=i_{0,k,k'}=0$ and, hence, 
$$
e_{k,k',0}=\tilde{e}_{1,k,k'}=T_p^{q-1}
$$
(cf. equation (21) in section 3.4). We deduce that $p_0=p_{0,k,k'}=X^{q-1}$, i.e. the sum over $h$ in equation (40) consists of only one term corresponding to $h=(q-1)L=(q-1)p^m$. Since furthermore $L=p^m\ge m$ we know that $h\ge L\ge m$ and, hence, $h(k-2)\ge m(k-2)\ge m$ and $h(k'-2)\ge m(k'-2)\ge m$; equation (42) thus yields
$$
\clubsuit\ge m
$$
(note that $i=i_{0,k,k'}=0$). Equation (41) therefore becomes
$$
m_j\omega^{-k}(\epsilon\delta)\,{\rm tr}\,(\zeta^\iota|L_{k-2})\equiv m_j\omega^{-k}(\epsilon\delta)\,{\rm tr}\,(\zeta^\iota|L_{k'-2})\pmod{{\mathfrak p}^{em}}\leqno(41^0)
$$
(note that $p_0=X^{q-1}$). Equations ($41^0$) and (40) taken together yield the claim in case $\alpha=0$. Thus, the proof of the Proposition is complete.

If $\alpha\in{\Bbb Q}$ with $\alpha>0$ we set 
$$
{{\sf C}_\alpha}=[N(\alpha)+\alpha^{-1}]+1,
$$ 
hence, ${{\sf C}_\alpha}\in{\Bbb N}$ and 
$$
{{\sf C}_\alpha}\ge N(\alpha)+\frac{1}{\alpha}.\leqno(45)
$$
In particular, ${\sf C}_\alpha\alpha\ge N(\alpha)\alpha+1\ge 1$ (cf. equation (32)) and we find that
$$
\left[\frac{m}{{\sf C}_\alpha\alpha}\right]\le\left[m\right]\le m.
$$
If $\alpha=0$ we set ${\sf C}_\alpha=0$. 

%[[We could as well have assigned any (finite) value to ${\sf C}_\alpha$ because multiplication of ${\sf C}_\alpha$ by $\alpha(=0)$ in the statements of the following Lemma will make it disappear anyway.]]

{\bf Lemma 2. } {\it 1.) The operator $S_{u,\delta}$, where $u\in{\Bbb N}\cup\{0\}$ with $p|u$, commutes with the operators $T_\ell$, where $\ell\not|N$ and with $e_\alpha$. 

2.) Let $A_1,\ldots,A_k,B\in {\rm End}_{\cal O}({\cal O}^n)$ be commuting operators and assume furthermore that all eigenvalues $\gamma$ of $B$ satisfy $v_p(\gamma)\ge c$ resp. all eigenvalues $\gamma$ of $B$ satisfy $v_p(\gamma-1)\ge c$. Then
$$
{\rm tr}\,A_1\cdots A_kB\equiv 0\pmod{p^c}
$$
resp.
$$
{\rm tr}\,A_1\cdots A_kB\equiv {\rm tr}\,A_1\cdots A_k\pmod{p^c}.
$$
}

{\it Proof. } 1.) Since $e_\alpha\in{\Bbb C}[T_p]$ it is sufficient to show that $S_{u,\delta}$, $p|u$ commutes with the Hecke operators $T_\ell$, $\ell\not|N$. Furthermore, since $\langle\delta\rangle$ commutes with the Hecke operators $T_\ell$ (cf. [D-S], p. 169), it is sufficient to show that $S_u=S_{u,1}$ ($p|u$) commutes with the Hecke operators $T_\ell$, $\ell\not|N$.
To verify this we use the adelic description of modular forms and Hecke operators. Let $F:\,{\rm GL}_2({\Bbb Q})\backslash {\rm GL}_2({\Bbb A})\rightarrow {\Bbb C}$ be the adelic function corresponding to the modular form $f\in{\cal M}_k(\Gamma)$, i.e. $f(z)=F(1,\ldots,1,g_\infty)$, where $g_\infty i=z$. We write
$$
\Gamma\Mat{1}{}{u}{1}\Gamma=\bigcup_i \Gamma\Mat{1}{}{u}{1}\gamma_i
$$
with $\gamma_i\in\Gamma$. We then obtain 
$$
S_u F(1,\ldots,1,g_\infty)=\sum_i F(\Mat{1}{}{-u}{1},\ldots,\Mat{1}{}{-u}{1},g_\infty).
$$
Using $K_1(Np)$-right invariance of $F$ and taking into account that $\Mat{1}{}{-u}{1}\in K_{1,\ell}(Np)={\rm GL}_2({\Bbb Z}_\ell)$ for all primes $\ell$, which do not divide $Np$ and even $\Mat{1}{}{-u}{1}\in K_{1,p}(Np)$ if $p|u$ we further obtain
$$
S_u F(1,\ldots,1,g_\infty)=\sum_i F(1,\ldots,1,\Mat{1}{}{-u}{1},\ldots,\Mat{1}{}{-u}{1},g_\infty).
$$
Here, the matrices $\Mat{1}{}{-u}{1}$ only appear in positions corresponding to primes $\ell$ dividing $N$. Since
$$
T_\ell F(1,\ldots,1,g_\infty)=\sum_x F(1,\ldots,1,x,g_\infty),
$$
where $x$ runs over certain matrices, which appear in the position corresponding to the prime $\ell$, it is obvious that $S_u$ with $p|u$ and the operators $T_\ell$, $\ell\not|N$ commute.

%[[Beachte, dass $\det K_{1,\ell}(Np)={\Bbb Z}_\ell^*$, d.h. es gibt nur eine Zusammenhangskomponente und wir mussen nicht $F(1,\ldots,1,\eta,g_\infty)$ betrachtem, wobei $\eta\in{\Bbb Q}^*\backslash {\Bbb A}^*/\det K_1(Np)$ die Zusammenhangskomponenten durchlauft, betrachten]]

2.) Since the operators $A_1,\ldots,A_k,B$ commute there is a basis ${\cal B}$ such that 
$$
{\cal D}_{\cal B}(A_h)=\left(
\begin{array}{ccc}
a_{h;1}&&*\\
&\ddots&\\
&&a_{h;n}\\
\end{array}
\right),\quad h=1,\ldots,k,
\qquad
{\cal D}_{\cal B}(B)=\left(
\begin{array}{ccc}
b_1&&*\\
&\ddots&\\
&&b_n\\
\end{array}
\right)
$$
all are upper triangular. Morover, by assumption we know that $v_p(b_i)\ge c$ for all $i=1,\ldots,n$ resp. $v_p(b_i-1)\ge c$ for all $i=1,\ldots,n$ and all $a_{h;i}$ are integral. This implies that
$$
{\rm tr}\,A_1\cdots A_kB=\sum_{i=1}^n \prod_{h=1}^k a_{h;i}b_i
$$
is congruent to $0$ resp. to ${\rm tr}\,A_1\cdots A_k$ modulo $p^c$. Thus, the lemma is proven.

Let $u\in{\Bbb N}$ be divisible by $p$; since the operator $S_{u,\delta}$ in particular commutes with the Hecke operator $T_p$ it leaves the slope subspaces invariant, i.e. we obtain a map
$$
S_{u,\delta}:\,{\cal M}_k(\Gamma,\chi\omega^{-k})^\alpha\rightarrow {\cal M}_k(\Gamma,\chi\omega^{-k})^\alpha.
$$

{\bf Proposition 2. }{\it Fix $\alpha\in{\Bbb Q}_{\ge 0}$. Assume that $k,k'\in{\Bbb N}$ satisfy 

$\bullet$ $k,k'\ge ({{\sf C}_\alpha}\alpha+1)^2+2 $

$\bullet$ $k\equiv k'\pmod{p^m}$ with $m\ge {{\sf C}_\alpha}\alpha+1$. 

Then, for any $u\in{\Bbb N}$, which is divisible by $p$, the following congruence holds true 
$$
\omega^{-k}(\delta)\,{\rm tr}\,S_{u,\delta} T_n |_{{\cal M}_k(\Gamma,\chi\omega^{-k})^\alpha}\equiv \omega^{-k'}(\delta)\,{\rm tr}\,S_{u,\delta} T_n |_{{\cal M}_{k'}(\Gamma,\chi\omega^{-k'})^\alpha}\pmod{p^{[\log_p [\frac{m}{{{\sf C}_\alpha}\alpha}]]-v_p(\varphi(N))}}
$$
if $\alpha>0$ and
$$
\omega^{-k}(\delta)\,{\rm tr}\,S_{u,\delta} T_n |_{{\cal M}_k(\Gamma,\chi\omega^{-k})^\alpha}\equiv \omega^{-k'}(\delta)\,{\rm tr}\,S_{u,\delta} T_n |_{{\cal M}_{k'}(\Gamma,\chi\omega^{-k'})^\alpha}\pmod{p^{m-v_p(\varphi(N))}}
$$
if $\alpha=0$. 

Moreover in case $\alpha>0$ it also holds that
$$
\omega^{-k}(\delta)\,{\rm tr}\,S_{u,\delta} T_n e_\alpha^{[\frac{m}{{{\sf C}_\alpha}\alpha}]}|_{{\cal M}_k(\Gamma,\chi\omega^{-k})^\alpha}\equiv \omega^{-k'}(\delta)\,{\rm tr}\,S_{u,\delta} T_n e_\alpha^{[\frac{m}{{{\sf C}_\alpha}\alpha}]}|_{{\cal M}_{k'}(\Gamma,\chi\omega^{-k'})^\alpha}\pmod{p^{[\frac{m}{{{\sf C}_\alpha}\alpha}]-v_p(\varphi(N))}}.
$$
}

{\it Proof. } We first assume $\alpha>0$ and prove the first congruence. As in Proposition 1 we set $L=\left[\frac{m}{{{\sf C}_\alpha}\alpha}\right]$. We first note that ${{\sf C}_\alpha}\ge N(\alpha)+\frac{1}{\alpha}$ implies $\frac{m}{{{\sf C}_\alpha}\alpha}\le (1-\frac{N(\alpha)}{{{\sf C}_\alpha}})m$. Since $L\le \frac{m}{{{\sf C}_\alpha}\alpha}$ we thus obtain
$$
L\le (1-\frac{N(\alpha)}{{{\sf C}_\alpha}})m.\leqno(46)
$$
Proposition 1.) in section 3.4 implies that all eigenvalues $\gamma$ of $e_\alpha^L$ acting on ${\cal M}_k(\Gamma,\chi\omega^{-k})^\beta$, $\beta\not=\alpha$, satisfy $v_p(\gamma)\ge L$. Hence, Lemma 2 shows that
$$
\omega^{-k}(\delta)\,{\rm tr}\,S_{u,\delta} T_n e_\alpha^L|_{{\cal M}_k(\Gamma,\chi\omega^{-k})^\beta}\equiv 0\pmod{p^L}
$$
for all $\beta\not=\alpha$. 
%
%[[Note that $e_\alpha^L|_{V^\beta}\in{\mathfrak p}_{k,k'}{\rm End}(V^\beta)$]]
%
Moreover, again by Proposition 1.) in section 3.4 we know that all eigenvalues $\gamma$ of $e_\alpha^L$ acting on ${\cal M}_k(\Gamma,\chi\omega^{-k})^\alpha$ satisfy $v_p(\gamma-1)\ge [\log_p L]$ (note that $(1+{\mathfrak p}^a)/(1+{\mathfrak p}^{a+1})\cong ({\cal O}/{\mathfrak p},+)$ and that multiplication with $p$ annihilates ${\cal O}/{\mathfrak p}$, i.e. $(1+{\mathfrak p}^a)^p\le 1+{\mathfrak p}^{a+1}$). Hence, Lemma 2 shows that
$$
\omega^{-k}(\delta)\,{\rm tr}\,S_{u,\delta} T_n e_\alpha^L|_{{\cal M}_k(\Gamma,\chi\omega^{-k})^\alpha}
\equiv 
\omega^{-k}(\delta)\,{\rm tr}\,S_{u,\delta} T_n |_{{\cal M}_k(\Gamma,\chi\omega^{-k})^\alpha}\pmod{p^{[\log_p L]}}.
$$
We thus obtain
$$
\omega^{-k}(\delta)\,{\rm tr}\,S_{u,\delta} T_n |_{{\cal M}_k(\Gamma,\chi\omega^{-k})^\alpha}
\equiv
\omega^{-k}(\delta)\,{\rm tr}\,S_{u,\delta} T_n e_\alpha^L|_{{\cal M}_k(\Gamma,\chi\omega^{-k})}\pmod{p^{[\log_p L]}}.\leqno(47)
$$
Of course, the same congruence holds true if we replace $k$ by $k'$. The claim now follows from (47) and the first congruence of Proposition 1 taking into account that 
\begin{eqnarray*}
[\log_p L]-v_p(\varphi(N))&\le& L-v_p(\varphi(N))\\
&\stackrel{(46)}{\le}& (1-\frac{N(\alpha)}{{{\sf C}_\alpha}})m-v_p(\varphi(N)).\\
\end{eqnarray*}
Hence, the Proposition is proven in case $\alpha>0$. 

We look at the case $\alpha=0$ and prove the second congruence. As in Proposition 1 we set $L=p^m$. Again, by using Proposition 1.) in section 3.4 and Lemma 2 we obtain
$$
\omega^{-k}(\delta)\,{\rm tr}\,S_{u,\delta} T_n e_\alpha^L|_{{\cal M}_k(\Gamma,\chi\omega^{-k})}\equiv \omega^{-k}(\delta)\,{\rm tr}\,S_{u,\delta} T_n |_{{\cal M}_k(\Gamma,\chi\omega^{-k})^\alpha}\pmod{p^m}.
$$
This congruence still holds if we replace $k$ by $k'$. Together with the second congruence in Proposition 1 this
yields the claim. 

Finally, the third congruence is proven in exactly the same way as the first congruence (using the first congruence of Proposition 1 and equation (46)) with one small modification: we omit the step, where we use that the eigenvalues of $e_\alpha$ on ${\cal M}_k(\Gamma,\chi\omega^{-k})^\alpha$ are congruent to $1$ modulo $p$. Thus, the proof of the Proposition is complete.

%{\bf Remark. } Since the slope decomposition is defined over ${\Bbb Q}$ we know that the spaces ${\cal M}_k(\Gamma,\chi)^\alpha$ are defined over ${\Bbb Q}(Np)$. Thus, ${\rm tr}\,T_n|_{{\cal M}_k(\Gamma,\chi\omega^{-k})^\alpha}\in{\Bbb Q}(Np)$.

%[[In our computation of ${\rm tr}\,T_{\ell^j}|_{{\cal M}_k(\Gamma,\chi\omega^{-k})^\alpha}$ in the proof of the preceeding Proposition we obscured this fact by passing to a splitting field for $T_p$, i.e. (Id) only shows that ${\rm tr}\,T_{\ell^j}|_{{\cal M}_k(\Gamma,\chi\omega^{-k})^\alpha}\in E$, where $E/{\Bbb Q}(Np)$ is any splitting field. If we had chosen the projectors $e_i$ algebraically, we would have found $p_i\in{\Bbb Q}[X]$, hence, (Id) would show that ${\rm tr}\,T_{\ell^j}|_{{\cal M}_k(\Gamma,\chi\omega^{-k})^\alpha}\in{\Bbb Q}(Np)$.]]

%2.) We note that the condition ${{\sf C}_\alpha}\ge N(\alpha)+\frac{1}{\alpha}$ implies that $\frac{1}{{{\sf C}_\alpha}\alpha}\le 1$. 

%[[Da wir schlussendlich nur eine Kongruenz modulo $p^{{\rm min}\,(L,(1-N(\alpha)/{{\sf C}_\alpha})m)}$ erhalten ist der Faktor $A$ in "wir erhalten Kongruenz modulo $p^{Am+B}$" immer kleiner als $1$ !!!]]

The basic trace identity now reads as follows.

{\bf Theorem. } {\it Let $\alpha\in{\Bbb Q}_{\ge 0}$. Assume that $k,k'\in{\Bbb N}$ satisfy 

$\bullet$ $k,k'\ge ({{\sf C}_\alpha}\alpha+1)^2+2$  

$\bullet$ $k\equiv k'\pmod{p^m}$ with $m\ge {{\sf C}_\alpha}\alpha+1$. 

Let $\ell_1,\ldots,\ell_s\in{\Bbb N}$ be prime numbers and let $M\in{\Bbb N}$ be any divisor of $N$; then 
$$
{\rm tr}\,\pi^{Np}_{Mp} T_{\ell_1}^{r_1}\cdot\ldots\cdot T_{\ell_s}^{r_s}|_{{\cal M}_k(\Gamma,\chi\omega^{-k})^\alpha}
\equiv{\rm tr}\,\pi^{Np}_{Mp} T_{\ell_1}^{r_1}\cdot\ldots\cdot T_{\ell_s}^{r_s}|_{{\cal M}_{k'}(\Gamma,\chi\omega^{-k'})^\alpha}\pmod{p^{[\log_p [\frac{m}{{{\sf C}_\alpha}\alpha}]]-v_p(\varphi(N))}}
$$
if $\alpha>0$ and
$$
{\rm tr}\,\pi^{Np}_{Mp} T_{\ell_1}^{r_1}\cdot\ldots\cdot T_{\ell_s}^{r_s}|_{{\cal M}_k(\Gamma,\chi\omega^{-k})^\alpha}
\equiv{\rm tr}\,\pi^{Np}_{Mp} T_{\ell_1}^{r_1}\cdot\ldots\cdot T_{\ell_s}^{r_s}|_{{\cal M}_{k'}(\Gamma,\chi\omega^{-k'})^\alpha}\pmod{p^{m-v_p(\varphi(N))}}.
$$
if $\alpha=0$.

In addition, in case $\alpha>0$ the following congruence holds
$$
{\rm tr}\,\pi^{Np}_{Mp} T_{\ell_1}^{r_1}\cdot\ldots\cdot T_{\ell_s}^{r_s}e_\alpha^{[\frac{m}{{\sf C}_\alpha\alpha}]}|_{{\cal M}_k(\Gamma,\chi\omega^{-k})^\alpha}
\equiv{\rm tr}\,\pi^{Np}_{Mp} T_{\ell_1}^{r_1}\cdot\ldots\cdot T_{\ell_s}^{r_s}e_\alpha^{[\frac{m}{{\sf C}_\alpha\alpha}]}|_{{\cal M}_{k'}(\Gamma,\chi\omega^{-k'})^\alpha}\pmod{p^{[\frac{m}{{{\sf C}_\alpha}\alpha}]-v_p(\varphi(N))}}.
$$

%[[ die folgende Kongruenz brauchen wir nicht mehr

%2.) Assume that $k,k'\ge 2\alpha+2$; then then the following congruence, which does not involve the slope, holds true
%$$
%{\rm tr}\,\pi^{Np}_{Mp} T_{\ell_1}^{r_1}\cdot\ldots\cdot T_{\ell_s}^{r_s}T_p^{\left[\frac{m}{2\alpha}\right]+1}|_{{\cal M}_k(\Gamma,\chi\omega^{-k})^\alpha}
%\equiv{\rm tr}\,\pi^{Np}_{Mp} T_{\ell_1}^{r_1}\cdot\ldots\cdot T_{\ell_s}^{r_s}T_p^{\left[\frac{m}{2\alpha}\right]+1}|_{{\cal M}_{k'}(\Gamma,\chi\omega^{-k'})^\alpha}\pmod{p^{m-v_p(\varphi(N))}}.
%$$
%]]
}

{\it Proof. } 1.) We first assume $\alpha>0$ and prove the first congruence. Let $n$ be any divisor of $\ell_1^{r_1}\cdot\ldots\cdot \ell_s^{r_s}$ and let $n'$ be relatively prime to $Np$. Since $T(n',n')$ acts on ${\cal M}_k(\Gamma,\chi\omega^{-k})^\alpha$ as multiplication by $\chi\omega^{-k}(n'){n'}^{k-2}$ (cf. (15)) we find
$$
{\rm tr}\, S_{u,\delta}\,T(n',n') T_n |_{{\cal M}_k(\Gamma,\chi\omega^{-k})^\alpha}=\chi\omega^{-k}(n'){n'}^{k-2}\,{\rm tr}\,S_{u,\delta} T_n |_{{\cal M}_k(\Gamma,\chi\omega^{-k})^\alpha}.\leqno(48)
$$
Since furthermore $k\equiv k'\pmod{p^m}$ implies 
$$
\chi\omega^{-k}(n'){n'}^{k-2}\equiv \chi\omega^{-k'}(n'){n'}^{k'-2}\pmod{p^m}\leqno(49)
$$
%
%[[denn: \begin{eqnarray*}
%\chi\omega^{-k}(n'){n'}^{k-2} - \chi\omega^{-k'}(n'){n'}^{k'-2}&=&\chi\omega^{-k}(n'){n'}^{k-2} (1- \omega^{-(k'-k)}(n'){n'}^{k'-k}\\
%&\equiv&1\pmod{p^m}
%\end{eqnarray*}
%da $\omega^{-(k'-k)}(n'){n'}^{k'-k}\equiv 1\pmod{p}$
%]]
%
and since
$$
\omega^{-k}(\delta)\,{\rm tr}\,S_{u,\delta}\,T_n\,|_{{\cal M}_k(\Gamma,\chi\omega^{-k})^\alpha}\equiv \omega^{-k}(\delta)\,{\rm tr}\,S_{u,\delta}\,T_n\,|_{{\cal M}_{k'}(\Gamma,\chi\omega^{-k'})^\alpha}\pmod{p^{[\log_p [\frac{m}{{{\sf C}_\alpha}\alpha}]]-v_p(\varphi(N))}}
$$
by the first congruence of Proposition 2, we deduce that 
\begin{eqnarray*}
(50)\quad&&\omega^{-k}(\delta)\,{\rm tr}\, S_{u,\delta}\,T(n',n') T_n\,|_{{\cal M}_k(\Gamma,\chi\omega^{-k})^\alpha}\\
&&\equiv \omega^{-k'}(\delta)\,{\rm tr}\,S_{u,\delta}\, T(n',n') T_n\,|_{{\cal M}_{k'}(\Gamma,\chi\omega^{-k'})^\alpha}\pmod{p^{[\log_p [\frac{m}{{{\sf C}_\alpha}\alpha}]]-v_p(\varphi(N))}}
\end{eqnarray*}
(note that $[\frac{m}{{{\sf C}_\alpha}\alpha}]\le m$, hence, $\log_p [\frac{m}{{{\sf C}_\alpha}\alpha}]<m$; cf. equation (45)).

On the other hand, since $p|Mp$ we know that $\pi^{Np}_{Mp}=\sum_{u\in{\bf M},\delta\in{\bf A}} S_{u,\delta}$, where any $u\in{\bf M}$ is divisible by $p$ and any $\delta\in{\bf A}$ is congruent to $1$ modulo $p$ (cf. the Corollary in section 3.2). Hence,
$$
\pi^{Np}_{Mp} T_{\ell_1}^{r_1}\cdot\ldots\cdot T_{\ell_s}^{r_s}\in\sum_{u,\delta} \sum_{n', n} {\Bbb Z}\,S_{u,\delta}  T(n',n') T_n \leqno(51)
$$
is a ${\Bbb Z}$-linear combination of terms $S_{u,\delta} T(n',n')T_n$ with all $u$ divisible by $p$ and all $\delta\equiv 1\pmod{p}$ (cf. the Remark in section 3.2). In particular, $\omega(\delta)=1$ and equations (50) and (51) together yield the claim in case $\alpha>0$. 

To prove the second congruence we remark that the same proof as in case $\alpha>0$ holds if we use the second congruence of Proposition 2 instead of the first one, which replaces the modulus of the congruence in equation (50) by $\pmod{p^{m-v_p(\varphi(N))}}$. 

Finally, the third congruence follows in exactly the same way as the first one if instead of the term ${\rm tr}\,S_{u,\delta}\,T(n',n') T_n |_{{\cal M}_k(\Gamma,\chi\omega^{-k})^\alpha}$ appearing in equation (48) we consider ${\rm tr}\,S_{u,\delta}\,T(n',n') T_n e_\alpha^{\left[\frac{m}{{\sf C}_\alpha \alpha}\right]}|_{{\cal M}_k(\Gamma,\chi\omega^{-k})^\alpha}$ and use the third congruence of Proposition 2 (note that $\left[\frac{m}{{\sf C}_\alpha \alpha}\right]\le m$). 
%
%2.) Equation (48) even holds on the space ${\cal M}_k(\Gamma,\chi\omega^{-k})$ and together with equation (49) and Proposition 2 we obtain \begin{eqnarray*}
%&&\omega^{-k}(\delta)\,{\rm tr}\,S_{u,\delta} T(n',n') T_n T_p^{\left[\frac{m}{2\alpha}\right]+1}|_{{\cal M}_k(\Gamma,\chi\omega^{-k})}\\
%&&\equiv \omega^{-k'}(\delta)\,{\rm tr}\,S_{u,\delta} T(n',n') T_n T_p^{\left[\frac{m}{2\alpha}\right]+1}|_{{\cal M}_{k'}(\Gamma,\chi\omega^{-k'})}\pmod{p^{m-v_p(\varphi(N))}}.\\
%\end{eqnarray*}
%In conjunction with equation (51) this immediately yields the claim. 
%
Hence, the proof of the Theorem is complete.

\newpage

\section{$p$-adic families of modular forms}

{\bf 4.1. Constance of dimension of the slope subspaces. } As in the previous sections, throughout section 4 we fix a tame level $N$ and a prime $p$ satisfying $(N,p)=1$ and we continue to set $\Gamma=\Gamma_1(Np)$.

We denote by $\Psi_{k,\chi\omega^{-k}}^\alpha$ the characteristic polynomial of $T_p$ acting on ${\cal M}_k(\Gamma,\chi\omega^{-k})^\alpha$. We have seen that $\Psi_{k,\chi\omega^{-k}}^\alpha$ is contained in ${\Bbb Z}[Np][X]$ (cf. the Lemma in section 3.3). We set 
$$
d_{k,\chi\omega^{-k}}^\alpha={\rm dim}\,{\cal M}_k(\Gamma,\chi\omega^{-k})^\alpha,
$$ 
hence, $d_{k,\chi\omega^{-k}}^\alpha$ is the degree of $\Psi_{k,\chi\omega^{-k}}^\alpha$. We write
$$
\Psi_{k,\chi\omega^{-k}}^\alpha=\sum_{j=0}^{d} (-1)^{j}a_{k,\chi\omega^{-k},j}^\alpha X^{d-j}
$$
where $d=d_{k,\chi\omega^{-k}}^\alpha$. This defines $a_{k,\chi\omega^{-k},j}^\alpha$ for $j=0,\ldots,d_{k,\chi\omega^{-k}}^\alpha$ and we set $a_{k,\chi\omega^{-k},j}^\alpha=0$ for $j>d_{k,\chi\omega^{-k}}^\alpha$. 

%[[beachte: die  $a_{k,\chi\omega^{-k},j}^\alpha$ mit $j>d$ korrespondieren zu negativen Potenzen von $X$ (?)]]

{\bf Lemma. }{\it Let $\alpha\in{\Bbb Q}_{\ge 0}$. Assume that $k,k'\in{\Bbb N}$ satisfy 

$\bullet$ $k,k'\ge({{\sf C}_\alpha}\alpha+1)^2+2$
 
$\bullet$ $k\equiv k'\pmod{p^m}$ with $m\ge {{\sf C}_\alpha}\alpha+1$.

Then for all $j\in{\Bbb N}\cup\{0\}$ we have
$$
a_{k,\chi\omega^{-k},j}^\alpha\equiv a_{k',\chi\omega^{-k'},j}^\alpha\pmod{p^{\left[\log_p [\frac{m}{{{\sf C}_\alpha}\alpha}]\right]-v_p(\varphi(N)j!)}}.
$$
if $\alpha>0$ and 
$$
a_{k,\chi\omega^{-k},j}^\alpha\equiv a_{k',\chi\omega^{-k'},j}^\alpha\pmod{p^{m-v_p(\varphi(N)j!)}}.
$$
if $\alpha=0$.
}

%[[We note that the statement of the Lemma is empty if $\left[\frac{m}{{{\sf C}_\alpha}\alpha}\right]-v_p(\varphi(N)j!)\le 0$. But since $j\le d(\alpha)$ is bounded by a constant only depending on $\alpha$ we see that $left[\frac{m}{{{\sf C}_\alpha}\alpha}\right]-v_p(\varphi(N)j!)\rightarrow+\infty$ for $m\rightarrow+\infty$ and fixed $\alpha$.
%]]

%[[da die Koeffizienten $a_{k,\chi\omega^{-k},j}^\alpha$ des char. Polynoms $\Psi_{k,\chi\omega^{-k}}^\alpha$ alle ganz sind nach Lemma in chapter 3.2 !]]

{\it Proof. } We first assume $\alpha>0$. We set
$$
\tau_{k,\chi,j}^\alpha={\rm tr}\,T_{p^j}|_{{\cal M}_k(\Gamma,\chi)^\alpha}.
$$ 
Since $T_{p^j}=T_p^j$ (note that $p$ divides the Level $Np$) we know that $\tau_{k,\chi\omega^{-k},j}^\alpha$ also is the trace of $T_p^j$ on ${\cal M}_k(\Gamma,\chi\omega^{-k})^\alpha$. We therefore obtain by a classical formula (cf. [Koe], 3.4.6 Satz, p. 117) 
$$
ja_{k,\chi\omega^{-k},j}^\alpha=\sum_{h=1}^j (-1)^{h+1} \tau_{k,\chi\omega^{-k},h}^\alpha a_{k,\chi\omega^{-k},j-h}^\alpha.\leqno(1)
$$ 
This formula is valid for all $j\in{\Bbb N}\cup\{0\}$ if we set $a_{k,\chi\omega^{-k},j}^\alpha=0$ for $j>d_{k,\chi\omega^{-k}}^\alpha$ as we have done. The lemma then follows by induction over $j=0,1,2,\ldots$ using (1) and the basic trace identity in section 3.5. In more detail, if $j=0$ then $a_{k,\chi\omega^{-k},0}^\alpha=a_{k',\chi\omega^{-k'},0}^\alpha=1$, whence the claim in this case (note that $a_{k,\chi\omega^{-k},0}^\alpha$ and $a_{k',\chi\omega^{-k'},0}^\alpha$ are the leading coefficients of the respective characteristic polynomials). Assume therefore that the claim holds for  $0,\ldots,j-1$; the first congruence of the Theorem in section 3.5 yields $\tau_{k,\chi\omega^{-k},j}^\alpha\equiv \tau_{k',\chi\omega^{-k'},j}^\alpha\pmod{p^{\left[\log_p [\frac{m}{{{\sf C}_\alpha}\alpha}]\right]-v_p(\varphi(N))}}$. Moreover, the induction hypotheses reads 
$$
a_{k,\chi\omega^{-k},j-h}^\alpha\equiv a_{k',\chi\omega^{-k'},j-h}^\alpha\pmod{p^{\left[\log_p [\frac{m}{{{\sf C}_\alpha}\alpha}]\right]-v_p(\varphi(N))-v_p((j-h)!)}}
$$ 
for all $h=j,\ldots,1$. We thus obtain 
\begin{eqnarray*}
ja_{k,\chi\omega^{-k},j}^\alpha&=&\sum_{h=1}^j (-1)^{h+1} \tau_{k,\chi\omega^{-k},h}^\alpha a_{k,\chi\omega^{-k},j-h}^\alpha\\
&\equiv& \sum_{h=1}^j (-1)^{h+1} \tau_{k',\chi\omega^{-k'},h}^\alpha a_{k',\chi\omega^{-k'},j-h}^\alpha\pmod{p^{\left[\log_p [\frac{m}{{{\sf C}_\alpha}\alpha}]\right]-v_p(\varphi(N))-v_p((j-1)!)}}\\
&=&ja_{k',\chi\omega^{-k'},j}^\alpha.\\
\end{eqnarray*}
%
%[[note that $j-h\le j-1$ $\Rightarrow (j-h)!|(j-1)!\Rightarrow v_p((j-h)!)\le v_p((j-1)!)$ ]]
%
Dividing by $j$ yields $a_{k,\chi\omega^{-k},j}^\alpha\equiv a_{k',\chi\omega^{-k'},j}^\alpha\pmod{p^{\left[\log_p [\frac{m}{{{\sf C}_\alpha}\alpha}]\right]-v_p(\varphi(N)j!)}}$. Thus, the proof of the lemma is complete in case $\alpha>0$. If $\alpha=0$ the same proof holds if we replace the summand "$\left[\log_p [\frac{m}{{{\sf C}_\alpha}\alpha}]\right]$" appearing in the modulus of the above congruences by "$m$" and use the second trace identity of the Theorem in section 3.5. Thus, the proof of the Lemma is complete.

Let $\alpha>0$ and assume that $d_{k,\chi\omega^{-k}}^\alpha\ge d_{k',\chi\omega^{-k'}}^\alpha$. The claim of the Lemma then is equivalent to the following two statements:

$\bullet$ $a_{k,\chi\omega^{-k},j}^\alpha\equiv a_{k',\chi\omega^{-k'},j}^\alpha\pmod{p^{\left[\log_p [\frac{m}{{{\sf C}_\alpha}\alpha}]\right]-v_p(\varphi(N)j!)}}\quad\mbox{for $j=0,\ldots,d_{k',\chi\omega^{-k'}}^\alpha$}$. 

$\bullet$ $a_{k,\chi\omega^{-k},j}^\alpha\equiv 0\pmod{p^{\left[\log_p [\frac{m}{{{\sf C}_\alpha}\alpha}]\right]-v_p(\varphi(N)j!)}}\quad\mbox{for $d_{k',\chi\omega^{-k'}}^\alpha<j\le d_{k,\chi\omega^{-k}}^\alpha$}$.

%[[Beachte, dass $a_{k',\chi\omega^{-k'},j}^\alpha=0$ fur $d_{k',\chi\omega^{-k'}}^\alpha<j\le d_{k,\chi\omega^{-k}}^\alpha$]]

(For $j>d_{k,\chi\omega^{-k}}^\alpha$ we obtain the congruence $0\equiv 0\pmod{p^{\left[\log_p [\frac{m}{{{\sf C}_\alpha}\alpha}]\right]-v_p(\varphi(N)j!)}}$, which is completely trivial.) In particular, taking into account that $j\le d(=d_{k,\chi\omega^{-k}}^\alpha)$ implies $j\le d(\alpha)$ (cf. section 1.5 for the definition of $d(\alpha)$) we obtain the congruences
$$
a_{k,\chi\omega^{-k},j}^\alpha\equiv a_{k',\chi\omega^{-k'},j}^\alpha\pmod{p^{\left[\log_p [\frac{m}{{{\sf C}_\alpha}\alpha}]\right]-v_p(\varphi(N)d(\alpha)!)}}\quad\mbox{for $j=0,\ldots,d_{k',\chi\omega^{-k'}}^\alpha$}\leqno(2a)
$$
$$
a_{k,\chi\omega^{-k},j}^\alpha\equiv 0\pmod{p^{\left[\log_p [\frac{m}{{{\sf C}_\alpha}\alpha}]\right]-v_p(\varphi(N)d(\alpha)!)}}\quad\mbox{for $d_{k',\chi\omega^{-k'}}^\alpha$}<j\le d_{k,\chi\omega^{-k}}^\alpha.\leqno(2b)
$$
We note that in case $\alpha=0$ the same congruences with "$\left[\log_p [\frac{m}{{{\sf C}_\alpha}\alpha}]\right]$" replaced by "$m$" hold true. We shall use equations (2a), (2b) to deduce that the function $k\mapsto d^\alpha_{k,\chi\omega^{-k}}$ is locally constant in the $p$-adic sense. To this end let $E/{\Bbb Q}(Np)$ be a splitting field of $\Psi_{k,\chi\omega^{-k}}^\alpha$. We continue to write $\Psi_{k,\chi\omega^{-k}}^\alpha=\sum_{j=0}^d (-1)^j a_{k,\chi\omega^{-k},j}^\alpha X^{d-j}$, where, again, we have set $d=d_{k,\chi\omega^{-k}}^\alpha$. Since $\Psi_{k,\chi\omega^{-k}}^\alpha=\prod_i (X-\xi_i)$ and the zeroes $\xi_i\in E$ of $\Psi_{k,\chi\omega^{-k}}^\alpha$ all have $p$-adic valuation $v_p(\xi_i)=\alpha$ we find 
$$
v_p(a_{k,\chi\omega^{-k},d}^\alpha)=d\alpha.\leqno(3)
$$

%[[Note: we can even determine the Newton polygon ${\cal N}$ of $\Psi_{k,\chi\omega^{-k}}^\alpha$: Since $v_p(a_{k,\chi\omega^{-k},j}^\alpha)=i$ for all $i$ we obtain $v_p(a_{k,\chi\omega^{-k},0}^\alpha)=0$, $v_p(a_{k,\chi\omega^{-k},d}^\alpha)=d\alpha$ and $v_p(a_{k,\chi\omega^{-k},j}^\alpha)\ge (d-j)\alpha$, $j=1,\ldots,d-1$. Hence, the Newton polygon ${\cal N}$ consists of one segment of length $d(=d_{k,\chi\omega^{-k}}^\alpha)$ and slope $-\alpha$, i.e. ${\cal N}$ is the line segment joining the points $(0,d\alpha)$ and $(d,0)$. Since the point $(0,v_p(a_{k,\chi\omega^{-k},d}^\alpha))$ corresponding to the constant term of $\Psi_{k,\chi\omega^{-k}}^\alpha$ always lies on the Newton polygon we again deduce that $v_p(a_{k,\chi\omega^{-k},d}^\alpha)=d\alpha$.
%]]

We set
$$
{{\sf B}_\alpha}={\rm max}\,\{{\sf C}_\alpha\alpha+1,{\sf C}_\alpha\alpha\big(\exp \big(\alpha d(\alpha)+v_p(\varphi(N)d(\alpha)!)+1\big)+1\big)\} .
$$
if $\alpha>0$ and ${\sf B}_\alpha=v_p(\varphi(N)d(\alpha)!)+1$ if $\alpha=0$.

{\bf Theorem B. }{\it Fix an arbitrary slope $\alpha\in{\Bbb Q}_{\ge 0}$. For all pairs of integers $k,k'\in{\Bbb N}$ satisfying 

$\bullet$ $k,k'\ge({{\sf C}_\alpha}\alpha+1)^2+2$ and 

$\bullet$ $k\equiv k'\pmod{p^m}$ with an integer $m>{\sf B}_\alpha$

it holds that
$$
d_{k,\chi\omega^{-k}}^\alpha=d_{k',\chi\omega^{-k'}}^\alpha.
$$
}

%[[$"... G(\alpha d(\alpha)+v_p(\varphi(N)d(\alpha)!)+1) ..." bedeutet genauer   "... G(\alpha d(\alpha)+v_p(\varphi(N) (d(\alpha)!) )+1) ..."]]

{\it Proof.} We first assume that $\alpha>0$. Straightforward calculation shows that $m>{\sf C}_\alpha\alpha\big(\exp (\alpha d(\alpha)+v_p(\varphi(N)d(\alpha)!)+1)+1\big)$ implies 
$$
\log_p (\frac{m}{{{\sf C}_\alpha}\alpha}-1)-v_p(\varphi(N)d(\alpha)!)>\alpha d(\alpha)+1.
$$
Since $\left[\frac{m}{{{\sf C}_\alpha}\alpha}\right]\ge \frac{m}{{{\sf C}_\alpha}\alpha}-1$ and $d(\alpha)\ge d^\alpha_{k,\chi\omega^{-k}}$ we obtain
$$
\log_p [\frac{m}{{{\sf C}_\alpha}\alpha}]-1-v_p(\varphi(N)d(\alpha)!)>\alpha d^\alpha_{k,\chi\omega^{-k}},
$$
hence,
$$
\left[\log_p [\frac{m}{{{\sf C}_\alpha}\alpha}]\right]-v_p(\varphi(N)d(\alpha)!)>\alpha d^\alpha_{k,\chi\omega^{-k}}.\leqno(4)
$$
We now assume there are $k,k'\ge({{\sf C}_\alpha}\alpha+1)^2+2$ such that $k\equiv k'\pmod{p^m}$ with $m>{\sf B}_\alpha$ and $d_{k,\chi\omega^{-k}}^\alpha\not=d_{k',\chi\omega^{-k'}}^\alpha$. Without loss of generality we may assume that $d_{k,\chi\omega^{-k}}^\alpha>d_{k',\chi\omega^{-k'}}^\alpha$. Equation (2b) then immediately yields 
$$
a_{k,\chi\omega^{-k},d}^\alpha\equiv 0\pmod{p^{\left[\log_p [\frac{m}{{{\sf C}_\alpha}\alpha}]\right]-v_p(\varphi(N)d(\alpha)!)}}
$$
($d=d_{k,\chi\omega^{-k}}^\alpha$), hence, by using (4), we find
$$
v_p(a_{k,\chi\omega^{-k},d}^\alpha)>\alpha d_{k,\chi\omega^{-k}}^\alpha.
$$
This last equation contradicts $(3)$, by which $v_p(a_{k,\chi\omega^{-k},d}^\alpha)=\alpha d_{k,\chi\omega^{-k}}^\alpha$. Thus, the assumption is false and the Theorem is proven in case $\alpha=0$. 

If $\alpha=0$ the same proof still holds: assume there are $k,k'\ge({{\sf C}_\alpha}\alpha+1)^2+2$ such that $k\equiv k'\pmod{p^m}$ with $m>{\sf B}_\alpha={\sf B}_0$ and $d_{k,\chi\omega^{-k}}^\alpha>d_{k',\chi\omega^{-k'}}^\alpha$; using the definition of ${\sf B}_\alpha={\sf B}_0$ we find that $v_p(a_{k,\chi\omega^{-k},d}^\alpha)>0$, while equation (3) yields $v_p(a_{k,\chi\omega^{-k},d}^\alpha)=0$, a contradiction. Thus, the proof of the theorem is complete.

%{\bf Remark.} Since the slope $\alpha$ occurs in ${\cal M}_k(\Gamma,\chi\omega^{-k})$ or in ${\cal M}_{k'}(\Gamma,\chi\omega^{-k'})$ we know that $\alpha\ge 0$. 

%[[wir brauchen $\alpha\ge 0$ um $\alpha d\ge 0$ folgern zu konnen]]

%Hence, $\left[\frac{m}{{{\sf C}_\alpha}\alpha}\right]-v_p(\varphi(N)d(\alpha)!)>\alpha d\ge 0$. Therefore the congruence $a_{k,\chi\omega^{-k},d}^\alpha\equiv 0\pmod{p^{\left[\frac{m}{{{\sf C}_\alpha}\alpha}\right]-v_p(\varphi(N)d(\alpha)!)}}$ occuring in the above proof is not an empty statement (note that $a_{k,\chi\omega^{-k},d}^\alpha\in{\Bbb Z}[Np]$).

%{\bf Remark. } The theorem in particular implies that the function $k\mapsto d_{k,\chi\omega^{-k}}^\alpha$ is locally constant, which is the first part of the Mazur-Gouvea conjecture.

{\bf 4.2. Transfer for modular forms from weight $k$ to weight $k'$. } As in section 3.1 we set ${\cal H}=\Gamma\backslash \Delta_1/\Gamma$ and recall that ${\cal H}$ acts on ${\cal M}_k(\Gamma,\chi\omega^{-k})^\alpha$. 
%
%[[because the operators $T_\ell$ commute with $T_p$ and also with $\langle \epsilon \rangle$]]
%
For any sequence $(\lambda_\ell)_{\ell\in P}$ with complex entries, where $\ell$ runs over the set $P$ of all rational primes, we set
$$
{\cal M}_{k,\chi}^\alpha(\lambda)=\{f\in {\cal M}_k(\Gamma,\chi)^\alpha:\,\mbox{for all $\ell\in P$ there is $k=k_\ell\in{\Bbb N}$}:\,\,(T_{\ell}-\lambda_\ell)^k f=0\}.
$$
We fix an (auxiliary) numbering $\ell_1,\ell_2,\ell_3,\ldots$ of the rational primes and we set
$$
{\cal M}_{k,\chi}^\alpha(\lambda_{\ell_1},\ldots,\lambda_{\ell_n})=\{f\in {\cal M}_k(\Gamma,\chi)^\alpha:\, \mbox{for all $i=1,\ldots,n$ there is $k=k_i\in{\Bbb N}$}:\,\,(T_{\ell_i}-\lambda_{\ell_i})^k f=0 \}.
$$
Since the operators $T_\ell$ commute, the finitely many Hecke operators $\{T_{\ell_1},\ldots,T_{\ell_s}\}$ can be simultaneously transformed into upper triangular form. 
%
%[[s. [Storch-Wiebe: Lehrbuch der Mathematik Bd 2: Lineare Algebra], das Kapitel uber diagonalisierbare und trigonalisierbare Operatoren]]
%
In particular, for any $n\in{\Bbb N}$ there is a basis ${\cal B}$ in ${\cal M}_{k,\chi\omega^{-k}}^\alpha(\lambda)$ such that for all $i\le n$ the Hecke operator $T_{\ell_i}$ is represented by the matrix
$$
{\cal D}_{\cal B}(T_{\ell_i}|_{{\cal M}_{k,\chi\omega^{-k}}^\alpha(\lambda)})
=\left(\begin{array}{ccc}
\lambda_{\ell_i}&&*_i\\
&\ddots&\\
&&\lambda_{\ell_i}\\
\end{array}
\right).\leqno(5)
$$
Let $\lambda=(\lambda_\ell)_{\ell\in P}$ be arbitrary. The spaces ${\cal M}_{k,\chi\omega^{-k}}^\alpha(\lambda_{\ell_1},\ldots,\lambda_{\ell_n})$ form a decreasing sequence
$$
{\cal M}_{k,\chi\omega^{-k}}^\alpha(\lambda_{\ell_1})\supseteq {\cal M}_{k,\chi\omega^{-k}}^\alpha(\lambda_{\ell_1},\lambda_{\ell_2})\supseteq {\cal M}_{k,\chi\omega^{-k}}^\alpha(\lambda_{\ell_1},\lambda_{\ell_2},\lambda_{\ell_3})\supseteq\cdots\leqno(6)
$$
such that
$$
\bigcap_{n\in{\Bbb N}} {\cal M}_{k,\chi\omega^{-k}}^\alpha(\lambda_{\ell_1},\ldots,\lambda_{\ell_n})={\cal M}_{k,\chi\omega^{-k}}^\alpha(\lambda).
$$
Since ${\rm dim}\,{\cal M}_k(\Gamma,\chi\omega^{-k})^\alpha$ is finite there are only finitely many $i\in{\Bbb N}$ such that 
${\cal M}_{k,\chi\omega^{-k}}^\alpha(\lambda_{\ell_1},\ldots,\lambda_{\ell_i})\not= {\cal M}_{k,\chi\omega^{-k}}^\alpha(\lambda_{\ell_1},\ldots,\lambda_{\ell_{i+1}})$ in (6). In particular, for every $\lambda=(\lambda_\ell)_{\ell\in P}$ there is an $n_\lambda\in{\Bbb N}$ such that
$$
{\cal M}_{k,\chi\omega^{-k}}^\alpha(\lambda_{\ell_1},\ldots,\lambda_{\ell_{n_\lambda}})= {\cal M}_{k,\chi\omega^{-k}}^\alpha(\lambda_{\ell_1},\ldots,\lambda_{\ell_{n_\lambda+1}})=\cdots={\cal M}_{k,\chi\omega^{-k}}^\alpha(\lambda).\leqno(7)
$$
In different words, the (finite) tuple $(\lambda_{\ell_1},\ldots,\lambda_{\ell_{n_\lambda}})$ already uniquely determines the subspace ${\cal M}_{k,\chi\omega^{-k}}^\alpha(\lambda)\le {\cal M}_k(\Gamma,\chi\omega^{-k})^\alpha$.

We now consider the set of all $\lambda$'s, i.e. we set
$$
\Lambda_{k,\chi}^\alpha=\{\lambda=(\lambda_\ell)_{\ell}:\,{\cal M}_{k,\chi\omega^{-k}}^\alpha(\lambda)\not=0\}
$$
and 
$$
\Lambda_{k,\chi,n}^\alpha=\{(\lambda_{\ell_1},\ldots,\lambda_{\ell_n}):\, {\cal M}_{k,\chi\omega^{-k}}^\alpha(\lambda_{\ell_1},\ldots,\lambda_{\ell_n})\not=0\}.
$$

Since ${\rm dim}\,{\cal M}_k(\Gamma,\chi\omega^{-k})^\alpha$ is finite, $\Lambda_{k,\chi}^\alpha$ is a finite set. We set 
$$
{\sf n}_k={\sf n}_{k,\chi}^\alpha={\rm max}\{n_\lambda,\,\lambda\in\Lambda_{k,\chi}^\alpha\}.
$$ 
Equation (7) then shows that for any two $\lambda=(\lambda_\ell)_{\ell\in P}$ and $(\lambda')=(\lambda_\ell')_{\ell\in P}$ with ${\cal M}_{k,\chi\omega^{-k}}^\alpha(\lambda)\not=0$, ${\cal M}_{k,\chi\omega^{-k}}^\alpha(\lambda')\not=0$ we have
$$
{\cal M}_{k,\chi\omega^{-k}}^\alpha(\lambda)={\cal M}_{k,\chi\omega^{-k}}^\alpha(\lambda')\Leftrightarrow (\lambda_1,\ldots,\lambda_{{\sf n}_k})=(\lambda_1',\ldots,\lambda_{{\sf n}_k}').
$$

%[[da ${\cal M}_{k,\chi\omega^{-k}}^\alpha(\lambda)={\cal M}_{k,\chi\omega^{-k}}^\alpha(\lambda_1,\ldots,\lambda_{{\sf n}_k})$ und genauso fur $\lambda'$]]

%[[trivialerweise gilt ${\cal M}_{k,\chi\omega^{-k}}^\alpha(\lambda)={\cal M}_{k,\chi\omega^{-k}}^\alpha(\lambda')\Leftrightarrow \lambda=\lambda'$ (d.h. Definition von ${\cal M}_{k,\chi\omega^{-k}}^\alpha(\lambda)$ benutzt alle Information von $\lambda=(\lambda_i)$)]]

In particular, for a given $(\mu_{\ell_1},\ldots,\mu_{\ell_{{\sf n}_k}})$ there is at most one $\lambda=(\lambda_\ell)\in\Lambda^\alpha_{k,\chi}$ such that $(\lambda_{\ell_1},\ldots,\lambda_{\ell_{{\sf n}_k}})=(\mu_{\ell_1},\ldots,\mu_{\ell_{{\sf n}_k}})$. Hence, for any $n\le {\sf n}_k$ we obtain  
$$
|\Lambda_{k,\chi,n}^\alpha|\le|\Lambda_{k,\chi,{\sf n}_k}^\alpha|=|\Lambda_{k,\chi}^\alpha|.\leqno(8)
$$
Since the operators $T_\ell$ split over ${\Bbb C}$ and commute with each other, we further obtain 

\begin{itemize}
\item (9a) ${\cal M}_k(\Gamma,\chi\omega^{-k})^\alpha=\bigoplus_{\lambda_{\ell_1}\in\Lambda_{k,\chi,1}^\alpha}{\cal M}_{k,\chi\omega^{-k}}^\alpha(\lambda_{\ell_1})$, 
\item (9b) ${\cal M}_{k,\chi\omega^{-k}}^\alpha(\lambda_{\ell_1})=\bigoplus_{\lambda_{\ell_2}:\,(\lambda_{\ell_1},\lambda_{\ell_2})\in\Lambda_{k,\chi,2}^\alpha} {\cal M}_{k,\chi\omega^{-k}}^\alpha(\lambda_{\ell_1},\lambda_{\ell_2})$, 
\item (9c) ${\cal M}_{k,\chi\omega^{-k}}^\alpha(\lambda_{\ell_1},\lambda_{\ell_2})=\bigoplus_{\lambda_{\ell_3}:\,(\lambda_{\ell_1},\lambda_{\ell_2},\lambda_{\ell_3})\in \Lambda_{k,\chi,3}^\alpha} {\cal M}_{k,\chi\omega^{-k}}^\alpha(\lambda_{\ell_1},\lambda_{\ell_2},\lambda_{\ell_3})$ 
\item ...
\end{itemize}

%[[da $T_{\ell_1}$ zerfallt gilt ${\cal M}_k(\Gamma,\chi\omega^{-k})^\alpha=\bigoplus_{\lambda_1\in\Lambda_{k,\chi,1}^\alpha}{\cal M}_{k,\chi\omega^{-k}}^\alpha(\lambda_{\ell_1})$; da $T_{\ell_1}, T_{\ell_2}$ kommutieren lasst $T_{\ell_2}$ ${\cal M}_{k,\chi\omega^{-k}}^\alpha(\lambda_{\ell_1})$ invariant (fur jedes $\lambda_1$). Da $T_{\ell_2}$ zerfallt gilt ${\cal M}_{k,\chi\omega^{-k}}^\alpha(\lambda_{\ell_1})=\bigoplus_{\lambda_{\ell_2}:\,(\lambda_{\ell_1},\lambda_{\ell_2})\in\Lambda_{k,\chi,2}^\alpha} {\cal M}_{k,\chi\omega^{-k}}^\alpha(\lambda_{\ell_1},\lambda_{\ell_2})$ ...
%]]

Thus, for any $n$ we obtain
$$
{\cal M}_k(\Gamma,\chi\omega^{-k})^\alpha=\bigoplus_{(\lambda_{\ell_1},\ldots,\lambda_{\ell_n})\in\Lambda^\alpha_{k,\chi,n}}{\cal M}_{k,\chi\omega^{-k}}^\alpha(\lambda_{\ell_1},\ldots,\lambda_{\ell_n}).
$$

In particular, using (7) and (8) we deduce that
\begin{eqnarray*}
(10)\qquad{\cal M}_k(\Gamma,\chi\omega^{-k})^\alpha&=&\bigoplus_{(\lambda_{\ell_1},\ldots,\lambda_{\ell_{{\sf n}_k}})\in\Lambda^\alpha_{k,\chi,{\sf n}_k}}{\cal M}_{k,\chi\omega^{-k}}^\alpha(\lambda_{\ell_1},\ldots,\lambda_{\ell_{{\sf n}_k}})\\
&=&\bigoplus_{\lambda=(\lambda_\ell)\in\Lambda^\alpha_{k,\chi}}{\cal M}_{k,\chi\omega^{-k}}^\alpha(\lambda_{\ell_1},\ldots,\lambda_{\ell_{{\sf n}_k}})\\
&=&\bigoplus_{\lambda=(\lambda_\ell)\in\Lambda_{k,\chi}^\alpha}{\cal M}_{k,\chi\omega^{-k}}^\alpha(\lambda).
\end{eqnarray*}

We describe an adelic version of the above decomposition of ${\cal M}_k(\Gamma,\chi\omega^{-k})^\alpha$. We set
$$
K_{1,\ell}(Np)=\{\Mat{a}{b}{c}{d}\in{\rm GL}_2({\Bbb Z}_\ell),\,c\equiv 0,\,d\equiv 1\pmod{Np}\}
$$
and
$$
K_1(Np)=\prod_{\ell\not=\infty}K_{1,\ell}(Np)\le \prod_{\ell\not=\infty}{\rm GL}_2({\Bbb Z}_\ell).
$$
For simplicity we put $K=K_1(Np)$ and $K_\ell=K_{1,\ell}(Np)$ and we note that $K_\ell={\rm GL}_2({\Bbb Z}_\ell)$ for all $\ell\not|Np$. We define local Hecke algebras
$$
{\cal H}_\ell=K_{\ell}\backslash D_{\ell}/K_{\ell},
$$
where
$$
D_{\ell}=\{\Mat{a}{b}{c}{d}\in M_2({\Bbb Z}_\ell),\,ad-bc\not=0, d\equiv 1, c\equiv 0\pmod{Np}\}.
$$

%[[cf. [M] p. 212/213]]

as well as the adelic Hecke algebra
$$
{\cal H}_{\Bbb A}=K\backslash D/K,
$$
where
$$
D=\prod_{\ell\not=\infty} D_{\ell}.
$$

%[[Beachte: ist $N$ eine Einheit in ${\Bbb Z}_\ell$, dan ist die Kongruenz $x\equiv y\pmod{N}$, $x,y\in{\Bbb Z}_\ell$ eine leere Bedingung]]

There are canonical isomorphisms
$$
\begin{array}{ccc}
{\cal H}&\cong&{\cal H}_{\Bbb A}\\
\Gamma\alpha\Gamma&\mapsto&K\alpha K\\
\end{array}
$$
and 
$$
\begin{array}{ccc}
{\cal H}_{\Bbb A}&\cong&\otimes_{\ell\not=\infty}{\cal H}_\ell\\
K(\alpha_\ell)_\ell K&\mapsto&\otimes_\ell K_\ell\alpha_\ell K_\ell.
\end{array}
$$
(cf. [M] Theorem 5.3.5, p. 214 and Theorem 4.5.18, p. 151).
%
%[[letztere Referenz wegen Unterschied von $\Gamma_0$ zu $\Gamma_1$]] 
%
In particular, any ${\cal H}$ module becomes canonically a ${\cal H}_{\Bbb A}$-module and even further an ${\cal H}_\ell$-module for any prime $\ell$. By $T_\ell$, $\langle \epsilon\rangle$ and $e_\chi$ we shall also denote the corresponding operators in ${\cal H}_{\Bbb A}$ and we set $T_\ell^{\rm loc}=K_\ell\Mat{1}{}{}{\ell}K_\ell\in{\cal H}_\ell$.

We denote be ${\cal R}_k$ the set of (finite parts of automorphic) representations $(\pi,V_\pi)$ of ${\rm GL}_2({\Bbb A}_f)$, which occur in $\lim_C H^1(S(C),L_{k,{\Bbb C}})$, $C\le{\rm GL}_2(\hat{\Bbb Z})$ running over all compact open subgroups, 
and which satisfy $V_\pi^K\not=0$ (recall that $K=K_1(Np)$). We then obtain
\begin{eqnarray*}
{\cal M}_k(\Gamma)&\cong& H^1(S(K),L_{k-2,{\Bbb C}})\\
&\cong&\bigoplus_{\pi\in {\cal R}_{k-2}} V_{\pi}^{K}.\\
\end{eqnarray*}
We note that these are isomorphisms of ${\cal H}_{\Bbb A}=\otimes_\ell {\cal H}_\ell$ modules. The representation $\pi$ decomposes $\pi=\otimes_\ell \pi_\ell$, where $(\pi_\ell,V_{\pi,\ell})$ is a representation of ${\rm GL}_2({\Bbb Q}_\ell)$ and we further obtain
$$
V^{K}_\pi=\otimes_\ell V_{\pi,\ell}^{K_\ell}.
$$
We look at the action of $T_\ell$ on $V_\pi^K$. If $\ell\not|Np$ the subspace of $K_\ell$ invariants $V_{\pi,\ell}^{K_\ell}=\langle v_{\pi,\ell}^0\rangle$ is $1$-dimensional (note that $K_\ell={\rm GL}_2({\Bbb Z}_\ell)$) and $T_\ell^{\rm loc}$ acts by multiplication with some scalar $\lambda_{\pi,\ell}\in{\Bbb C}$, i.e.
$$
T_\ell^{\rm loc} v_{\pi,\ell}^0=\lambda_{\pi,\ell} \,v_{\pi,\ell}^0.
$$
If $\ell|Np$, let $\{v_{\ell,t}\}_t$ be a basis of $V_{\pi,\ell}^{K_\ell}$; ${\cal B}=\{\otimes_{\ell\not|Np} v_{\pi,\ell}^0\otimes\otimes_{\ell|Np}v_{\ell,t_\ell}\}_{(t_\ell)_{\ell|Np}}$ then is a basis of $V_\pi^K$ with respect to which all operators $T_\ell$ with $\ell\not|Np$ are diagonal and even central
$$
{\cal D}_{\cal B}(T_\ell|_{V_\pi^K})=\left(
\begin{array}{ccc}
\lambda_{\pi,\ell}&&\\
&\ddots&\\
&&\lambda_{\pi,\ell}\\
\end{array}
\right)=\lambda_{\pi,\ell}\,{\rm Id}.\leqno(11)
$$

%[[Beachte: $T_\ell$ operiert nur in der $\ell$-Komponente von $\otimes_\ell v_\ell$ und dort als $T^{\rm loc}_{\ell}$]]

If $\ell|Np$ we denote by $\lambda_{\pi,\ell;h}$, $h=1,\ldots,n_{\pi,\ell}$ the eigenvalues of $T_\ell$ acting on $V_{\pi,\ell}^{K_\ell}$. 
Furthermore, we denote by $V_{\pi,\ell;h}^{K_\ell}\subset V_{\pi,\ell}^{K_\ell}$ the generalized eigenspace attached to the eigenvalue $\lambda_{\pi,\ell;h}$ of $T_\ell$. Hence, 
there is a basis ${\cal B}_h$ of $V_{\pi,\ell;h}^{K_\ell}$ such that $T_\ell^{\rm loc}|_{V_{\pi,\ell;h}^{K_\ell}}$ is represented by the matrix
$$
{\cal D}_{{\cal B}_h}(T_\ell^{\rm loc}|_{V_{\pi,\ell;h}^{K_\ell}})
=\left(
\begin{array}{ccc}
\lambda_{\pi,\ell;h}&&*\\
&\ddots&\\
&&\lambda_{\pi,\ell;h}\\
\end{array}
\right).
$$
Clearly, $V_{\pi,\ell}^{K_\ell}$ decomposes
$$
V_{\pi,\ell}^{K_\ell}=\bigoplus_{h} V_{\pi,\ell;h}^{K_\ell}.
$$

%[[Beachte $V_{\pi,\ell,h}$ ist i. allg. kein irreduzibler Jordanblock !]]

For $\pi\in {\cal R}_k$ and any tuple $h=(h_\ell)_{\ell|Np}$ we define $\lambda_{\pi;h}=(\lambda_{\pi;h,\ell})_\ell\in{\Bbb C}^P$ by letting $\lambda_{\pi;h,\ell}=\lambda_{\pi,\ell}$ if $\ell\not|Np$ and $\lambda_{\pi;h,\ell}=\lambda_{\pi,\ell;h_\ell}$ if $\ell|Np$. Moreover, we set 
$$
V_{\pi;h}^K=\otimes_{\ell\not|Np} V_{\pi,\ell}^{K_\ell}\otimes\otimes_{\ell|Np} V_{\pi,\ell;h_\ell}^{K_\ell}.
$$
Clearly,
$$
V_{\pi;h}^K={\cal M}_{k,\chi\omega^{-k}}^\alpha(\lambda_{\pi;h})
$$
and
$$
{\cal M}_k(\Gamma)\cong\bigoplus_{\pi\in {\cal R}_{k-2}}\bigoplus_{h=(h_\ell)_{\ell|Np}} V_{\pi;h}^K.
$$

Finally, we obtain the adelic version of the decomposition (10) of ${\cal M}_k(\Gamma,\chi\omega^{-k})^\alpha$ as follows. We have $V_{\pi;h}^K\le {\cal M}_k(\Gamma,\chi\omega^{-k})^\alpha$ precisely if
$$
\pi\in{\cal R}_{k-2},\quad \omega_\pi|_{\hat{\Bbb Z}}=\chi\quad\mbox{ and } \quad v_p(\lambda_{\pi,p;h_p})=\alpha.
$$
Here, $\omega_\pi$ is the central character of $\pi$; note that since ${\Bbb Q}^*\backslash {\Bbb A}^*\cong \prod_{\ell\not=\infty}{\Bbb Z}_\ell^*\times{\Bbb R}^*_{>0}$, the idele class character $\omega_\pi:\,{\Bbb Q}^*\backslash {\Bbb A}^*\rightarrow{\Bbb C}^*$ is uniquely determined by the pair $(\tilde{\omega}_\pi,k)$, where $\tilde{\omega}_\pi=\omega_\pi|_{\hat{\Bbb Z}}$ and $\omega_\pi|_{{\Bbb R}^*_{>0}}=|\cdot|^{2-k}_\infty$ (cf. the definition of the representation $I_k$ in section 1.3). 
%
%[[Beachte: es gilt 
%$$
%e_\alpha V_{\pi,h}=\otimes_{\ell\not|Np} V_{\pi,\ell}\otimes_{\ell|N} V_{\pi,\ell,h_\ell}\otimes \{e_\alpha V_{\pi,p,h_p}\}
%$$
%]]
%
We deduce that ${\cal M}_{k,\chi\omega^{-k}}^\alpha(\lambda)$ vanishes unless there are $\pi$ and $h$ such that the above conditions hold and  $\lambda=\lambda_{\pi,h}$. Moreover, in this case we obtain
$$
{\cal M}_{k,\chi\omega^{-k}}^\alpha(\lambda)=V_{\pi;h}^K.
$$
Notice that $v_p(\lambda_{\pi,p;h_p})=\alpha$ implies $V_{\pi;h}^K=V_{\pi;h}^{K,\alpha}$. 
%
%[[Beachte, dass unter der bedingung $v_p(\lambda_{\pi,p,h_p})=\alpha$ gilt $V_{\pi,h}=V_{\pi,h}^\alpha$, denn $V_{\pi,h}$ ist ja schon ein (verallgemeinerter) $T_p$-Eigenraum.]]
%
In particular, we obtain the adelic analogoue of (10)
\begin{eqnarray*}
(tl:12)\qquad{\cal M}_k(\Gamma,\chi\omega^{-k})^\alpha&=&\bigoplus_{\pi\in {\cal R}_{k-2}\atop \omega_\pi|_{\hat{\Bbb Z}}=\chi\omega^{-k}}\bigoplus_{h=(h_\ell)_{\ell|Np}} V_{\pi;h}^\alpha\\
&=&\bigoplus_{\pi\in {\cal R}_{k-2}\atop \omega_\pi|_{\hat{\Bbb Z}}=\chi\omega^{-k}}\bigoplus_{h=(h_\ell)_{\ell|Np}\atop v_p(\lambda_{\pi;h,p})=\alpha} V_{\pi;h}.\\
\end{eqnarray*}

%[[Beachte $V_{\pi,h}$ ist $T_\ell$-invariant fur all $\ell$ (nach Definition von $V_{\pi,h}$); insbesondere ist $V_{\pi,h}$ $T_p$-invariant
%]]

%[[Bemerkung. 1.) Der slope $\alpha$ einer Modulform $f$ ist eine invariante der Modulform, nicht aber der zu $f$ gehorigen Darstellung $\pi_f$ !! Beweis: die beiden Altformen $f_1,f_2$ zum level $Np$, die man aus der Neuform $f$ zum Level $N$ 
%erhalt haben i. allg. verschiedenen slope (vgl. [M-G], p. 796, der Abschnitt "$p$-oldforms").
%2.) Aber: auf ${\cal M}_{k,\chi\omega^{-k}}^\alpha(\lambda)$ ist der slope konstant, d.h. alle $f\in {\cal M}_{k,\chi\omega^{-k}}^\alpha(\lambda)$ haben den gleichen slope ! $\Rightarrow$ slope ist eine Invariante von ${\cal M}_{k,\chi\omega^{-k}}^\alpha(\lambda)$.
%]]

We obtain analogous decompositions of the space of modular forms if we only use primes $\ell$, which do not divide a certain fixed integer. More precisely, let $M\in{\Bbb N}$ and let $\lambda=(\lambda_\ell)_{\ell\not|M}$ be a sequence of complex numbers, whose entries are defined for all primes not dividing $M$. We set
$$
{\cal M}_{k,\chi\omega^{-k}}^\alpha(\lambda)=\{f\in{\cal M}_k(\Gamma,\chi\omega^{-k})^\alpha:\mbox{for all $\ell\not|M$ there is $n=n_\lambda\in{\Bbb N}$: $(T_\ell-\lambda_\ell)^nf=0$}\}.
$$
We define
$$
\Lambda_{k,\chi,M}^\alpha=\{\lambda=(\lambda_\ell)_{\ell\not|M}:\,{\cal M}_{k,\chi\omega^{-k}}^\alpha(\lambda)\not=0\}.
$$
Quite similar to equations (9a,b,c,...) and (10) we obtain
$$
{\cal M}_k(\Gamma,\chi\omega^{-k})^\alpha=\bigoplus_{\lambda\in\Lambda_{k,\chi,M}^\alpha} {\cal M}_{k,\chi\omega^{-k}}^\alpha(\lambda).\leqno(13^\alpha)
$$
We also define the spaces ${\cal M}_{k,\chi\omega^{-k}}(\lambda)$ and the set $\Lambda_{k,\chi,M}$: we set
$$
{\cal M}_{k,\chi\omega^{-k}}(\lambda)=\{f\in{\cal M}_k(\Gamma,\chi\omega^{-k}):\mbox{for all $\ell\not|M$ there is $n\in{\Bbb N}$: $(T_\ell-\lambda_\ell)^nf=0$}\}.
$$
and 
$$
\Lambda_{k,\chi,M}=\{\lambda=(\lambda_\ell)_{\ell\not|M}:\,{\cal M}_{k,\chi\omega^{-k}}(\lambda)\not=0\}.
$$ 
Again, 
$$
{\cal M}_k(\Gamma,\chi\omega^{-k})=\bigoplus_{\lambda\in\Lambda_{k,\chi,M}} {\cal M}_{k,\chi\omega^{-k}}(\lambda).\leqno(13)
$$
Clearly, in case $p\not|M$ we have
$$
{\cal M}_{k,\chi\omega^{-k},M}^\alpha(\lambda)=\left\{
\begin{array}{ccc}
0&{\rm if}&v_p(\lambda_p)\not=\alpha\\
{\cal M}_{k,\chi\omega^{-k},M}(\lambda)&{\rm if}&v_p(\lambda_p)=\alpha.
\end{array}
\right.
$$
We note that in what follows we shall mostly consider the case $M=Np$. 

It is immediate by the definitions that
$$
{\cal M}_{k,\chi\omega^{-k}}^\alpha(\lambda)={\cal M}_k(\Gamma,\chi\omega^{-k})^\alpha\cap {\cal M}_{k,\chi\omega^{-k}}(\lambda).
$$
We claim that 
$$
{\cal M}_{k,\chi\omega^{-k}}^\alpha(\lambda)={\cal M}_{k,\chi\omega^{-k}}(\lambda)^\alpha
$$
i.e. ${\cal M}_{k,\chi\omega^{-k}}^\alpha(\lambda)$ is the slope $\alpha$ subspace of ${\cal M}_{k,\chi\omega^{-k}}(\lambda)$.  In order to prove this we verify both inclusions. "$\subseteq$": Clearly, ${\cal M}_{k,\chi\omega^{-k}}^\alpha(\lambda)\subseteq{\cal M}_{k,\chi\omega^{-k}}(\lambda)$. On the other hand, all eigenvalues of $T_p$ on ${\cal M}_{k,\chi\omega^{-k}}^\alpha(\lambda)$ have $p$-adic value equal to $\alpha$ because ${\cal M}_{k,\chi\omega^{-k}}^\alpha(\lambda)\subseteq{\cal M}_k(\Gamma,\chi\omega^{-k})^\alpha$, hence, ${\cal M}_{k,\chi\omega^{-k}}^\alpha(\lambda)$ must be contained in the slope $\alpha$ subspace of ${\cal M}_{k,\chi\omega^{-k}}(\lambda)$ (use that $V^\alpha=\bigoplus_{v_p(\gamma)=\alpha}V_\gamma$, where $V_\gamma\subset V$ is the generalized eigenspace attached to $\gamma$). "$\supseteq$": Clearly, ${\cal M}_{k,\chi\omega^{-k}}(\lambda)^\alpha\subseteq {\cal M}_{k,\chi\omega^{-k}}(\lambda)$. Moreover, since all eigenvalues of $T_p$ on ${\cal M}_{k,\chi\omega^{-k}}(\lambda)^\alpha$ have $p$-adic value equal to $\alpha$ we deduce that ${\cal M}_{k,\chi\omega^{-k}}(\lambda)^\alpha\subseteq{\cal M}_k(\Gamma,\chi\omega^{-k})^\alpha$. Hence, altogether we see that ${\cal M}_{k,\chi\omega^{-k}}(\lambda)^\alpha$ is contained in the intersection, i.e. ${\cal M}_{k,\chi\omega^{-k}}(\lambda)^\alpha$ is contained in ${\cal M}_{k,\chi\omega^{-k}}^\alpha(\lambda)$. Thus, the proof is complete.

%The space ${\cal M}_{k,\chi\omega^{-k}}(\lambda)$ decomposes
%$$
%{\cal M}_{k,\chi\omega^{-k}}(\lambda)=\bigoplus_{\lambda_p} {\cal M}_{k,\chi\omega^{-k}}(\lambda,\lambda_p),
%$$
%where ${\cal M}_{k,\chi\omega^{-k}}(\lambda,\lambda_p)$ is the generalized eigenspace in ${\cal M}_{k,\chi\omega^{-k}}(\lambda)$ attached to the operator $T_p$ and the eigenvalue $\lambda_p$; we obtain
%$$
%{\cal M}_{k,\chi\omega^{-k}}^\alpha(\lambda)=\bigoplus_{\lambda_p,\,v_p(\lambda_p)=\alpha} {\cal M}_{k,\chi\omega^{-k}}(\lambda,\lambda_p),
%$$
%i.e. ${\cal M}_{k,\chi\omega^{-k}}^\alpha(\lambda)$ is the slope-$\alpha$ subspace of ${\cal M}_{k,\chi\omega^{-k}}(\lambda)$. 

Let $(\pi,V_\pi)$ be a (finite part of a cuspidal automorphic) representation of ${\rm GL}_2({\Bbb A}_f)$ such that $V_\pi^{K}$ occurs in ${\cal M}_k(\Gamma,\chi\omega^{-k})$. $V_\pi^K$ decomposes $V_\pi^K=\otimes_\ell V_{\pi,\ell}^{K_\ell}$, where 
$$
V_{\pi,\ell}^{K_\ell}=\langle v_{\pi,\ell}^0 \rangle
$$ 
if $\ell\not|Np$ and $T_\ell$ acts on $V_{\pi,\ell}^{K_\ell}$, $\ell\not|Np$, as multiplication by some scalar $\lambda_{\pi,\ell}$. We set
$$
\lambda_{\pi}=(\lambda_{\pi,\ell})_{\ell\not|Np}\in\Lambda_{k,\chi,Np}.
$$ 
We denote by $V_{\pi}^{K,\alpha}$ resp. $V_{\pi,\ell}^{K_\ell,\alpha}$ the slope $\alpha$ subspace of $V_\pi^K$ resp. $V_{\pi,\ell}^{K_\ell}$ with respect to the operator $T_p$ resp. $T_p^{\rm loc}$. 
${\cal M}_k(\Gamma,\chi\omega^{-k})$ is the direct sum
$$
{\cal M}_k(\Gamma,\chi\omega^{-k})=\bigoplus_{\pi\in{\cal R}_{k-2}\atop \omega_\pi|_{\hat{\Bbb Z}}=\chi\omega^{-k}} V_\pi^{K}.\leqno(14)
$$
As an immediate consequence, we obtain
$$
{\cal M}_k(\Gamma,\chi\omega^{-k})^\alpha=\bigoplus_{\pi\in{\cal R}_{k-2}\atop \omega_\pi|_{\hat{\Bbb Z}}=\chi\omega^{-k}} V_\pi^{K,\alpha}.\leqno(14^\alpha)
$$
%
%[[Beachte: ${\cal M}_k(\Gamma,\chi\omega^{-k})^\alpha$ ist die Summe der Unterraume ${\cal M}_k(\Gamma,\chi\omega^{-k})^\alpha(\lambda_p)$
%mit $v_p(\lambda_p)=\alpha$ und jedes $V_\pi^K$ ist $T_p$-invariant !
%]]
%
Clearly, $V_\pi^K$ is contained in ${\cal M}_{k,\chi\omega^{-k}}(\lambda_{\pi})$. Moreover, the strong multiplicity-1 Theorem implies that  $\pi$ is the only representation in ${\cal R}_{k-2}$ such that $V_\pi$ is contained in ${\cal M}_{k,\chi\omega^{-k}}(\lambda_{\pi})$,
%
%[[$V_{\pi'}\subset {\cal M}_{k,\chi}^\alpha(\lambda_{\pi})$ $\Rightarrow \lambda_{\pi',\ell}=\lambda_{\pi,\ell}$ for all $\ell\not|Np$
%$\Rightarrow$ $\pi=\pi'$ by strong multiplicity-1 theorem.
%]]
%
i.e. the assignment $\pi\mapsto\lambda_\pi$ is injective. Thus, comparing (13) and (14) we obtain
$$
{\cal M}_{k,\chi\omega^{-k}}(\lambda_{\pi})=V_\pi^K.\leqno(15)
$$ 
Taking into account that ${\cal M}_{k,\chi\omega^{-k}}^\alpha(\lambda)={\cal M}_{k,\chi\omega^{-k}}(\lambda)^\alpha$ we obtain as an immediate consequence
$$
V_\pi^{K,\alpha}={\cal M}_{k,\chi\omega^{-k}}^\alpha(\lambda_\pi).\leqno(15^\alpha)
$$
In view of (15) we define the conductor of $\lambda=(\lambda_\ell)_{\ell\not|Np}$ as the conductor of $\pi$, where $\pi$ is the representation of ${\rm GL}_2({\Bbb A}_f)$ corresponding to $\lambda$ (i.e. $\lambda_\pi=\lambda$). For later purpose, we note that $V_\pi^{K,\alpha}$ decomposes
$$
V_\pi^{K,\alpha}=\otimes_{\ell\not=p} V_{\pi,\ell}^{K_\ell}\otimes V_{\pi,p}^{K_p,\alpha}.\leqno(16)
$$

Finally, let $\lambda=(\lambda_\ell)_{\ell\not|Np}$ and $\lambda'=(\lambda'_\ell)_{\ell\not|Np}$ be contained in $\Lambda_{k,\chi,Np}^\alpha$ and assume that $\lambda\not=\lambda'$. The strong multiplicity-1 Theorem implies that there is a prime $\ell(\lambda,\lambda')$ such that $\ell(\lambda,\lambda')$ does not divide the level $Np$ and $\lambda_{\ell(\lambda,\lambda')}\not=\lambda_{\ell(\lambda,\lambda')}'$. In particular, there is an integer $A_{\lambda,\lambda'}$ such that 
$$
\lambda_{\ell(\lambda,\lambda')}\not\equiv\lambda_{\ell(\lambda,\lambda')}'\pmod{p^{A_{\lambda,\lambda'}}}.\leqno(17)
$$ 

We set $P_k=P_{k,\chi}^\alpha=\{\ell(\lambda,\lambda'),\,\lambda,\lambda'\in\Lambda_{k,\chi,Np}^\alpha\}$.

{\bf Notation. } We denote by $A_k=A_{k,\chi}^\alpha$ any integer which is greater than all $A_{\lambda,\lambda'}$, $\lambda,\lambda'\in\Lambda_{k,\chi,Np}^\alpha$. In different words, $A_k$ is any integer such that 
$$
\lambda=\lambda'\Leftrightarrow \lambda_\ell\equiv\lambda_\ell'\pmod{p^{A_k}}\quad\mbox{for all}\;\ell\in P_k.\leqno(18)
$$

For example, we may choose $A_k={\rm max}_{\lambda,\lambda'\in\Lambda_{k,\chi,Np}^\alpha}\,\{A_{\lambda,\lambda'}\}$. Notice that $P_k$ is a finite set and that the maximum exists, both because the set $\Lambda_{k,\chi,Np}^\alpha$ is finite. 

%
%[[We note that any element $\lambda=(\lambda_\ell)_{\ell\not|M}\in\Lambda_{k,\chi,M}$ with $M|Np$ (e.g. $M=N$) defines via restriction an element $\lambda^*=(\lambda_\ell)_{\ell\not|Np}$ in $\Lambda_{k,\chi,Np}$. In particular, the conductor of $\lambda$ is defined as the conductor of $\lambda^*$ and for any pair of elements $\lambda=(\lambda_\ell)_{\ell\not|N}$ and $\lambda'=(\lambda'_\ell)_{\ell\not|N}$ the primes $\ell(\lambda,\lambda')$ are defined as those attached to the pair of restrictions $\lambda^*$ and ${\lambda'}^*$.
%]]

Let $\lambda=(\lambda_\ell)_{\ell\not|N}\in\Lambda_{k,\chi,Np}^\alpha$. Since the projection operator $\pi^{Np}_{Fp}$ ($F|N$) commutes with the Hecke operators $T_\ell$, $\ell\not|Np$ as well as with $T_p$ (use the Corollary in section 3.2 as well as the Lemma 2, part 1.) in section 3.5) we deduce that $\pi^{Np}_{Fp}$ leaves the space ${\cal M}_{k,\chi\omega^{-k}}^\alpha(\lambda)$ invariant, i.e. for any $F|N$ we obtain a map
$$
\pi^{Np}_{Fp}:\,{\cal M}_{k,\chi\omega^{-k}}^\alpha(\lambda)\rightarrow {\cal M}_{k,\chi\omega^{-k}}^\alpha(\lambda)^{\Gamma_1(Fp)}.
$$

{\bf Theorem. }(Transfer from weight $k$ to weight $k'$). {\it Let $\alpha\in{\Bbb Q}_{\ge 0}$ and assume that $k,k'\in{\Bbb N}$ satisfy

$\bullet$ $k,k'\ge ({{\sf C}_\alpha}\alpha+1)^2+2$ 

$\bullet$ $k\equiv k'\pmod{p^m}$ where $m\ge {{\sf B}_\alpha}$. 

Let $\lambda=(\lambda_\ell)_{\ell\not|Np}\in\Lambda_{k,\chi,Np}^\alpha$. 

1.) There is $\lambda'=(\lambda_\ell')_{\ell\not|Np}\in\Lambda_{k',\chi,Np}^\alpha$ such that

$\bullet$ ${\cal M}_{k',\chi\omega^{-k'}}^\alpha(\lambda')^{\Gamma_1(Fp)}\not=(0)$, i.e. the conductor of $\lambda'$ is a divisor of $Fp$

$\bullet$ $\lambda_\ell'\equiv\lambda_\ell\pmod{p^{\sf D}}$ for all $\ell$ which do not divide $Np$.

Here, $F$ is the prime to $p$ part of the conductor of $\lambda$, i.e. ${\rm cond}\,(\lambda)=F\wp$, where $(F,p)=1$ and $\wp|p$ and 
$$
{\sf D}={\sf D}_{\alpha,m}=\frac{\left[\frac{m}{{{\sf C}_\alpha}\alpha}\right]-v_p(\varphi(N)[{\rm GL}_2({\Bbb Z}):\Gamma_1(Np)]d_{k,\chi\omega^{-k}}^\alpha!)-1}{d_{k,\chi\omega^{-k}}^\alpha}-A_k\leqno(19)
$$
if $\alpha>0$ and
$$
{\sf D}={\sf D}_{\alpha,m}=\frac{m-v_p(\varphi(N)[{\rm GL}_2({\Bbb Z}):\Gamma_1(Np)]d_{k,\chi\omega^{-k}}^\alpha!)-1}{d_{k,\chi\omega^{-k}}^\alpha}-A_k\leqno(19^0)
$$
if $\alpha=0$.
%
%[[Note that ${\sf D}\in{\Bbb Q}$, i.e. congruences modulo $p^{\sf D}$ make sense]]

2.) Assume in addition that ${\rm dim}\,{\cal M}_k(\Gamma,\chi\omega^{-k})^\alpha=1$. Then there is  $\lambda'=(\lambda_\ell')_{\ell\not|Np}\in\Lambda_{k',\chi,Np}^\alpha$ such that

$\bullet$ ${\cal M}_{k',\chi\omega^{-k'}}^\alpha(\lambda')^{\Gamma_1(Fp)}\not=(0)$, i.e. the conductor of $\lambda'$ is a divisor of $Fp$

$\bullet$ $\lambda_\ell'\equiv\lambda_\ell\pmod{p^{\sf D}}$ for all $\ell$ which do not divide $Np$.

where
$$
{\sf D}={\sf D}_{\alpha,m}=\frac{\left[\frac{m}{{{\sf C}_\alpha}\alpha}\right]-v_p(\varphi(N)[{\rm GL}_2({\Bbb Z}):\Gamma_1(Np)])-1}{d_{k,\chi\omega^{-k}}^\alpha}\leqno(19)
$$
if $\alpha>0$ and
$$
{\sf D}={\sf D}_{\alpha,m}=\frac{m-v_p(\varphi(N)[{\rm GL}_2({\Bbb Z}):\Gamma_1(Np)])-1}{d_{k,\chi\omega^{-k}}^\alpha}\leqno(19^0)
$$
if $\alpha=0$.

}

{\it Proof. } 1.) We fix an arbitrary element $\lambda=(\lambda_\ell)_{\ell\not|Np}$ in $\Lambda_{k,\chi,Np}^\alpha$ and we assume that there is no $\lambda'=(\lambda_\ell')_{\ell\not|Np}\in \Lambda_{k',\chi,Np}^\alpha$ satisfying the conditions of the Theorem.

We denote by $\Lambda_{k',1}^\alpha$ the set of $\lambda'\in \Lambda_{k',\chi,Np}^\alpha$ such that there is a prime $\ell(\lambda')$, which does not divide $Np$ and which satisfies 
$$
\lambda_{\ell(\lambda')}'\not\equiv \lambda_{\ell(\lambda')}\pmod{p^{\sf D}}.
$$ 
Moreover, by $\Lambda_{k',2}^\alpha$ we denote the set of elements $\lambda'\in \Lambda_{k',\chi,Np}^\alpha$, which are not contained in $\Lambda_{k',1}^\alpha$ and which satisfy 
$$
{\cal M}_{k',\chi\omega^{-k'}}^\alpha(\lambda')^{\Gamma_1(Fp)}=(0). 
$$

Our assumption that there is no $\lambda'$ as in the Theorem then is equivalent to 
$$
\Lambda_{k',\chi,Np}^\alpha=\Lambda_{k',1}^\alpha\dot{\cup}\Lambda_{k',2}^\alpha\qquad\mbox{(disjoint union).}
$$

We set
$$
e=\prod_{\mu\in \Lambda_{k,\chi,Np}^\alpha\atop \mu\not=\lambda}(T_{\ell(\mu,\lambda)}-\mu_{\ell(\mu,\lambda)})
\prod_{\mu'\in\Lambda_{k',1}^\alpha}(T_{\ell(\mu')}-\mu'_{\ell(\mu')})\in{\cal H}
$$
and 
$$
\tilde{e}=\pi^{Np}_{Fp}\,e\,e_\alpha^{[\frac{m}{{\sf C}_\alpha\alpha}]}.
$$
Here, $\ell(\mu,\lambda)$ is the prime defined before equation (17). Note that $\ell(\mu,\lambda)$ does not divide $Np$, i.e. the term $\mu_{\ell(\mu,\lambda)}$ makes sense.

We determine ${\rm tr}\,\tilde{e}|_{{\cal M}_k(\Gamma,\chi\omega^{-k})^\alpha}$ and ${\rm tr}\,\tilde{e}|_{{\cal M}_{k'}(\Gamma,\chi\omega^{-k'})^\alpha}$ and we start with a general observation. Let $V={\cal M}_{k,\chi\omega^{-k}}^\alpha(\gamma)$, $\gamma\in\Lambda_{k,\chi,Np}^\alpha$ or $V={\cal M}_{k',\chi\omega^{-k'}}^\alpha(\gamma')$, $\gamma'\in\Lambda_{k',\chi,Np}^\alpha$. Let $v_1,\ldots,v_s$ be a basis of $V^{\Gamma_1(Fp)}$ and extend to a basis ${\cal B}=\{v_1,\ldots,v_s,v_{s+1},\ldots,v_m\}$ of $V$. With respect to ${\cal B}$ we obtain
$$
{\cal D}_{\cal B}(\pi^{Np}_{Fp}|_{V})=\left(
\begin{array}{cccccc}
[\Gamma_1(Fp):\Gamma_1(Np)]&&&&*&\\
&\ddots&&&&\\
&&[\Gamma_1(Fp):\Gamma_1(Np)]&&&\\
&&&0&\cdots&0\\
&&&&\ddots&\vdots\\
&&&&&0\\
\end{array}
\right).\leqno(20)
$$
Thus $\pi^{Np}_{Fp}|_{V}$ has eigenvalue $0$ with multiplicity ${\rm dim}\,V-{\rm dim}\,V^{\Gamma_1(Fp)}$ and eigenvalue $[\Gamma_1(Fp):\Gamma_1(Np)]$ with multiplicity ${\rm dim}\,V^{\Gamma_1(Fp)}$. Since $\pi^{Np}_{Fp}$ and $T_\ell$, $\ell\not|Np$ and $e_\alpha$ commute, they are simultaneously trigonalizable, i.e. there is a basis ${\cal B}$ of $V$ such that
$$
{\cal D}_{\cal B}(T_\ell|_{V})=\left(
\begin{array}{ccc}
\kappa&&*\\
&\ddots&\\
&&\kappa\\
\end{array}
\right)\leqno(21)
$$
for all $\ell\not|Np$, where $\kappa=\gamma_\ell$ or $\kappa=\gamma'_\ell$ (note that $T_\ell$ has only eigenvalue $\gamma_\ell$ resp. $\gamma_\ell'$ on ${\cal M}_{k,\chi\omega^{-k}}^\alpha(\gamma)$ resp. ${\cal M}_{k',\chi\omega^{-k'}}^\alpha(\gamma')$),
$$
{\cal D}_{\cal B}(\pi^{Np}_{Fp}|_V)=\left(
\begin{array}{ccc}
\kappa_1&&*\\
&\ddots&\\
&&\kappa_m\\
\end{array}
\right)\leqno(22)
$$
and
$$
{\cal D}_{\cal B}(e_\alpha^{[\frac{m}{{\sf C}_\alpha\alpha}]}|_{V})=\left(
\begin{array}{ccc}
\epsilon&&*\\
&\ddots&\\
&&\epsilon\\
\end{array}
\right),\leqno(23)
$$
where $\epsilon\equiv 1\pmod{p}$ (cf. Proposition 1.) in section 3.4). Moreover, (20) implies that precisely ${\rm dim}\,V^{\Gamma_1(Np)}$ of the diagonal entries of (22) equal $[\Gamma_1(Fp):\Gamma_1(Np)]$, while the remaining ones are equal to $0$. Using (21), (22) and (23) we see that
$$
{\rm tr}\,\tilde{e}|_V=\kappa \epsilon [\Gamma_1(Fp):\Gamma_1(Np)] {\rm dim}\,V^{\Gamma_1(Fp)}.
$$
We shall use this in the following computation of ${\rm tr}\,\tilde{e}|_{{\cal M}_k(\Gamma,\chi\omega^{-k})^\alpha}$ and ${\rm tr}\,\tilde{e}|_{{\cal M}_{k'}(\Gamma,\chi\omega^{-k'})^\alpha}$. To compute these traces we look at the spaces ${\cal M}_{k,\chi\omega^{-k}}^\alpha(\gamma)$, $\gamma\in\Lambda_{k,\chi,Np}^\alpha$ and ${\cal M}_{k',\chi\omega^{-k'}}^\alpha(\gamma')$, $\gamma'\in\Lambda_{k',\chi,Np}^\alpha$ individually.

$\bullet$ The spaces ${\cal M}_{k',\chi\omega^{-k'}}^\alpha(\gamma')$, $\gamma'\in\Lambda_{k',1}^\alpha$. The definition of $e$ and equation (21) imply that
$$
{\cal D}_{\cal B}(e|_{{\cal M}_{k',\chi\omega^{-k'}}^\alpha(\gamma)})=\left(
\begin{array}{ccc}
\zeta&&*\\
&\ddots&\\
&&\zeta\\
\end{array}
\right),
$$
where
$$
\zeta=\prod_{\mu\in \Lambda_{k,\chi,Np}^\alpha\atop \mu\not=\lambda}(\gamma'_{\ell(\mu,\lambda)}-\mu_{\ell(\mu,\lambda)})
\prod_{\mu'\in\Lambda_{k',1}^\alpha}(\gamma'_{\ell(\mu')}-\mu'_{\ell(\mu')}).
$$
Since the factor corresponding to $\mu'=\gamma'$ vanishes, $\zeta$ equals $0$ and we deduce that ${e}|_{{\cal M}_{k',\chi\omega^{-k'}}^\alpha(\gamma')}=0$. Hence, since $e_\alpha\in{\Bbb C}[T_p]$ leaves the space ${\cal M}_{k',\chi\omega^{-k'}}^\alpha(\gamma')$ invariant we obtain $\tilde{e}|_{{\cal M}_{k',\chi\omega^{-k'}}^\alpha(\gamma')}=0$.

$\bullet$ The spaces ${\cal M}_{k',\chi\omega^{-k'}}^\alpha(\gamma')$, $\gamma'\in\Lambda_{k',2}^\alpha$.  Since 
${\cal M}_{k',\chi\omega^{-k'}}^\alpha(\gamma')$ is stable under all Hecke operators $T_\ell$, $\ell\not|Np$, it is also stable under $e$. Moreover, since ${\cal M}_{k',\chi\omega^{-k'}}^\alpha(\gamma')^{\Gamma_1(Fp)}=(0)$ by the definition of $\Lambda_{k',2}^\alpha$, the operator $\pi^{Np}_{Fp}$ annihilates 
${\cal M}_{k',\chi\omega^{-k'}}^\alpha(\gamma')$. Hence, $\tilde{e}|_{{\cal M}_{k',\chi\omega^{-k'}}^\alpha}(\gamma')=0$. 

Thus, altogether we have seen that $\tilde{e}|_{{\cal M}_{k'}(\Gamma,\chi\omega^{-k'})^\alpha}=0$ and therefore
$$
{\rm tr}\,\tilde{e}|_{{\cal M}_{k'}(\Gamma,\chi\omega^{-k'})^\alpha}=0.\leqno(24)
$$

$\bullet$ The spaces ${\cal M}_{k,\chi\omega^{-k}}^\alpha(\gamma)$, $\gamma\in\Lambda_{k,\chi,Np}^\alpha$, $\gamma\not=\lambda$. Again, the definition of $e$ and equation (21) imply that
$$
{\cal D}_{\cal B}(e|_{{\cal M}_{k,\chi\omega^{-k}}^\alpha(\gamma)})=\left(
\begin{array}{ccc}
\zeta&&*\\
&\ddots&\\
&&\zeta\\
\end{array}
\right),
$$
where
$$
\zeta=\prod_{\mu\in \Lambda_{k,\chi,Np}^\alpha\atop \mu\not=\lambda}(\gamma_{\ell(\mu,\lambda)}-\mu_{\ell(\mu,\lambda)})
\prod_{\mu'\in\Lambda_{k',1}^\alpha}(\gamma_{\ell(\mu')}-\mu'_{\ell(\mu')}).
$$
Since the factor corresponding to $\mu=\gamma$ vanishes $\zeta$ equals $0$ and we deduce that
$e|_{{\cal M}_{k,\chi\omega^{-k}}^\alpha(\gamma)}=0$. Hence, we obtain $\tilde{e}|_{{\cal M}_{k,\chi\omega^{-k}}^\alpha(\gamma)}=0$ and therefore ${\rm tr}\,\tilde{e}|_{{\cal M}_{k,\chi\omega^{-k}}^\alpha(\gamma)}=0$.

$\bullet$ The space ${\cal M}_{k,\chi\omega^{-k}}^\alpha(\lambda)$ (i.e. $\gamma=\lambda$). Using equations (21), (22) and (23) we find
$$
{\rm tr}\,\pi^{Np}_{Fp}ee_\alpha^{\left[\frac{m}{{\sf C}_\alpha\alpha}\right]}|_{{\cal M}_{k,\chi\omega^{-k}}^\alpha(\lambda)}=\zeta \epsilon\,[\Gamma_1(Fp):\Gamma_1(Np)]\,d,
$$
where $d={\rm dim}\,{{\cal M}_{k,\chi\omega^{-k}}^\alpha}(\lambda)^{\Gamma_1(Fp)}$, $v_p(\epsilon)=0$ and
$$
\zeta=\prod_{\mu\in \Lambda_{k,\chi,Np}^\alpha\atop \mu\not=\lambda}(\lambda_{\ell(\mu,\lambda)}-\mu_{\ell(\mu,\lambda)})
\prod_{\mu'\in\Lambda_{k',1}^\alpha}(\lambda_{\ell(\mu')}-\mu'_{\ell(\mu')}).\leqno(25)
$$
Taking into account that 

- $|\Lambda_{k,\chi,Np}^\alpha|\le {\rm dim}\,{\cal M}_k(\Gamma,\chi\omega^{-k})^\alpha= d_{k,\chi\omega^{-k}}^\alpha$ 

- $|\Lambda_{k',1}^\alpha|\le {\rm dim}\,{\cal M}_{k'}(\Gamma,\chi\omega^{-k'})^\alpha= d_{k',\chi\omega^{-k'}}^\alpha=d_{k,\chi\omega^{-k}}^\alpha$, because $d_{k',\chi\omega^{-k'}}^\alpha=d_{k,\chi\omega^{-k}}^\alpha$ by Theorem B 

%[[Hier ist der neuralgische Punkt: Nur fur slope $\alpha$ Raume ${\cal M}_k(\Gamma,\chi\omega^{-k})$ wissen wir dass die Dimension UNABHANGIG VOM GEWICHT $k$ beschrankt ist, d.h. es gilt $|\Lambda_{k,\chi,Np}^\alpha|\le d(\alpha)$. Die dimension von ganz ${\cal M}_k(\Gamma,\chi\omega^{-k})$ ist nicht unabhangig von $k$ beschrankt, d.h. eine ungleichung $|\Lambda_{k,\chi,Np}|\le d(\alpha)$ fur alle $k$ KANN NICHT GELTEN !!
%]]

- $\lambda_{\ell(\lambda,\mu)}\not\equiv \mu_{\ell(\lambda,\mu)}\pmod{p^{A_k}}$ (cf. equation (17)) 

and 

- $\lambda_{\ell(\mu')}\not\equiv\mu'_{\ell(\mu')}\pmod{p^{\sf D}}$ because $\mu'$ runs over $\Lambda_{k',1}^\alpha$

we deduce from (25) that
$$
v_p(\zeta)\le d_{k,\chi\omega^{-k}}^\alpha A_k+ d_{k,\chi\omega^{-k}}^\alpha {\sf D}= d_{k,\chi\omega^{-k}}^\alpha (A_k+{\sf D}).
$$
Hence, taking into account that $v_p(d)\le v_p(d_{k,\chi\omega^{-k}}^\alpha!)$ (note that $d\le d_{k,\chi\omega^{-k}}^\alpha$),
%
%[[da $d\le d_{k,\chi\omega^{-k}}^\alpha \Rightarrow v_p(d)\le v_p(d_{k,\chi\omega^{-k}}^\alpha!$]]
%
$v_p(\epsilon)=0$ and $[\Gamma_1(Fp):\Gamma_1(Np)]\big| [{\rm GL}_2({\Bbb Z}):\Gamma_1(Np)]$ we obtain

\begin{eqnarray*}
v_p({\rm tr}\,(\tilde{e}|_{{\cal M}_{k,\chi\omega^{-k}}^\alpha(\lambda)}))&=&v_p([\Gamma_1(Fp):\Gamma_1(Np)])+v_p(\zeta)+v_p(d) \\
&\le&v_p([{\rm GL}_2({\Bbb Z}):\Gamma_1(Np)])+d_{k,\chi\omega^{-k}}^\alpha(A_k+{\sf D})+v_p(d_{k,\chi\omega^{-k}}^\alpha!).\\
\end{eqnarray*}

Thus, we obtain altogether
\begin{eqnarray*}
(26)\quad v_p({\rm tr} (\tilde{e}|_{{\cal M}_k(\Gamma,\chi\omega^{-k})^\alpha}))&=&v_p({\rm tr}\,\tilde{e}|_{{\cal M}_{k,\chi\omega^{-k}}^\alpha(\lambda)})\\
&\le& v_p([{\rm GL}_2({\Bbb Z}):\Gamma_1(Np)])+d_{k,\chi\omega^{-k}}^\alpha(A_k+{\sf D})+v_p(d_{k,\chi\omega^{-k}}^\alpha!).\\
\end{eqnarray*}

To proceed we distinguish cases and first assume $\alpha>0$. Since ${\sf D}$ is such that 
$$
v_p([{\rm GL}_2({\Bbb Z}):\Gamma_1(Np)])+d_{k,\chi\omega^{-k}}^\alpha(A_k+{\sf D})+v_p(d_{k,\chi\omega^{-k}}^\alpha!)=\left[\frac{m}{{\sf C}_\alpha\alpha}\right]-v_p(\varphi(N))-1
$$ 
equations (24) and (26) yield
$$
{\rm tr} (\tilde{e}|_{{\cal M}_k(\Gamma,\chi\omega^{-k})^\alpha})
\not\equiv 
{\rm tr}\,(\tilde{e}|_{{\cal M}_{k'}(\Gamma,\chi\omega^{-k'})^\alpha})\pmod{p^{\left[\frac{m}{{\sf C}_\alpha\alpha}\right]-v_p(\varphi(N))}}.\leqno(27)
$$
On the other hand, since all $\lambda_{\ell(\lambda,\mu)}$ and $\lambda_{\ell(\mu')}$ are integral elements in some finite extension $E/{\Bbb Q}$ we deduce that $\tilde{e}$ is an ${\cal O}_E$-linear combination of terms of the form $\pi^{Np}_{Fp} T_{\ell_1}\cdot\ldots\cdot T_{\ell_s}e_\alpha^{\left[\frac{m}{{\sf C}_\alpha\alpha}\right]}$, where $\ell_i=\ell(\lambda,\mu)$ for some $\mu\in\Lambda_{k,\chi,Np}^\alpha-\{\lambda\}$ or $\ell_i=\ell(\mu')$ for some $\mu'\in\Lambda_{k',1}^\alpha$. The third congruence of the Theorem in section 3.5 therefore yields (note that $m\ge{\sf B}_\alpha\ge {\sf C}_\alpha\alpha+1$)
$$
{\rm tr} (\tilde{e}|_{{\cal M}_k(\Gamma,\chi\omega^{-k})^\alpha})
\equiv 
{\rm tr}\,(\tilde{e}|_{{\cal M}_{k'}(\Gamma,\chi\omega^{-k'})^\alpha})\pmod{p^{\left[\frac{m}{{\sf C}_\alpha\alpha}\right]-v_p(\varphi(N))}},
$$ 
which contradicts (27). Thus, our assumption is wrong and the Theorem therefore is proven in case $\alpha>0$. 

Finally, we look at the case $\alpha=0$. Here, the definition of ${\sf D}$ is such that
$$
v_p([{\rm GL}_2({\Bbb Z}):\Gamma_1(Np)])+d_{k,\chi\omega^{-k}}^\alpha(A_k+{\sf D})+v_p(d_{k,\chi\omega^{-k}}^\alpha!)= m-v_p(\varphi(N))-1,
$$
hence, equations (24) and (26) yield
$$
{\rm tr} (\tilde{e}|_{{\cal M}_k(\Gamma,\chi\omega^{-k})^\alpha})
\not\equiv 
{\rm tr}\,(\tilde{e}|_{{\cal M}_{k'}(\Gamma,\chi\omega^{-k'})^\alpha})\pmod{p^{m-v_p(\varphi(N))}}.\leqno(27^0)
$$
On the other hand, the second congruence appearing in Theorem in section 3.5 yields 
$$
{\rm tr} (\tilde{e}|_{{\cal M}_k(\Gamma,\chi\omega^{-k})^\alpha})
\equiv 
{\rm tr}\,(\tilde{e}|_{{\cal M}_{k'}(\Gamma,\chi\omega^{-k'})^\alpha})\pmod{p^{m-v_p(\varphi(N))}},
$$ 
which again, contradicts equation ($27^0$). Hence our assumption is wrong and we obtain the claim in case $\alpha=0$. Thus, the proof of part 1.) is complete.

2.) Since ${\rm dim}\,{\cal M}_k(\Gamma,\chi\omega^{-k})^\alpha=1$ we know in particular that ${\cal M}_k(\Gamma,\chi\omega^{-k})^\alpha={\cal M}_{k,\chi\omega^{-k}}^\alpha(\lambda)$. We therefore define the simpler (im comparison with the element defined in part 1.)) element
$$
e=\prod_{\mu'\in\Lambda_{k',1}^\alpha}(T_{\ell(\mu')}-\mu'_{\ell(\mu')})\in{\cal H}
$$
and 
$$
\tilde{e}=\pi^{Np}_{Fp}\,e\,e_\alpha^{[\frac{m}{{\sf C}_\alpha\alpha}]}.
$$
As in the proof of part 1.) we find

- $\tilde{e}|_{{\cal M}_{k'}(\Gamma,\chi\omega^{-k'})^\alpha}=0$, hence, ${\rm tr}\,\tilde{e}|_{{\cal M}_{k'}(\Gamma,\chi\omega^{-k'})^\alpha}=0$

- ${\rm tr}\, \tilde{e}|_{{\cal M}_{k,\chi\omega^{-k}}^\alpha(\lambda)}=\zeta\epsilon\,[\Gamma_1(Fp):\Gamma_1(Np)])$ (note that $d={\rm dim}\,({\cal M}_{k,\chi\omega^{-k}}^\alpha(\lambda)^{\Gamma_1(Fp)})=1$, because $\lambda$ has conductor $Fp$), where $v_p(\epsilon)=0$ and (simpler now) $\zeta=\prod_{\mu'\in\Lambda_{k',1}^\alpha}(\lambda_{\ell(\mu')}-\mu'_{\ell(\mu')})$. Hence, as in part 1.) we obtain
$$
v_p({\rm tr}\, \tilde{e}|_{{\cal M}_k(\Gamma,\chi\omega^{-k})^\alpha})\le v_p([{\rm GL}_2({\Bbb Z}):\Gamma_1(Np)])+d_{k,\chi\omega^{-k}}^\alpha{\sf D}.
$$
Using the definition of ${\sf D}$ these facts imply that
$$
{\rm tr} (\tilde{e}|_{{\cal M}_k(\Gamma,\chi\omega^{-k})^\alpha})
\not\equiv 
{\rm tr}\,(\tilde{e}|_{{\cal M}_{k'}(\Gamma,\chi\omega^{-k'})^\alpha})\pmod{p^{\left[\frac{m}{{\sf C}_\alpha\alpha}\right]-v_p(\varphi(N))}},
$$
in case $\alpha>0$ and
$$
{\rm tr} (\tilde{e}|_{{\cal M}_k(\Gamma,\chi\omega^{-k})^\alpha})
\not\equiv 
{\rm tr}\,(\tilde{e}|_{{\cal M}_{k'}(\Gamma,\chi\omega^{-k'})^\alpha})\pmod{p^{m-v_p(\varphi(N))}},
$$
in case $\alpha=0$. Again, this contradicts the Theorem in section 3.5 and the proof of the Theorem therefore is complete.

We set
$$
{\rm b}(\alpha,k)=-\frac{v_p(\varphi(N)[{\rm GL}_2({\Bbb Z}):\Gamma_1(Np)]d_{k,\chi\omega^{-k}}^\alpha!)+2}{d_{k,\chi\omega^{-k}}^\alpha}
$$
if $\alpha>0$ and
$$
{\rm b}(\alpha,k)=-\frac{v_p(\varphi(N)[{\rm GL}_2({\Bbb Z}):\Gamma_1(Np)]d_{k,\chi\omega^{-k}}^\alpha!)+1}{d_{k,\chi\omega^{-k}}^\alpha}
$$
if $\alpha=0$ and
$$
{\sf a}(\alpha,k)=\left\{
\begin{array}{ccc}
\frac{1}{d_{k,\chi\omega^{-k}}^\alpha{\sf C}_\alpha\alpha}&\mbox{if}&\alpha>0\\
\frac{1}{d_{k,\chi\omega^{-k}}^\alpha}&\mbox{if}&\alpha=0.\\
\end{array}
\right.
$$
We then find ${\sf a}(\alpha,k)m+{\rm b}(\alpha,k)-A_k\le {\sf D}$ for all $\alpha\ge 0$ and even ${\sf a}(\alpha,k)m+{\rm b}(\alpha,k)\le {\sf D}$ for all $\alpha\ge 0$ if in addition ${\rm dim}\,{\cal M}_k(\Gamma,\chi\omega^{-k})^\alpha=1$ holds. In case $\alpha=0$ these inequalities are obvious, in case $\alpha>0$ they follow from $[\frac{m}{d_{k,\chi\omega^{-k}}^\alpha{\sf C}_\alpha\alpha}]\ge\frac{m}{d_{k,\chi\omega^{-k}}^\alpha{\sf C}_\alpha\alpha}-1$. Thus, the element $\lambda'$ appearing in the above Theorem in particular satisfies the congruences
$$
\lambda'_\ell\equiv\lambda_\ell\pmod{p^{{\sf a}(\alpha,k)m+{\sf b}(\alpha,k)-A_k}}
$$
for all $\ell\not|Np$ if $\alpha\ge 0$ and even
$$
\lambda'_\ell\equiv\lambda_\ell\pmod{p^{{\sf a}(\alpha,k)m+{\sf b}(\alpha,k)}}
$$
for all $\ell\not|Np$ if $\alpha=0$. Thus, if we set
$$
{\sf a}(\alpha)={\rm min}\,\{{\sf a}(\alpha,k),\,k\ge ({\sf C}_\alpha\alpha+1)^2+1\}
$$
and
$$
{\sf b}(\alpha)={\rm min}\,\{{\sf b}(\alpha,k),\,k\ge({\sf C}_\alpha\alpha+1)^2+1\}
$$
we obtain the following

{\bf Theorem'. }{\it Let the assumptions and notations be as in the preceeding Theorem. Then, for any $\lambda=(\lambda_\ell)_{\ell\not|Np}\in\Lambda_{k,¸\chi}^\alpha$ there is a $\lambda'=(\lambda_\ell')_{\ell\not|Np}\in\Lambda_{k',\chi}^\alpha$ such that 

$\bullet$ ${\cal M}_{k',\chi\omega^{-k'}}^\alpha(\lambda')^{\Gamma_1(Fp)}\not=(0)$, i.e. the conductor of $\lambda'$ is a divisor of $Fp$

$\bullet$ $$
\lambda_\ell'\equiv\lambda_\ell\pmod{p^{{\sf a}m+{\sf b}-A_k}}
$$ 
for all $\ell\not|Np$. 

If in addition ${\rm dim}\,{\cal M}_k(\Gamma,\chi\omega^{-k})^\alpha=1$ holds we even obtain the congruences
$$
\lambda_\ell'\equiv\lambda_\ell\pmod{p^{{\sf a}m+{\sf b}}}.
$$
Here, ${\sf a}={\sf a}(\alpha$ and ${\sf b}={\sf b}(\alpha)$ only depend on $\alpha$ and ${\sf a}(\alpha)$ satisfies the inequalities 
$$
{\sf a}(\alpha)\le \frac{1}{{\rm dim}\,{\cal M}_k(\Gamma,\chi\omega^{-k})^\alpha}
$$ 
for all $k\ge ({\sf C}_\alpha\alpha+1)^2+1$. Moreover, in case $\alpha=0$ it even holds that
$$
{\sf a}(\alpha)={\rm min}\, \{\frac{1}{{\rm dim}\,{\cal M}_k(\Gamma,\chi\omega^{-k})^\alpha},\,k\ge ({\sf C}_\alpha\alpha+1)^2+1\}.
$$

}

{\it Proof. } Everything has been proven above except for the claimed inequality satisfied by ${\sf a}$, which follows immediately from ${\sf C}_\alpha\alpha\ge 1$ (cf. equation (45) in section 3.5) and the definition of ${\sf a}(\alpha,k)$. Thus, Theorem' is proven.

{\bf Remark. } Since ${\sf a}(\alpha,k)$ and ${\sf b}(\alpha,k)$ only depend on $k$ modulo $p^{{\sf B}_\alpha}$ we find that
$$
{\rm min}\,\{{\sf a}(\alpha,k),\,k\ge ({\sf C}_\alpha\alpha+1)^2+1\}={\rm min}\,\{{\sf a}(\alpha,k),\,({\sf C}_\alpha\alpha+1)^2+1\le k\le ({\sf C}_\alpha\alpha+1)^2+1+{\sf B}_\alpha\},
$$
i.e. ${\sf a}(\alpha)$ is the minimum of a finite set. In particular, since all ${\sf a}(\alpha,k)$ are strictly positive we see that
${\sf a}(\alpha)>0$.

{\bf Remark. } We set ${\sf a}={\sf a}(\alpha)$ and ${\sf b}={\sf b}(\alpha)$ and in what follows we shall always use the (slightly weaker) congruence  "$\pmod{p^{{\sf a}m+{\sf b}}}$" of Theorem' instead of the congruence "$\pmod{p^{\sf D}}$" of the Theorem.

{\bf 4.3. $p$-adic families of modular forms. } To begin with we note three more corollaries to the above Theorem'. To this end, we let $\lambda=(\lambda_\ell)_{\ell\not|Np}\in\Lambda_{k,\chi,Np}^\alpha$ and $\lambda'=(\lambda_\ell')_{\ell\not|Np}\in\Lambda_{k',\chi,Np}^\alpha$ be as in Theorem', i.e. $\lambda_\ell\equiv \lambda_\ell'\pmod{p^{{\sf a}m+{\sf b}-A_k}}$ for all $\ell\not|Np$ and $F'|F$, where we denote by $F$ resp. $F'$ the prime-to- $p$ part of the conductor of $\lambda$ resp. of $\lambda'$. Let $\pi=\pi_\lambda$ resp. $\pi'=\pi_{\lambda'}$ be the representation of ${\rm GL}_2({\Bbb A}_f)$ corresponding to $\lambda$ resp. $\lambda'$, i.e. ${\cal M}_{k,\chi\omega^{-k}}^\alpha(\lambda)=V_\pi^{K,\alpha}$ and ${\cal M}_{k',\chi\omega^{-k'}}^\alpha(\lambda')=V_{\pi'}^{K,\alpha}$ (cf. ($15^\alpha$)).

{\bf Corollary C1.} {\it We fix $\alpha\in{\Bbb Q}_{\ge 0}$ and assume that $k,k'$ satisfy $k,k'\ge ({{\sf C}_\alpha}\alpha+1)^2+2$ and $k\equiv k'\pmod{p^m}$ with $m\ge {\sf B}_\alpha$. Then, with the above notations we have
$$
{\rm dim}\,V_{\pi,\ell}^{K_{\ell}}\le {\rm dim}\,V_{\pi',\ell}^{K_\ell}\leqno(28a)
$$
for all primes $\ell\not=p$ and
$$
{\rm dim}\,V_{\pi,p}^{K_{p},\alpha}\le {\rm dim}\,V_{\pi',p}^{K_p,\alpha}.\leqno(28b)
$$ 
In particular,
$$
{\rm dim}\,V_{\pi'}^{K,\alpha}\ge{\rm dim}\,V_{\pi}^{K,\alpha}\leqno(28c)
$$
(cf. equation (16)).
}

{\it Proof. } Let $N=\prod_i \ell_i^{n_i}$, $F=\prod_i\ell_i^{e_i}$, $F'=\prod_i\ell_i^{e_i'}$ be the prime decompositions of $N$, $F$ and $F'$. We distinguish cases. 

In case $\ell$ does not divide $Np$ the claim is obvious since both spaces in question are $1$-dimensional. 

In case $\ell$ divides $N$ we will apply the following result, which is part of Casselman's representation theoretic reformulation of  Atkin-Lehner theory: if $(\rho,V_\rho)$ is an admissible ${\rm GL}_2({\Bbb Q}_\ell)$-module with conductor $\ell^e$, then ${\rm dim}\,V_\rho^{K_1(\ell^a)}= a-e+1$ for any $a\ge e$ (cf. [Ca], Corollary to the Proof, p. 306). 
Since $F'|F|N$ we obtain $e_i'\le e_i\le n_i$. We then compute using Casselman's formula (recall that $K=K_1(Np)$)
$$
{\rm dim}\,V_{\pi,\ell}^{K_\ell}={\rm dim}\, V_{\pi,\ell}^{K_{1,\ell}(\ell^{n_i})}= n_i-e_i+1.\leqno(29)
$$
In the same way we find
$$
{\rm dim}\,V_{\pi',\ell}^{K_\ell}= n_i-e_i'+1\leqno(30)
$$
and since $e_i'\le e_i$ equations (29) and (30) yield the assertion in case $\ell|N$.

Finally we look at the case $\ell=p$. Since $V_\pi$ and $V_{\pi'}$ both appear in the space of $\Gamma_1(Np)$-invariant forms, the conductor
of $V_{\pi,p}$ and of $V_{\pi',p}$ equals $1$ or $p$, or, equivalently, ${\rm dim}\,V_{\pi,p}^{K_p}$ and ${\rm dim}\, V_{\pi',p}^{K_p}$ both equal $1$ or $2$. Thus, the claim follows if we can show that ${\rm dim}\,V_{\pi,p}^{K_p,\alpha}\le 1$ (note that $V_{\pi',p}^{K_p,\alpha}\not=0$, hence, ${\rm dim}\,V_{\pi',p}^{K_p,\alpha}\ge 1$ because $\lambda'\in\Lambda_{k',\chi,Np}^\alpha$, i.e. $0\not={\cal M}_{k',\chi\omega^{-k'}}^\alpha(\lambda')\subset V_{\pi'}^{K,\alpha}$). Clearly, ${\rm dim}\,V_{\pi,p}^{K_p,\alpha}\le 1$ holds trivially if ${\rm cond}\,V_{\pi,p}=p$, hence, we may assume ${\rm cond}\,V_{\pi,p}=1$ or equivalently, ${\rm dim}\,V_{\pi,p}^{K_p}=2$. We choose a basis $\{f_1,f_2\}$ of $V_{\pi,p}^{K_p}$ such that $T_p$ has triangular form $\Mat{\lambda_1}{*}{}{\lambda_2}$ with repect to this basis. We set $\alpha_i=v_p(\lambda_i)$, $i=1,2$. If $\alpha_1\not=\alpha_2$ then again we trivially have ${\rm dim }\,V_{\pi,p}^{K_p,\alpha}\le 1$. Thus we may assume $\alpha_1=\alpha_2$. Since $\alpha_1+\alpha_2=k-1$ (cf. [M-G], p. 796) we obtain 
$$
\alpha=\alpha_1=\alpha_2=(k-1)/2.
$$
We distinguish cases and first assume $\alpha>0$. By our assumption on $k$ we know that
$$
k\ge ({\sf C}_{\alpha}\alpha+1)^2\ge ({\sf C}_{\alpha}\alpha)^2.
$$
Recalling that ${\sf C}_{\alpha}\ge N(\alpha)+\frac{1}{\alpha}$ (cf. equation (45) in section 3) we further obtain $k\ge (N(\alpha)\alpha+1)^2$ and since $N(\alpha)\ge 1$ (cf. equation (32) in section 3) this implies 
$$
k\ge N(\alpha)^2\alpha^2+2N(\alpha)\alpha+1\ge \alpha^2+1=(k-1)^2/4+1\ge (k-1)^2/4-1/4.
$$ 
Altogether, we obtain $0\ge k^2-6k$, which yields that $0\le k\le 6$. Hence, $k$ is one of the finitely many weights $2,3,4,5,6$. On the other hand, since ${\sf C}_\alpha\alpha\ge N(\alpha)\alpha+1>1$ (cf. the definition of ${\sf C}_\alpha$ in section 3.5 and equation (32) in section 3.4) we find $k\ge({\sf C}_\alpha\alpha+1)^2+2>6$, which is a contradiction. Thus, our assumption $\alpha_1=\alpha_2$ is false and we deduce that ${\rm dim}\,V_{\pi,p}^{K_p,\alpha}\le 1$ (of course, since we assume $\lambda\in\Lambda_{k,\chi,Np}^\alpha$ the space ${\rm dim}\,V_{\pi,p}^{K_p,\alpha}$ is not trivial, i.e. we have ${\rm dim}\,V_{\pi,p}^{K_p,\alpha}= 1$). Finally, in case $\alpha=0$, i.e. ${\sf C}_\alpha=0$, the assumption $\alpha_1=\alpha_2=\alpha$ leads to $0=v_p(\alpha)=(k-1)/2$, i.e. $k=1$, which is not possible. This completes the proof of the Corollary.

We introduce two more pieces of notation. Let $k,k'$ be as in Theorem' in section 4.2. For  $\lambda=(\lambda_\ell)_{\ell\not|Np}\in\Lambda_{k,\chi,Np}^\alpha$ we denote by $\lambda(k')\in\Lambda_{k',\chi,Np}^\alpha$ any element satisfying the conditions
$$
{\cal M}_{k',\chi\omega^{-k}}^\alpha(\lambda)^{\Gamma_1(Fp)}\not=(0)\leqno(31)
$$ 
and 
$$
\lambda'_\ell\equiv\lambda_\ell\pmod{p^{{\sf a}m+{\sf b}-A_k}}\leqno(32)
$$ 
for all $\ell\not|Np$, where $k\equiv k'\pmod{p^m}$. We note that in case ${\rm dim}\,{\cal M}_k(\Gamma,\chi\omega^{-k})^\alpha=1$ Theorem' in section 4.2 shows that stronger congruences than those of equation (32) hold, but we will not need this (but see also the Remark following Corollary C2). Theorem' in section 4.2 guaranties the existence of an element $\lambda(k')$ satisfying (31) and (32) if 
$$
k,k'\ge({{\sf C}_\alpha}\alpha+1)^2+2\leqno(33)
$$ 
and 
$$
k\equiv k'\pmod{p^m}\quad\mbox{with}\quad m>{\sf B}_\alpha.\leqno(34)
$$ 
Moreover, we denote by $\pi_\lambda$ the representation of ${\rm GL}_2({\Bbb A}_f)$ attached to $\lambda$ and by $\pi_{\lambda(k')}$ the ${\rm GL}_2({\Bbb A}_f)$ representation attached to $\lambda(k')$. Finally, if $\lambda,\mu\in\Lambda_{k,\chi,Np}^\alpha$ we write 
$$
\lambda\equiv \mu\pmod{p^{{\sf a}m+{\sf b}-A_k}},
$$ 
to denote that $\lambda_\ell\equiv \mu_\ell\pmod{p^{{\sf a}m+{\sf b}-A_k}}$ for all $\ell\not|Np$. Equation (18) then may be reformulated as follows
$$
\lambda,\mu\in\Lambda_{k,\chi,Np}^\alpha,\,\lambda\not=\mu\;\Rightarrow\;\lambda\not\equiv\mu\pmod{p^{A_k}}.\leqno(35)
$$

{\bf Corollary C2. }{\it Let $\alpha\in{\Bbb Q}_{\ge 0}$ and let $k\ge({{\sf C}_\alpha}\alpha+1)^2+2$ be any integer. Set 
$$
{\sf L}(\alpha,k)=[{\rm max}\,\{{{\sf B}_\alpha},(2A_k-{\sf b})/{\sf a}\}]+1\in{\Bbb N}.
$$ 
For all $k'\ge ({{\sf C}_\alpha}\alpha+1)^2+2$ satisfying $k'\equiv k\pmod{p^{m}}$ with $m\ge {\sf L}(\alpha,k)$ the following properties hold.

1.) ${\cal M}_{k'}(\Gamma,\chi\omega^{-k'})^\alpha=\bigoplus_{\lambda\in\Lambda_{k,\chi,Np}^\alpha} {\cal M}_{k',\chi\omega^{-k'}}^\alpha(\lambda(k'))\Big(=\bigoplus_{\lambda\in\Lambda_{k,\chi,Np}^\alpha} V_{\pi_{\lambda(k')}}^{K,\alpha}\Big)$.

2.) ${\rm dim}\,V_{\pi_\lambda}^{K,\alpha}={\rm dim}\,V_{\pi_{\lambda(k')}}^{K,\alpha}$ ($\lambda\in\Lambda_{k,\chi,Np}^\alpha$ arbitrary). In particular, ${\rm cond}\,\lambda={\rm cond}\,\lambda(k')$ (because we know that ${\rm cond}\,\lambda(k')|{\rm cond}\,\lambda$).

%, hence, ${\rm dim}\,V_{\pi_\lambda,\ell}^{K_\ell}={\rm dim}\,V_{\pi_{\lambda(k')},\ell}^{K_\ell}$ for all primes $\ell\not=p$ and ${\rm dim}\,V_{\pi_\lambda,p}^{K_p,\alpha}={\rm dim}\,V_{\pi_{\lambda(k')},p}^{K_p,\alpha}$. 

3.) For any $\lambda\in\Lambda_{k,\chi,Np}^\alpha$ there is precisely one $\lambda'\in\Lambda_{k',\chi,Np}^\alpha$ satisfying $$
\lambda'\equiv\lambda\pmod{p^{{\sf a}m+{\sf b}-A_k}}.
$$ 
Since $\lambda(k')$ also satisfies (32) we have $\lambda'=\lambda(k')$ and 2.) implies that $\lambda'$ also satisfies ${\rm cond}\,\lambda={\rm cond}\,\lambda(k')$.

4.) We denote by 
$$
\varphi:\,\Lambda_{k,\chi,Np}^\alpha\rightarrow\Lambda_{k',\chi,Np}^\alpha
$$
the map which sends $\lambda$ to the uniquely determined element $\lambda(k')$ satisfying $\lambda(k')\equiv\lambda\pmod{p^{{\sf a}m+{\sf b}-A_k}}$. Then, $\varphi$ is a bijection.

}

{\it Proof. } Let $\lambda,\mu\in\Lambda_{k,\chi,Np}^\alpha$, $\mu\not=\lambda$. We set $\lambda'=\lambda(k')$ and $\mu'=\mu(k')$ and we denote by $(\pi_{\lambda'},V_{\pi_{\lambda'}})$ resp. $(\pi_{\mu'},V_{\pi_{\mu'}})$ the representation of ${\rm GL}_2({\Bbb A}_f)$ corresponding to $\lambda'$ resp. $\mu'$. Since the numbers $k,k'$ satisfy (31) and (32) the existence of $\lambda'$ and $\mu'$ is guarantied by Theorem' in section 4.2. Equation (18) implies that 
$$
\lambda_{\ell(\mu,\lambda)}\not\equiv\mu_{\ell(\mu,\lambda)}\pmod{p^{A_k}}.
$$ 
Since $m\ge {\sf L}(\alpha,k)$, equation (32) in particular implies 
$$
\lambda_{\ell(\mu,\lambda)}\equiv\lambda'_{\ell(\mu,\lambda)}\pmod{p^{{\sf a}{\sf L}(\alpha,k)+{\sf b}-A_k}}\leqno(36)
$$ 
and 
$$
\mu_{\ell(\mu,\lambda)}\equiv\mu'_{\ell(\mu,\lambda)}\pmod{p^{{\sf a}{\sf L}(\alpha,k)+{\sf b}-A_k}}.\leqno(37)
$$ 
Since ${\sf a}{\sf L}(\alpha,k)+{\sf b}-A_k\ge A_k$, we altogether deduce that 
$$
\lambda'_{\ell(\mu,\lambda)}\not\equiv\mu'_{\ell(\mu,\lambda)}\pmod{p^{A_k}}.\leqno(38)
$$
Thus, we have proven that the $\lambda(k')$, $\lambda\in\Lambda_{k,\chi,Np}^\alpha$, are pairwise non-congruent modulo $p^{A_k}$. In particular, the representations $\pi_{\lambda(k')}$, $\lambda\in\Lambda_{k,\chi,Np}^\alpha$, are pairwise non-isomorphic, hence, the spaces $V_{\pi_{\lambda(k')}}$, $\lambda\in\Lambda_{k,\chi,Np}^\alpha$, are pairwise different. Since the $V_{\pi_{\lambda(k')}}^K$, $\lambda\in\Lambda_{k,\chi,Np}^\alpha$, are irreducible $K\backslash {\rm GL}_2({\Bbb A}_f)/K$-modules, we deduce that the sum $\bigoplus_{\lambda\in\Lambda_{k,\chi,Np}^\alpha} V_{\pi_{\lambda(k')}}^K$ is direct.
%
%[[Note that $K\backslash {\rm GL}_2({\Bbb A}_f)/K$ is the big hecke algebra, i.e. it is bigger than the Hecke algebra generated by the classical hecke operators $T_\ell$, $\ell$ prime. $V_{\pi_{\lambda(k')}}$ need not be irreducible under the Hecke algebra generated by the $T_\ell$, $\ell$ prime.]]
%
In particular, looking at slope $\alpha$ subspaces we obtain
$$
\bigoplus_{\lambda\in\Lambda_{k,\chi,Np}^\alpha} V_{\pi_{\lambda(k')}}^{K,\alpha}\subset {\cal M}_{k'}(\Gamma,\chi\omega^{-k'})^{\alpha}.\leqno(39)
$$

Since furthermore $k,k'$ are bigger than $({{\sf C}_\alpha}\alpha+1)^2+2$ and $k\equiv k'\pmod{p^{{\sf L}(\alpha,k)}}$ where ${\sf L}(\alpha,k)>{{\sf B}_\alpha}$, Theorem B implies ${\rm dim}\,{\cal M}_{k}(\Gamma,\chi\omega^{-k})^\alpha={\rm dim}\,{\cal M}_{k'}(\Gamma,\chi\omega^{-k'})^\alpha$. Also, ${\rm dim}\,V_{\pi_{\lambda(k')}}^K\ge{\rm dim}\,V_{\pi_\lambda}^K$ by (28c) and altogether we obtain the chain of 
inequalities

\begin{eqnarray*}
(40)\qquad{\rm dim}\,{\cal M}_{k'}(\Gamma,\chi\omega^{-k'})^\alpha&\stackrel{(39)}{\ge}& {\rm dim}\,\bigoplus_{\lambda\in\Lambda_{k,\chi,Np}^\alpha} V_{\pi_{\lambda(k')}}^{K,\alpha}\\
&=&\sum_{\lambda\in\Lambda_{k,\chi,Np}^\alpha} {\rm dim}\,V_{\pi_{\lambda(k')}}^{K,\alpha}\\
&\stackrel{(28c)}{\ge}& \sum_{\lambda\in\Lambda_{k,\chi,Np}^\alpha} {\rm dim}\,V_{\pi_{\lambda}}^{K,\alpha}\\
&=&{\rm dim}\,{\cal M}_{k}(\Gamma,\chi\omega^{-k})^\alpha\\
&\stackrel{\rm (Thm B)}{=}&{\rm dim}\,{\cal M}_{k'}(\Gamma,\chi\omega^{-k'})^\alpha.
\end{eqnarray*}

Thus, we have equality everywhere and, hence, we obtain
$$
{\rm dim}\,{\cal M}_{k'}(\Gamma,\chi\omega^{-k'})^\alpha= {\rm dim}\,\bigoplus_{\lambda\in\Lambda_{k,\chi,Np}^\alpha} V_{\pi_{\lambda(k')}}^{K,\alpha},
$$
which yields the first claim of the Corollary. Moreover, equation (40) further yields
$$
\sum_{\lambda\in\Lambda_{k,\chi,Np}^\alpha} {\rm dim}\,V_{\pi_{\lambda(k')}}^{K,\alpha}={\rm dim}\,{\cal M}_{k'}(\Gamma,\chi\omega^{-k'})^\alpha=\sum_{\lambda\in\Lambda_{k,\chi,Np}^\alpha} {\rm dim}\,V_{\pi_{\lambda}}^{K,\alpha}.
$$
In conjunction with equation (28c) this implies  ${\rm dim}\,V_{\pi_\lambda}^{K,\alpha}={\rm dim}\,V_{\pi_{\lambda(k')}}^{K,\alpha}$, which is the second claim of the Corollary. Furthermore, if there were two different elements $\lambda'$ and $\lambda''$ in $\Lambda_{k',\chi,Np}^\alpha$ satisfying equation (32) then one of them would equal $\lambda(k')$, i.e. it satisfies in addition (31). Without loss of generality we may assume $\lambda'=\lambda(k')$ and we would obtain 
$$
\bigoplus_{\lambda\in\Lambda_{k,\chi,Np}^\alpha} V_{\pi_{\lambda(k')}}^{K,\alpha}\oplus V_{\pi_{\lambda''}}^{K,\alpha}\subseteq {\cal M}_{k'}(\Gamma,\chi\omega^{-k'})^\alpha,
$$
which contradicts ${\cal M}_{k'}(\Gamma,\chi\omega^{-k'})^\alpha=\bigoplus_{\lambda\in\Lambda_{k,\chi,Np}^\alpha} V_{\pi_{\lambda(k')}}^{K,\alpha}$. Thus, there is only one element $\lambda(k')$ satisfying (32), which is the third claim. Finally, to prove the last claim we assume there are $\lambda,\mu\in\Lambda_{k,\chi,Np}^\alpha$ such that $\varphi(\lambda)=\varphi(\mu)$, i.e. $\lambda(k')=\mu(k')$. We then obtain $\lambda\equiv\lambda(k')=\mu(k')\equiv\mu\pmod{p^{{\sf a}m+{\sf b}-A_k}}$. Since $m\ge {\sf L}(\alpha,k)$ we find ${\sf a}m+{\sf b}-A_k\ge A_k$. Hence, equation (35) implies that $\lambda=\mu$, which is the injectivity of $\varphi$. To prove surjectivity, let $\lambda'\in\Lambda_{k',\chi,Np}^\alpha$. By 1.) there is $\lambda\in\Lambda_{k,\chi,Np}^\alpha$ such that ${\cal M}_{k',\chi\omega^{-k'}}^\alpha(\lambda')={\cal M}_{k',\chi\omega^{-k'}}^\alpha(\lambda(k'))$, which immediately implies that $\lambda'=\lambda(k')=\varphi(\lambda)$. Hence, $\varphi$ also is surjective and therefore a bijection, which is the fourth claim. Thus, the proof of the Corollary is complete.

{\bf Remark. } Let $\alpha\in{\Bbb Q}_{\ge 0}$, let $k\ge({{\sf C}_\alpha}\alpha+1)^2+2$ be any integer and assume that in addition ${\rm dim}\,{\cal M}_k(\Gamma,\chi\omega^{-k})^\alpha=1$. Hence, $\Lambda_{k,\chi,Np}^\alpha=\{\lambda\}$ and ${\cal M}_k(\Gamma,\chi\omega^{-k})^\alpha=V_{\pi_\lambda}^{K,\alpha}$. Let $k'\ge({{\sf C}_\alpha}\alpha+1)^2+2$ be any integer satisfying $k'\equiv k\pmod{p^{{\sf B}_\alpha}}$. Theorem' in section 4.2 shows that there is an element $\lambda(k')\in\Lambda_{k',\chi,Np}^\alpha$ such that $\lambda\equiv\lambda(k')\pmod{p^{{\sf a}m+{\sf b}}}$. On the other hand Theorem B shows that ${\rm dim}\,{\cal M}_k(\Gamma,\chi\omega^{-k})={\rm dim}\,{\cal M}_{k'}(\Gamma,\chi\omega^{-k'})$, hence, $\Lambda_{k',\chi,Np}^\alpha=\{\lambda(k')\}$ and ${\cal M}_{k'}(\Gamma,\chi\omega^{-k'})=V_{\pi_{\lambda(k')}}^{K,\alpha}$. In particular, 
$$
\varphi:\,\Lambda_{k,\chi,Np}^\alpha\rightarrow\Lambda_{k',\chi,Np}^\alpha,\quad\lambda\mapsto\lambda(k')
$$
is a bijection,
$$
{\cal M}_{k'}(\Gamma,\chi\omega^{-k})^\alpha=V_{\pi_{\lambda(k')}}^{K,\alpha}
$$
and 
$$
{\rm dim}\,V_{\pi_\lambda}^{K,\alpha}={\rm dim}\,V_{\pi_{\lambda(k')}}^{K,\alpha}(=1).
$$
Thus, we obtain in a trivial way all the statements of Corollary C2. On the other hand, we note that in case ${\rm dim}\,{\cal M}_k(\Gamma,\chi\omega^{-k})=1$ these statements hold for all $k'$ which satisfy the congruence $k\equiv k'\pmod{p^{{\sf B}_\alpha}}$, i.e. we replaced ${\sf L}(\alpha,k_0)$ by ${\sf B}_\alpha$, which no longer depends on $k_0$.

We note one more consequence of Theorem':

%[[Im folgenden Corollar C3 taucht zum ersten mal $k_0$ auf]]

{\bf Corollary C3. } {\it Let $\alpha\in{\Bbb Q}_{\ge 0}$ and let $k_0\ge({\sf C}_\alpha\alpha+1)^2+2$. Then for any integer $k$ satisfying $k\equiv k_0\pmod{p^{{\sf L}(\alpha,k_0)}}$ and $k\ge({\sf C}_\alpha\alpha+1)^2+2$ the following holds: let $\lambda,\mu\in\Lambda_{k,\chi,Np}^\alpha$; then 
$$
\lambda\equiv \mu\pmod{p^{A_{k_0}}} \Leftrightarrow \lambda=\mu.
$$ 
In different words, as long as we consider integers satisfying $k\equiv k_0\pmod{p^{{\sf L}(\alpha,k_0)}}$ (and $k\ge({\sf C}_\alpha\alpha+1)^2+2$) we may choose $A_k=A_{k_0}$ independently of $k$ (cf. the Notation in section 4.2 for the definition of $A_k$).
}

%[[Does this imply that we can choose ${\sf L}(\alpha,k)$ such that it only depends on $k$ moduko $p^{{\sf L}(\alpha,k)}$ ???]]

{\it Proof. } Only the implication "$\Rightarrow$" needs proof. To this end we let $\lambda,\mu\in\Lambda_{k,\chi,Np}^\alpha$ and assume $\lambda\not=\mu$. We write $k\equiv k_0\pmod{p^m}$; by Corollary C2 4.) there are elements $\lambda_0,\mu_0\in\Lambda_{k_0,\chi,Np}^\alpha$ such that $\lambda=\lambda_0(k)$ and $\mu=\mu_0(k)$, i.e. $\lambda_0$ and $\mu_0$ satisfy the congruences $\lambda\equiv \lambda_0\pmod{p^{{\sf a}m+{\sf b}-A_{k_0}}}$ and $\mu\equiv \mu_0\pmod{p^{{\sf a}m+{\sf b}-A_{k_0}}}$. Since $m\ge {\sf L}(\alpha,k_0)$ we know that ${\sf a}m+{\sf b}-A_{k_0}\ge A_{k_0}$. Furthermore we know that $\mu_0\not=\lambda_0$ because $\mu\not=\lambda$, hence, $\mu\not\equiv \lambda\pmod{p^{A_{k_0}}}$ by  equation (35) and we obtain
$$
\lambda\equiv\lambda_0\not\equiv\mu_0\equiv \mu\pmod{p^{A_{k_0}}}.
$$
Thus, $\lambda\not\equiv \mu\pmod{p^{A_{k_0}}}$, which completes the proof of the Corollary.

We are now ready to prove that any $\lambda\in\Lambda_{k_0,\chi,Np}^\alpha$, $k_0\ge ({{\sf C}_\alpha}\alpha+1)^2+2$ can be placed into a $p$-adic family of systems of eigenvalues.

%[[Im folgenden theorem benutzen wir erstmals $k_0$ ]]

{\bf Theorem D1. } {\it Let $\alpha\in{\Bbb Q}_{\ge 0}$ and let $k_0\ge ({{\sf C}_\alpha}\alpha+1)^2+2$. Then, for any  $\lambda\in\Lambda_{k_0,\chi,Np}^\alpha$ there is a family $(\lambda(k))_k$, where $k$ runs over all integers satisfying $k\ge({{\sf C}_\alpha}\alpha+1)^2+2$ and $k\equiv k_0\pmod{p^{{\sf L}(\alpha,k_0)}}$, which satisfies the following properties

$\bullet$ $\lambda(k)\in\Lambda_{k,\chi,Np}^\alpha$

$\bullet$ $\lambda(k_0)=\lambda$

$\bullet$ if $k\equiv k'\pmod{p^m}$ then $\lambda(k)\equiv \lambda(k')\pmod{p^{{\sf a}m+{\sf b}-A_{k_0}}}$.

Here, ${\sf a}={\sf a}(\alpha)$ and ${\sf b}={\sf b}(\alpha)$ are the integers defined in section 4.2, hence, 
$$
{\sf a}\le \frac{1}{{\rm dim}\,{\cal M}_{k}(\Gamma,\chi\omega^{-k})^\alpha}
$$ 
for all $k\ge({\sf C}_\alpha\alpha+1)^2+1$ and even
$$
{\sf a}={\rm min}\,\{\frac{1}{{\rm dim}\,{\cal M}_{k}(\Gamma,\chi\omega^{-k})^\alpha},\;k\ge({\sf C}_\alpha\alpha+1)^2+1\}
$$ 
if $\alpha=0$. 

Moreover, if $(\mu(k))_k$, $\mu(k)\in\Lambda_{k,\chi,Np}^\alpha$ is another family satisfying the above three conditions, then $\mu(k)=\lambda(k)$ for all $k$. 

}

%[[Beachte: aus dem Beweis des Theorems zu dem "Notation" gehort, folgt sogar ${\sf a}\ge {\rm dim}\,{\cal M}_{k_0,\chi}^\alpha(\lambda_0)$ !!!! ]]

{\it Proof. } Let $\lambda\in\Lambda_{k_0,\chi,Np}^\alpha$, $k_0\ge ({{\sf C}_\alpha}\alpha+1)^2+2$. We define the family $(\lambda(k))_k$ as follows. We have seen (cf. Corollary C2) that for any weight $k$ satisfying $k\ge ({{\sf C}_\alpha}\alpha+1)^2+2$ and $k\equiv k_0\pmod{p^{{\sf L}(\alpha,k_0)}}$ there is a (unique) element $\lambda(k)\in\Lambda_{k,\chi,Np}^\alpha$ such that $\lambda(k)\equiv\lambda\pmod{p^{{\sf a}m+{\sf b}-A_{k_0}}}$. In this way we obtain a family $(\lambda(k))_k$, which clearly satisfies the first two conditions. We will show that it also satisfies the third condition. 
To this end, let $k,k'$ be integers satisfying $k,k'\ge ({{\sf C}_\alpha}\alpha+1)^2+2$ and $k,k'\equiv k_0\pmod{p^{{\sf L}(\alpha,k_0)}}$. We write $k\equiv k_0\pmod{p^x}$, $k'\equiv k_0\pmod{p^{y}}$ and $k\equiv k'\pmod{p^{m}}$. We obtain

$$
\lambda(k)\equiv\lambda(k_0)\pmod{p^{{\sf a}x+{\sf b}-A_{k_0}}}.
$$
and 
$$
\lambda(k')\equiv\lambda(k_0)\pmod{p^{{\sf a}y+{\sf b}-A_{k_0}}}.
$$
Clearly, $x,y\ge {\sf L}(\alpha,k_0)$ and since ${\sf a}{\sf L}(\alpha,k_0)+{\sf b}-A_{k_0}\ge A_{k_0}$ we obtain
$$
\lambda(k)\equiv\lambda(k_0)\equiv\lambda(k')\pmod{p^{A_{k_0}}}.\leqno(42)
$$
On the other hand, since $k,k'\equiv k_0\pmod{p^{{\sf L}(\alpha,k_0)}}$ we also find that $m\ge {\sf L}(\alpha,k_0)$ and since ${\sf L}(\alpha,k_0)\ge {\sf B}_\alpha$ Theorem' in section 4.2 shows that attached to $\lambda(k)\in\Lambda_{k,\chi,Np}^\alpha$ there is an element $\mu'\in\Lambda_{k',\chi,Np}^\alpha$ satisfying
$$
\mu'\equiv\lambda(k)\pmod{p^{{\sf a}m+{\sf b}-A_{k_0}}}.
$$
Note that we may choose $A_k=A_{k_0}$ by Corollary C3; this is why "$A_{k_0}$" may appear in the modulus of the above congruence in place of 
"$A_k$". Since $m\ge {\sf L}(\alpha,k_0)$ we obtain 
$$
\mu'\equiv\lambda(k)\pmod{p^{A_{k_0}}}.\leqno(43)
$$
Equations (42) and (43) yield $\lambda(k')\equiv \mu'\pmod{p^{A_{k_0}}}$, which by Corollary C3 (and the definition of $A_{k_0}$; cf. the Notation in section 4.2) implies that $\lambda(k')=\mu'$. By the choice of $\mu'$ we then deduce that $\lambda(k')\equiv \lambda(k)\pmod{p^{{\sf a}m+{\sf b}-A_{k_0}}}$. 

%[[Klar: $\lambda(k_0)=\lambda]]

The uniqueness of the family $(\lambda(k))_k$ is a direct consequence of Corollary C2 3.). Thus, the proof of the theorem is complete.

If we know in addition that ${\rm dim}\,{\cal M}_k(\Gamma,\chi\omega^{-k})^\alpha=1$ we can derive a stronger result than Theorem D. In fact, any $\lambda\in\Lambda_{k_0,\chi,Np}^\alpha$ can be placed into a $p$-adic family, whose base and modulus no longer depend on the initial weight $k_0$:

{\bf Theorem D2. } {\it Let $\alpha\in{\Bbb Q}_{\ge 0}$, let $k_0\ge ({{\sf C}_\alpha}\alpha+1)^2+2$ and assume in addition that ${\rm dim}\,{\cal M}_{k_0}(\Gamma,\chi\omega^{-k_0})^\alpha=1$. Then, for any $\lambda\in\Lambda_{k_0,\chi,Np}^\alpha$ there is a family $(\lambda(k))_k$, where $k$ runs over all integers satisfying $k\ge({{\sf C}_\alpha}\alpha+1)^2+2$ as well as $k\equiv k_0\pmod{p^{{\sf B}_\alpha}}$, which satisfies the properties

$\bullet$ $\lambda(k)\in\Lambda_{k,\chi,Np}^\alpha$ 

$\bullet$ $\lambda(k_0)=\lambda$

$\bullet$ if $k\equiv k'\pmod{p^m}$ then $\lambda(k)\equiv \lambda(k')\pmod{p^{{\sf a}m+{\sf b}}}$.

Here, as in Theorem D1, ${\sf a}={\sf a}(\alpha)$ and ${\sf b}={\sf b}(\alpha)$ are the integers defined in section 4.2.

}

{\it Proof. } Theorem B implies that ${\rm dim}\,{\cal M}_k(\Gamma,\chi\omega^{-k})^\alpha={\rm dim}\,{\cal M}_{k_0}(\Gamma,\chi\omega^{-k_0})^\alpha=1$ for all $k$ satisfying the conditions in Theorem D2, hence, $\Lambda_{k,\chi,Np}^\alpha=\{\lambda_k\}$ for some element $\lambda_k$ for all those weights $k$. We will show that the family $(\lambda_k)_k$, where $k$ runs over all integers satisfying $k\ge({{\sf C}_\alpha}\alpha+1)^2+2$ and $k\equiv k_0\pmod{p^{{\sf B}_\alpha}}$, satisfies the claimed  congruences. In fact, let $k,k'$ be integers satisfying $k,k'\ge ({{\sf C}_\alpha}\alpha+1)^2+2$ and $k,k'\equiv k_0\pmod{p^{{\sf B}_\alpha}}$. We write $k\equiv k'\pmod{p^{m}}$. Using Theorem' in section 4.2 there is an element $\lambda'\in\lambda_{k',\chi,Np}^\alpha$ such that $\lambda'\equiv \lambda_k\pmod{p^{{\sf a}m+{\sf b}}}$. Since $\Lambda_{k',\chi,Np}^\alpha=\{\lambda_{k'}\}$ we deduce that $\lambda_{k'}=\lambda'$, hence, we obtain the requested congruence $\lambda_{k'}\equiv \lambda_k\pmod{p^{{\sf a}m+{\sf b}}}$. This finishes the proof of Theorem D2.

%{\bf Remark. } We note that the base of the family $(\lambda(k))_k$ in Theorem D' as well as the modulus of the congruences satisfied by this family no longer depend on the initial weight $k_0$.

{\bf 4.4. $p$-adic interpolation of the non-cuspidal part of the spectrum. } Let $\lambda=(\lambda_\ell)_{\ell\not|Np}\in\Lambda_{k,\chi,Np}^\alpha$ be a system of Hecke eigenvalues. Equation (15) shows that ${\cal M}_{k,\chi\omega^{-k}}(\lambda)=V_\pi^{K}$ for some irreducible automorphic representation $\pi$ and we say that $\lambda$ is cuspidal if $\pi$ is cuspidal and we call $\lambda$ non-cuspidal if $\pi$ is not cuspidal. Thus, if $\lambda$ is cuspidal resp. non-cuspidal the space ${\cal M}_{k,\chi\omega^{-k}}(\lambda)$ consists entirely of cusp forms resp. of non-cusp forms.

{\bf Proposition E. }{\it Fix $\alpha\in{\Bbb Q}_{\ge 0}$ and let $k_0\in{\Bbb N}$, $k_0>2$. Let furthermore $\lambda\in\Lambda_{k_0,\chi,Np}^\alpha$ be non-cuspidal. Then the following holds.

1.) The slope $\alpha$ of $\lambda$ equals $0$ or $k_0$. In particular, if $(\lambda_k)_k$ is any family of non-cuspidal elements $\lambda_k\in\Lambda_{k,\chi,Np}^\alpha$ of constant slope $\alpha$ (i.e. any $\lambda_k$ has slope $\alpha$) then $\alpha=0$.

2.) If the slope $\alpha$ of $\lambda$ equals $0$, then there is a family $(\lambda_k)_k$, $k\in k_0+(p-1){\Bbb Z}$, where $\lambda_k=(\lambda_{k,\ell})_{\ell\not|Np}\in \Lambda_{k,\chi,Np}^\alpha$, such that 

$\bullet$ $\lambda_k$ is non-cuspidal for all $k\in k_0+(p-1){\Bbb Z}$

$\bullet$ $k\equiv k'\pmod{p^{m}}$ implies $\lambda_{k,\ell}\equiv \lambda_{k',\ell}\pmod{p^{m+1}}$ for all $\ell\not|Np$.

$\bullet$ $\lambda_{k_0}=\lambda$

%3.) If $(\lambda(k))_k$ ($\lambda(k)\in\Lambda_{k,\chi,Np}^\alpha$, $k\equiv k_0\pmod{(p-1){\Bbb Z}}$) is any non-cuspidal family of constant slope $\alpha$ (i.e. any $\lambda(k)$ is non-cuspidal and has slope $\alpha$) then $\alpha=0$.
}

{\it Proof. } Let $\Pi$ be the automorphic representation such that ${\cal M}_{k_0,\chi\omega^{-k_0}}(\lambda)=V_\Pi^K$. Hence, $V_\Pi^K$ occurs in ${\cal M}_{k_0}(\Gamma,\chi\omega^{-k_0})\subset H^1(\Gamma,L_{k_0,{\Bbb C}})$ and equation (9) in section 1 shows that
$$
\Pi\cong{\rm Ind}(\chi_{1,f},\chi_{2,f}),
$$
where $\chi_1,\chi_2$ satisfy $\chi_{1,\infty}=|\cdot|_\infty^{1/2}$, $\chi_{2,\infty}={\rm sgn}^\epsilon|\cdot|_\infty^{3/2-k_0}$ with $\epsilon=(-1)^{k_0}$ and $\chi_1\chi_2=\chi\omega^{-k_0}$. Moreover, ${\rm cond}\,\chi_1,\,{\rm cond}\,\chi_2$ both divide $Np$. 

We denote by $\tilde{\omega}:\,{\Bbb Q}^*\backslash{\Bbb A}^*\rightarrow{\Bbb C}^*$ the adelic Teichmuller character, i.e. $\tilde{\omega}$ has conductor $p$ and $\omega$ and $\tilde{\omega}$ are related by $\tilde{\omega}|_{({\Bbb Z}_p/p{\Bbb Z}_p)^*}=\omega^{-1}$, or, equivalently, $\tilde{\omega}_\ell(\ell)=\omega(\ell)$ for all primes $\ell$ different from $p$. We define a family of induced representations by
$$
\Pi_k={\rm Ind}(\chi_{1,f},\tilde{\omega}_f^{k_0-k}|\cdot|_f^{k_0-k}\chi_{2,f})
$$
and we write $\Pi_k=\otimes_\ell\Pi_{k,\ell}$, where
$$
\Pi_{k,\ell}\cong{\rm Ind}(\chi_{1,\ell},\tilde{\omega}^{k_0-k}_\ell|\cdot|_\ell^{k_0-k}\chi_{2,\ell}).
$$
Since the Teichmuller character $\tilde{\omega}$ only ramifies at $p$ and $\chi_1,\chi_2$ are unramified outside $Np$ we see that the same is true for $\Pi_k$, i.e. $\Pi_{k,\ell}$ is unramified for all $\ell\not|Np$. We let $\psi_{k,\ell}\in\Pi_{k,\ell}$, $\ell\not|Np$ be the spherical function and we set 
$$
\psi_k=\otimes_\ell\psi_{k,\ell},\leqno(46)
$$ 
where we have chosen $\psi_{k,\ell}\in\Pi_{k,\ell}$ arbitrarily if $\ell|Np$. Then, for all $\ell\not|Np$ the vector $\psi_k$ is an eigenvector for the Hecke operator $T_\ell$ with eigenvalue
$$
\lambda_{k,\ell}=\ell^{1/2}\,(\chi_{1,\ell}(\ell)+\chi_{2,\ell}(\ell)\tilde{\omega}_\ell^{k_0-k}(\ell)\ell^{k-k_0}).\leqno(47)
$$
%[[cf. [Bump], Proposition 4.6.6, p. 496 ]]
We set $\lambda_k=(\lambda_{k,\ell})_{\ell\not|Np}$. As an immediate consequence of (47) we obtain
$$
\lambda_{k,\ell}-\lambda_{k',\ell}=\chi_{2,\ell}\tilde{\omega}_\ell^{k_0}(\ell)\ell^{-k_0+1/2}   
\big(   (\tilde{\omega}_\ell^{-1}(\ell)\ell)^k    -(\tilde{\omega}_\ell^{-1}(\ell)\ell)^{k'})     \big).\leqno(48)
$$
Since $\tilde{\omega}_\ell^{-1}(\ell)\ell=(\omega^{-1}(\ell)\ell)\equiv 1\pmod{p}$ and $\ell^{-k_0+1/2}$ is a $p$-adic unit because $\ell\not|Np$, equation (48) shows that $k\equiv k'\pmod{p^m}$ implies 
$$
\lambda_{k,\ell}\equiv\lambda_{k',\ell}\pmod{p^{m+1}}.\leqno(49)
$$

%[[(Im folgenden Erorterung sind $\omega$ und $\tilde{\omega}$ vertauscht !) Sei $\tilde{\omega}:\,{\Bbb Z}/p{\Bbb Z}\rightarrow\langle\zeta_{p-1}\rangle$ der teichmuller charachter und $\omega:\,{\Bbb Q}^*\backslash {\Bbb A}^*\rightarrow \langle\zeta_{p-1}\rangle$ der zugehorige Idele klassencharacter. Da der teichmuller character den Fuhrer $p$ hat, gilt also $\tilde{\omega}=\omega|_{{\Bbb Z}/{\Bbb Z}}$. Weiter gilt $1=\omega(-1,\ldots,-1)=\omega_\infty(-1)\omega_p(-1)$ und da $-1\equiv -1\pmod{p}$ und $-1$ eine $(p-1)$-te Einheitswurzel ist gilt auch $\omega_p(-1)=-1$ $\Rightarrow$ $\omega_\infty(-1)=-1$. Der teichmuller character ist also immmer ungerade. Weiter gilt $\omega_\ell(x)=\omega(1,\ldots,1,\ell,1,\ldots,1)=\omega(\ell^{-1},\ldots,\ell^{-1},1,\ell^{-1},\ldots,\ell^{-1})
%=\omega_\infty(\ell^{-1})\omega_p(\ell^{-1})=\omega_p(\ell)^{-1}=\tilde{\omega}^{-1}(\ell)$ (beachte: $\ell^{-1}>0$ ist positiv).
%Mit fast derselben Rechnung beweist man $\tilde{\omega}_p(p)=1$
%]]

To prove the Proposition we still have to look at the place $p$ and to determine the slope subspaces of $V_{\Pi_k,p}^{K_p}$. Since $\chi_{1,\infty}=|\cdot|_\infty^{1/2}$ and $\chi_{2,\infty}=|\cdot|_\infty^{3/2-k_0}$ we may write
$\chi_1=\Psi_1|\cdot|^{1/2}$ and $\chi_2=\Psi_2|\cdot|^{3/2-k_0}$ with $\Psi_1,\Psi_2$ idele class characters of finite order. In particular, $\Psi_1,\Psi_2$ are unitary and $\Pi_p\cong{\rm Ind}(|\cdot|_p^{1/2}\psi_{1,p},|\cdot|_p^{3/2-k_0}\psi_{2,p})$ is irreducible because $1/2-(3/2-k_0)=k_0-1>1$ (note that by our assumption $k_0>2$). Hence,
$$
{\rm Ind}(\chi_{1,p},\chi_{2,p})\cong {\rm Ind}(\chi_{2,p},\chi_{1,p}).\leqno(50)
$$
Since $\Pi$ occurs in ${\cal M}_{k_0}(\Gamma,\chi\omega^{-k_0})$ we know that ${\rm cond}\,\Pi|Np$. Hence, ${\rm cond}\,\Pi_p=0$ or $=p$.

We distinguish cases and first asssume ${\rm cond}\,\Pi=0$. In this case the characters $\chi_{1,p}$ and $\chi_{2,p}$ are unramified and since $\omega^{k-k_0}=1$ if $k\in k_0+(p-1){\Bbb Z}$ all representations $\Pi_{k,p}$, $k\in k_0+(p-1){\Bbb Z}$ then are unramified. We want to determine the effect of the (ramified !) Hecke operator $T_p$ on the space $\Pi_{k,p}^{K_p}$, where $K_p=K_{1,p}(Np)=\{\Mat{a}{b}{c}{d}\in{\rm GL}_2({\Bbb Z}_p):\,c\equiv 0, d\equiv 1\pmod{p}\}$. We know that ${\rm dim}\,\Pi_{k,p}^{K_p}=2$ and we construct a basis of $\Pi_{k,p}^{K_p}$ as follows. For any $\psi\in \Pi_{k,p}^{K_p}$ and $\Mat{\alpha}{\beta}{}{\delta}\in B_2({\Bbb Q}_p)$, $\Mat{a}{b}{c}{d}\in K_p$  we know that
$$
\psi(\Mat{\alpha}{\beta}{}{\delta}g\Mat{a}{b}{c}{d})=\chi_{1,p}(\alpha)|\alpha|_p^{1/2}\chi_{2,p}(\delta)|\delta|_p^{k_0-k-1/2}
\,\chi_{1,p}\chi_{2,p}(d)\leqno(51)
$$
(note that $\tilde{\omega}^{k-k_0}=1$ and $d$ is a $p$-adic unit, hence, $|d|_p^{k_0-k}=1$). Hence, the decomposition

\begin{eqnarray*}
(52)\qquad{\rm GL}_2({\Bbb Q}_p)&=&B_2({\Bbb Q}_p)\Mat{1}{}{1}{1}K_0(p)\dot{\cup} B_2({\Bbb Q}_p)\Mat{1}{}{p}{1}K_0(p)\\
&=&B_2({\Bbb Q}_p)\Mat{1}{}{1}{1}K_0(p)\dot{\cup} B_2({\Bbb Q}_p)K_0(p)\\
\end{eqnarray*}
shows that any element in $\Pi_{k,p}^{K_p}$ is uniquely determined by the values $\psi(\Mat{1}{}{1}{1})$ and $\psi(1)$. We denote by ${\bf 1}_{k,1}$ rep. ${\bf 1}_{k,2}$ the uniquely determined element in $\Pi_{k,p}^{K_p}$, which satisfies ${\bf 1}_{k,1}(\Mat{1}{}{1}{1})=1$ and ${\bf 1}_{k,1}(1)=0$ resp. ${\bf 1}_{k,2}(\Mat{1}{}{1}{1})=0$ and ${\bf 1}_{k,2}(1)=1$. Clearly, ${\bf 1}_{k,1}$ and ${\bf 1}_{k,2}$ span $\Pi_{k,p}^{K_1(p)}$.

We compute the effect of the Hecke operator $T_p$ on ${\bf 1}_{k,1}$and ${\bf 1}_{k,2}$ and we begin with ${\bf 1}_{k,1}$. By definition, 
$$
T_p\,{\bf 1}_{k,1}(\Mat{1}{}{1}{1})=\sum_{u=0}^{p-1} {\bf 1}_{k,1}(\Mat{1}{}{1}{1}\Mat{1}{u}{}{p}).
$$
Using the decomposition
$$
\Mat{1}{}{1}{1}\Mat{1}{u}{}{p}=\Mat{\frac{p}{u+p}}{\frac{u}{u+p}}{}{1}\Mat{1}{}{1}{1}\Mat{1}{}{}{u+p}\leqno(53)
$$
and
$$
\Mat{1}{}{1}{1}\Mat{1}{0}{}{p}=\Mat{p}{1-p}{}{1}\Mat{1}{}{1}{1}\Mat{1}{p-1}{}{1}\leqno(54)
$$
we obtain

\begin{eqnarray*}
T_p{\bf 1}_{k,1} (\Mat{1}{}{1}{1})&=&\sum_{u=1}^{p-1}{\bf 1}_{k,1}(\Mat{\frac{p}{u+p}}{\frac{u}{u+p}}{}{1}\Mat{1}{}{1}{1}\Mat{1}{}{}{u+p})
+{\bf 1}_{k,1}(\Mat{p}{1-p}{}{1}\Mat{1}{}{1}{1}\Mat{1}{p-1}{}{1})\\
&=&\sum_{u=1}^{p-1} \chi_{1,p}(\frac{p}{u+p})|\frac{p}{u+p}|_p^{1/2}  \chi_{1,p}\chi_{2,p}(u+p)
+\chi_{1,p}(p)|p|_p^{1/2}\\
&=&\sum_{u=1}^{p-1} \chi_{1,p}(p)|p|_p^{1/2}+\chi_{1,p}(p)|p|_p^{1/2}\\
&=&p\,\chi_{1,p}(p)|p|_p^{1/2}\\
&=&p^{1/2}\chi_{1,p}(p).
\end{eqnarray*}
The second equality holds because $\tilde{\omega}_p^{k_0-k}=1$ and $u+p$ is a $p$-adic unit if $u\not=0$ and the third equality holds because $\chi_{1,p}$ and $\chi_{2,p}$ are unramified and $u+p$ is a $p$-adic unit if $u\not=0$. Hence, we obtain
$$
T_p {\bf 1}_{k,1} (\Mat{1}{}{1}{1})=p^{1/2}\chi_{1,p}(p).\leqno(55)
$$
Similarly, we find
$$
T_p{\bf 1}_{k,1}(1)=\sum_{u=0}^{p-1} {\bf 1}_{k,1}(\Mat{1}{u}{}{p})\leqno(56)
$$
and this equals $0$ because $\Mat{1}{u}{}{p}\in B_2({\Bbb Q}_p)K_0(p)$. Thus, (55) and (56) yield
$$
T_p{\bf 1}_{k,1}=p^{1/2}\,\chi_{1,p}(p){\bf 1}_{k,1}.
$$ 
We note that 
$$
v_p(p^{1/2}\,\chi_{1,p}(p))=0\leqno(57)
$$
i.e. $p^{1/2}\,\chi_{1,p}(p)$ is a $p$-adic unit, because $\chi_1=|\cdot|^{1/2}\psi_1$ with $\psi_1$ of finite order. Thus, we have seen that ${\bf 1}_{k,1}$ is a $T_p$-eigenform with eigenvalue $p^{1/2}\chi_{1,p}(p)$, in particular, ${\bf 1}_{k,1}$ has slope $0$. 

We compute the effect of $T_p$ on ${\bf 1}_{k,2}$. By definition,
$$
T_p\,{\bf 1}_{k,2}(\Mat{1}{}{1}{1})=\sum_{u=0}^{p-1} {\bf 1}_{k,2}(\Mat{1}{}{1}{1}\Mat{1}{u}{}{p}).
$$
Since equations (53) and (54) imply that 
$$
\Mat{1}{}{1}{1}\Mat{1}{u}{}{p}\in B_2({\Bbb Q})\Mat{1}{}{1}{1}K_0(p)
$$
we obtain 
$$
T_p\,{\bf 1}_{k,2}(\Mat{1}{}{1}{1})=0.\leqno(58)
$$ 
Moreover, 
\begin{eqnarray*}
(59)\qquad T_p\,{\bf 1}_{k,2}(1)&=&\sum_{u=0}^{p-1}{\bf 1}_{k,2}(\Mat{1}{u}{}{p})\\
&=&\sum_{u=0}^{p-1}\chi_{2,p}(p)|p|_p^{k_0-k-1/2}\\
&=&p \chi_{2,p}(p)|p|_p^{k_0-k-1/2}\\
&=&p^{3/2} \chi_{2,p}(p)|p|_p^{k_0-k}.\\
\end{eqnarray*}
Equations (58) and (59) together yield
$$
T_p\,{\bf 1}_{k,2}=p^{3/2} \chi_{2,p}(p)|p|_p^{k_0-k}\,{\bf 1}_{k,2}.
$$
Since $\chi_2=|\cdot|^{3/2-k_0}\psi_2$ with $\psi_2$ of finite order we obtain
$$
v_p(p^{3/2} \chi_{2,p}(p)|p|_p^{k_0-k})=k.\leqno(60)
$$
Thus, we have seen that ${\bf 1}_{k,2}$ is a $T_p$-eigenform with eigenvalue $p^{3/2}\chi_{2,p}(p)|p|_p^{k_0-k}$, in particular, ${\bf 1}_{k,2}$ has slope $k$. Equations (57) and (60) imply that the eigenvalues $p^{-1/2}\,\chi_{1,p}(p)$ and $p^{3/2} \chi_{2,p}(p)|p|^{k_0-k}$ are different ($k$ is an integer $>0$). Hence, $T_p$ is diagonalizable on $\Pi_{k,p}^{K_p}$ and the only slopes occuring in $\Pi_{k,p}^{K_p}$ are $\alpha=0$ and $\alpha=k$ with corresponding slope spaces
$$
\Pi_{k,p}^{K_p,0}={\Bbb C}\,{\bf 1}_{k,1}\qquad \Pi_{k,p}^{K_p,k}={\Bbb C}{\bf 1}_{k,2}.
$$
This holds for all $k\in k_0+(p-1){\Bbb Z}$ and proves the first claim of the Proposition. In particular, we obtain a family of constant slope only we select $\psi_{k,p}={\bf 1}_{k,1}$ and the family then has constant slope $\alpha=0$. On the other hand if $\lambda$ has slope $0$ we specify the $p$-component of the family $(\psi_k)_k$ as $\psi_{k,p}={\bf 1}_{k,1}$ (we leave the components $\psi_{k,\ell}$, $\ell |N$ unspecified); is then immediate that $(\psi_k)_k$ satisfies the following properties

$\bullet$ $\psi_k$ is an eigenform for all Hecke operators $T_\ell$, $\ell\not|Np$ with eigenvalue $\lambda_{k,\ell}$ defined in (47); in particular, the eigenvalues $\lambda_{k,\ell}$, $\ell\not|Np$ satisfy the congruences (49)

$\bullet$ $\psi_k$ is an eigenform for $T_p$ with eigenvalue $p^{1/2}\chi_{1,p}(p)$, i.e. $\psi_k$ has slope $0$

$\bullet$ Since, in particular, $k\equiv k_0\pmod{2}$ we find $(-1)^k=(-1)^{k_0}$ and since furthermore $\tilde{\omega}^{k_0-k}=1$, equation (9) in section 1 implies that $\Pi_{k,p}^{K_p}\hookrightarrow {\cal M}_k(\Gamma,\chi\omega^{-k})$. Hence, $\psi_k\in {\cal M}_k(\Gamma,\chi\omega^{-k})^0$, i.e. the system of Hecke eigenvalues $\lambda_k=(\lambda_{k,\ell})_{\ell\not|Np}$ is contained in $\Lambda_{k,\chi,Np}^\alpha$.

$\bullet$ Since $\Pi_k$ is not a cuspidal representation, the system of Hecke eigenvalues $\lambda_k$ is non-cuspidal.

Thus, $\lambda$ can be placed into a $p$-adic family of constant slope $0$. This proves second claim and, hence, the proposition is proven in case $\Pi_p$ is unramified.

We look at the second case, i.e. we assume that ${\rm cond}\,\Pi_p=p$. In this case either $({\rm cond}\,\chi_{1,p},{\rm cond}\,\chi_{2,p})=(1,p)$ or $({\rm cond}\,\chi_{1,p},{\rm cond}\,\chi_{2,p})=(p,1)$ and in view of (50) we may assume that $({\rm cond}\,\chi_{1,p},{\rm cond}\,\chi_{2,p})=(1,p)$. In particular, $\chi_{1,p}$ is unramified. Again, since $k$ is congruent to $k_0$ modulo $(p-1){\Bbb Z}$ the factor $\tilde{\omega}_p^{k_0-k}$ vanishes and all characters $\tilde{\omega}^{k_0-k}_p\chi_{2,p}=\chi_{2,p}$ have conductor equal to $p$; in particular, the conductor of $\Pi_{k,p}$ equals $p$ for all $k\in k_0+(p-1){\Bbb Z}$. Hence, the space of $K_p$-invariants $\Pi_{k,p}^{K_p}$ in $\Pi_{k,p}$ is $1$-dimensional and the discussion on top of p. 306 in [Ca] shows that $\Pi_{k,p}^{K_p}$ is spanned by the function ${\bf 1}_{k,1}$:
$$
\Pi_{k,p}^{K_p}={\Bbb C}\,{\bf 1}_{k,1}.
$$ 
Here, ${\bf 1}_{k,1}$ is the unquely determined function ${\rm GL}_2({\Bbb Q}_p)\rightarrow{\Bbb C}$ satisfying (51) as well as ${\bf 1}_{k,1}(\Mat{1}{}{1}{1})=1$ and ${\bf 1}_{k,1}(1)=0$ (compare equation (52); note that the analogoue of the function ${\bf 1}_{k,2}$ is not well defined because $\chi_{2,p}$ ramifies). We compute the effect of $T_p$ on  ${\bf 1}_{k,1}$. To this end it suffices to determine $T_p\,{\bf 1}_{k,1}(\Mat{1}{}{1}{1})$ because ${\bf 1}_{k,1}$ is a $T_p$-eigenform. A computation completely analogous to the computation of $T_p\,{\bf 1}_{k,1}(\Mat{1}{}{1}{1})$ in the unramified case yields
$$
T_p\,{\bf 1}_{k,1}(\Mat{1}{}{1}{1})=p^{1/2}\chi_{1,p}(p).
$$
Again, this holds for all $k\in k_0+(p-1){\Bbb Z}$. In particular, ${\bf 1}_{k,1}$ is a $T_p$-eigenvector with eigenvalue $p^{1/2}\chi_{1,p}(p)$. Since $\chi_1=|\cdot|^{1/2}\Psi_1$, with $\Psi_1$ of finite order we obtain
$$
v_p(p^{1/2}\chi_{1,p}(p))=0.
$$
Hence, the only slope occuring in $\Pi_{k,p}^{K_p}$ is $\alpha=0$:
$$
\Pi_p^{K_p}=\Pi_p^{K_p,0}={\Bbb C} {\bf 1}_{k,1},
$$
which proves the first claim of the Proposition. In particular, we have to specify the $p$-component $\psi_{k,p}$ of $\psi_k$ to be ${\bf 1}_{k,1}$; $\psi_k$ then has slope $0$ for all $k\in k_0+(p-1){\Bbb Z}$. Altogether we see that the family $(\psi_k)_k$ satisfies the following properties, which are completely analogous to the unramified case:

$\bullet$ $\psi_k$ is an eigenform for all Hecke operators $T_\ell$, $\ell\not|Np$ with eigenvalue $\lambda_{k,\ell}$ defined in (47); in particular, the eigenvalues $\lambda_{k,\ell}$ satisfy the congruences (49)

$\bullet$ $\psi_k$ is an eigenform for $T_p$ with eigenvalue $p^{1/2}\chi_{1,p}(p)$, i.e. $\psi_k$ has slope $0$

$\bullet$ As in the first case, equation (9) in section 1 shows that $\Pi_{k,p}^{K_p}\hookrightarrow {\cal M}_k(\Gamma,\omega^{-k}\chi)$, hence, $\psi_k\in {\cal M}_k(\Gamma,\omega^{-k}\chi)^\alpha$, i.e. the system of Hecke eigenvalues $\lambda_k=(\lambda_{k,\ell})_{\ell\not|Np}$ is contained in $\Lambda_{k,\chi,Np}^\alpha$.

$\bullet$ Since $\Pi_k$ is not a cuspidal representation, the system of Hecke eigenvalues $\lambda_k$ is non-cuspidal.

This proves the second claim and concludes the proof of the Proposition in case $\Pi_p$ ramifies. Thus the proposition is proven.

{\bf Corollary. }{\it Fix $\alpha\in{\Bbb Q}_{\ge 0}$ and let $\lambda\in\Lambda_{k_0,\chi,Np}^\alpha$ be non-cuspidal. Assume that  $k_0\ge({\sf C}_\alpha\alpha+1)^2+2$. Then, $\alpha=0$ and, hence, there is a family $(\lambda_k)_k$, $k\in k_0+(p-1){\Bbb Z}$, where $\lambda_k=(\lambda_{k,\ell})_{\ell\not|Np}\in \Lambda_{k,\chi,Np}^0$ such that 

$\bullet$ $\lambda_{k_0}=\lambda$

$\bullet$ $\lambda_k$ is non-cuspidal for all $k\in k_0+(p-1){\Bbb Z}$

$\bullet$ $k\equiv k'\pmod{p^{m}}$ implies $\lambda_{k,\ell}\equiv \lambda_{k',\ell}\pmod{p^{m+1}}$ for all $\ell\not|Np$.

}

{\it Proof. } By the above Proposition $\lambda$ has slope $\alpha=0$ or $\alpha=k_0$ and the claim follows from the proposition if we can show that $\alpha=0$. To this end assume that $\alpha=k_0$. Since $\psi$ has weight $k_0$, by our assumption on $k_0$ we obtain
$$
k_0\ge ({\sf C}_{k_0}k_0+1)^2+2\ge ({\sf C}_{k_0}k_0)^2.
$$
Recalling that ${\sf C}_{k_0}\ge N(k_0)+\frac{1}{k_0}$ (cf. equation (45)) we further obtain $k_0\ge (N(k_0)k_0+1)^2$ and since $N(k_0)\ge 1$ (cf. equation (32) in section 3) this implies $k_0\ge N(k_0)^2k_0^2+2N(k_0)k_0+1\ge k_0^2+1$, a contradiction ($k_0$ is a positive integer). Thus, the slope $\alpha=k_0$ does not satisfy $k_0\ge ({\sf C}_\alpha\alpha+1)^2$ and we deduce that $\alpha=0$. This finishes the proof of the Corollary.

From now on we asssume $k>2$. We denote by $E_{k,\chi,Np}^\alpha$ the set of all irreducible representations $\Pi$ of ${\rm Gl}_2({\Bbb A}_f)$ such that $\Pi^{K,\alpha}\not=0$, $\Pi^{,\alpha}K$ occurs in ${\cal M}_k(\Gamma,\chi\omega^{-k})^\alpha$ and $\Pi$ is non-cuspidal, i.e. $\Pi^K$ does not occur in the subspace ${\cal S}_k(\Gamma,\chi\omega^{-k})$ of cusp forms. Proposition D shows that $E_{k,\chi,Np}^\alpha\not=\emptyset$ implies $\alpha=0$ or $\alpha=k$. We recall from equation (9) in section 1.3 that ${\cal M}_k(\Gamma,\chi\omega^{-k})^\alpha$ is the direct sum
$$
{\cal M}_k(\Gamma,\chi\omega^{-k})^\alpha={\cal E}_k(\Gamma,\chi\omega^{-k})^\alpha\oplus  {\cal S}_k(\Gamma,\chi\omega^{-k})^\alpha,
$$
where
$$
{\cal E}_k(\Gamma,\chi\omega^{-k})^\alpha=\bigoplus_{\Pi\in E_{k,\chi,Np}^\alpha}\Pi^{K,\alpha}
$$
(note that $k>2$). Moreover, equation (9) in section 1.3 shows that any $\Pi\in E_{k,\chi,Np}^\alpha$ is of the form
$$
\Pi={\rm Ind}(\chi_{1,f},\chi_{2,f})
$$
for characters $\chi_1,\chi_2:\,{\Bbb Q}^*\backslash{\Bbb A}^*\rightarrow{\Bbb C}^*$ of conductor $Np$, which satisfy the conditions $\chi_{1,\infty}=|\cdot|_\infty^{1/2}$, $\chi_{2,\infty}={\rm sgn}^\epsilon|\cdot|_\infty^{3/2-k_0}$ with $\epsilon=(-1)^{k_0}$ and $\chi_1\chi_2=\chi\omega^{-k_0}$. In the proof of Proposition D we have seen that the slope $0$ subspace of any non-cuspidal representation $\Pi={\rm Ind}(\chi_{1,f}\chi_{2,f})$ occuring in ${\cal M}_k(\Gamma,\chi\omega^{-k})$ is non-trivial (in fact we have seen that ${\rm dim}\,\Pi^{K,0}=1$); we thus obtain a map
$$
\begin{array}{cccc}
\Phi_{k,k'}:&E_{k,\chi,Np}^0&\rightarrow&E_{k',\chi,Np}^0\\
&{\rm Ind}(\chi_{1,f},\chi_{2,f})&\mapsto&{\rm Ind}(\chi_{1,f},\tilde{\omega}_f^{k'-k}|\cdot|_f^{k-k'}\chi_{2,f}).\\
\end{array}
$$
Clearly, $\Phi_{k,k'}$ is a bijection with inverse $\Phi_{k',k}$.

{\bf Lemma. }{\it Let $\Pi\in E_{k,\chi,Np}^0$ and $\Pi'\in E_{k',\chi,Np}^0$ correspond to each other under $\Phi_{k,k'}$. Assume that $k\equiv  k'\pmod{(p-1)p^m}$. Then, for all primes $\ell_1,\ldots,\ell_s$, which do not divide $Np$, we have

%1.) ${\rm dim}\,\Pi^K={\rm dim}\,{\Pi'}^K$ and equality of dimension also holds for the slope $0$ subspaces: ${\rm dim}\,\Pi^{K,0}={\rm dim}\,{\Pi'}^{K,0}$

$$
{\rm tr}\,\pi^{Np}_{Mp} T_{\ell_1}\cdots T_{\ell_s}|_{\Pi} 
\equiv {\rm tr}\,\pi^{Np}_{Mp} T_{\ell_1}\cdots T_{\ell_s}|_{\Pi'}\pmod{p^{m+1}}\qquad(M|N).
$$

}

{\it Proof. } %1.) Since $k\equiv k'\pmod{p-1}$ we know that $\tilde{\omega}^{k'-k}=1$, hence, $\Pi$ and $\Pi'$ share the same conductor, which by the Corollary to proof in [Ca], p. 306 implies that ${\rm dim}\,\Pi^K={\rm dim}\,{\Pi'}^K$. On the other hand, in the proof of Proposition D we have seen that ${\rm dim}\,\Pi_p^{K,0}={\rm dim}\,{\Pi'_p}^{K,0}=1$, which yields the second claim.
$T_{\ell_i}$, $\ell_i\not|Np$ acts on $\Pi^{K,0}$ resp. on ${\Pi'}^{K,0}$ as multiplication with some number $\lambda_{\ell_i}$ resp. $\lambda_{\ell_i}'$ and in the proof of Proposition E we have seen that $\lambda_{\ell_i}\equiv \lambda_{\ell_i}'\pmod{p^{m+1}}$ if $k\equiv k'\pmod{(p-1)p^m}$. 
%(note that $\Pi^{K,0}={\cal M}_{k,\chi\omega^{-k}}^0(\lambda_\Pi)$). 
In particular, $T_{\ell_i}$ acting on $\Pi^{K,0}$ resp. on ${\Pi'}^{K,0}$ is represented by the matrix $\lambda_{\ell_i}{\rm Id}$ resp. $\lambda_{\ell_i}'{\rm Id}$ with respect to any basis of $\Pi^{K,0}$ resp. of ${\Pi'}^{K,0}$. On the other hand, in the proof of the Theorem in section 4.2 we have seen that $\pi^{Np}_{Mp}|_{\Pi^{K,0}}$ resp. $\pi^{Np}_{Mp}|_{{\Pi'}^{K,0}}$
with respect to a suitable basis is represented by the matrix
$$
{\cal D}_{\cal B}(\pi^{Np}_{Fp})=\left(
\begin{array}{cccccc}
[\Gamma_1(Mp):\Gamma_1(Np)]&&&&*&\\
&\ddots&&&&\\
&&[\Gamma_1(Mp):\Gamma_1(Np)]&&&\\
&&&0&\cdots&0\\
&&&&\ddots&\vdots\\
&&&&&0\\
\end{array}
\right).
$$
Here, the entry $[\Gamma_1(Mp):\Gamma_1(Np)]$ appears $d={\rm dim}\,({\Pi^{K,0}})^{K_1(Mp)}$-times resp. $d'={\rm dim}\,({{\Pi'}^{K,0}})^{K_1(Mp)}$-times on the diagonal. Since $\lambda_{\ell_i}\equiv \lambda_{\ell_i}'\pmod{p^{m+1}}$ the claim follows if we can show that $d=d'$. To prove this write $({\Pi^{K,0}})^{K_1(Mp)}=\Pi^{K_1(Mp),0}=\otimes_{\ell\not=p,\infty}\Pi_\ell^{K_{1,\ell}(Mp)}\otimes\Pi_p^{K_{1,p}(Mp),0}$ and analogously for $({{\Pi'}^{K,0}})^{K_1(Mp)}$. The Corollary to proof in [Ca], p. 306 shows that ${\rm dim}\,\Pi_\ell^{K_{1,\ell}(Mp)}={\rm dim}\,{\Pi_\ell'}^{K_{1,\ell}(Mp)}$ for all $\ell\not=p,\infty$ because $\Pi$ and $\Pi'$ share the same conductor. In the proof of Proposition E we have seen that ${\rm dim}\,\Pi_p^{K_{1,p}(Mp),0}={\rm dim}\,{\Pi_p'}^{K_{1,p}(Mp),0}=1$. Hence, $d=d'$ and the proof of the lemma is complete.

Substracting the identity of the Lemma from the Basic trace identity we obtain

{\bf Proposition }(Basic Trace identity for cusp froms). {\it Let $\alpha\in{\Bbb Q}_{\ge 0}$. Assume that $k,k'$ satisfy

$\bullet$  $k,k'\ge ({\sf C}_\alpha\alpha+1)^2+2$

$\bullet$  $k\equiv k'\pmod{(p-1)p^m}$ with $m>{\sf B}_\alpha$.

and let $\ell_1,\ldots,\ell_s$ be primes, which do not divide $Np$. Then, the same congruences as in the Theorem in section 3.5 with "${\cal M}$" replaced by "${\cal S}$" hold true, i.e.
$$
{\rm tr}\,\pi^{Np}_{Mp} T_{\ell_1}^{r_1}\cdot\ldots\cdot T_{\ell_s}^{r_s}|_{{\cal S}_k(\Gamma,\chi\omega^{-k})^\alpha}
\equiv{\rm tr}\,\pi^{Np}_{Mp} T_{\ell_1}^{r_1}\cdot\ldots\cdot T_{\ell_s}^{r_s}|_{{\cal S}_{k'}(\Gamma,\chi\omega^{-k'})^\alpha}\pmod{p^{m-v_p(\varphi(N))}}.
$$
if $\alpha=0$ and
$$
{\rm tr}\,\pi^{Np}_{Mp} T_{\ell_1}^{r_1}\cdot\ldots\cdot T_{\ell_s}^{r_s}e_\alpha^{[\frac{m}{{\sf C}_\alpha\alpha}]}|_{{\cal S}_k(\Gamma,\chi\omega^{-k})^\alpha}
\equiv{\rm tr}\,\pi^{Np}_{Mp} T_{\ell_1}^{r_1}\cdot\ldots\cdot T_{\ell_s}^{r_s}e_\alpha^{[\frac{m}{{\sf C}_\alpha\alpha}]}|_{{\cal S}_{k'}(\Gamma,\chi\omega^{-k'})^\alpha}\pmod{p^{[\frac{m}{{{\sf C}_\alpha}\alpha}]-v_p(\varphi(N))}}.
$$
if $\alpha>0$.

}

{\it Proof. } We note that
$$
{\cal M}_k(\Gamma,\chi\omega^{-k})^\alpha={\cal E}_k(\Gamma,\chi\omega^{-k})^\alpha\oplus {\cal S}_k(\Gamma,\chi^{-k})^\alpha
$$
Thus, to prove the proposition it remains to compare ${\rm tr}\,\pi^{Np}_{Mp} T_{\ell_1}^{r_1}\cdot\ldots\cdot T_{\ell_s}^{r_s}e_\alpha^{[\frac{m}{{\sf C}_\alpha\alpha}]}|_{{\cal E}_k(\Gamma,\chi\omega^{-k})^\alpha}$ and ${\rm tr}\,\pi^{Np}_{Mp} T_{\ell_1}^{r_1}\cdot\ldots\cdot T_{\ell_s}^{r_s}e_\alpha^{[\frac{m}{{\sf C}_\alpha\alpha}]}|_{{\cal E}_{k'}(\Gamma,\chi\omega^{-k'})^\alpha}$.
We distinguish: 

Case 1: $\alpha>0$. The above Corollary shows that $E_{k,\chi,Np}^\alpha=\emptyset$, hence, ${\cal E}_k(\Gamma,\chi\omega^{-k})^\alpha={\cal E}_{k'}(\Gamma,\chi\omega^{-k'})^\alpha=0$ and the Proposition is immediate by the Theorem in section 3.5.

Case 2: $\alpha=0$. Since $\Phi_{k,k'}:\,E_{k,\chi,Np}^0\rightarrow E_{k',\chi,Np}^0$ is a bijection such that 
$$
{\rm tr}\,\pi^{Np}_{Mp} T_{\ell_1}\cdots T_{\ell_s}|_{\Pi} 
\equiv {\rm tr}\,\pi^{Np}_{Mp} T_{\ell_1}\cdots T_{\ell_s}|_{\Phi_{k,k'}(\Pi)}\pmod{p^{m+1}}\qquad(M|N)
$$
(cf. the above Lemma) we obtain
$$
{\rm tr}\,\pi^{Np}_{Mp} T_{\ell_1}\cdots T_{\ell_s}|_{{\cal E}_k(\Gamma,\chi\omega^{-k})^0} 
\equiv {\rm tr}\,\pi^{Np}_{Mp} T_{\ell_1}\cdots T_{\ell_s}|_{{\cal E}_{k'}(\Gamma,\chi\omega^{-k'})^0}\pmod{p^{m+1}}\qquad(M|N).
$$
Again, the claim follows immediately from this and the Theorem in section 3.5. Thus, the proposition is proven.

Following the exposition in section 4.2 while replacing "${\cal M}$" by "${\cal S}$" throughout we obtain the following result from Corollary C3. We denote by $\Lambda_{k_0,\chi,Np}^{\alpha,cusp}$ the set of all $\lambda\in \Lambda_{k_0,\chi,Np}^{\alpha}$ such that 
${\cal M}_ {k,\chi\omega^{-k}}^\alpha(\lambda)\subset {\cal S}_k(\Gamma,\chi\omega^{-k})^\alpha$, i.e. ${\cal M}_ {k,\chi\omega^{-k}}^\alpha(\lambda)$ is cuspidal.

{\bf Addendum to Theorem D1, D2. } {\it Let $\alpha\in{\Bbb Q}_{\ge 0}$, let $k_0\ge ({{\sf C}_\alpha}\alpha+1)^2+2$. Then, any $\lambda\in\Lambda_{k_0,\chi,Np}^{\alpha,cusp}$ can be placed into a cuspidal $p$-adic family. More precisely, there is a family $(\lambda(k))_k$, where $k$ runs over all integers satisfying $k\ge({{\sf C}_\alpha}\alpha+1)^2+2$ as well as 
$k\equiv k_0\pmod{p^{{\sf L}(\alpha,k_0)}}$, which satisfies the properties

$\bullet$ $\lambda(k)\in\Lambda_{k,\chi,Np}^{\alpha,cusp}$

$\bullet$ $\lambda(k_0)=\lambda$

$\bullet$ $k\equiv k'\pmod{p^m}$ implies $\lambda(k)\equiv \lambda(k')\pmod{p^{{\sf a}m+{\sf b}-A_{k_0}}}$.

Here, ${\sf a}={\sf a}(\alpha)$ and ${\sf b}={\sf b}(\alpha)$ only depend on $\alpha$ and are defined as the same quantities in section 4.2 with "${\cal M}$" replaced by "${\cal S}$", hence, 
$$
{\sf a}\le {\rm dim}\,{\cal S}_k(\Gamma,\chi\omega^{-k})^\alpha
$$ 
for all $k\ge{\sf C}_\alpha\alpha+1)^2+1$ and even
$$
{\sf a}={\rm min}\,\{ {\rm dim}\,{\cal S}_k(\Gamma,\chi\omega^{-k})^\alpha,\,k\ge {\sf C}_\alpha\alpha+1)^2+1\}
$$
if $\alpha=0$.

If in addition ${\rm dim}\,{\cal S}_{k_0}(\Gamma,\chi\omega^{-k_0})^\alpha=1$ holds then any $\lambda\in\Lambda_{k_0,\chi,Np}^{\alpha,cusp}$ can be placed into a cuspidal $p$-adic family, which exists for all $k$ such that $k\ge({{\sf C}_\alpha}\alpha+1)^2+2$ and
$k\equiv k_0\pmod{p^{{\sf B}_\alpha}}$ and which satisfies the congruences
$$
\lambda(k)\equiv \lambda(k')\pmod{p^{{\sf a}m+{\sf b}}}\quad {\rm if} \quad k\equiv k'\pmod{p^m}.
$$

}

{\bf Remark. } Under the assumption $k>({\sf C}_\alpha\alpha+1)^2+2$ we know that ${\cal E}_k(\Gamma,\chi\omega^{-k})^\alpha$ does not vanish only for $\alpha=0$. Thus, for all $\alpha>0$ the space ${\cal M}_k(\Gamma,\chi\omega^{-k})^\alpha$ appearing in Theorem D2 is spanned by  cusp forms, i.e. 
$$
{\cal M}_k(\Gamma,\chi\omega^{-k})^\alpha={\cal S}_k(\Gamma,\chi\omega^{-k})^\alpha
$$
and we are in the situation of the above Addendum to Theorem D1, D2. In different words, if $\alpha>0$ the statements of Theorems D1, D2 and of the Addendum are equivalent.

{\bf 4.5. Eigenvalues of Hecke operators at ramified places. } Let $k_0\ge ({\sf C}_\alpha\alpha+1)^2+2$ and $\lambda\in\Lambda_{k_0,\chi}^\alpha$. Theorem D1 shows that there is a family $(\lambda(k))_k$, where $k$ runs over all integers larger than $({\sf C}_\alpha\alpha+1)^2+2$ and congruent to $k_0\pmod{p^{{\sf L}(\alpha,k_0)}}$ and where  $\lambda(k)=(\lambda(k)_\ell)_{\ell\not|Np}\in\Lambda_{k,\chi,Np}^\alpha$ such that the following congruence holds: $\lambda(k)\equiv \lambda(k')\pmod{p^{{\sf a}m+{\sf b}-A_{k_0}}}$ if $k\equiv k'\pmod{p^m}$. We denote by $(\pi(k),V(k))$ the represntation of ${\rm Gl}_2({\Bbb A}_f)$ attached to $\lambda(k)$ and we write $V(k)=\otimes_{\ell\not=\infty}V(k)_\ell$. In this last section we shall show that for any $k$ as above and any $\ell$ dividing $Np$ we can pick an eigenvalue $\lambda(k)_\ell$ of $T_\ell$ acting on $V(k)_\ell$ such that the same type of congruences as for the places $\ell\not|Np$ also hold at the ramified places $\ell|Np$. 
Hence, any $\lambda\in\Lambda_{k_0,\chi}^\alpha$ can be placed into a $p$-adic family $(\lambda(k))_k$, where now $\lambda(k)\in \Lambda_{k,\chi}^\alpha$, i.e. the elements of the family are defined for all primes $\ell$ and not just the primes $\ell$, which do not divide $Np$.

We recall that we defined a $p$-adic valuation $v_p$ on ${\Bbb C}$ by using an isomorphism $i:\,{\Bbb C}_p\cong {\Bbb C}$; for any polynomial $f=\sum_{h=0}^d a_h X^h\in{\Bbb C}[X]$ we therefore may set
$$
v_p(f)={\rm min}\,\{v_p(a_h),\, h=0,\ldots,d\}.
$$

{\bf Lemma. }{\it Assume that $f=\sum_{h=0}^d a_h X^h\in\bar{\Bbb Q}[X]$ and $f'=\sum_{h=0}^d a'_h X^h\in\bar{\Bbb Q}[X]$ satisfy the following properties

$\bullet$ $f,f'$ have $p$-adically integral coefficients and leading term equal to $1$, i.e. $v_p(a_h),v_p(a_h')\ge 0$ for all $h=1,\ldots,d$  and $a_d=a_d'=1$

$\bullet$ $v_p(f-f')\ge m$ with $m\ge 0$. 

Then, for any zero $\lambda$ of $f$ there is a zero $\lambda'$ of $f'$ such that $v_p(\lambda-\lambda')\ge m/d$.

}

{\it Proof. } We denote by ${\cal O}$ the ring of integers of ${\Bbb Q}_p(a_1,\ldots,a_d)$. Since $v_p(a_h)\ge 0$ for all $h=1,\ldots,d$, the coefficients $a_h$ of $f$ are contained in ${\cal O}$ and since furthermore the leading coefficient of $f$ equals $1$, any root $\lambda$ of $f$ is integral over ${\cal O}$. Hence, $v_p(\lambda)\ge 0$ for any root $\lambda$ of $f$. Clearly, the same is true of the roots of $f'$ and we deduce that the Lemma holds trivially in case $m=0$. We therefore assume $m>0$ and we fix any root $\lambda$ of $f$. Let $\lambda_1',\ldots,\lambda_d'$ denote the roots of $f'$, each root appearing as often as its multiplicity. Since $v_p(f-f')\ge m$ we obtain $v_p(f(\lambda)-f'(\lambda))\ge m$. Taking into account that $f(\lambda)=0$ we thus find
$$
v_p(f'(\lambda))=v_p(\prod_{h=0}^d ({\lambda-\lambda_h'}))\ge m.\leqno(63)
$$
Now, if $v_p(\lambda-\lambda_i')$ was strictly less than $m/d$ for all zeroes $\lambda_i'$ of $f'$ we obtained $v_p(\prod_{h=0}^d ({\lambda-\lambda_i'}))<m$, a contradiction to (63). Thus, the lemma is proven.
%
%[[Note that if $m>0$ then $a_d=1$ implies that $a_d'$ is congruent $1$ modulo $p^m$, hence $a_d'$ too is a $p$ adic unit and in particular does not vanish]]

{\bf Proposition. } {\it Let $f_k=\sum_{h=0}^d a_{k,h}X^h\in {\Bbb C}[X]$, $k\in {\Bbb N}$ be a family of polynomials such that the following holds.

$\bullet$ Each $f_k$ has $p$-adically integral coefficients and leading coefficient equal to $1$

$\bullet$ There are $a,b\in{\Bbb Q}$, $a>0$ such that the following congruence holds. 
$$
\mbox{$k\equiv k'\pmod{p^m}$ implies $v_p(f_k- f_{k'})\ge{am+b}$.}\leqno(64) 
$$
We denote by ${\cal Z}_k$ the set of roots of $f_k$. Then, if $\lambda\in {\cal Z}_{k_0}$ is any root of $f_{k_0}$, there is a family $(\lambda_k)_k$, $\lambda_k\in{\cal Z}_k$ such that $\lambda_{k_0}=\lambda$ and the following congruence holds: 
$$
\mbox{$k\equiv k'\pmod{p^m}$ implies $v_p(\lambda_k-\lambda_{k'})\ge \frac{am+b}{d}$}.\leqno(65)
$$

%[[
%Dass alle $f_k$ alle denselben grad haben kann man O:E: annehmen: setze einfach fuhrende Koeffizienten $=0$
%]]
}

{\it Proof.} We define the sequence $\lambda_k$ inductively. Of course, we set $\lambda_{k_0}=\lambda$. Next, we assume that zeroes $\lambda_{k_0},\lambda_{k_1},\ldots,\lambda_{k_n}$ have been chosen such that $\lambda_{k_h}\in{\cal Z}_{k_h}$ and the congruences (64) hold 
for $k,k'\in\{k_0,\ldots,k_n\}$. We select any $k=k_{n+1}\in{\Bbb N}-\{k_0,\ldots,k_n\}$. Let $k_i\in\{k_0,\ldots,k_n\}$ be closest to $k$, i.e. 
$$
v_p(k-k_i)\ge v_p(k-k_j)
$$
for all $j=0,\ldots,n$. We set $m=v_p(k-k_i)$. By using the above Lemma we find a zero $\lambda_{k_{n+1}}\in{\cal Z}_{k_{n+1}}$ such that
$$
v_p(\lambda_{k_{n+1}}-\lambda_{k_i})\ge \frac{am+b}{d}.
$$
We will show that $\lambda_{k_0},\lambda_{k_1},\ldots,\lambda_{k_n},\lambda_{k_{n+1}}$ satisfies the congruence (64) thereby completing the induction and, hence, proving the Proposition. Clearly, only the congruences involving $\lambda=\lambda_{k_{n+1}}$ still need to be verified. To this end we select an arbitrary  $\lambda_{k_j}\in\{\lambda_{k_0},\lambda_{k_1},\ldots,\lambda_{k_n}\}$. We set $m'=v_p({k_i}-{k_j})$ 
and distinguish:

Case 1: $v_p(k-k_i)>v_p(k_i-k_j)$, i.e. $m>m'$. We obtain, $v_p(k-k_j)=v_p(k-k_i+k_i-k_j)={\rm min}(v_p(k-k_i),v_p(k_i-k_j))=m'$ by the $p$-adic triangle inequality. On the other hand, $v_p(\lambda-\lambda_{k_j})=v_p(\lambda-\lambda_{k_i}+\lambda_{k_i}-\lambda_{k_j}))={\min}(\frac{am+b}{d},\frac{am'+b}{d})=\frac{am'+b}{d}$. 
%
%[[Hier brauchen wir $a>0$ damit aus $m>m'$ folgt $am+b>am'+b$ ]]
%
Hence, the claimed congruence (65) holds in this case.

Case 2: $v_p(k-k_i)<v_p(k_i-k_j)$, i.e. $m<m'$. Here, the congruence (65) follows in the same way as in case 1.

Case 3: $v_p(k-k_i)=v_p(k_i-k_j)$, i.e. $m=m'$. Here, $v_p(k-k_j)=v_p(k-k_i+k_i-k_j)\ge m$ and since $k_i$ was closest to $k$ we obtain $v_p(k-k_j)=m$.  On the other hand, $v_p(\lambda-\lambda_{k_j})=v_p(\lambda-\lambda_{k_i}+\lambda_{k_i}-\lambda_{k_j}))\ge{\min}(\frac{am+b}{d},\frac{am'+b}{d})=\frac{am+b}{d}$. Hence, the claimed congruence (65) holds in this case too.

This completes the proof of the Proposition.

Before we turn to the proof of our last theorem we make the following observation concerning the Theorem and Theorem' in section 4.2. The Theorem and Theorem' in section 4.2 still hold if we replace the set $\Lambda_{k,\chi,Np}^\alpha$ by $\Lambda_{k,\chi,N}^\alpha$ and $\Lambda_{k',\chi,Np}^\alpha$ by $\Lambda_{k',\chi,N}^\alpha$, i.e. for any system of eigenvalues $\lambda=(\lambda_\ell)_{\ell\not|N}$ in $\Lambda_{k,\chi,N}^\alpha$ there is a system of eigenvalues $\lambda'=(\lambda'_\ell)_{\ell\not|N}$ in $\Lambda_{k',\chi,N}^\alpha$ such that the following holds

$\bullet$ $\lambda_\ell\equiv\lambda_\ell'\pmod{p^{{\sf a}m+{\sf b}-A_k}}$ for all $\ell\not|N$

$\bullet$ $F'|F$, where $F$ resp. $F'$ is the prime to $p$ part of the conductor of $\lambda$ resp. of $\lambda'$.

The reason for this is that in addition to the Hecke operators $T_\ell$, $\ell\not|Np$, the Hecke operator $T_p$ too commmutes with $\pi^{Np}_{Mp}$. The proof of the Theorem in section 4.2 therefore carries over word for word to yield the above stronger result, which includes the prime $p$. As a consequence, Corollary C2 and Theorems D1, D2 also hold if we replace $\Lambda_{k,\chi,Np}^\alpha$ by $\Lambda_{k,\chi,N}^\alpha$ and $\Lambda_{k',\chi,Np}^\alpha$ by $\Lambda_{k',\chi,N}^\alpha$, i.e. for any $\lambda\in\Lambda_{k_0,\chi,N}^\alpha$ there is a family $(\lambda_k)_k$, where $k$ runs as in Theorem D1, D2, $\lambda_k\in\Lambda_{k,\chi,N}^\alpha$ and $k\equiv k'\pmod{p^m}$ implies that
$$
\lambda(k)_\ell\equiv \lambda(k')_\ell\pmod{p^{{\sf a}m+{\sf b}-A_{k_0}}}
$$
for all $\ell\not|N$.

Let $k_0>({\sf C}_\alpha\alpha+1)^2+2$ be any weight and let $\lambda=(\lambda_{\ell})_{\ell\not|Np}\in\Lambda^\alpha_{k_0,\chi,N}$. In particular, for any weight $k$ satisfying 
$$
k>({\sf C}_\alpha\alpha+1)^2+2\leqno(66)
$$
and
$$
k\equiv k_0\pmod{p^m}\quad\mbox{with}\quad m>{{\sf L}(\alpha,k_0)}\leqno(67)
$$
there is a uniquely determined element $\lambda(k)=(\lambda_{k,\ell})_{\ell\not|N}\in\Lambda^\alpha_{k,\chi,N}$, such that the family $(\lambda(k))_k$ satisfies $\lambda(k_0)=\lambda$ as well as the congruences
$$
\lambda(k)_\ell\equiv\lambda(k')_\ell\pmod{p^{{\sf a}m+{\sf b}-A_{k_0}}}\quad\mbox{if $k\equiv k'\pmod{p^m}$}
$$
for all $\ell\not|N$. Moreover, the family $(\lambda(k))_k$ is uniquely determined by $\lambda$. 
We denote by $(\pi(k),V(k))$ the representation of ${\rm GL}_2({\Bbb A}_f)$ corresponding to $\lambda(k)$. For any prime number $\ell$ dividing $N$, we denote by $\Psi_\ell=\Psi^\alpha_{k,\chi,\ell}$ the characteristic polynomial of $T_\ell$ acting on ${\cal M}_{k,\chi\omega^{-k}}^\alpha(\lambda(k))$ (cf. section 4.2 for the definition of ${\cal M}_{k,\chi\omega^{-k}}^\alpha(\lambda(k))$ with $\lambda\in\Lambda_{k,\chi,N}^\alpha$). We note that the extension of Corollary C2 2.) implies that ${\rm dim}\,{\cal M}_{k,\chi\omega^{-k}}^\alpha(\lambda(k))={\rm dim}\,V(k)^{K,\alpha}={\rm dim}\,V(k_0)^{K,\alpha}={\rm dim}\,{\cal M}_{k_0,\chi\omega^{-k_0}}^\alpha(\lambda(k'))$ as long as $k\equiv k_0\pmod{p^{{\sf L}(\alpha,k_0)}}$. Hence, the degree of $\Psi_\ell$ is the same for all $k$ satisfying (66) and (67) and we may write
$$
\Psi_\ell=\sum_{h=0}^{d_\ell} a_{k,h} X^h
$$
with $d_\ell={\rm dim}\,{\cal M}_{k,\chi\omega^{-k}}^\alpha(\lambda(k))={\rm dim}\,{\cal M}_{k_0,\chi\omega^{-k_0}}^\alpha(\lambda)$.

{\bf Theorem F. }{\it Let $k',k\ge ({\sf C}_\alpha\alpha+1)^2+1$ be weights which satisfy equation (67). Then, with the above notations, for any $\ell|N$ we have
$$
v_p(\Psi_{k,\chi,\ell}^\alpha-\Psi_{k',\chi,\ell}^\alpha)
\ge 
{\sf a}m+{\sf b}-(d_{k_0,\chi}^\alpha+1)A_{k_0})-v_p(d_{k_0,\chi}^\alpha!)).
$$
where $d_{k_0,\chi}^\alpha={\rm dim}\,{\cal M}_{k_0}(\Gamma,\chi\omega^{-k_0})^\alpha$.

}

{\it Proof.} We write $k\equiv k'\pmod{p^m}$, hence, $m>{\sf L}(\alpha,k_0)$. Let $\lambda\in\Lambda_{k_0,\chi,N}^\alpha$. Since elements, which are distinct, already are distinct modulo $p^{A_{k_0}}$ (cf. equation (35)) we find for any $\mu\in\Lambda_{k_0,\chi,N}^\alpha$, $\mu\not=\lambda$ a prime $\ell(\lambda,\mu)\in{\Bbb N}$, which does not divide $Np$ and which satisfies $\lambda_{\ell(\lambda,\mu)}\not\equiv \mu_{\ell(\lambda,\mu)}\pmod{p^{A_{k_0}}}$. We define for any $k$ satisfying (66) and (67), any $a\in{\Bbb N}$ and any $\ell|N$ an element in the Hecke algebra
$$
e_{k}=e_{k,a,\lambda}=\prod_{\mu\in\Lambda_{k_0,\chi,Np}^\alpha\atop \mu\not=\lambda}(T_{\ell(\mu,\lambda)}-\mu(k)_{\ell(\mu,\lambda)})\cdot T_\ell^a\cdot e_\alpha^{[\frac{m}{{\sf C}_\alpha\alpha}]}.
$$
We recall that $\mu(k)\in\Lambda_{k,\chi,N}^\alpha$ is the image of $\mu\in\Lambda_{k_0,\chi,N}^\alpha$ under the map $\varphi$ (cf. the extension of Corollary C2 4.)). We note that the above extension of Theorem D1 implies that the family $(\mu(k))_k$ satisfies the congruences
$$
\mu(k)_\ell\equiv\mu(k')_\ell\pmod{p^{{\sf a}m+{\sf b}-A_{k_0}}}\quad\mbox{for all $\ell\not|N$}.
$$

We compute the trace of $e_{k}$ acting on ${\cal M}_k(\Gamma,\chi\omega^{-k})^\alpha=\bigoplus_{\mu\in\Lambda_{k_0,\chi,N}^\alpha}{\cal M}_{k,\chi\omega^{-k}}^\alpha(\mu(k))$ (cf. the extension of Corollary C2).

The definition of $e_{k}$ immediately yields $e_{k}|_{{\cal M}_{k,\chi\omega^{-k}}^\alpha(\mu(k))}=0$ if $\mu\not=\lambda$.

To determine $e_{k}|_{{\cal M}_{k,\chi\omega^{-k}}^\alpha(\lambda(k))}$ we choose a basis ${\cal B}$ of ${\cal M}_{k,\chi\omega^{-k}}^\alpha(\lambda(k))$ with respect to which all Hecke operators $T_\ell$, $\ell$ prime, are represented by an upper triangular matrix
$$
{\cal D}_{\cal B}(T_\ell|_{{\cal M}_{k,\chi\omega^{-k}}^\alpha(\lambda(k))})=\left(\begin{array}{ccc}
\lambda(k)_\ell&&*\\
&\ddots&\\
&&\lambda(k)_\ell
\end{array}
\right).
$$
We then obtain
$$
{\cal D}_{\cal B}(e_{k}|_{{\cal M}_{k,\chi\omega^{-k}}^\alpha(\lambda(k))})
=\left(
\begin{array}{ccc}
\zeta_k&&*\\
&\ddots&\\
&&\zeta_k
\end{array}
\right){\cal D}_{\cal B}(T_\ell^a|_{{\cal M}_{k,\chi\omega^{-k}}^\alpha(\lambda(k))})\,{\cal D}_{\cal B}(e_\alpha^{[\frac{m}{{\sf C}_\alpha\alpha}]}|_{{\cal M}_{k,\chi\omega^{-k}}^\alpha(\lambda(k))}),
$$
where
$$
\zeta_k=\prod_{\mu\in\Lambda_{k_0,\chi,Np}^\alpha\atop\mu\not=\lambda}(\lambda(k)_{\ell(\lambda,\mu)}-\mu(k)_{\ell(\lambda,\mu)}).
$$
Since ${\cal D}_{\cal B}(T_\ell^a|_{{\cal M}_{k,\chi}^\alpha(\lambda(k))})$ too is upper triangular we find
$$
{\rm tr}\,e_{k}|_{{\cal M}_k(\Gamma,\chi\omega^{-k})^\alpha}={\rm tr}\,e_{k}|_{{\cal M}_{k,\chi\omega^{-k}}^\alpha(\lambda(k))}= \zeta_k\cdot {\rm tr}\;T_\ell^a|_{{\cal M}_{k,\chi\omega^{-k}}^\alpha(\lambda(k))}\;{\rm tr}\,e_\alpha^{[\frac{m}{{\sf C}_\alpha\alpha}]}|_{{\cal M}_{k,\chi\omega^{-k}}^\alpha(\lambda(k))}.\leqno(68)
$$

On the other hand the Hecke algebra $\langle T_\ell,\,\ell\,{\rm prime}\rangle$ is isomorphic to the polynomial algebra ${\Bbb C}[T_\ell,\,\ell\,{\rm prime}]$; for any $f=\sum_{\underline{m}}a_{\underline{m}} \prod_{\ell\,{\rm prime}} T_\ell^{m_\ell}\in \langle T_\ell,\,\ell\,{\rm prime}\rangle$ we define 
$$
v_p(f):={\rm min}_{\underline{m}}\,v_p(a_{\underline{m}}),
$$
i.e. $v_p(f-f')\ge r$ precisely if $f'=f+g$ with $g\in p^r\langle T_\ell,\,\ell\;{\rm prime}\rangle$.

Since $\mu(k)\equiv \mu(k')\pmod{p^{{\sf a}m+{\sf b}-A_{k_0}}}$, the definition of $e_k$, $e_k'$ implies
$$
v_p(e_k-e_{k'})\ge {\sf a}m+{\sf b}-A_{k_0}.
$$ 
In particular, if $V$ is any ${\Bbb C}$-vector space with basis ${\cal C}$, on which the Hecke algebra acts we obtain 
$$
{\cal D}_{\cal C}(e_{k}|V)\equiv{\cal D}_{\cal C}(e_{k'}|V)\pmod{p^{{\sf a}m+{\sf b}-A_{k_0}}}.
$$
(a congruence between matrices is to be understood entry-wise) and therefore
$$
{\rm tr}\,e_{k}|_V\equiv {\rm tr}\,e_{k'}|_V\pmod{p^{{\sf a}m+{\sf b}-A_{k_0}}}.\leqno(69)
$$
On the other hand, the Theorem in section 3.5 yields
$$
{\rm tr}\,e_{k}|_{{\cal M}_k(\Gamma,\chi\omega^{-k})^\alpha}\equiv {\rm tr}\,e_{k}|_{{\cal M}_{k'}(\Gamma,\chi\omega^{-k'})^\alpha}
\pmod{p^{\left[\frac{m}{{\sf C}_\alpha\alpha}\right]-v_p(\varphi(N))}}\leqno(70)
$$
and since ${\sf a}m+{\sf b}\le{\sf D}\le\left[\frac{m}{{\sf C}_\alpha\alpha}\right]-v_p(\varphi(N))$ (cf. the definition of ${\sf a,b}$ (and ${\sf D}={\sf D}_{\alpha,m}$) in section 4.2) equations (70) and (69) yield
$$
{\rm tr}\,e_{k}|_{{\cal M}_k(\Gamma,\chi\omega^{-k})^\alpha}\equiv {\rm tr}\,e_{k'}|_{{\cal M}_{k'}(\Gamma,\chi\omega^{-k'})^\alpha}
\pmod{p^{{\sf a}m+{\sf b}-A_{k_0}}}.
$$
Together with equation (68) we obtain
\begin{eqnarray*}
(71)&&\zeta_k\,{\rm tr}\,T_\ell^a|_{{\cal M}_{k,\chi\omega^{-k}}^\alpha(\lambda(k))}\,{\rm tr}\,e_\alpha^{[\frac{m}{{\sf C}_\alpha\alpha}]}|_{{\cal M}_{k,\chi\omega^{-k}}^\alpha(\lambda(k))}\\
&\equiv&
\zeta_{k'}{\rm tr}\,T_\ell^a|_{{\cal M}_{k',\chi\omega^{-k'}}^\alpha(\lambda(k'))}\,{\rm tr}\, e_\alpha^{[\frac{m}{{\sf C}_\alpha\alpha}]}|_{{\cal M}_{k',\chi\omega^{-k'}}^\alpha(\lambda(k'))}
\pmod{p^{{\sf a}m+{\sf b}-A_{k_0}}}.\\
\end{eqnarray*}
Again, since $\mu(k)\equiv \mu(k')$ and $\lambda(k)\equiv \lambda(k')$ $\pmod{p^{{\sf a}m+{\sf b}-A_{k_0}}}$ we find 
$$
\zeta_k\equiv \zeta_{k'}\pmod{p^{{\sf a}m+{\sf b}-A_{k_0}}}.
$$ 
We know that $\lambda_{\ell(\lambda,\mu)}\not\equiv \mu_{\ell(\lambda,\mu)}\pmod{p^{A_{k_0}}}$ for all $\mu\not=\lambda$. Moreover, since ${\sf a}{\sf L}(\alpha,k_0)+{\sf b}\ge A_{k_0}$ (by the definition of ${\sf L}(\alpha,k_0)$) we also obtain $\mu(k)\equiv\mu$ and $\lambda(k)\equiv\lambda$ $\pmod{p^{A_{k_0}}}$. Thus, we obtain altogether
$$
\lambda(k)_{\ell(\lambda,\mu)}\not\equiv  \mu(k)_{\ell(\lambda,\mu)}\pmod{p^{A_{k_0}}}.
$$
In particular,
$$
v_p(\zeta_k)\le |\Lambda_{k,\chi,N}^\alpha| A_{k_0} \le{\rm dim}\, {\cal M}_k(\Gamma,\chi\omega^{-k})^\alpha A_{k_0}={\rm dim}\, {\cal M}_{k_0}(\Gamma,\chi\omega^{-k_0})^\alpha A_{k_0}.
$$
The last equality is a consequence of Theorem B. We set $d_{k_0,\chi}^\alpha={\rm dim}\, {\cal M}_{k_0}(\Gamma,\chi\omega^{-k_0})^\alpha$ and thus obtain from (71)
\begin{eqnarray*}
&&{\rm tr}\,T_\ell^a|_{{\cal M}_{k,\chi\omega^{-k}}^\alpha(\lambda(k))}\,{\rm tr}\,e_\alpha^{[\frac{m}{{\sf C}_\alpha\alpha}]}|_{{\cal M}_{k,\chi\omega^{-k}}^\alpha(\lambda(k))}\\
&\equiv&  
{\rm tr}\,T_\ell^a|_{{\cal M}_{k',\chi\omega^{-k'}}^\alpha(\lambda(k'))}\,{\rm tr}\,e_\alpha^{[\frac{m}{{\sf C}_\alpha\alpha}]}|_{{\cal M}_{k',\chi\omega^{-k'}}^\alpha(\lambda(k'))} 
\pmod{p^{{\sf a}m+{\sf b}-A_{k_0}-d_{k_0,\chi}^\alpha A_{k_0}}}.\\
\end{eqnarray*}
%
%[[beachte: modulo $p^{{\sf a}m+{\sf b}-A_{k_0}}$ konnen wir $\zeta_k$ durch $\zeta_{k_0}$ ersetzen]]
%
Moreover, since $\lambda(k)_p\equiv \lambda(k')_p$ and $e_\alpha$ is a polynomial in $T_p$ we deduce that ${\rm tr}\,e_\alpha^{[\frac{m}{{\sf C}_\alpha\alpha}]}|_{{\cal M}_{k,\chi\omega^{-k}}^\alpha(\lambda(k))}\equiv {\rm tr}\,e_\alpha^{[\frac{m}{{\sf C}_\alpha\alpha}]}|_{{\cal M}_{k',\chi\omega^{-k'}}^\alpha(\lambda(k'))}\pmod{p^{{\sf a}m+{\sf b}-A_{k_0}}}$. Since furthermore ${\rm tr}\,e_\alpha^{[\frac{m}{{\sf C}_\alpha\alpha}]}|_{{\cal M}_{k,\chi\omega^{-k}}^\alpha(\lambda(k))}$ and ${\rm tr}\,e_\alpha^{[\frac{m}{{\sf C}_\alpha\alpha}]}|_{{\cal M}_{k',\chi\omega^{-k'}}^\alpha(\lambda(k'))}$ are $p$-adic units by Proposition 1.) in section 3.4, we obtain
$$
{\rm tr}\,T_\ell^a|_{{\cal M}_{k,\chi\omega^{-k}}^\alpha(\lambda(k))}\equiv  
{\rm tr}\,T_\ell^a|_{{\cal M}_{k',\chi\omega^{-k'}}^\alpha(\lambda(k'))} 
\pmod{p^{{\sf a}m+{\sf b}-A_{k_0}-d_{k_0,\chi}^\alpha A_{k_0}}}.\leqno(72)
$$

We write
$$
\Psi_{k,\chi,\ell}^\alpha=\sum_{j=1}^d (-1)^j a_{k,\chi,j}^\alpha X^{d-j}
$$
and
$$
\Psi_{k',\chi,\ell}^\alpha=\sum_{j=1}^d (-1)^j a_{k',\chi,j}^\alpha X^{d-j}.
$$
We set $\sigma_i={\rm tr}\,T_\ell^i|_{{\cal M}_{k,\chi\omega^{-k}}^\alpha}(\lambda(k))$. Since 
$$
a_{k,\chi,j}^\alpha=\frac{1}{j!}\det \left(
\begin{array}{ccccc}
\sigma_1&1&&&\\
\sigma_2&\sigma_1&2&&\\
&\cdots&&\cdots&\\
&\cdots&&\cdots&j-1\\
\sigma_j&\sigma_{j-1}&\cdots&\cdots&\sigma_1
\end{array}
\right)
$$
and the same formula holds if we replace $k$ by $k'$ (cf. [Koe], 3.4.6. Korollar, p. 117) we deduce from (72)
$$
a_{k,\chi,j}^\alpha\equiv a_{k',\chi,j}^\alpha\pmod{p^{{\sf a}m+{\sf b}-(d_{k_0,\chi}^\alpha+1)A_{k_0}-v_p(j!)}}.
$$
Since $j\le {\rm dim}\,{\cal M}_{k,\chi\omega^{-k}}^\alpha(\lambda(k))\le {\rm dim}\,{\cal M}_k(\Gamma,\chi\omega^{-k})^\alpha\stackrel{\rm Thm B}{=}d_{k_0,\chi}^\alpha$ we obtain 
$$
\Psi_{k,\chi,\ell}^\alpha\equiv \Psi_{k',\ell,\chi}^\alpha\pmod{p^{{\sf a}m+{\sf b}-(d_{k_0,\chi}^\alpha+1)A_{k_0}-v_p(d_{k_0,\chi}^\alpha!)}}.
$$ 
This completes the proof of the Theorem.

Let $\lambda\in\Lambda_{k_0,\chi}^\alpha$. Then there is a family $(\lambda(k))_k$, where $k\ge({\sf C}_\alpha\alpha+1)^2+2$ and $k\equiv k_0\pmod{p^m}$ with $m\ge {\sf L}(\alpha,k_0)$ such that 
$$
\lambda(k)_\ell\equiv\lambda(k')_\ell\pmod{p^{{\sf a}m+{\sf b}-A_{k_0}}}\quad \mbox{if}\quad k\equiv k'\pmod{p^{{\sf a}m+{\sf b}-A_{k_0}}}
$$
for all $\ell\not|N$ (cf. the extension of Theorem D1). Applying the above Proposition to the family of characteristic polynomials $\Psi_{k,\chi}^\alpha$ of $T_\ell$ acting on ${\cal M}_k(\Gamma,\chi\omega^{-k})^\alpha$, $\ell|N$ (which is possible by Theorem F) we see that for any $k$ we can choose eigenvalues $\lambda(k)_\ell$ of $T_\ell$ acting on ${\cal M}_k(\Gamma,\chi\omega^{-k})^\alpha$ such that the family $(\lambda(k)_\ell)_k$ satisfies the congruences
$$
v_p(\lambda(k)_\ell-\lambda(k')_\ell)\ge \frac{{\sf a}m+{\sf b}-(d_{k_0,\chi}^\alpha+1)A_{k_0}-v_p(d_{k_0,\chi}^\alpha!)}{d_{\lambda,\ell}}\quad \mbox{if}\quad k\equiv k'\pmod{p^m},
$$
where $d_{\lambda,\ell}={\rm deg}\,\Psi_{k,\chi,\ell}^\alpha={\rm dim}\,{\cal M}_{k_0,\chi\omega^{-k_0}}^\alpha(\lambda)$. In particular, for any $k$ as above we obtain elements $(\lambda(k)_\ell)_\ell\in\Lambda_{k,\chi}^\alpha$, where now $\ell$ runs over {\it all} rational primes $\ell$. Hence, any $\lambda\in\Lambda_{k_0,\chi}^\alpha$ can be placed into a $p$-adic family which satisfies the requested type of congruences also at the ramified places. Moreover, since in particular, $d_{\lambda,\ell}\le {\rm dim}\,{\cal M}_{k_0}(\Gamma,\chi\omega^{-k_0})^\alpha=d_{\chi,k_0}^\alpha$ in conjunction with the extension of Theorems D1 and D2 we obtain the following final result.

{\bf Theorem G. }{\it Let $\alpha\in{\Bbb Q}_{\ge 0}$ and let $k_0\ge({\sf C}_\alpha\alpha+1)^2+2$ be any weight. Let $\lambda=(\lambda_\ell)_\ell\in\Lambda_{k_0,\chi}^\alpha$. Then, for any weight $k$ satisfying $k>({\sf C}_\alpha\alpha+1)^2+2$ and $k\equiv k'\pmod{p^m}$ with $m>{\sf L}(\alpha,k_0)$ there is a system of eigenvalues $\lambda_k=(\lambda_{k,\ell})$, $\ell$ running over all rational primes, such that the following conditions hold:

$\bullet$ $\lambda(k)\in\Lambda_{k,\chi}^\alpha$

$\bullet$ $\lambda(k_0)=\lambda$

$\bullet$ $k\equiv k'\pmod{p^m}$ implies 
$$
\lambda(k)_\ell\equiv\lambda(k')_\ell\pmod{p^{{\sf a}m+{\sf b}-A_{k_0}}}
$$ 
for all $\ell\not|N$ and 
$$
v_p(\lambda(k)_\ell-\lambda(k')_\ell)\ge \frac{{\sf a}m+{\sf b}-(d_{k_0,\chi}^\alpha+1)A_{k_0}-v_p(d_{k_0,\chi}^\alpha!)}{d_{k_0,\chi}^\alpha}
$$ 
for all $\ell|N$. Here, $d_{k_0,\chi}^\alpha={\rm dim}\,{\cal M}_{k_0}(\Gamma,\chi\omega^{-k_0})^\alpha$.

$\bullet$ if $\lambda$ is cuspidal then any $\lambda(k)$ is cuspidal.

Here, ${\sf a}={\sf a}(\alpha)$ and ${\sf b}={\sf b}(\alpha)$ only depend on $\alpha$. Moreover, ${\sf a}$ satisfies
$$
{\sf a}\le\frac{1}{{\rm dim}\,{\cal M}_{k}(\Gamma,\chi\omega^{-k})^\alpha}
$$
for all $k\ge ({\sf C}_\alpha\alpha+1)^2+1$ and even
$$
{\sf a}={\rm min}\,\{ \frac{1}{{\rm dim}\,{\cal M}_{k}(\Gamma,\chi\omega^{-k})^\alpha},\;k\ge ({\sf C}_\alpha\alpha+1)^2+1\}
$$
if $\alpha=0$.

}

We note that Theorem G has an obvious reformulation in terms of modular forms.

%Note: We do not obtain a statement on spaces of cusp forms because Proposition E does not include the prime $p$

\newpage

{\Large\bf References}

[B] Bewersdorff, J., Eine Lefschetzsche Fixpunktformel fur Hecke Operatoren, Bonner Math. Schriften {\bf 164}, Bonn, 1985

[Bu]  Buzzard, K., Families of modular forms. 21st Journ{\'e}es arithm{\'e}tique (Rome 2001), J. Th{\'e}or. Nombres Bordeaux 13 (2001), no. 1, 43-52

[B-C] Buzzard, K., Calegari, F., A counterexample to the Mazur-Gouvea Conjecture, C. R. Math. Acad. Sci. Paris 338 (2004), no. 10 ,751-753

[Ca] Casselman, W., On some results of Atkin and Lehner, Math. Ann. {\bf 201}, 301-314 (1973)

[C] Coleman, R., $p$-adic Banach spaces and families of modular forms, Inv. Math. {\bf 127}, 417 - 479 (1997)

%[C2] Coleman, R., Classical and overconvergent modular forms, Inv. Math. {\bf 124}, 215 - 241 (1996)

[D-S] Diamond, F., Shurman, J., A first course in modular forms, Springer, GTM {\bf 228}, 2005

[H] Harder, G., Eisenstein cohomology of arithmetic groups. The case ${\rm GL}_2$, Inv. Math. {\bf 89}, 37 - 118 (1987)

[Hi] Hida, H., Elementary theory of $L$-functions and Eisenstein series, Cambridge Univ. Press, 1993

[Koe] Koecher, M., Lineare Algebra und analytische Geometrie, 2. Auflage, Springer, 1985

[G-MacPh] Goresky, M., MacPherson, R., The topological trace formula, J. reine angew. Math. {\bf 560}, 77 - 150 (2003) 

[M] Miyake, T., Modular forms, Springer 2006

[M-G] Mazur, B., Gouvea, F., Families of modular eigenforms, Math. Comp. {\bf 58} (198), 793 - 805 (1992)

[N] Neukirch J., Algebraische Zahlentheorie, Springer 1992
 
[Wa] Wan, D., Dimension variation of spaces classical and $p$-adic modular forms, Inv. Math 133 (1998), 449-463

\end{document}